\documentclass[12pt, leqno, openany]{book}

\usepackage{makeidx}
\makeindex

\def\Rad{\mathop{\rm Rad}}
\def\Re{\mathop{\rm Re}}
\def\span{\mathop{\rm span}}
\def\Int{\mathop{\rm Int}}
\def\Order{\mathop{\rm Order}}

\newtheorem{theorem}{Theorem}
\newtheorem{lemma}[theorem]{Lemma}
\newtheorem{proposition}[theorem]{Proposition}
\newtheorem{sublemma}[theorem]{Sublemma}
\newtheorem{definition}[theorem]{Definition}
\newtheorem{corollary}[theorem]{Corollary}
\newtheorem{problem}[theorem]{Problem}
\newtheorem{remark}[theorem]{Remark}
\newtheorem{claim}[theorem]{Claim}
\newtheorem{assumptions}[theorem]{Assumptions}
\newtheorem{sassumptions}[theorem]{Standing Assumptions}
\newtheorem{examples}[theorem]{Examples}
\newtheorem{exercise}[theorem]{Exercise}

\newcommand{\begintheorem}{\addtocounter{equation}{1}\begin{theorem}}
\newcommand{\beginlemma}{\addtocounter{equation}{1}\begin{lemma}}
\newcommand{\beginproposition}{\addtocounter{equation}{1}\begin{proposition}}
\newcommand{\beginsublemma}{\addtocounter{equation}{1}\begin{sublemma}}
\newcommand{\begindefinition}{\addtocounter{equation}{1}\begin{definition}}
\newcommand{\begincorollary}{\addtocounter{equation}{1}\begin{corollary}}
\newcommand{\beginproblem}{\addtocounter{equation}{1}\begin{problem}}
\newcommand{\beginremark}{\addtocounter{equation}{1}\begin{remark}}
\newcommand{\beginclaim}{\addtocounter{equation}{1}\begin{claim}}
\newcommand{\beginassumptions}{\addtocounter{equation}{1}\begin{assumptions}}
\newcommand{\beginsassumptions}{\addtocounter{equation}{1}\begin{sassumptions}}
\newcommand{\beginexamples}{\addtocounter{equation}{1}\begin{examples}}
\newcommand{\beginexercise}{\addtocounter{equation}{1}\begin{exercise}}

\begin{document}

\frontmatter

\title{Elements of Linear and Real Analysis}

\author{Stephen Semmes \\
	Rice University \\
	Houston, Texas}

\date{}

\maketitle

\chapter{Preface}

	This book deals with some basic themes in mathematical
analysis along the lines of classical norms on functions and
sequences, general normed vector spaces, inner product spaces, linear
operators, some maximal and square-function operators, interpolation
of operators, and quasisymmetric mappings between metric spaces.
Aspects of the broad area of harmonic analysis are entailed in
particular, involving famous work of M.~Riesz, Hardy, Littlewood,
Paley, Calder\'on, and Zygmund.

	However, instead of working with arbitrary continuous or
integrable functions, we shall often be ready to use only step
functions on an interval, i.e., functions which are
piecewise-constant.  Similarly, instead of infinite-dimensional
Hilbert or Banach spaces, we shall frequently restrict our attention
to finite-dimensional inner product or normed vector spaces.  We
shall, however, be interested in quantitative matters.

	We do not attempt to be exhaustive in any way, and there are
many related and very interesting subjects that are not addressed.
The bibliography lists a number of books and articles with further
information.

	The formal prerequisites for this book are quite limited.
Much of what we do is connected to the notion of integration, but for
step functions ordinary integrals reduce to finite sums.  A sufficient
background should be provided by standard linear algebra of real and
complex finite-dimensional vector spaces and some knowledge of
beginning analysis, as in the first few chapters of Rudin's celebrated
\emph{Principles of Mathematical Analysis} \cite{Ru1}.  This is not to
say that the present monograph would necessarily be easy to read with
this background, as the types of issues considered may be unfamiliar.
On the other hand, it is hoped that this monograph can be helpful to
readers with a variety of perspectives.

\tableofcontents

\mainmatter

\chapter{Notation and conventions}
\label{Notation and conventions}

	If $a$ and $b$ are real numbers with $a \le b$, then the
following are the intervals\index{intervals} in the real line ${\bf
R}$\index{$R$@${\bf R}$} with endpoints $a$ and $b$:
\begin{eqnarray*}
	[a,b] & = & \{x \in {\bf R} : a \le x \le b\};		\\
	(a,b) & = & \{x \in {\bf R} : a < x < b\};		\\  \
	[a,b) & = & \{x \in {\bf R} : a \le x < b\};		\\
	(a,b] & = & \{x \in {\bf R} : a < x \le b\}.		
\end{eqnarray*}
All but the first is the empty set when $a = b$, while $[a,b]$
consists of the one point $a = b$.  In general, the first of these
intervals is called the closed interval with endpoints $a$ and $b$,
and the second is the open interval with endpoints $a$ and $b$.  The
third and fourth are half-open, half-closed intervals, with the third
being left-closed and right-open, and the fourth left-open and
right-closed.

	The \emph{length}\index{length of an interval} of each of
these intervals is defined to be $b-a$.  If an interval is denoted
$I$, we may write $|I|$ \index{$"|I"|$ ($I$ an interval)} for the
length of $I$.

	For the record, see Chapter 1 in \cite{Ru1} concerning
detailed properties of the real numbers (as well as the complex
numbers ${\bf C}$).\index{$C$@${\bf C}$} In particular, let us recall the
``least upper bound'' or ``completeness'' property, to the effect that
a nonempty set $E$ of real numbers which has an upper bound has a
least upper bound.  The least upper bound is also called the supremum
of $E$, and is denoted $\sup E$.  Similarly, if $F$ is a nonempty set
of real numbers which has a lower bound, then $F$ has a greatest lower
bound, or infimum, which is denoted $\inf F$.  We shall sometimes use
the \emph{extended real numbers} (as in \cite{Ru1}), with $\infty$ and
$-\infty$ added to the real numbers, and write $\sup E = \infty$ and
$\inf F = -\infty$ if $E$ and $F$ are nonempty sets of real numbers
such that $E$ does not have an upper bound and $F$ does not have a
lower bound.

	If $A$ is a subset of some specified set $X$ (like the real
line), we let ${\bf 1}_A(x)$ denote the \emph{indicator
function}\index{indicator functions} of $A$ on $X$ (sometimes called
the \emph{characteristic function}\index{characteristic functions}
associated to $A$, although in other contexts this name can be used
for something quite different).  This is the function which is equal
to $1$ when $x \in A$, and is equal to $0$ when $x \in X \backslash
A$.

\begindefinition [Step functions]
\label{Step functions}
\index{step functions} A function on the real line, or on an interval
in the real line, is called a \emph{step function} if it is a finite
linear combination of indicator functions of intervals.
\end{definition}

	This is equivalent to saying that there is a partition of the
domain into intervals on which the function is constant.

	In these notes, one is normally welcome to assume that a given
function on the real line, or on an interval in the real line, is a
step function.  In fact, one is normally welcome to assume that a
given function is a dyadic step function, as defined in the next
chapter.  For step functions, it is very easy to define the integral
over an interval in the domain of definition, by reducing to linear
combinations of lengths of intervals.

	An exception to this convention occurs when we consider convex
or monotone functions, which we do not necessarily wish to ask to be
step functions.  When dealing with integrals, typically the function
being integrated can be taken to be a step function.  (This function
might be the composition of a non-step function with a step function,
which is still a step function.)

\chapter{Dyadic intervals}

\section{The unit interval and dyadic subintervals}
\label{The unit interval and dyadic subintervals}\index{unit interval}

	Normally, a reference to ``the unit interval'' might suggest
the interval $[0,1]$ in the real line.  It will be convenient to use
$[0,1)$ instead, for minor technical reasons (and one could easily
work around this anyway).

\begindefinition [Dyadic intervals in $[0,1)$]
\label{Dyadic intervals in [0,1)}
\index{dyadic intervals}\index{intervals!dyadic}
The \emph{dyadic subintervals} of the unit interval $[0,1)$ are the
intervals of the form $[j \, 2^{-k}, (j+1) \, 2^{-k})$, where $j$ and
$k$ are nonnegative integers, and $j+1 \le 2^k$.  (Thus the length of
such an interval is of the form $2^{-k}$, where $k$ is a nonnegative
integer.)
\end{definition}

	In general one can define the dyadic intervals in ${\bf R}$
to be the intervals of the same form, except that $j$ and $k$ are
allowed to be arbitrary integers.

	The half-open, half-closed condition leads to nice properties
in terms of disjointness, as in the following lemmas.  (With closed
intervals one could get disjointness of interiors in similar
circumstances.  This would be fine in terms of integrals, measures,
etc.)

\beginlemma [Partitions of $[0,1)$]
\label{Partitions of [0,1)}
For each nonnegative integer $k$, $[0,1)$ is the union of the dyadic
subintervals of itself of length $2^{-k}$, and these intervals are
pairwise disjoint.
\end{lemma}

\beginlemma [Comparing pairs of intervals]
\label{Comparing pairs of intervals}
If $J_1$ and $J_2$ are two dyadic subintervals of $[0,1)$, then either
$J_1 \subseteq J_2$, or $J_2 \subseteq J_1$, or $J_1 \cap J_2 =
\emptyset$.  (The first two possibilities are not mutually exclusive,
as one could have $J_1 = J_2$.)
\end{lemma}

	These two lemmas are not hard to verify, just from the
definitions.  (Exercise.)  For the second lemma, one might be a bit
more precise and say the following.  Suppose that $J_1$ and $J_2$ are
dyadic subintervals of $[0,1)$, and that $J_1$ has length $2^{-k_1}$,
and $J_2$ has length $2^{-k_2}$.  If $k_1 \le k_2$ (so that $2^{-k_1}
\ge 2^{-k_2}$), then either $J_2 \subseteq J_1$ or $J_1 \cap J_2 =
\emptyset$.  

\beginlemma [Partitions of dyadic intervals]
\label{Partitions of dyadic intervals}
If $J$ is a dyadic subinterval of $[0,1)$ of length $2^{-k}$, and if
$m$ is an integer greater than $k$, then $J$ is the union of the
dyadic subintervals of itself of length $2^{-m}$ (and these intervals
are pairwise disjoint).  Also, every dyadic subinterval of $[0,1)$ of
length $2^{-m}$ is contained in a dyadic subinterval of $[0,1)$ of
length $2^{-k}$ when $m \ge k$.
\end{lemma}

	(Exercise.)

\beginlemma [Structure of unions of dyadic intervals]
\label{Structure of unions of dyadic intervals}
Let $\mathcal{F}$ be an arbitrary nonempty collection of dyadic
subintervals of $[0,1)$.  Then there is a subcollection
$\mathcal{F}_0$ of $\mathcal{F}$ such that
\begin{equation}
\label{bigcup_{J in mathcal{F}_0} J = bigcup_{J in mathcal{F}} J}
	\bigcup_{J \in \mathcal{F}_0} J = \bigcup_{J \in \mathcal{F}} J
\end{equation}
and the elements of $\mathcal{F}_0$ are pairwise disjoint, i.e., if
$J_1$, $J_2$ are distinct elements of $\mathcal{F}_0$, then $J_1 \cap
J_2 = \emptyset$.
\end{lemma}

	To prove this, we take $\mathcal{F}_0$ to be the set of
\emph{maximal}\index{maximal dyadic intervals} elements of
$\mathcal{F}$.  Here ``maximal'' means maximal with respect to the
ordering that comes from set-theoretic inclusion, so that $J \in
\mathcal{F}$ is maximal if there is no element $K$ of $\mathcal{F}$
such that $J \subseteq K$ and $J \ne K$.

	If $J$ is any dyadic subinterval of $[0,1)$, then there are
only finitely many dyadic subintervals of $[0,1)$ which contain $J$ as
a subset.  Indeed, if $J$ has length $2^{-j}$, then there is exactly
one such dyadic interval of length $2^{-k}$ for each nonnegative
integer $k$, $k \le j$.

	Using this, it is easy to see that every interval in
$\mathcal{F}$ is contained in a maximal interval in $\mathcal{F}$.  In
particular, the set $\mathcal{F}_0$ of maximal elements of
$\mathcal{F}$ is nonempty, since $\mathcal{F}$ is nonempty, by
assumption.  We also obtain (\ref{bigcup_{J in mathcal{F}_0} J =
bigcup_{J in mathcal{F}} J}).

	Any two maximal elements of $\mathcal{F}$ which are distinct
are disjoint.  This follows from Lemma \ref{Comparing pairs of
intervals}.  This proves the second property of $\mathcal{F}_0$ in
Lemma \ref{Structure of unions of dyadic intervals}.

\section{Functions on the unit interval}
\label{Functions on the unit interval}

	Suppose that $f$ is a function on the unit interval $[0,1)$
(real or complex-valued).  We assume that $f$ is at least mildly-well
behaved, so that it makes sense to talk about integrals of $f$ over
subintervals of $[0,1)$ (and with the customary basic properties).

	Assume that a nonnegative integer $k$ is also given.  Let us
define a function $E_k(f)$\index{$E_k$} on $[0,1)$ as follows.  If $x$
is any element of $[0,1)$, then there is a unique dyadic subinterval
$J$ of $[0,1)$ of length $2^{-k}$ which contains $x$, as in Lemma
\ref{Dyadic intervals in [0,1)}.  Then we set
\begin{equation}
\label{def of E_k(f)}
	E_k(f)(x) = 2^{-k} \int_J f(y) \, dy.
\end{equation}
Note that $E_k(f)$ is linear in $f$, so that
\begin{equation}
	E_k(a \, f_1 + b \, f_2) = a \, E_k(f_1) + b \, E_k(f_2)
\end{equation}
when $a$, $b$ are constants and $f_1$, $f_2$ are functions on $[0,1)$
(which are at least mildly-well behaved, as before).

\beginlemma [Some properties of $E_k(f)$]
\label{Some properties of E_k(f)}
(a) For each $f$, $E_k(f)$ is constant on the dyadic subintervals of
$[0,1)$ of length $2^{-k}$.

(b) If $f$ is constant on the dyadic subintervals of $[0,1)$ of length
$2^{-k}$, then $E_k(f) = f$.  (In particular, note that $E_k(1) = 1$
for all $k$.)

(c) For all $f$, if $j$ is an integer such that $j \ge k$, then
$E_j(E_k(f)) = E_k(f)$.  Also, $E_k(E_j(f)) = E_k(f)$ in this case.

(d) If $g$ is a function on $[0,1)$ which is constant on the dyadic
subintervals of $[0,1)$ of length $2^{-k}$, then $E_k(g \, f) = g \,
E_k(f)$ for all $f$.
\end{lemma}

	(Exercise.  Note that the first part of (c) holds simply
because $E_k(f)$ is constant on dyadic subintervals of $[0,1)$ of
length $2^{-j}$ when $j \ge k$, while in the second part one is first
averaging $f$ on the (smaller) dyadic intervals of length $2^{-j}$ to
get $E_j(f)$, and then averaging the result on the larger dyadic
intervals of size $2^{-k}$ to get $E_k(E_j(f))$, and the conclusion is
that this is the same as averaging over the dyadic intervals of length
$2^{-k}$ directly.  As in Lemma \ref{Partitions of dyadic intervals},
dyadic subintervals of $[0,1)$ of size $2^{-k}$ are disjoint unions of
dyadic intervals of size $2^{-j}$ when $j \ge k$.)

\begindefinition [Dyadic step functions]
\label{Dyadic step functions}
\index{dyadic step functions}\index{step functions!dyadic} A function
$f$ on $[0,1)$ is called a \emph{dyadic step function} if it is a
finite linear combination of indicator functions of dyadic
subintervals of $[0,1)$.
\end{definition}

\beginlemma 
\label{characterizations of dyadic step functions}
Let $f$ be a function on $[0,1)$.  The following are equivalent:

	(a) $f$ is a dyadic step function;

	(b) There is a nonnegative integer $k$ such that $f$ is
constant on every dyadic subinterval of $[0,1)$ of length $2^{-k}$;

	(c) $E_k(f) = f$ for some nonnegative integer $k$ (and hence
for all sufficiently large integers $k$).
\end{lemma}

	(Exercise, using Lemma \ref{Some properties of E_k(f)}.)

\section{Haar functions}
\label{Haar functions}\index{Haar functions}

	Let $I$ be a dyadic subinterval of $[0,1)$.  Thus $I$ is the
disjoint union of two dyadic subintervals $I_l$ and $I_r$ of half the
size of $I$, corresponding to the left and right halves of $I$.
Define the Haar function $h_I(x)$ on $[0,1)$ associated to the
interval $I$ by
\begin{eqnarray}
\label{def of h_I(x)}
	h_I(x) & = & - |I|^{1/2} \qquad\hbox{when } x \in I_l
								\\
	       & = & |I|^{1/2} \qquad\enspace\, \hbox{when } x \in I_r
							\nonumber \\
	       & = & 0 \qquad\qquad\enspace\hbox{when } 
					x \in [0,1) \backslash I.
							\nonumber
\end{eqnarray}
Notice that
\begin{equation}
\label{int_{[0,1)} h_I(x) dx = 0}
	\int_{[0,1)} h_I(x) \, dx = 0
\end{equation}
and
\begin{equation}
\label{int_{[0,1)} h_I(x)^2 dx = 1}
	\int_{[0,1)} h_I(x)^2 \, dx = 1.
\end{equation}

	In addition to these functions, we define a special Haar 
function $h_0(x)$ on $[0,1)$ by simply taking $h_0(x) = 1$ for all
$x \in [0,1)$.  For this function we also have that
\begin{equation}
	\int_{[0,1)} h_0(x)^2 \, dx = 1.
\end{equation}

	If $I$ and $J$ are distinct dyadic subintervals of $[0,1)$,
then $h_I$ and $h_J$ satisfy the orthogonality condition
\begin{equation}
\label{int_{[0,1)} h_I(x) h_J(x) dx = 0}
	\int_{[0,1)} h_I(x) \, h_J(x) \, dx = 0.
\end{equation}
To see this, it is helpful to consider some cases separately.  If $I$
and $J$ are disjoint, then the product $h_I(x) \, h_J(x)$ is equal to
$0$ for all $x \in [0,1)$, and the integral vanishes automatically.
Otherwise one of the intervals $I$ and $J$ is contained in the other.
We may as well assume that $J \subseteq I$, since the two cases are
completely symmetric.  Since $J \ne I$, $J$ is either a subinterval of
$I_l$ or of $I_r$.  In both situations we have that $h_I$ is constant
on $J$, where $h_J$ is concentrated.  As a result, (\ref{int_{[0,1)}
h_I(x) h_J(x) dx = 0}) follows from (\ref{int_{[0,1)} h_I(x) dx = 0}).

	In addition, if $I$ is any dyadic subinterval of $[0,1)$,
we have the orthogonality relation
\begin{equation}
	\int_{[0,1)} h_0(x) \, h_I(x) \, dx = 0,
\end{equation}
again by (\ref{int_{[0,1)} h_I(x) dx = 0}).

	Given a nonnegative integer $k$, consider the space of dyadic
step functions on $[0,1)$ which are constant on dyadic intervals of
length $2^{-k}$.  This space contains the Haar functions $h_0$ and
$h_I$ for all dyadic subintervals $I$ of $[0,1)$ such that $|I| \ge
2^{-k+1}$.  In fact, every dyadic step function on $[0,1)$ which is
constant on dyadic intervals of length $2^{-k}$ is a linear
combination of these Haar functions.  This can be verified using
induction on $k$.  Alternatively, one can note that the space of these
functions has dimension $2^k$, and the total number of Haar functions
being used is also $2^k$.  More precisely, for each nonnegative
integer $j < k$ there are $2^j$ Haar functions $h_J$ associated to
dyadic subintervals of $[0,1)$ of length $2^{-j}$, and if we sum over
$j$ then the total number of these is $\sum_{j=0}^{k-1} 2^j = 2^k -
1$.  If we add the Haar function $h_0$, then we obtain a total of
$2^k$ Haar functions in this space.

	In particular, the space of all dyadic step functions on
$[0,1)$ is equal to the space of functions which can be given as
finite linear combinations of Haar functions (associated to arbitrary
dyadic subintervals of the unit interval, and also the Haar function
$h_0$).

	If $f$ is a dyadic step function on $[0,1)$, then 
\begin{equation}
	f = \langle f, h_0 \rangle \, h_0 
		+ \sum_I \langle f, h_I \rangle \, h_I,
\end{equation}
where $\langle f, h_0 \rangle = \int_{[0,1)} f(x) \, h_0(x) \, dx$,
and analogously for $h_I$ instead of $h_0$, and where the sum is taken
over all dyadic subintervals $I$ of $[0,1)$.  The sum is in fact a
finite sum, in the sense that all but finitely many terms are equal to
$0$.  This expression for $f$ follows from the fact that $f$ is equal
to a finite linear combination of the Haar functions, and from the
orthonormality conditions for the Haar functions given above.

	Similarly, if $k$ is a nonnegative integer, and $E_k(f)$ is as
defined in Section \ref{Functions on the unit interval}, then
\begin{equation}
	E_k(f) = \langle f, h_0 \rangle \, h_0 
		+ \sum_{|I| \ge 2^{-k+1}} \langle f, h_I \rangle \, h_I.
\end{equation}
Now the sum is taken over all dyadic subintervals $I$ of $[0,1)$ whose
length is at least $2^{-k+1}$, and this is interpreted as being $0$
when $k = 0$.

\section{Binary sequences}
\label{Binary sequences}\index{binary sequences}

	Let $\mathcal{B}$ denote the set of all \emph{binary
sequences}, i.e., sequences $\{x_j\}_{j=1}^\infty$ such that $x_j = 0$
or $1$ for all $j$.  There is a natural correspondence between binary
sequences and real numbers in the interval $[0,1]$, through standard
binary expansions.  More precisely, every binary sequence defines an
element of $[0,1]$, every element of $[0,1]$ has such a binary
expansion, and the binary expansion is unique except for rational
numbers of the form $j 2^{-k}$, where $j$ and $k$ are positive
integers.  If we restrict ourselves to the interval $[0,1)$, then the
binary sequence with all $1$'s does not correspond to a point in
$[0,1)$, but this does not cause serious trouble.  In particular,
these various exceptions do not cause trouble for the integrals that
we consider.

	If one accounts for these exceptions in a suitable (and
simple) way, then dyadic intervals in $[0,1)$ correspond in a nice
manner to subsets of $\mathcal{B}$ defined by prescribing the values
of a binary sequence $\{x_j\}_{j=1}^\infty$ for $j = 1, \ldots, n$ for
some $n \ge 1$, and leaving the other $x_j$'s free.  Similarly, dyadic
step functions on $[0,1)$ correspond to functions on $\mathcal{B}$
which depend only on the first $n$ terms of a binary sequence for some
positive integer $n$ (with suitable allowances at exceptional points,
as above).

\chapter{Convexity and some basic inequalities}
\label{Convexity and some basic inequalities}

\section{Convex functions}
\label{Convex functions}
\index{convex functions}

	Let $I$ be a subset of the real line which is either an open
interval, or an open half-line, or the whole real line ${\bf R}$
itself.  We shall always assume that $I$ is of this type in this
section.

	A real-valued function $\phi(x)$ on $I$ is said to be
\emph{convex} if
\begin{equation}
\label{convexity inequality for phi}
	\phi(\lambda \, x + (1-\lambda) \, y)
		\le \lambda \, \phi(x) + (1-\lambda) \, \phi(y)
\end{equation}
for all $x, y \in I$ and $\lambda \in [0,1]$.  (Note that $\lambda x +
(1-\lambda) \, y \in I$ when $x, y \in I$, so that $\phi(\lambda \, x
+ (1-\lambda) \, y)$ is defined.)

	If $\phi(x)$ is affine, $\phi(x) = a \, x + b$ for some real numbers
$a$ and $b$, then $\phi(x)$ is convex on the whole real line, and in fact
one has equality in (\ref{convexity inequality for phi}) for all $x$, $y$,
and $\lambda$.  This equality for all $x$, $y$, and $\lambda$ characterizes
affine functions, as one can easily check.

	For $\phi(x) = |x|$ we have that
\begin{equation}
	|x+y| \le |x| + |y| \quad\hbox{and}\quad 
				|\lambda \, x| = |\lambda| \, |x|
\end{equation}
for all $x, y, \lambda \in {\bf R}$, and the convexity condition
(\ref{convexity inequality for phi}) follows easily from this.

	If $\phi(x)$ is an arbitrary convex function on $I$, as above,
and if $a$ is a real number, then the translation $\phi(x-c)$ of
$\phi(x)$ is a convex function on 
\begin{equation}
	I + c = \{x + c : x \in I\}.
\end{equation}
In particular, for each real number $c$, $|x-c|$ defines a convex function
on ${\bf R}$.  

\beginlemma
\label{convexity and increasing difference quotients}
A function $\phi(x)$ on $I$ is convex if and only if the following
is true: for any points $s, t, u \in I$ such that $s < t < u$,
\begin{equation}
\label{increasing difference quotients inequalities}
	\frac{\phi(t) - \phi(s)}{t-s} \le \frac{\phi(u) - \phi(s)}{u-s}
					\le \frac{\phi(u) - \phi(t)}{u-t}.
\end{equation}
\end{lemma}

	Roughly speaking, this condition says that the difference quotients
of $\phi$ increase, or remain the same, as the points used in the difference
quotient become larger in ${\bf R}$.

	If $s$, $t$, $u$ are as in the lemma,
then
\begin{equation}
\label{t = frac{t-s}{u-s} u + frac{u-t}{u-s} s}
	t = \frac{t-s}{u-s} \, u + \frac{u-t}{u-s} \, s.
\end{equation}
This is easy to check, and we have that
\begin{equation}
	\frac{t-s}{u-s} \in (0,1), \quad \frac{u-t}{u-s} = 1 - \frac{t-s}{u-s}.
\end{equation}
Thus, if $\phi(x)$ is convex, then
\begin{equation}
\label{convexity inequality rewritten}
	\phi(t) \le \frac{t-s}{u-s} \, \phi(u) + \frac{u-t}{u-s} \, \phi(s).
\end{equation}
One can rewrite this in two different ways to get (\ref{increasing
difference quotients inequalities}).  To get the converse, one can
work backwards.  That is, either one of the inequalities in
(\ref{increasing difference quotients inequalities}) can be rewritten
to give (\ref{convexity inequality rewritten}), and this gives
(\ref{convexity inequality for phi}) when $s$, $t$, and $u$ are
obtained in the natural way from $x$, $y$, and $\lambda$, following
(\ref{t = frac{t-s}{u-s} u + frac{u-t}{u-s} s}).

\beginlemma
\label{convexity and comparison with affine functions}
A function $\phi(x)$ on $I$ is convex if and only if the following
condition holds: for each $t \in I$ there is a real-valued affine
function $A(x)$ on ${\bf R}$ such that $A(t) = \phi(t)$ and $A(x) \le
\phi(x)$ for all $x \in I$.
\end{lemma}

	To see that this condition is sufficient for $\phi$ to be
convex, one can apply the hypothesis with $t = \lambda \, x + (1 -
\lambda) \, y$ (given $x$, $y$, and $\lambda$ in the usual manner) to
get an affine function $A(x)$ as above, and then observe that
\begin{eqnarray}
\label{using A to get convexity inequality}
	\phi(\lambda \, x + (1 - \lambda) \, y)
		& = & A(\lambda \, x + (1 - \lambda) \, y)		\\
		& = &  \lambda \, A(x) + (1 - \lambda) \, A(y)	
								\nonumber \\
		& \le &  \lambda \, \phi(x) + (1 - \lambda) \, \phi(y).
								\nonumber
\end{eqnarray}

	Conversely, suppose that $\phi(x)$ is convex, and let $t \in I$
be given.  We want to choose a real number $a$ so that $A(x) = \phi(t)
+ a \, (x-t)$ satisfies $A(x) \le \phi(x)$ for all $x \in I$.  In other
words, the latter condition asks that
\begin{equation}
	a \, (x-t) \le \phi(x) - \phi(t)
\end{equation}
for all $x \in I$.  This is automatic when $x = t$, and we can rewrite
the inequality as
\begin{equation}
\label{a le frac{phi(x) - phi(t)}{x-t}}
	a \le \frac{\phi(x) - \phi(t)}{x-t}
\end{equation}
when $x > t$, and as
\begin{equation}
\label{frac{phi(t) - phi(x)}{t-x} le a}
	\frac{\phi(t) - \phi(x)}{t-x} \le a
\end{equation}
when $x < t$.

	From Lemma \ref{convexity and increasing difference
quotients} we have that
\begin{equation}
	\frac{\phi(t) - \phi(s)}{t - s} 
			\le \frac{\phi(u) - \phi(t)}{u - t}
\end{equation}
for all $s, u \in I$ such that $s < t < u$.  This implies that
\begin{equation}
	D_l = \sup \biggl\{\frac{\phi(t) - \phi(s)}{t-s} 
						: s \in I, s < t \biggr\}
\end{equation}
and
\begin{equation}
	D_r = \inf \biggl\{\frac{\phi(u) - \phi(t)}{u-t}
						: u \in I, u > t \biggr\}
\end{equation}
are well-defined (i.e., that the set of numbers of which the supremum is
taken is bounded from above, and that the set of numbers of which the
infimum is taken is bounded from below), and that
\begin{equation}
	D_l \le D_r.
\end{equation}
To get (\ref{a le frac{phi(x) - phi(t)}{x-t}}) and (\ref{frac{phi(t) -
phi(x)}{t-x} le a}), we want to choose $a \in {\bf R}$ so that
\begin{equation}
	D_l \le a \le D_r,
\end{equation}
and this we can do.  This completes the proof of Lemma \ref{convexity
and comparison with affine functions}.

\section{Jensen's inequality}
\label{Jensen's inequality}
\index{Jensen's inequality}

	Let $I$ be as in the previous section, and suppose that
$\phi(x)$ is a convex function on $I$.  Also let $J$ be an interval in
${\bf R}$ (with nonzero length) and let $f(x)$ be a function on $J$
such that $f(x) \in I$ for all $x \in J$.  Then $|J|^{-1} \int_J f(x)
\, dx \in I$, and
\begin{equation}
\label{phi(|J|^{-1} int_J f(x) dx) le |J|^{-1} int_J phi(f(x)) dx}
	\phi\Bigl(|J|^{-1} \int_J f(x) \, dx \Bigr)
		\le |J|^{-1} \int_J \phi(f(x)) \, dx.
\end{equation}
This is known as \emph{Jensen's inequality}.

	Let us first consider a version of this for sums.  Suppose
that $x_1, x_2, \ldots, x_n$ are elements of $I$, and that $\lambda_1,
\lambda_2, \ldots, \lambda_n$ are nonnegative real numbers such that
\begin{equation}
	\sum_{i=1}^n \lambda_i = 1.
\end{equation}
Then $\sum_{i=1}^n \lambda_i \, x_i \in I$, and
\begin{equation}
\label{generalized convexity inequality}
	\phi\Bigl(\sum_{i=1}^n \lambda_i \, x_i \Bigr)
		\le \sum_{i=1}^n \lambda_i \, \phi(x_i).
\end{equation}
This is equivalent to (\ref{phi(|J|^{-1} int_J f(x) dx) le |J|^{-1}
int_J phi(f(x)) dx}) for step functions $f(x)$.  (For more general
functions, one should be a bit more careful.)

	The definition of convexity (\ref{convexity inequality for
phi}) is the same as (\ref{generalized convexity inequality}) when $n
= 2$, and one can use it repeatedly to get the general case (via
induction).  Alternatively, one can use the characterization of
convexity in Lemma \ref{convexity and comparison with affine
functions}, and make an argument very similar to the one used in
(\ref{using A to get convexity inequality}) in order to get
(\ref{generalized convexity inequality}).

	Let us note that
\begin{eqnarray}
\label{phi(t) = t^p is convex when t ge 1}
	&& \hbox{$\phi(t) = t^p$ is a convex function on $[0,\infty)$}	\\
	&& \hbox{when $p$ is a real number such that $p \ge 1$.}	
								\nonumber
\end{eqnarray}
This is well known.  Although our current convention is to use open
intervals and half-lines for the domains of convex functions, it is
easy to allow for the endpoint $0$ here, with properties like those
that have been discussed.  In particular, we have that
\begin{equation}
\label{(|J|^{-1} int_J f(x) dx)^p le |J|^{-1} int_J f(x)^p dx}
	\Bigl(|J|^{-1} \int_J f(x) \, dx\Bigr)^p
		\le |J|^{-1} \int_J f(x)^p \, dx
\end{equation}
for nonnegative functions $f(x)$ on an interval $J$.

\section{H\"older's inequality}
\label{section on Holder's inequality}

	Let $p$, $q$ be real numbers such that $p, q \ge 1$ and
\begin{equation}
\label{1/p + 1/q = 1}
	\frac{1}{p} + \frac{1}{q} = 1.
\end{equation}
In this case we say that $p$ and $q$ are \emph{conjugate exponents}.
\index{conjugate exponents}

	Fix an interval $J$ in ${\bf R}$.  If $f(x)$ and $g(x)$ are
nonnegative functions on $J$, then \emph{H\"older's
inequality}\index{H\"older's inequality} states that
\begin{equation}
\label{Holder's inequality written out}
	\int_J f(x) \, g(x) \, dx
		\le \Bigl(\int_J f(y)^p \, dy\Bigr)^{1/p}
			\Bigl(\int_J g(z)^q \, dz\Bigr)^{1/q}.
\end{equation}
(As usual, one should make some modest assumptions so that the
integrals make sense, and by our conventions one is free to restrict
one's attention to step functions.)

	We can allow one of $p$ and $q$ to be $1$ and the other to be
$\infty$ (which is consistent with (\ref{1/p + 1/q = 1})).  If $p = 1$
and $q = \infty$, then the analogue of (\ref{Holder's inequality
written out}) is that
\begin{equation}
\label{Holder's inequality, p=1, q=infty}
	\int_J f(x) \, g(x) \, dx
		\le \Bigl(\int_J f(y) \, dy \Bigr)
			\Bigl(\sup_{z \in J} g(z) \Bigr)
\end{equation}
(which is very simple).

	Let us prove H\"older's inequality in the case where $p, q >
1$.  We begin with some preliminary observations.  The inequality is
trivial if $f$ or $g$ is identically $0$, since the left side of
(\ref{Holder's inequality written out}) is then $0$.  Thus we assume
that neither is identically $0$, so that
\begin{equation}
\label{(int_J f(y)^p dy)^{1/p}}
	\Bigl(\int_J f(y)^p \, dy\Bigr)^{1/p}
\end{equation}
and
\begin{equation}
\label{(int_J g(z)^q dz)^{1/q}}
	\Bigl(\int_J g(z)^q \, dz\Bigr)^{1/q}
\end{equation}
are both nonzero.

	Next, we may assume that (\ref{(int_J f(y)^p dy)^{1/p}}) and
(\ref{(int_J g(z)^q dz)^{1/q}}) are both equal to $1$.  In other
words, if we can prove (\ref{Holder's inequality written out}) in this
special case, then it follows in general, by multiplying $f$ and $g$
by constants.

	Our basic starting point is that
\begin{equation}
	s \, t \le \frac{s^p}{p} + \frac{t^q}{q}
\end{equation}
for all nonnegative real numbers $s$ and $t$.  This is a version of
the geometric-arithmetic mean inequalities, and it is the same as the
convexity of the exponential function.

	By integrating we are lead to
\begin{equation}
	\int_J f(x) \, g(x) \, dx
		\le \frac{1}{p} \int_J f(x)^p \, dx
			+ \frac{1}{q} \int_J g(x)^q \, dx.
\end{equation}
This gives (\ref{Holder's inequality written out}) when (\ref{(int_J
f(y)^p dy)^{1/p}}) and (\ref{(int_J g(z)^q dz)^{1/q}}) are equal to
$1$, which is what we wanted.

	Let us note the version of H\"older's inequality for sums.
Namely,
\begin{equation}
\label{Holder's inequality for sums}
\index{H\"older's inequality}
	\sum_i a_i \, b_i \le \Bigl(\sum_j a_j^p\Bigr)^{1/p}
				\Bigl(\sum_k b_k^q\Bigr)^{1/q}
\end{equation}
holds for arbitrary nonnegative real numbers $a_i$, $b_i$ when $p$ and
$q$ conjugate exponents.  If $p = 1$ and $q = \infty$, then this 
should be interpreted as
\begin{equation}
	\sum_i a_i b_i \le \Bigl(\sum_j a_j\Bigr) \Bigl(\sup_k b_k\Bigr),
\end{equation}
which is very simple, as before.  One can prove (\ref{Holder's
inequality for sums}) in the same manner as before, or derive it
from the previous version.

	The $p = q = 2$ case of (\ref{Holder's inequality for sums})
is the Cauchy-Schwarz inequality.\index{Cauchy-Schwarz inequality} The
$p = q = 2$ case of (\ref{Holder's inequality written out}) is also
called a Cauchy-Schwarz inequality.

	Observe that (\ref{(|J|^{-1} int_J f(x) dx)^p le |J|^{-1}
int_J f(x)^p dx}) can be derived from H\"older's inequality
(\ref{Holder's inequality written out}) by taking $g \equiv 1$.

\section{Minkowski's inequality}
\label{Minkowski's inequality}

	Fix an interval $J$ in ${\bf R}$.  Let $f$ and $g$ be
nonnegative functions on $J$, and let $p$ be a real number, $p \ge
1$.  Then
\begin{eqnarray}
\label{Minkowski's inequality, written out}
\index{Minkowski's inequality}\index{triangle inequality}
\lefteqn{\Bigl(\int_J (f(x) + g(x))^p \, dx \Bigr)^{1/p}}		\\
	&& \le \Bigl(\int_J f(x)^p \, dx \Bigr)^{1/p}
		+ \Bigl(\int_J g(x)^p \, dx \Bigr)^{1/p}.	
								\nonumber
\end{eqnarray}
This is \emph{Minkowski's inequality}.

	The analogue of this for $p = \infty$ is
\begin{equation}
\index{triangle inequality}
	\sup_{x \in J} (f(x) + g(x))
		\le \sup_{x \in J} f(x) + \sup_{x \in J} g(x).
\end{equation}
This inequality is easy to verify.

	The $p=1$ case of (\ref{Minkowski's inequality, written out})
is clear (and with equality), and so we focus now on $1 < p < \infty$.
We shall describe two arguments, the first employing H\"older's
inequality.

	One begins by writing
\begin{equation}
	(f(x) + g(x))^p = f(x) \, (f(x) + g(x))^{p-1} 
				   + g(x) \, (f(x) + g(x))^p.
\end{equation}
Thus
\begin{eqnarray}
\lefteqn{\quad\int_J (f(x) + g(x))^p \, dx}	\\
	&& = \int_J f(x) \, (f(x) + g(x))^{p-1} \, dx
			+ \int_J g(x) \, (f(x) + g(x))^{p-1} \, dx.
								\nonumber
\end{eqnarray}
If $q > 1$ is the conjugate exponent of $p$, $1/p + 1/q = 1$, then
H\"older's inequality implies that
\begin{eqnarray}
\lefteqn{\int_J f(x) \, (f(x) + g(x))^{p-1} \, dx}	\\
	&& \le \Bigl(\int_J f(y)^p \, dy \Bigr)^{1/p}
	        \Bigl(\int_J (f(z) + g(z))^{q(p-1)} \, dz \Bigr)^{1/q}.
								\nonumber
\end{eqnarray}
Because of the choice of $q$, we have that $q (p-1) = p$, and hence
\begin{eqnarray}
\lefteqn{\int_J f(x) \, (f(x) + g(x))^{p-1} \, dx}	\\
	&& \le \Bigl(\int_J f(y)^p \, dy \Bigr)^{1/p}
	        \Bigl(\int_J (f(z) + g(z))^p \, dz \Bigr)^{1 - 1/p}.
								\nonumber
\end{eqnarray}
Similarly,
\begin{eqnarray}
\lefteqn{\int_J g(x) \, (f(x) + g(x))^{p-1} \, dx}	\\
	&& \le \Bigl(\int_J g(y)^p \, dy \Bigr)^{1/p}
	         \Bigl(\int_J (f(z) + g(z))^p \, dz \Bigr)^{1 - 1/p},
								\nonumber
\end{eqnarray}
and therefore
\begin{eqnarray}
\lefteqn{\qquad\qquad\int_J (f(x) + g(x))^p \, dx}	\\
	&& \le \biggl\{\Bigl(\int_J f(y)^p \, dy \Bigr)^{1/p}
			+ \Bigl(\int_J g(y)^p \, dy \Bigr)^{1/p}\biggr\}
	        \Bigl(\int_J (f(z) + g(z))^p \, dz \Bigr)^{1 - 1/p}.
								\nonumber
\end{eqnarray}
It is easy to derive (\ref{Minkowski's inequality, written out}) from this.

	Before discussing the second argument, let us mention that the
version of Minkowski's inequality\index{Minkowski's inequality} for
sums is
\begin{equation}
\label{Minkowski's inequality for sums}\index{triangle inequality}
	\Bigl(\sum_i (a_i + b_i)^p \Bigr)^{1/p}
		\le \Bigl(\sum_i a_i^p \Bigr)^{1/p}
			+ \Bigl(\sum_i b_i^p \Bigr)^{1/p},
\end{equation}
where $a_i, b_i \ge 0$ for all $i$ and $p \ge 1$.  This can be proved
in the same manner as above.  (One can also derive the version for
sums from the version for integrals, as well as the other way around.)
The analogue for $p = \infty$ is
\begin{equation}
\index{triangle inequality}
	\sup_i (a_i + b_i) \le \sup_i a_i + \sup_i b_i.
\end{equation}

	Let us go through the second argument in the case of sums.
Fix $p$, $1 < p < \infty$, and assume for the moment that
\begin{equation}
\label{(sum_i a_i^p)^{1/p} = (sum_i b_i^p)^{1/p} = 1}
	\Bigl(\sum_i a_i^p \Bigr)^{1/p}
			= \Bigl(\sum_i b_i^p \Bigr)^{1/p} = 1.
\end{equation}
Suppose that $t$ is a real number such that $0 \le t \le 1$.
We would like to show that
\begin{equation}
\label{(sum_i (t a_i + (1-t) b_i)^p)^{1/p} le 1}
	\Bigl(\sum_i (t \, a_i + (1-t) \, b_i)^p \Bigr)^{1/p} \le 1.
\end{equation}
This would follow from Minkowski's inequality, and it is not hard to
verify that one can derive Minkowski's inequality from this version.

	We can rewrite (\ref{(sum_i a_i^p)^{1/p} = (sum_i b_i^p)^{1/p}
= 1}) and (\ref{(sum_i (t a_i + (1-t) b_i)^p)^{1/p} le 1}) as
\begin{equation}
\label{sum_i a_i^p = sum_i b_i^p = 1}
	\sum_i a_i^p = \sum_i b_i^p = 1
\end{equation}
and
\begin{equation}
\label{sum_i (t a_i + (1-t) b_i)^p le 1}
	\sum_i (t \, a_i + (1-t) \, b_i)^p \le 1.
\end{equation}
To go from (\ref{sum_i a_i^p = sum_i b_i^p = 1}) to (\ref{sum_i (t a_i
+ (1-t) b_i)^p le 1}) it suffices to know that
\begin{equation}
	(t \, a_i + (1-t) \, b_i)^p \le t \, a_i^p + (1-t) \, b_i^p
\end{equation}
for each $i$.  This inequality follows from the convexity of the
function $t^p$, as in (\ref{phi(t) = t^p is convex when t ge 1}) in
Section \ref{Jensen's inequality}.

\section{$p < 1$}
\label{p < 1}

	If $0 < p < 1$, then (\ref{Minkowski's inequality, written out})
and (\ref{Minkowski's inequality for sums}) no longer hold in general.
One has the alternatives
\begin{equation}
	\int_J (f(x) + g(x))^p \, dx 
	  \le \int_J f(x)^p \, dx + \int_J g(x)^p \, dx 
\end{equation}
and
\begin{equation}
	\sum_i (a_i + b_i)^p \le \sum_i a_i^p + \sum_i b_i^p.
\end{equation}
These alternatives do not hold in general when $p > 1$ (and become
equations when $p = 1$).

	The underlying building block for these inequalities is the fact
that
\begin{equation}
\label{(c + d)^p le c^p + d^p}
	(c + d)^p \le c^p + d^p
\end{equation}
when $c$ and $d$ are nonnegative real numbers, and $p$ is a real number
such that $0 < p \le 1$.  Once one has (\ref{(c + d)^p le c^p + d^p}),
the previous inequalities follow simply by integrating or summing.

\chapter{Normed vector spaces}

	In this book, all vector spaces use either the real or complex
numbers as their underlying scalar field.  We may sometimes wish to
restrict ourselves to one or the other, but in many cases either is
fine.  Throughout this chapter, we make the standing assumption that
all vector spaces are finite-dimensional.  Given a real or complex
number $t$, we let $|t|$ denote its absolute value or modulus.

\section{Definitions and basic properties}
\label{definitions, etc. (normed vector spaces)}

	Let $V$ be a vector space (real or complex).  By a
\emph{norm}\index{norm (on a vector space)} we mean a nonnegative
real-valued function $\| \cdot \|$ on $V$ which satisfies the
following properties: $\|v\| = 0$ if and only if $v$ is the zero
vector in $V$; $\| t \, v \| = |t| \|v\|$ for all vectors $v \in V$
and all scalars $t$; and $\|v+w\| \le \|v\| + \|w\|$ for all $v, w \in
V$.  The last property is called the \emph{triangle
inequality}\index{triangle inequality} for $\|\cdot \|$.

	A vector space equipped with a choice of norm is called a
\emph{normed vector space}.\index{normed vector space}

	If $B_1 = \{v \in V : \|v\| \le 1\}$ is the (closed) unit
ball\index{unit ball (in a normed vector space)}\index{closed unit
ball (in a normed vector space)} corresponding to the norm $\|\cdot
\|$, then $B_1$ is a convex\index{convex subset of a vector space}
subset of $V$.  In other words, if $v$, $w$ are elements of $B_1$ and
$t$ is a real number, $0 \le t \le 1$, then $t \, v + (1-t) \, w \in
B_1$.  This is easy to derive from the homogeneity property and
triangle inequality for $\|\cdot \|$.  Conversely, if one assumes the
homogeneity property for $\|\cdot \|$ and that the unit ball $B_1$ is
convex, then it is easy to show that the triangle inequality holds for
$\|\cdot \|$.

	As a basic family of examples, take $V$ to be ${\bf R}^n$ or
${\bf C}^n$, and, for a given $1 \le p \le \infty$, set
\begin{equation}
\label{def of ||v||_p}
	\|v\|_p = \Bigl(\sum_{j=1}^n |v_j|^p\Bigr)^{1/p}
\end{equation}
when $p < \infty$, and
\begin{equation}
	\|v\|_{\infty} = \max_{1 \le j \le n} |v_j|.
\end{equation}
Here $v_j$, $1 \le j \le n$, denote the components of $v$ in ${\bf
R}^n$ or ${\bf C}^n$.  That these are indeed norms is easy to check,
using Section \ref{Minkowski's inequality} for the triangle inequality.

	Let $V$ be a vector space, and let $\|\cdot\|$ be a norm on
$V$.  Notice that
\begin{equation}
\label{| ||v|| - ||w|| | le ||v-w||}
	\Bigl|\|v\| - \|w\| \Bigr| \le \|v-w\|
\end{equation}
for all $v, w \in V$.  This follows from
\begin{equation}
\label{||v|| le ||w|| + ||v-w||}
	\|v\| \le \|w\| + \|v-w\|
\end{equation}
and the analogous inequality with $v$ and $w$ interchanged, which are
instances of the triangle inequality.

	Suppose that $V = {\bf R}^n$ or ${\bf C}^n$, which is not a
real restriction, since every vector space is isomorphic to one of
these.  Let $|x|$ denote the standard Euclidean norm on ${\bf R}^n$ or
${\bf C}^n$ (which is the same as the norm $\|x\|_2$ in (\ref{def of
||v||_p})).  Then there is a positive constant $C$, depending on $n$
and the given norm $\|\cdot \|$ such that
\begin{equation}
\label{||v|| le C |v|}
	\|v\| \le C |v|
\end{equation}
for all $v \in V$.  This is not hard to check.

	From (\ref{| ||v|| - ||w|| | le ||v-w||}) and (\ref{||v|| le C
|v|}) we obtain in particular that $\|v\|$ is a continuous real-valued
function on $V$, with respect to the usual Euclidean metric and
topology.  As a consequence, there is a positive real number $C'$ such
that
\begin{equation}
\label{|v| le C' ||v||}
	|v| \le C' \|v\|
\end{equation}
for all $v \in V$.  More precisely, because of homogeneity, it
suffices to check this for $v$ in the standard unit sphere, $|v| = 1$.
We want to show that $\|v\|$ is bounded from below by a positive
constant on this set.  The definition of a norm ensures that $\|v\| >
0$ for all $v \ne 0$, and the compactness of the unit sphere $\{v \in
V : |v| = 1\}$ and the continuity of $\|v\|$ then imply that $\|v\|$
is bounded from below by a positive real number, as desired.

\section{Dual spaces and norms}
\label{Dual spaces and norms}

	Let $V$ be a vector space, real or complex.  By a \emph{linear
functional}\index{linear functional} on $V$ we mean a linear mapping
from $V$ into the field of scalars (the real or complex numbers, as
appropriate).  The \emph{dual of $V$}\index{dual vector space} is
the vector space of linear functionals on $V$, a vector space over the
same field of scalars.  The dual space of $V$ is denoted
$V^*$,\index{$V^*$} and it has the same dimension as $V$.

	Now suppose that $\|\cdot \|$ is a norm on $V$.  The
\emph{dual norm}\index{dual norm} $\|\cdot \|^*$\index{$"\"|\cdot "\"|^*$}
on $V^*$ is defined as follows.  If $\lambda$ is an element of $V^*$,
and thus a linear functional on $V$, then
\begin{equation}
	\|\lambda\|^* = \sup \{|\lambda(v)| : v \in V, \|v\| = 1\}.
\end{equation}
The remarks near the end of the Section \ref{definitions, etc. (normed
vector spaces)} show that this supremum is finite.  This definition of
$\|\lambda\|^*$ is equivalent to saying that
\begin{equation}
	|\lambda(v)| \le \|\lambda\|^* \, \|v\|
\end{equation}
for all $v \in V$, and that $\|\lambda\|^*$ is the smallest
nonnegative real number with this property.  It is not hard to verify
that $\| \cdot \|^*$ does indeed define a norm on $V^*$.

	Let us look at this in the context of the examples mentioned
in Section \ref{definitions, etc. (normed vector spaces)}.  That is,
suppose that $V = {\bf R}^n$ or ${\bf C}^n$, and take $\|\cdot \|_p$
for the norm on $V$, for some $p$ which satisfies $1 \le p \le
\infty$.

	We can identify $V^*$ with ${\bf R}^n$ or ${\bf C}^n$
(respectively) by associating to each $w$ in ${\bf R}^n$ or ${\bf
C}^n$
the linear functional $\lambda_w$ on $V$ given by
\begin{equation}
	\lambda_w(v) = \sum_{j=1}^n w_j \, v_j.
\end{equation}
Let $q$, $1 \le q \le \infty$ be the exponent conjugate to $p$, so
that $1/p + 1/q = 1$.  If $\|\cdot \| = \|\cdot \|_p$, let us check
that $\|\lambda_w \|^* = \|w\|_q$ for all $w$.

	First, we have that
\begin{equation}
	|\lambda_w(v)| \le \|w\|_q \, \|v\|_p
\end{equation}
for all $w$ and $v$ by H\"older's inequality.\index{H\"older's
inequality} To show that $\|\lambda_w\|^* = \|w\|_q$, we would like to
check that for each $w$ there is a nonzero $v$ so that
\begin{equation}
\label{|lambda_w(v)| = ||w||_q ||v||_p}
	|\lambda_w(v)| = \|w\|_q \, \|v\|_p.
\end{equation}
Let $w$ be given.  We may as well assume that $w \ne 0$, since
otherwise any $v$ would do.  Let us also assume for the moment that $p
> 1$, so that $q < \infty$.  Under these conditions, define $v$ by
\begin{equation}
	v_j = \overline{w_j} \, |w_j|^{q-2}
\end{equation}
when $w_j \ne 0$, and by $v_j = 0$ when $w_j = 0$.  Here
$\overline{w_j}$ denotes the complex conjugate of $w_j$, which is not
needed when we are working with real numbers instead of complex
numbers.  With this choice of $v$, we have that $\lambda_w(v) = \sum
|w_j|^q = \|w\|_q^q$.  It remains to check that
\begin{equation}
\label{||v||_p = ||w||_q^{q-1}}
	\|v\|_p = \|w\|_q^{q-1}
\end{equation}
in this case.  If $q = 1$, then $p = \infty$, and this equation is the
same as saying that $\max |v_j| = 1$.  Indeed, if $q=1$, then $|v_j| =
1$ whenever $v_j \ne 0$, and this happens for at least one $j$ because
$w \ne 0$.  Thus $\max |v_j| = 1$.  If $q > 1$, then $|v_j| =
|w_j|^{q-1}$ for all $j$, and one can verify (\ref{||v||_p =
||w||_q^{q-1}}) using the fact that $p$ and $q$ are conjugate
exponents.

	Finally, if $p = 1$, so that $q = \infty$, then choose (a
single) $k$, $1 \le k \le n$, so that $|w_k| = \max |w_j| =
\|w\|_{\infty}$.  Define $v$ by $v_k = \overline{w_k} \, |w_k|^{-1}$
and $v_j = 0$ when $j \ne k$.  Then $\lambda_w(v) = |w_k| =
\|w\|_{\infty}$ and $\|v\|_1 = 1$, and (\ref{|lambda_w(v)| = ||w||_q
||v||_p}) holds in this situation as well.  This completes the proof
that $\|\lambda_w\|^* = \|w\|_q$.

\section{Second duals}
\label{Second duals}

	Let $V$ be a vector space, and $V^*$ its dual space.  Since
$V^*$ is a vector space in its own right, we can take the dual of it
to get the second dual $V^{**}$\index{$V^{**}$} of $V$.

	There is a canonical isomorphism from $V$ onto $V^{**}$, which
is defined as follows.  Let $v \in V$ be given.  For each $\lambda \in
V^*$, we get a scalar by taking $\lambda(v)$.  The mapping $\lambda
\mapsto \lambda(v)$ defines a linear functional on $V^*$, and hence an
element of $V^{**}$.  Since we can do this for every $v \in V$, we
get a mapping from $V$ into $V^{**}$ which one can check is linear and
an isomorphism.

	Now suppose that we also have a norm $\|\cdot \|$ on $V$.
This leads to a dual norm $\|\cdot \|^*$ on $V^*$, as in the previous
section.  By the same token, we get a double dual norm $\|\cdot
\|^{**}$ on $V^{**}$.  Using the canonical isomorphism between $V$ and
$V^{**}$ just described, we can think of $\| \cdot \|^{**}$ as
defining a norm on $V$.  Let us show that
\begin{equation}
\label{||v||^{**} = ||v|| for all v in V}
	\|v\|^{**} = \|v\|   \qquad\hbox{for all } v \in V.
\end{equation}
Note that this holds for the $p$-norms $\| \cdot \|_p$ on ${\bf R}^n$
and ${\bf C}^n$ by the analysis of their duals in the preceding
section.

	Let $v \in V$ be given.  To prove (\ref{||v||^{**} = ||v|| for
all v in V}) for this choice of $v$, it suffices to establish the
following two statements.  First,
\begin{equation}
	|\lambda(v)| \le \|\lambda\|^* \, \|v\| 
			\qquad\hbox{for all } \lambda \in V^*.
\end{equation}
Second, there is a nonzero $\lambda_0 \in V^*$ such that
\begin{equation}
\label{lambda_0(v) = ||lambda_0||^* ||v||}
	\lambda_0(v) = \|\lambda_0\|^* \, \|v\|.
\end{equation}
The first statement comes from the definition of $\|\lambda\|^*$, and
so we only need to prove the second one.  For this we may as well
suppose that $v \ne 0$.  We shall apply the next extension result.
(Alternatively, one could approach this as in Section \ref{Separation
of convex sets}.)

\begintheorem
\label{extension theorem}
Let $V$ be a vector space (real or complex), and let $\|\cdot \|$ be
a norm on $V$.  Suppose that $W$ is a vector subspace of $V$ and that
$\mu$ is a linear functional on $W$ such that
\begin{equation}
	|\mu(w)| \le \|w\|  \qquad\hbox{for all } w \in W.
\end{equation}
Then there is a linear functional $\widehat{\mu}$ on $V$ which is an
extension of $\mu$ and has norm less than or equal to $1$.
\end{theorem}

	The existence of a nonzero $\lambda_0 \in V^*$ satisfying
(\ref{lambda_0(v) = ||lambda_0||^* ||v||}) follows easily from this,
by first defining $\lambda_0$ on the span of $v$ so that $\lambda_0(v)
= \|v\|$, and then extending to a linear functional on $V$ with norm
$1$.

	To prove the theorem, we first assume that $V$ is a
\emph{real} vector space.  Afterwards we shall discuss the complex
case.

	Let $W$ and $\mu$ be given as in the theorem.  For each
nonnegative integer $j$ less than or equal to the codimension of $W$
in $V$, we would like to show that there is a vector subspace $W_j$ of
$V$ and a linear functional $\mu_j$ on $W_j$ such that $W \subseteq
W_j$, the codimension of $W_j$ in $V$ is equal to the codimension of
$W$ in $V$ minus $j$, $\mu_j = \mu$ on $W$, and
\begin{equation}
\label{|mu_j(w)| le ||w|| for all w in W_j}
	|\mu_j(w)| \le \|w\|  \qquad\hbox{for all } w \in W_j.
\end{equation}
If we can do this with $j$ equal to the codimension of $W$ in $V$,
then $W_j$ would be equal to $V$, and this would give a linear
functional on $V$ with the required properties.

	Let us show that we can do this for each $j$ less than or
equal to the codimension of $W$ in $V$ by induction.  For the base
case of $j = 0$ we simply take $W_0 = W$ and $\mu_0 = \mu$.

	Suppose that $j$ is a nonnegative integer strictly less than
the codimension of $W$ in $V$ such that $W_j$ and $\mu_j$ exist as
above.  We would like to choose $W_{j+1}$ and $\mu_{j+1}$ with the
analogous properties for $j+1$ instead of $j$.

	Under these conditions, the codimension of $W_j$ in $V$ is
positive, so that there is a vector $z$ in $V$ which does not lie in
$W_j$.  Fix any such $z$, and take $W_{j+1}$ to be the span of $W_j$
and $z$.  Thus the codimension of $W_{j+1}$ in $V$ is the codimension
of $W_j$ in $V$ minus $1$, and hence is equal to the codimension of
$W$ in $V$ minus $j+1$, as desired.

	We shall choose $\mu_{j+1}$ to be an extension of $\mu_j$ from
$W_j$ to $W_{j+1}$.  Let $\alpha$ be a real number, to be chosen later
in the argument.  If we set $\mu_{j+1}(z) = \alpha$, then $\mu_{j+1}$
is determined on all of $W_{j+1}$ by linearity and the condition that
$\mu_{j+1}$ be an extension of $\mu_j$.  Note that $W \subseteq
W_{j+1}$ and $\mu_{j+1} = \mu$ on $W$, because of the corresponding
statements for $W_j$ and $\mu_j$.

	The remaining point is to show that $\mu_{j+1}$ satisfies the
counterpart of (\ref{|mu_j(w)| le ||w|| for all w in W_j}) for $j+1$,
i.e., that
\begin{equation}
	|\mu_{j+1}(w)| \le \|w\|  
		\qquad\hbox{for all } w \in W_{j+1}.
\end{equation}
This is the same as saying that
\begin{equation}
\label{|mu_{j+1}(u + t z)| le ||u + t z||, u in W_j, t in {bf R}}
	|\mu_{j+1}(u + t z)| \le \|u + t z\|
		\qquad\hbox{for all } u \in W_j
		\hbox{ and } t \in {\bf R},
\end{equation}
since $W_{j+1}$ is the span of $W_j$ and $z$.  To establish this, it
is enough to show that
\begin{equation}
\label{|mu_{j+1}(u+z)| le ||u+z|| for all u in W_j}
	|\mu_{j+1}(u+z)| \le \|u+z\| 
			\qquad\hbox{for all } u \in W_j.
\end{equation}
Indeed, the case where $t = 0$ in the previous statement corresponds
exactly to our induction hypothesis (\ref{|mu_j(w)| le ||w|| for
all w in W_j}) for $\mu_j$ and $W_j$.  If $t \ne 0$, then one can
eliminate it using linearity of $\mu_{j+1}$, homogeneity of the
norm $\|\cdot \|$, and the fact that $W_j$ is invariant under scalar
multiplication (since it is a vector subspace).

	Let us rewrite (\ref{|mu_{j+1}(u+z)| le ||u+z|| for all u
in W_j}) as
\begin{equation}
\label{|mu_j(u) + alpha| le ||u + z|| for all u in W_j}
	|\mu_j(u) + \alpha| \le \|u + z\|
			   \qquad\hbox{for all } u \in W_j.
\end{equation}
In other words, we want to choose a $\alpha \in {\bf R}$ so that
(\ref{|mu_j(u) + alpha| le ||u + z|| for all u in W_j}) holds.  If
we can do this, then the rest works, and we get $\mu_{j+1}$ and
$W_{j+1}$ with the desired features.

	Because we are in the real case, $\mu_j(u)$ is a real
number for all $u \in W_j$, and (\ref{|mu_j(u) + alpha| le ||u +
z|| for all u in W_j}) is equivalent to
\begin{equation}
  - \mu_j(u) -\|u + z\| \le \alpha \le -\mu_j(u) + \|u + z\|
		\qquad\hbox{for all } u \in W_j.
\end{equation}
To show that there is an $\alpha \in {\bf R}$ which satisfies this 
property, it is enough to establish that
\begin{equation}
	\quad
   -\mu_j(u_1) - \|u_1 + z\| \le - \mu_j(u_2) + \|u_2 + z\|
		\qquad\hbox{for all } u_1, u_2 \in W_j.
\end{equation}
This condition is the same as 
\begin{equation}
	\mu_j(u_2 - u_1) \le \|u_2 + z\| + \|u_1 + z\|
		\qquad\hbox{for all } u_1, u_2 \in W_j.
\end{equation}
By the triangle inequality, it is enough to know that
\begin{equation}
	\mu_j(u_2 - u_1) \le \|u_2 - u_1\|
		\qquad\hbox{for all } u_1, u_2 \in W_j.
\end{equation}
This last is a consequence of (\ref{|mu_j(w)| le ||w|| for all w in
W_j}).  Thus we can choose $\alpha \in {\bf R}$ with the required
properties, and the induction argument is now finished.  This
completes the proof of Theorem \ref{extension theorem} when $V$
is a real vector space.

	To handle the case of complex vector spaces, let us make the
following remarks.  Given a complex linear functional on a complex
vector space, we can take its real part to get a linear functional on
the vector space now viewed as a real vector space, i.e., where we
forget about multiplication by $i$.  Conversely, given a real-valued
function $h$ on a complex vector space which is linear with respect to
real scalars (and vector addition), there is a complex linear
functional on the vector space whose real part is $h$.  Specifically,
one can take $h(w) - i h(iw)$ for this complex linear functional.  If
our complex vector space is equipped with a norm, then the norm of a
(complex) linear functional on the space is equal to the norm of the
real part of the functional, viewed as a linear functional on the real
version of the vector space (forgetting about $i$).  This is not hard
to verify.  The main point is that for each complex linear functional
$\lambda$ on the space and each vector $v$ in the space, there is a
complex number $\beta$ such that $|\beta| = 1$ and $\lambda(\beta \,
v) = \beta \, \lambda(v)$ is real.  Also, the norm of $\beta \, v$ is
equal to the norm of $v$, by the properties of a norm on a complex
vector space.

	Using these remarks, it is not hard to derive the version of
Theorem \ref{extension theorem} for complex vector spaces from the
version for real vector spaces.

\section{Linear transformations and norms}
\label{Linear transformations and norms}

	Suppose that $V_1$ and $V_2$ are vector spaces (both real or
both complex) equipped with norms $\|\cdot \|_1$ and $\|\cdot \|_2$
(where the subscripts are now simply labels for arbitrary norms,
rather than referring to the $p$-norms from Section \ref{definitions,
etc. (normed vector spaces)}).  If $T : V_1 \to V_2$ is a linear
transformation, then the \emph{operator norm}\index{operator norm}
$\|T\|_{op}$\index{$"\"|\cdot "\"|_{op}$} of $T$ (with respect to the
given norms on $V_1$ and $V_2$) is defined by
\begin{equation}
\label{defn of ||T||_{op}}
	\|T\|_{op} = \sup \{\|T(v)\|_2 : v \in V_1, \|v\|_1 = 1\}.
\end{equation}
This is equivalent to saying that
\begin{equation}
\label{bound from ||T||_{op}}
	\|T(v)\|_2 \le \|T\|_{op} \, \|v\|_1
		\qquad\hbox{for all } v \in V_1
\end{equation}
and that $\|T\|_{op}$ is the smallest nonnegative real number with
this property.

	The finiteness of $\|T\|_{op}$ is easy to establish using the
remarks near the end of Section \ref{definitions, etc. (normed vector
spaces)}.  Notice that the dual norm on the space of linear
functionals on a vector space is a special case of the operator norm,
in which one takes $V_2$ to be the $1$-dimensional vector space of
scalars, with the standard norm.

	The set of all linear transformations from $V_1$ to $V_2$ is a
vector space in a natural way, using addition and scalar
multiplication of linear transformations (with the same scalar field
as for $V_1$, $V_2$).  Again, the dual space of a vector space is a
special case of this.  It is not hard to check that the operator norm
$\| \cdot \|_{op}$ on the vector space of linear transformations from
$V_1$ to $V_2$ is indeed a norm in the sense described in Section
\ref{definitions, etc. (normed vector spaces)}.

	Suppose that $V_3$ is another vector space, with the same
field of scalars as for $V_1$ and $V_2$, and that $\|\cdot \|_3$ is a
norm on it.  Suppose that in addition to $T: V_1 \to V_2$ we have a
linear mapping $U : V_2 \to V_3$.  Then the composition $U \, T$ is a
linear transformation from $V_1$ to $V_3$, and
\begin{equation}
\label{op norm UT le op norm U x op norm T}
	\|U \, T\|_{op, 13} \le \|U\|_{op, 23} \, \|T\|_{op, 12},
\end{equation}
where the subscripts in the operator norms indicate which vector
spaces and norms on them are being used.  This inequality is easy to
check.

\section{Linear transformations and duals}

	Let $V_1$ and $V_2$ be vector spaces (both real or both
complex), and let $T : V_1 \to V_2$ be a linear mapping between them.
Associated to $T$ is a canonical linear transformation $T' : V_2^* \to
V_1^*$ \index{$T'$} called the \emph{transpose}\index{transpose (of a
linear transformation)} of $T$ or \emph{dual linear
transformation},\index{dual linear transformation} and it is defined
by
\begin{equation}
	T'(\mu) = \mu \circ T  \qquad\hbox{for all } \mu \in V_2^*.
\end{equation}
In other words, if $\mu$ is a linear functional on $V_2$, then $\mu
\circ T$ is a linear functional on $V_1$, and this linear functional
is $T'(\mu)$.

	Sometimes one might call $T'$ the adjoint of $T$, or denote it
$T^*$.\index{$T^*$} We prefer not to do that, to avoid confusions with
similar but distinct objects in the setting of inner product spaces.

	If $S : V_1 \to V_2$ is another linear transformation, and if
$a$, $b$ are scalars, then $a \, S + b \, T$ is a linear mapping from
$V_1$ to $V_2$ whose dual $(a \, S + b \, T)'$ is $a \, S' + b\, T'$.
If $V_3$ is another vector space with the same field of scalars as
$V_1$, $V_2$, and if $U : V_2 \to V_3$ is a linear mapping, then we
can consider the composition $U \, T : V_1 \to V_3$.  The dual $(U \,
T)' : V_3^* \to V_1^*$ of $U \, T$ is given by $T' \, U'$, as one can
easily check.

	The dual of the identity transformation on a vector space $V$
is the identity transformation on the dual space $V^*$.  A linear
mapping $T : V_1 \to V_2$ is invertible if and only if the dual
transformation $T' : V_2^* \to V_1^*$ is invertible.  

	Now suppose that $\|\cdot \|_1$ and $\|\cdot \|_2$ are norms
on $V_1$ and $V_2$.  Associated to these are the operator norm
$\|\cdot \|_{op}$ for linear mappings from $V_1$ to $V_2$, and the
dual norms $\|\cdot \|_1^*$ and $\|\cdot \|_2^*$ on the dual spaces
$V_1^*$ and $V_2^*$.  Furthermore, we can use the dual norms $\|\cdot
\|_2^*$, $\|\cdot \|_1^*$ to get an operator norm $\|\cdot \|_{op*}$
for linear mappings from $V_2^*$ to $V_1^*$.  (This should not be
confused with the dual of $\|\cdot \|_{op}$, viewed as a norm in its
own right.)  With this notation,
\begin{equation}
\label{||T||_{op} = ||T'||_{op*}}
	\|T\|_{op} = \|T'\|_{op*}
\end{equation}
for every linear mapping $T : V_1 \to V_2$.

	To see this, one can begin by noticing that $\|T'\|_{op*} \le
\|T\|_{op}$.  This is not hard to verify directly from the
definitions.  The opposite inequality can be derived using linear
functionals as in (\ref{lambda_0(v) = ||lambda_0||^* ||v||}).  Another
way to look at this is to apply the first inequality to $T'$ instead
of $T$, and with $V_2^*$ instead of $V_1$, $V_1^*$ instead of $V_2$,
etc.  We have already seen that $V_1^{**}$ and $V_2^{**}$ are
canonically isomorphic to $V_1$ and $V_2$, respectively, and it is
easy to check that $T''$ corresponds to $T$ under this isomorphism.
The second dual norms $\|\cdot \|_1^{**}$ and $\|\cdot \|_2^{**}$
correspond to $\|\cdot \|_1$ and $\|\cdot \|_2$, as in Section
\ref{Second duals}, and the operator norm $\|\cdot \|_{op**}$
associated to the second duals thus reduces to $\|\cdot \|_{op}$.  In
this way the first inequality applied to $T'$ instead of $T$ leads to
$\|T\|_{op} = \|T''\|_{op**} \le \|T'\|_{op*}$, and this together with
the first inequality (as it is, applied to $T$ directly) implies
(\ref{||T||_{op} = ||T'||_{op*}}).

\section{Inner product spaces}
\label{Inner product spaces}

	Let $V$ be a vector space.  An \emph{inner
product}\index{inner product} on $V$ is a scalar-valued function
$\langle \cdot, \cdot \rangle$\index{$\langle \cdot, \cdot \rangle$}
on $V \times V$ with the following properties: (1) for each $w \in V$,
the function $v \mapsto \langle v, w \rangle$ is a linear functional
on $V$; (2) if $V$ is a real vector space, then
\begin{equation}
\label{langle w, v rangle = langle v, w rangle for all v, w in V}
	\langle w, v \rangle = \langle v, w \rangle 
				  \qquad\hbox{for all } v, w \in V,
\end{equation}
and if $V$ is a complex vector space, then
\begin{equation}
\label{langle w, v rangle = overline{langle v, w rangle} for all v, w in V}
	\langle w, v \rangle = \overline{\langle v, w \rangle}
				  \qquad\hbox{for all } v, w \in V
\end{equation}
(where $\overline{a}$ denotes the complex conjugate of the complex
number $a$); (3) the inner product is positive definite, in the sense
that $\langle v, v \rangle$ is a positive real number (whether $V$ is
real or complex) for all $v \in V$ such that $v \ne 0$.  Note that
$\rangle v, w \langle = 0$ whenever either $v$ or $w$ is the zero
vector by the first two properties of an inner product.

	The first two properties of an inner product imply that for
each $v \in V$, the function $w \mapsto \langle v, w \rangle$ is
linear when $V$ is a real vector space, and is ``conjugate linear''
when $V$ is a complex vector space (which means that complex
conjugations are applied to scalars at appropriate moments).

	A vector space equipped with an inner product is called an
\emph{inner product space}.\index{inner product space} The inner
product $\langle \cdot, \cdot \rangle$ leads to a norm on the vector
space, by setting
\begin{equation}
\index{norm associated to an inner product}
	\|v\| = \langle v, v \rangle^{1/2}.
\end{equation}
The inner product and norm satisfy the \emph{Cauchy-Schwarz 
inequality},\index{Cauchy-Schwarz inequality} given by
\begin{equation}
\label{Cauchy-Schwarz inequality in inner product spaces}
	|\langle v, w \rangle| \le \|v\| \, \|w\|
\end{equation}
for all $v$ and $w$ in the vector space.  This is well known, and the
fact that $\|\cdot \|$ satisfies the triangle inequality is normally
derived from this.  

	For each positive integer $n$, the standard inner products
on ${\bf R}^n$ and ${\bf C}^n$ are given by
\begin{equation}
	\langle x, y \rangle = \sum_{j=1}^n x_j \, y_j
\end{equation}
on ${\bf R}^n$ and
\begin{equation}
	\langle v, w \rangle = \sum_{j=1}^n v_j \, \overline{w_j}
\end{equation}
on ${\bf C}^n$.  The associated norms are the standard Euclidean norms
on ${\bf R}^n$ and ${\bf C}^n$, and are the same as $\|\cdot \|_2$
in Section \ref{definitions, etc. (normed vector spaces)}.

	If $(V, \langle \cdot, \cdot \rangle)$ is an inner product
space, and if $v$, $w$ are two vectors in $V$, then $v$ and $w$
are said to be \emph{orthogonal}\index{orthogonal} if
\begin{equation}
	\langle v, w \rangle = 0.
\end{equation}
This condition is symmetric in $v$ and $w$.  When $v$ and $w$ are
orthogonal, we have that
\begin{equation}
\label{||v+w|| = (||v||^2 + ||w||^2)^{1/2}, v,w orthogonal}
	\|v+w\| = (\|v\|^2 + \|w\|^2)^{1/2}.
\end{equation}

	A collection of vectors is said to be orthogonal if any two
vectors in the collection are orthogonal.  A collection of vectors is
said to be \emph{orthonormal}\index{orthonormal} if it is an
orthogonal collection of vectors, and all of the vectors have norm
$1$.

	Any collection of nonzero orthogonal vectors is automatically
linearly independent (as one can easily verify).  A collection of
vectors in $V$ is said to be an \emph{orthogonal
basis}\index{orthogonal basis} in $V$ if it is orthogonal and a basis
in the usual sense.  Since nonzero orthogonal vectors are
automatically linearly independent, this amounts to saying that the
collection if orthogonal, and that the vectors in the collection are
nonzero and span $V$.  Similarly, an \emph{orthonormal
basis}\index{orthonormal basis} is an orthonormal collection of
vectors which is also a basis, which is equivalent to saying that it
is an orthonormal collection which spans $V$.

	In ${\bf R}^n$ and ${\bf C}^n$ one has the \emph{standard
bases}\index{standard bases (in ${\bf R}^n$ or ${\bf C}^n$)}
consisting of the $n$ vectors with $1$ component equal to $1$, and the
rest equal to $0$.  These bases are orthonormal with respect to the
standard inner products.

	Every inner product space admits an orthonormal basis.  This
famous result is often established by starting with any basis of the
vector space and converting it to one which is orthonormal by applying
the \emph{Gram-Schmidt process}.  As a consequence, one obtains that
every inner product space is isomorphic to ${\bf R}^n$ or ${\bf C}^n$
with the standard inner product, according to whether the vector space
is real or complex, where $n$ is the dimension of the initial vector
space.  

\beginremark
\label{polarization and the parallelogram law}
{\rm If $(V, \langle \cdot, \cdot \rangle)$ is an inner product space,
and if $\|\cdot\|$ is the norm associated to $\langle \cdot, \cdot
\rangle$, then one can give a simple formula for the inner product in
terms of the norm, through \emph{polarization}.\index{polarization} In
particular, the inner product is uniquely determined by the norm, and
a linear mapping on $V$ which preserves the norm also preserves the
inner product.

	For the norm by itself, one has the 
\emph{parallelogram law}\index{parallelogram law}
\begin{equation}
\label{the parallelogram law}
	\|v+w\|^2 + \|v-w\|^2 = 2 \, (\|v\|^2 + \|w\|^2)
		\qquad\hbox{for all } v, w \in V.
\end{equation}
Conversely, if $V$ is a vector space and $\|\cdot\|$ is a norm on $V$
which satisfies (\ref{the parallelogram law}), then there is an inner
product on $V$ such that $\|\cdot\|$ is the associated norm.  This is
a well-known fact.  To establish this, one can start by defining
$\langle v, w \rangle$ as a function using the formula for the inner
product in terms of the norm when the former exists.  The point is then
to show that this does indeed define an inner product if the norm
satisfies (\ref{the parallelogram law}).
}
\end{remark}

\section{Inner product spaces, continued}
\label{Inner product spaces, continued}

	Let $(V, \langle \cdot, \cdot \rangle)$ be an inner product
space.  For each $w \in V$, one can define a linear functional 
$L_w$ on $V$ by 
\begin{equation}
	L_w(v) = \langle v, w \rangle.
\end{equation}

	The mapping $w \mapsto L_w$ defines a mapping from $V$
into its dual space $V^*$.  This mapping is linear when $V$ is a
real vector space, and it is conjugate-linear when $V$ is a complex
vector space.  In either case, this mapping is one-to-one and sends
$V$ onto $V^*$, as one can check.

	Using the norm $\|\cdot \|$ on $V$ associated to the inner
product, we get a dual norm on $V^*$. The dual norm of $L_w$ is then
equal to $\|w\|$.  Indeed, that the dual norm is less than or equal to
$\|w\|$ follows from the Cauchy-Schwarz inequality
(\ref{Cauchy-Schwarz inequality in inner product spaces}).  To get
the opposite inequality, one can observe that $L_w(w) = \|w\|^2$.

	Now let $S$ be a vector subspace of $V$.  The \emph{orthogonal
complement}\index{orthogonal complements}
$S^\perp$\index{$S \perp@$S^\perp$ (orthogonal complement)} of $S$ is
defined by
\begin{equation}
	S^\perp = \{v \in V : \langle v, w \rangle = 0
				\hbox{ for all } w \in V\}.
\end{equation}
This is also a vector subspace of $V$.

	If $S$ is not all of $V$, then $S^\perp$ contains a nonzero 
vector.  Indeed, if $S$ is not all of $V$, then there is a linear
functional on $V$ which is equal to $0$ on $S$ but which is not equal
to $0$ on all of $V$.  This linear functional can be represented as
$L_w$ for some $w \in V$, and $w \ne 0$ since the functional is not
$0$ on all of $V$.  On the other hand, $w \in S^\perp$ because the
functional vanishes on $S$.

	The span of $S$ and $S^\perp$ is equal to all of $V$.  To see
this, we can apply the previous observation.  Namely, the span of $S$
and $S^\perp$ is automatically a subspace of $V$, and if it is not all
of $V$, then there is a nonzero vector $w$ in $V$ such that $w$ is
orthogonal to all elements of the span of $S$ and $S^\perp$.  In
particular, $w$ should be orthogonal to all elements of $S$, which
means that $w \in S^\perp$.  Since $w$ is also supposed to be
orthogonal to all elements of $S^\perp$, we conclude that $w$ is
orthogonal to itself.  This implies that $w = 0$, a contradiction.

	The statement that the span of $S$ and $S^\perp$ is equal to
$V$ can be reformulated as saying that every vector $v$ in $V$ can be
written as $u + z$, where $u \in S$ and $z \in S^\perp$.  This
decomposition of $v$ is unique, because $S \cap S^\perp = \{0\}$.

	We can express this decomposition in the form of a linear
mapping $P_S : V \to V$ such that $P_S(v) \in S$ and $v - P_S(v) \in
S^\perp$ for all $v \in V$.  This is called the \emph{orthogonal
projection of $V$ onto $S$}\index{orthogonal projection} from $V$ onto
$S$.  It is also characterized by the property that $P_S(v) = v$ for
all $v \in S$ and $P_S(v) = 0$ for all $v \in S^\perp$.

	Another important feature of the orthogonal projection is that
it has norm $1$\index{projections with operator norm $1$} in the
operator norm associated to the inner product norm on $V$ (unless $S$
is the trivial subspace consisting of only the zero vector $0$).  In
other words,
\begin{equation}
\label{||P_S(v)|| le ||v|| for all v in V}
	\|P_S(v)\| \le \|v\| \qquad\hbox{for all } v \in V,
\end{equation}
and equality holds when $v \in S$.  This inequality follows easily
from (\ref{||v+w|| = (||v||^2 + ||w||^2)^{1/2}, v,w orthogonal}) in
Section \ref{Inner product spaces}.  Conversely, a projection on an
inner product space is an orthogonal projection if it has norm $1$.
See Lemma \ref{a projection is orthogonal iff it has norm 1} in
Section \ref{Projections}.

\section{Separation of convex sets}
\label{Separation of convex sets}
\index{separation of convex sets}

	Let $E$ be a nonempty closed convex subset\index{convex subset
of a vector space} of ${\bf R}^n$ for some positive integer $n$, and
let $p$ be a point in ${\bf R}^n$ which does not lie in $E$.  A famous
\emph{separation theorem} states that there is a hyperplane $H$ in
${\bf R}^n$ such that $E$ lies on one side of $H$ and $p$ lies on the
other.

	This theorem can be proved in the following manner.  First,
let $q$ be an element of $E$ such that the distance from $p$ to $q$ in
the usual Euclidean norm is as small as possible.  The existence of
$q$ follows from standard considerations of continuity and
compactness.  Although $E$ may not be compact, because it may not be
bounded, only a bounded subset of $E$ is needed for this minimization
(namely, a part of $E$ which is not too far from $p$).

	Of course $p \ne q$, because $p$ does not lie in $E$.
Consider the hyperplane $H_0$ in ${\bf R}^n$ which passes through $q$
and which is orthogonal to $p-q$.  It is not too hard to show that $E$
lies in one of the \emph{closed} half-spaces bounded by $H_0$ (the one
that does not contain $p$).  (Exercise.)  To get a hyperplane $H$
which \emph{strictly} separates $p$ from $E$, one can slide $H_0$ over
towards $p$.

	Now consider a slightly different situation, where $p$ is an
element of $E$, but $p$ does not lie in the interior of $E$.  In other
words, $p$ is in the boundary of $E$.  In this case a related result
is that there is a hyperplane $H$ in ${\bf R}^n$ which passes through
$p$ and for which $E$ is contained in one of the closed half-spaces
bounded by $H$.

	If $E$ has smooth boundary, then $H$ should exactly be taken
to be the tangent hyperplane to $E$ at $p$.  In general, however, $E$
can have corners at the boundary.

	At any rate, one can deal with this case using the previous
assertion.  Specifically, because $p$ does not lie in the interior of
$E$, there are points $p_1$ in ${\bf R}^n$ which are not in $E$ and
which are as close as one likes to $p$.  This leads to hyperplanes
that pass near $p$ and for which $E$ is contained in one of the
corresponding closed half-spaces.  To get a hyperplane that passes
through $p$, one can take a limit of preliminary hyperplanes of this
type.  More precisely, one can start with a sequence of these
preliminary hyperplanes for which the associated points $p_1$ converge
to $p$, and then use the compactness of the unit sphere in ${\bf R}^n$
to pass to a subsequence of hyperplanes with unit normal vectors
converging to a unit vector in ${\bf R}^n$.  The hyperplane in ${\bf
R}^n$ that passes through $p$ and has this limiting vector as a normal
vector also has $E$ contained in one of the closed half-spaces of
which it is the boundary, as one can verify.

	These separation properties can be rephrased in terms of
linear functionals on ${\bf R}^n$.  Indeed, given any hyperplane $H$
in ${\bf R}^n$, there is a nonzero linear functional $\lambda$ on
${\bf R}^n$ and a real number $c$ such that
\begin{equation}
	H = \{x \in {\bf R}^n : \lambda(x) = c\}.
\end{equation}
The separation properties of $E$ and $p$ can be described in terms of
the values of $\lambda$ on $E$ and at $p$.

	One can also look at separations of pairs of convex sets,
rather than a convex set and a single point.  Specifically, let $E_1$
and $E_2$ be disjoint nonempty convex subsets of ${\bf R}^n$.  A
useful trick is to consider the set of differences $E_1 - E_2$,
\begin{equation}
	E_1 - E_2 = \{x-y : x \in E_1, y \in E_2\}.
\end{equation}
It is easy to see that $E_1 - E_2$ is also convex, and that $0$ does
not lie in $E_1 - E_2$, since $E_1$ and $E_2$ are disjoint.  However,
$E_1 - E_2$ is not closed in general, even if $E_1$ and $E_2$ are.
(This does hold if one of $E_1$, $E_2$ is compact and the other is
closed.)

	Let $F$ denote the closure of $E_1 - E_2$.  Thus $F$ is a
nonempty closed convex subset of ${\bf R}^n$.  Because $0$ does not
lie in $E_1 - E_2$, one can verify that $0$ does not lie in the
interior of $F$, i.e., it either does not lie in $F$ at all, or it
lies in the boundary of $F$.  Hence there is a hyperplane $H$ that
passes through $0$ such that $F$ is contained in one of the two closed
half-spaces bounded by $H$, as before.  This is equivalent to saying
that there is a linear functional $\lambda$ on ${\bf R}^n$ which is
nonnegative everywhere on $E_1 - E_2$.  As a consequence, the supremum
of $\lambda$ on $E_2$ is bounded by the infimum of $\lambda$ on $E_1$.

	The first result mentioned in this section has a simple and
well-known corollary, which is that a set $E$ in ${\bf R}^n$ is closed
and convex if and only if it can be realized as the intersection of a
family of closed half-spaces in ${\bf R}^n$.  The ``if'' part of this
assertion is automatic, since closed half-spaces are closed and
convex, and these properties are preserved by intersections.  For the
converse, if $E$ is a closed and convex subset of ${\bf R}^n$, then
one can compare $E$ to the intersection of all closed half-spaces in
${\bf R}^n$ that contain $E$.  This intersection contains $E$ by
definition, and it contains no other element of ${\bf R}^n$ when $E$
is closed and convex because of the separation theorem.

\section{Some variations}

	In this section, $V$ will be a finite-dimensional real vector
space.

\begindefinition [Sublinear functions]
\label{def of sublinear functions}
\index{sublinear functions} A 
real-valued function $p(v)$ on $V$ is said to be \emph{sublinear} if
$p(v+w) \le p(v) + p(w)$ for all $v, w \in V$ and $p(t \, v) = t \,
p(v)$ for all $v \in V$ and all nonnegative real numbers $t$.  In
particular, $p(0) = 0$.
\end{definition}

	Unlike a norm, a sublinear function $p(v)$ is not asked to
satisfy $p(-v) = p(v)$.

\begintheorem
\label{extension theorem, 2}
Let $V$ be a finite-dimensional real vector space, and let $p$ be a
sublinear function on $V$.  Suppose that $W$ is a vector subspace of
$V$ and that $\mu$ is a linear functional on $W$ such that
\begin{equation}
\label{mu(w) le p(w) for all w in W}
	\mu(w) \le p(w)  \qquad\hbox{for all } w \in W.
\end{equation}
Then there is a linear functional $\widehat{\mu}$ on $V$ which is an
extension of $\mu$ and satisfies
\begin{equation}
\label{widehat{mu}(v) le p(v) for all v in V}
	\widehat{\mu}(v) \le p(v) \qquad\hbox{for all } v \in V.
\end{equation}
\end{theorem}

	Note that we simply have $\mu(w)$ and $\widehat{\mu}(v)$ in
(\ref{mu(w) le p(w) for all w in W}) and (\ref{widehat{mu}(v) le p(v)
for all v in V}), rather than their absolute values.  If $p(-v) =
p(v)$ for all $v \in V$, then this would not make a real difference.

	Theorem \ref{extension theorem, 2} can be proved in a similar
manner as Theorem \ref{extension theorem} in Section \ref{Second
duals}.  The first part of the proof is essentially the same except
for simple changes following the differences in the statements of the
theorems.  This works up to (\ref{|mu_{j+1}(u + t z)| le ||u + t z||,
u in W_j, t in {bf R}}), whose counterpart in the present setting is
\begin{equation}
\label{mu_{j+1}(u + t z) le p(u + t z), u in W_j, t in {bf R}}
	\mu_{j+1}(u + t z) \le p(u + t z)
		\qquad\hbox{for all } u \in W_j
		\hbox{ and } t \in {\bf R}.
\end{equation}
To be more precise, for the proof we want to choose the real number
$\alpha$ (which $\mu_{j+1}(z)$ is assigned to be) so that
(\ref{mu_{j+1}(u + t z) le p(u + t z), u in W_j, t in {bf R}})
holds, and if we can do that, then the theorem will follow.

	For the next step we need a modification, which is to say
that (\ref{mu_{j+1}(u + t z) le p(u + t z), u in W_j, t in {bf R}})
can be reduced to
\begin{eqnarray}
\label{mu_{j+1}(u+z) le p(u+z), mu_{j+1}(u-z) le p(u-z) for all u in W_j}
	&& \mu_{j+1}(u+z) \le p(u+z)
		\hbox{ and } \mu_{j+1}(u-z) \le p(u-z)		\\
	&& \hbox{for all } u \in W_j.			\nonumber
\end{eqnarray}
As before, the $t = 0$ case of (\ref{mu_{j+1}(u + t z) le p(u + t z),
u in W_j, t in {bf R}}) comes from the induction hypothesis (for
$\mu_j$ and $W_j$), and for $t \ne 0$ one can convert (\ref{mu_{j+1}(u
+ t z) le p(u + t z), u in W_j, t in {bf R}}) into (\ref{mu_{j+1}(u+z)
le p(u+z), mu_{j+1}(u-z) le p(u-z) for all u in W_j}) by employing the
homogeneity of the sublinear function $p(v)$.  Now we have to treat
$t > 0$ and $t < 0$ separately, since we are only assuming homogeneity
for positive real numbers.

	We can rewrite (\ref{mu_{j+1}(u+z) le p(u+z), mu_{j+1}(u-z) le
p(u-z) for all u in W_j}) as
\begin{eqnarray}
\label{mu_j(u) + alpha le p(u+z), mu_j(u) - alpha le p(u-z) for all u in W_j}
	&& \mu_j(u) + \alpha \le p(u+z)
		\hbox{ and } \mu_j(u) - \alpha \le p(u-z)		\\
	&& \hbox{for all } u \in W_j,				\nonumber
\end{eqnarray}
using the definition of $\mu_{j+1}$, as in the previous situation.
The problem is to verify that $\alpha \in {\bf R}$ can be chosen so
that this condition holds.

	This requirement can be rephrased as 
\begin{equation}
\label{mu_j(u) - p(u-z) le alpha le p(u+z) - mu_j(u) for all u in W_j}
	\mu_j(u) - p(u-z) \le \alpha \le p(u+z) - \mu_j(u)
	\quad\hbox{for all } u \in W_j.
\end{equation}
The rest of the argument is nearly the same as before.  To show that
there is an $\alpha \in {\bf R}$ for which (\ref{mu_j(u) - p(u-z) le
alpha le p(u+z) - mu_j(u) for all u in W_j}) is satisfied, it suffices
to verify that
\begin{equation}
\label{mu_j(u_1) - p(u_1-z) le  p(u_2 +z) - mu_j(u_2) for all u_1, u_2 in W_j}
	\mu_j(u_1) - p(u_1 - z) \le p(u_2 + z) - \mu_j(u_2)
	\quad\hbox{for all } u_1, u_2 \in W_j.
\end{equation}
This is equivalent to saying that
\begin{equation}
	\mu_j(u_1 + u_2) \le p(u_1 - z) + p(u_2 + z) 
	\quad\hbox{for all } u_1, u_2 \in W_j.
\end{equation}
The subadditivity property of $p(v)$ shows that this condition holds
if
\begin{equation}
	\mu_j(u_1 + u_2) \le p(u_1 + u_2) 
	\quad\hbox{for all } u_1, u_2 \in W_j,
\end{equation}
which is the same as
\begin{equation}
	\mu_j(u) \le p(u) \qquad\hbox{for all } u \in W_j.
\end{equation}
This is exactly part of the ``induction hypothesis'' for $W_j$ and
$\mu_j$ in this setting, i.e., the counterpart of (\ref{mu(w) le p(w)
for all w in W}) for $W_j$ and $\mu_j$.  Thus (\ref{mu_j(u_1) -
p(u_1-z) le p(u_2 +z) - mu_j(u_2) for all u_1, u_2 in W_j}) is valid,
and it is possible to choose $\alpha \in {\bf R}$ with the desired
feature.  This completes the proof of Theorem \ref{extension theorem, 2}.

	Now let us turn to convex cones.  A subset $C$ of $V$ is
called a \emph{convex cone}\index{convex cone} if 
\begin{equation}
\label{x + y in C whenever x, y in C}
	x + y \in C \hbox{ whenever } x, y \in C
\end{equation}
and
\begin{equation}
\label{t x in C whenever x in C, t > 0}
	t \, x \in C \hbox{ whenever } x \in C, \, t > 0.
\end{equation}
In the second condition, $t$ is a real number.

	Clearly convex cones are convex sets.  Conversely, $C$ is a
convex cone if it is convex and satisfies (\ref{t x in C whenever x in
C, t > 0}).

	A subset $C$ of $V$ is called an \emph{open convex
cone},\index{open convex cone} or a \emph{closed convex
cone},\index{closed convex cone} if it is a convex cone and if it is
open or closed, respectively, as a subset of $V$.  For the latter we
use a linear isomorphism between $V$ and ${\bf R}^n$, where $n$ is
the dimension of $V$, to define open and closed sets in $V$.  Which
isomorphism ones employs does not matter, because invertible linear
mappings on ${\bf R}^n$ are homeomorphisms.

	For example, if $C$ is the set of $x = (x_1, \ldots, x_n)$ in
${\bf R}^n$ such that the coordinates $x_i$ of $x$ are all positive,
then $C$ is an open convex cone in ${\bf R}^n$.  If $C$ is the set of
$x \in {\bf R}^n$ such that the coordinates $x_i$ are all nonnegative,
then $C$ is a closed convex cone.

	Here is a simple recipe for producing convex cones.  If $C_0$
is a convex set in the vector space $V$, then consider the set
\begin{equation}
\label{{t x : t > 0, x in C_0}}
	\{t \, x : t > 0, \ x \in C_0\}.
\end{equation}
One can check that this is a convex cone in $V$, and in fact it is an
open convex cone if $C_0$ is an open convex subset of $V$.  If $C_0$
is a compact convex subset of $V$ which does not contain $0$, then
\begin{equation}
\label{{t x : t > 0, x in C_0} cup {0}}
	\{t \, x : t > 0, \ x \in C_0\} \cup \{0\}
\end{equation}
is a closed convex cone in $V$.

	If $P$ is an affine plane in $V$ that does not pass through
the origin, and $C_0$ is contained in $P$, then (\ref{{t x : t > 0, x
in C_0}}) and (\ref{{t x : t > 0, x in C_0} cup {0}}) define convex
cones whose intersection with $P$ is exactly $C_0$.

	To avoid degeneracies, one sometimes asks that a convex cone
$C$ in $V$ not be contained in a proper vector subspace of
$V$.  One also sometimes asks that a convex cone $C$ not
contain any lines.  Let us make the standing convention that
\begin{equation}
\label{if C is a closed convex cone in {bf R}^n, then 0 in C}
\hbox{if $C$ is a closed convex cone, then $0 \in C$}.
\end{equation}
In other words, a closed convex cone should not be the empty set.

	Now suppose that a norm $\|\cdot \|$ on $V$ has been chosen.
If $C$ is a nonempty convex cone in $V$, define $p(v)$ on $V$ by
\begin{equation}
	p(v) = \inf \{ \|v-x\| : x \in C \}.
\end{equation}
In other words, $p(v)$ is the distance from $v$ to $C$ relative
to the norm $\|\cdot \|$.  It is not hard to verify that $p(v)$ is
a sublinear function when $C$ is a nonempty convex cone.  Note that
the closure of $C$ is a closed convex cone which leads to the same
function $p(v)$.

	If $C$ is a convex cone in $V$, then we can associate to it a
subset $C^*$\index{$C^*$ (dual cone)} of the dual space $V^*$ of $V$,
namely, the set of linear functionals on $V$ which are nonnegative on
$C$.  Actually, this makes sense for any subset $C$ of $V$, but if $C$
is not already a convex cone, then one can take the convex cone that
it generates, and this convex cone will define the same set in $V^*$.
Similarly, the closure of $C$ leads to the same set in the dual, and
adding $0$ to $C$ if it is not already an element of $C$ does not
affect $C^*$.

	It is not hard to check that $C^*$ is always a closed convex
cone in $V^*$, containing $0$ in particular.  This is called the
\emph{dual cone}\index{dual cone} to $C$.

	If one takes $C$ to be the degenerate cone consisting of only
$0$, then the dual cone is all of $V^*$.  Similarly, if $C$ is all of
$V$, then the dual cone consists of only the zero linear functional.

	Suppose that $V = {\bf R}^n$, and for $C$ let us consider the
closed convex cone of points $x$ all of whose coordinates are
nonnegative.  A linear function $\lambda$ on ${\bf R}^n$ can be given
as
\begin{equation}
	\lambda(x) = \sum_{i=1}^n a_i \, x_i
\end{equation}
for some real numbers $a_1, \ldots, a_n$, and this linear functional
is nonnegative on $C$ exactly when $a_i \ge 0$ for each $i$.

	Let $C$ be a closed convex cone in $V$, and let $C^*$ be the
dual cone in $V^*$.  Associated to it is the second dual cone $C^{**}$
in $V^{**}$.  As in Section \ref{Second duals}, we can identify
$V^{**}$ with $V$ in a simple way, where each element of $V$ defines
a linear functional on $V^*$ by evaluation.  It is easy to see that
$C \subseteq C^{**}$.  In fact,
\begin{equation}
	C = C^{**}.
\end{equation}
To verify this, one wants to show that an element $v$ of $V$ does lie
in $C$ if every linear functional in $C^*$ is nonnegative at $v$.  In
other words, if $v \in V$ does not lie in $C$, then one would like to
say that there is a linear functional on $V$ which is nonnegative on
$C$, but negative at $v$.  This can be established as in Section
\ref{Separation of convex sets}, or using Theorem \ref{extension
theorem, 2}.  For the latter, one starts with a linear functional on
the $1$-dimensional subspace spanned by $v$, and one uses a sublinear
function $p$ associated to $C$ as above.  There is also a
multiplication by $-1$ involved, for switching between nonnegativity
and nonpositivity on $C$.

	See Chapter III of \cite{SW-book} for some very interesting
topics concerning convexity, duality, Fourier transforms, and
holomorphic functions on certain regions in ${\bf C}^n$.

\chapter{Strict convexity}

\section{Functions of one real variable}

	Let $I$ be a nonempty subset of the real line ${\bf R}$ which
is either an open interval, an open half-line, or the whole real line.
A real-valued function $\phi(x)$ on $I$ is said to be \emph{strictly
convex}\index{strictly convex (function of one real variable)} if
\begin{equation}
\label{strict convexity inequality}
	\phi(\lambda \, x + (1-\lambda) \, y) <
		\lambda \, \phi(x) + (1-\lambda) \, \phi(y)
\end{equation}
for all $x, y \in I$ such that $x \ne y$ and all $\lambda \in (0,1)$.
Thus strictly convex functions are convex in particular.

\beginlemma
\label{strict convexity and comparison with affine functions}
A real-valued function $\phi$ on $I$ is strictly convex if and only if
for every point $t \in I$ there is a real-valued affine function
$A(x)$ on ${\bf R}$ such that $A(t) = \phi(t)$ and $A(x) < \phi(x)$
for all $x \in I \backslash \{t\}$.
\end{lemma}

	This is the analogue of Lemma \ref{convexity and comparison
with affine functions} in Section \ref{Convex functions} for strictly
convex functions.

	To prove the ``if'' part of Lemma \ref{strict convexity and
comparison with affine functions}, let $x, y \in I$ and $\lambda \in
(0,1)$ be given, and set $t = \lambda \, x + (1-\lambda) \, y$.  Then
$t$ lies in $I$, and by hypothesis there is a real-valued affine
function $A(z)$ on ${\bf R}$ such that $A(t) = \phi(t)$ and $A(z) <
\phi(z)$ for all $z \in I$, $z \ne t$.  If $x \ne y$, then $x \ne t$
and $y \ne t$, so that $A(x) < \phi(x)$ and $A(y) < \phi(y)$.  Because
$A(z)$ is affine, we have that
\begin{equation}
	A(t) = \lambda \, A(x) + (1-\lambda) \, A(y),
\end{equation}
and hence
\begin{equation}
	\phi(t) = A(t) < \lambda \, \phi(x) + (1-\lambda) \, \phi(y).
\end{equation}
This is exactly the inequality that we want.

	Conversely, suppose that $\phi$ is a strictly convex function
on $I$, and let $t$ be an element of $I$.  By Lemma \ref{convexity and
comparison with affine functions} we know that there is a real-valued
affine function $A(x)$ on ${\bf R}$ such that $A(t) = \phi(t)$ and
$A(x) \le \phi(x)$ for all $x \in I$.  We want to show that $A(x) <
\phi(x)$ when $x \in I$, $x \ne t$.

	Suppose to the contrary that $\phi(x) = A(x)$ for some $x \in
I$, $x \ne t$.  For each $\lambda \in (0,1)$, we have that
\begin{equation}
	A(\lambda \, t + (1-\lambda) \, x) 
		\le \phi(\lambda \, t + (1-\lambda) \, x)
\end{equation}
by the conditions on $A(y)$, and
\begin{equation}
	A(\lambda \, t + (1-\lambda) \, x)
		= \lambda \, A(t) + (1-\lambda) \, A(x),
\end{equation}
since $A(y)$ is affine.  Because $A(t) = \phi(t)$ and $A(x) =
\phi(x)$, we obtain that
\begin{equation}
	\lambda \, \phi(t) + (1-\lambda) \, \phi(x)
		\le \phi(\lambda \, t + (1-\lambda) \, x).
\end{equation}
This contradicts the strict convexity of $\phi$, and the lemma
follows.

	As a basic family of examples (from ``calculus''), the
functions $t \mapsto t^p$ on $(0,\infty)$ are strictly convex when $p
> 1$.  One might as well say that these functions are strictly convex
on $[0,\infty)$, even if this is not an open half-line, because they
satisfy similar properties there (in terms of (\ref{strict convexity
inequality}) and comparison with affine functions).

\section{The unit ball in a normed vector space}

	Let $(V, \|\cdot \|)$ be a normed vector space (real or
complex).  As in Section \ref{definitions, etc. (normed vector
spaces)}, define the \emph{unit ball}\index{unit ball (in a normed
vector space)} in $V$ to be the set
\begin{equation}
	B_1 = \{v \in V : \|v\| \le 1\}.
\end{equation}
To be more precise, this is sometimes called the ``closed unit
ball'',\index{closed unit ball (in a normed vector space)} to
distinguish it from the open unit ball,\index{open unit ball (in a
normed vector space)} which is defined in the same manner except that
the inequality $\le$ is replaced with $<$.

	Let us say that the unit ball $B_1$ in a normed vector space
$(V, \|\cdot \|)$ is \emph{strictly convex}\index{strictly convex
(unit ball in a normed vector space)} if for every $x, y \in V$ such
that $\|x\| = \|y\| = 1$ and $x \ne y$ and every real number $\lambda
\in (0,1)$ we have that
\begin{equation}
\label{strict convexity of the unit ball}
	\|\lambda \, x + (1-\lambda) \, y\| < 1.
\end{equation}

	As basic examples, one can start with the fields of scalars
${\bf R}$ and ${\bf C}$, and their standard norms.  The unit balls for
these are strictly convex in this sense.  More generally, the unit
ball in any inner product space is strictly convex.  This is not
difficult to compute, using the standard analysis of the case of
equality for the Cauchy-Schwarz inequality.

	Now suppose that $V$ is ${\bf R}^n$ or ${\bf C}^n$, $n \ge 2$,
equipped with a norm $\|\cdot \|_p$ as in Section \ref{definitions,
etc. (normed vector spaces)}, $1 \le p \le \infty$.  If $p = 1$ or $p
= \infty$, then it is not hard to show that the unit ball is not
strictly convex.  However, if $1 < p < \infty$, then the unit ball is
strictly convex.  This assertion can be derived from the strict
convexity of the function $t \mapsto t^p$ on $[0,\infty)$, mentioned
in the previous section.  For this it is convenient to reformulate the
strict convexity of the unit ball as follows: for each $x$ and $y$ in
${\bf R}^n$ or ${\bf C}^n$ such that 
\begin{equation}
	\sum_{j=1}^n |x_j|^p = \sum_{j=1}^n |y_j|^p = 1
\end{equation}
and $x \ne y$, and for each $\lambda \in (0,1)$, we have that
\begin{equation}
	\sum_{j=1}^n |\lambda \, x_j + (1-\lambda) \, y_j|^p < 1.
\end{equation}

\section{Linear functionals}

	As in earlier chapters, we shall make the standing assumption
that vector spaces have finite dimension in this section.

\beginproposition
\label{strict convexity of the unit ball and linear functionals}
Let $(V, \|\cdot \|)$ be a normed vector space.  The unit ball $B_1$
of $V$ is strictly convex if and only if for every nonzero linear
functional $\mu$ on $V$, there is a unique point $x \in B_1$ such that
\begin{equation}
\label{mu(x) = sup {mu(z) : z in B_1}}
	\mu(x) = \sup \{\mu(z) : z \in B_1\}
\end{equation}
when $V$ is a real vector space, and
\begin{equation}
\label{Re mu(x) = sup {Re mu(z) : z in B_1}}
	\Re \mu(x) = \sup \{\Re \mu(z) : z \in B_1\}
\end{equation}
when $V$ is a complex vector space.  Here $\Re a$ denotes the real
part of a complex number $a$.
\end{proposition}

	In other words, the last part says that there is a unique
point in $B_1$ at which $\mu$ or $\Re \mu$ attains its maximum.  We
can in fact rewrite (\ref{mu(x) = sup {mu(z) : z in B_1}}) and
(\ref{Re mu(x) = sup {Re mu(z) : z in B_1}}) as
\begin{equation}
	\mu(x) = \|\mu\|^*
\end{equation}
for both real and complex vector spaces, where $\|\mu\|^*$ denotes the
dual norm of $\mu$ (as in Section \ref{Dual spaces and norms}).  In
other words, $\mu(z)$ or $\Re \mu(z)$ should be equal to $|\mu(z)|$ in
order to be as large as possible, since they are always less than or
equal to $|\mu(z)|$, and equality can be arranged on $B_1$ by
multiplying $z$ by a scalar of magnitude $1$.

	Every finite-dimensional vector space is linearly isomorphic
to a Euclidean space, and this permits one to apply standard results
on Euclidean spaces pertaining to continuity, compactness, and so on.
It also does not matter which linear isomorphism one uses for this,
i.e., the results would be the same.  The general fact that continuous
real-valued functions on compact sets attain their maxima can be
employed to conclude that for any normed vector space $(V, \|\cdot
\|)$, every linear functional on $V$ attains its maximum on the unit
ball $B_1$ in $V$.  More precisely, linear functionals and their real
parts are continuous with respect to this topology on $V$, coming from
the usual one on a Euclidean space, and the ball $B_1$ is compact
because it is closed and bounded.  The comments near the end of
Section \ref{definitions, etc. (normed vector spaces)} are relevant
for the latter.

	Proposition \ref{strict convexity of the unit ball and linear
functionals} is concerned with the uniqueness of points at which a
maximum is attained.  Let us begin with the ``only if'' part.  Suppose
that $B_1$ is strictly convex, and that $\mu$ is a nonzero linear
functional on $V$.  Assume, for the sake of a contradiction, that
there are two distinct points $x, y \in B_1$ at which the maximum is
attained.  As above, this means that
\begin{equation}
\label{mu(x) = mu(y) = ||mu||^*}
	\mu(x) = \mu(y) = \|\mu \|^*
\end{equation}
(whether $V$ is a real or complex vector space).  We also have that
$\|x\| = \|y\| = 1$, since 
\begin{equation}
\label{|mu(w)| le ||mu||^* ||w|| < ||mu||^* when ||w|| < 1}
	|\mu(w)| \le \|\mu \|^* \, \|w\| < \|\mu\|^*
			\qquad\hbox{when } \|w\| < 1.
\end{equation}

	For each real number $t \in (0,1)$, consider $t \, x + (1-t)
\, y \in V$.  We have that
\begin{equation}
	\|t \, x + (1-t) \, y \| \le 1
\end{equation}
by the usual properties of the norm.  On the other hand,
\begin{equation}
	\mu(t \, x + (1-t) \, y) = \|\mu \|^*,
\end{equation}
because of (\ref{mu(x) = mu(y) = ||mu||^*}).  Thus
\begin{equation}
	\|t \, x + (1-t) \, y \| = 1,
\end{equation}
as in (\ref{|mu(w)| le ||mu||^* ||w|| < ||mu||^* when ||w|| < 1}).
This contradicts the strict convexity of $B_1$.  Therefore we obtain
that there is only one point in $B_1$ at which the maximum is attained
when $B_1$ is strictly convex.

	For the ``if'' part of the proposition, let us assume for the
sake of a contradiction that there are points $x, y \in V$ such that
$\|x\| = \|y\| = 1$, $x \ne y$, and 
\begin{equation}
	\|t \, x + (1-t) \, y \| = 1
\end{equation}
for some real number $t \in (0,1)$.  Set $u = t \, x + (1-t) \, y$.

	As in (\ref{lambda_0(v) = ||lambda_0||^* ||v||}) in Section
\ref{Second duals}, there is a nonzero linear functional $\mu$ on $V$
such that $\mu(u) = \|\mu \|^*$.  Consider $\mu(x)$, $\mu(y)$.  On
the one hand, 
\begin{equation}
\label{mu(u) = mu(t x + (1-t) y) = t mu(x) + (1-t) mu(y)}
	\mu(u) = \mu(t \, x + (1-t) \, y) 
		    = t \, \mu(x) + (1-t) \, \mu(y).
\end{equation}
On the other hand,
\begin{equation}
	|\mu(x)|, |\mu(y)| \le \|\mu \|^*,
\end{equation}
since $\|x\| = \|y\| = 1$.  It follows that
\begin{equation}
	\mu(x) = \mu(y) = \|\mu \|^*.
\end{equation}
Thus we have more than one point in $B_1$ at which the maximum is
attained.  This completes the proof of Proposition \ref{strict
convexity of the unit ball and linear functionals}.

\beginproposition
\label{strict convexity of the dual unit ball and linear functionals}
Let $(V, \|\cdot \|)$ be a normed vector space, and let $\|\cdot \|^*$
be the dual norm of $\|\cdot \|$ on the dual space $V^*$.  The unit
ball $B_1^*$ of $V^*$ is strictly convex if and only if for every
point $x \in V$ with $\|x\| = 1$ there is a unique linear functional
$\mu$ on $V$ such that $\|\mu \|^* = 1$ and $\mu(x) = 1$.
\end{proposition}

	Here we are interchanging the roles of $x$ and $\mu$, compared
to the situation in Proposition \ref{strict convexity of the unit ball
and linear functionals}.  As before, it is the uniqueness of $\mu$
which is at issue, since (\ref{lambda_0(v) = ||lambda_0||^* ||v||}) in
Section \ref{Second duals} implies that there is a $\mu \in V^*$ such
that $\|\mu \|^* = 1$ and $\mu(x) = 1$ for a given $x \in V$ with
$\|x\| = 1$.

	One can derive Proposition \ref{strict convexity of the dual
unit ball and linear functionals} from Proposition \ref{strict
convexity of the unit ball and linear functionals} by applying the
latter to $(V^*, \|\cdot \|^*)$ instead of $(V, \|\cdot \|)$.  The
second dual of $V$ is now in the position of the dual space before,
and we know from Section \ref{Second duals} that the second dual of
$V$ is equivalent to $V$, both as a vector space and for the norm.

\section{Uniqueness of points of minimal distance}
\index{minimal distances}
\label{Uniqueness of points of minimal distance}

	Let us mention a variant of the ``only if'' part of
Proposition \ref{strict convexity of the unit ball and linear
functionals}.

\beginproposition
\label{strict convexity of the unit ball in minimizing distances}
Suppose that $(V, \|\cdot\|)$ is a (finite-dimensional)
normed vector space in which the unit ball $B_1$ is strictly convex.
Let $A$ be a closed convex subset of $V$ (where ``closed'' can be
defined in terms of a linear isomorphism with a Euclidean space, as
before), and let $w$ be any element of $V$.  Then there is a unique
point $z \in A$ such that
\begin{equation}
\label{||w-z|| = inf {||w-u|| : u in A}}
	\|w-z\| = \inf \{\|w-u\| : u \in A\}.
\end{equation}
\end{proposition}

	For instance, $A$ could be a vector subspace of $V$.

	Let $A$, $w$, etc., be given as above.  The existence of $z
\in A$ which satisfies (\ref{||w-z|| = inf {||w-u|| : u in A}})
follows from general considerations of continuity and compactness
again.  More precisely, $A$ may not be bounded, and hence not compact,
but one can reduce to that situation because only a bounded subset of
$A$ is needed for the infimum.  In other words, points which are too
far away will not have a chance of minimizing the distance to $w$.
Strictly speaking, this uses remarks in Section \ref{definitions,
etc. (normed vector spaces)}, as does the continuity of $\|w-u\|$ as a
function of $u$.  For this part of the argument one does not need the
assumption that the unit ball $B_1$ in $V$ is strictly convex, nor
does one need the convexity of $A$.

	Suppose for the sake of a contradiction that there is point $y
\in A$, $y \ne z$, which also satisfies
\begin{equation}
	\|w-y\| = \inf \{\|w-u\| : u \in A\}.
\end{equation}
In particular, $\|w-y\| = \|w-z\|$.  Strict convexity of the unit ball
$B_1$ in $V$ implies that any nontrivial convex combination of $w-z$
and $w-y$ has norm strictly less than $\|w-z\| = \|w-y\|$.  On the
other hand, a convex combination of $w-z$ and $w-y$ can be written as
$w-p$, where $p$ is a convex combination of $z$ and $y$.  This point
$p$ lies in $A$, because $A$ is convex.  Hence $\|w-p\| \ge \|w-z\|$,
since $z$ was chosen so that $\|w-z\|$ is minimal.  This contradicts
the previous assertion to the effect that $\|w-p\| < \|w-z\|$, and
Proposition \ref{strict convexity of the unit ball in minimizing
distances} follows.

\beginremark 
{\rm Assume that $(V, \langle \cdot, \cdot \rangle)$ is a
(finite-dimensional) inner product space and $\|\cdot \|$ is the norm
that comes from the inner product.  Let $S$ be a vector subspace of
$V$, and let $w$ be any element of $V$.  As in Section \ref{Inner
product spaces, continued}, there is a point $z \in S$ such that $w-z
\in S^\perp$, i.e., $w-z$ is orthogonal to every element of $S$.  This
point $z$ is exactly the one for which $\|w-z\|$ is minimal.  Indeed,
if $u$ is any other element of $S$, then we can write $w-u$ as $(w-z)
+ (z-u)$, and $w-z$ and $z-u$ are orthogonal by the choice of $z$.
Hence
\begin{equation}
	\|w-u\| = (\|w-z\|^2 + \|z-u\|^2)^{1/2} = \|w-z\|.
\end{equation}
}
\end{remark}

\section{Clarkson's inequalities}
\label{Clarkson's inequalities}
\index{Clarkson's inequalities}

	Let $n$ be a positive integer, and $p$ a real number such that
$1 < p < \infty$.  Consider ${\bf C}^n$ equipped with the norm
$\|\cdot\|_p$ defined in (\ref{def of ||v||_p}).  If $p \ge 2$, then
\begin{equation}
   \|\textstyle\frac{1}{2}(x+y)\|_p^p + \|\textstyle\frac{1}{2}(x-y)\|_p^p
	\le \textstyle\frac{1}{2} \|x\|_p^p + \textstyle\frac{1}{2} \|y\|_p^p
\end{equation}
for all $x, y \in {\bf C}^n$.  If $p \le 2$, and if $p'$ is the exponent
conjugate to $p$ (so that $1/p + 1/p' = 1$), then
\begin{equation}
  \|\textstyle\frac{1}{2}(x+y)\|_p^{p'} + \|\textstyle\frac{1}{2}(x-y)\|_p^{p'}
	\le \Bigl(\textstyle\frac{1}{2} \|x\|_p^p + 
		\textstyle\frac{1}{2} \|y\|_p^p\Bigr)^{\frac{1}{p-1}}
\end{equation}
for all $x, y \in {\bf C}^n$.  These are \emph{Clarkson's inequalities}.
See (15.5) on p225 and (15.8) on p227 of \cite{HS}. 

	Thus, for instance, if $\|x\|_p = \|y\|_p = 1$, then one gets
an upper bound for $\|x-y\|_p$ in terms of how close $\|(x+y)/2\|_p$
is to $1$.

\chapter{Spectral theory}

	In this chapter, vector spaces are assumed to be
finite-dimensional.  In Sections \ref{The spectrum and spectral
radius} -- \ref{Spectral radius and norms, 2}, vector spaces are
assumed to be complex as well, while real and complex vector spaces
are considered in Sections \ref{Inner product spaces (and spectral
theory)} -- \ref{Commuting families of operators}.

\section{The spectrum and spectral radius}
\label{The spectrum and spectral radius}

	Let $V$ be a vector space, and let $T$ be a linear operator
from $V$ to $V$.  The \emph{spectrum}\index{spectrum (of a linear
transformation)} is the set of complex numbers $\alpha$ such that
$\alpha$ is an \emph{eigenvalue}\index{eigenvalue} of $T$, i.e., there
is a nonzero vector $v$ in $V$ such that $T(v) = \alpha v$.  In this
case one says that $v$ is an \emph{eigenvector}\index{eigenvector} of
$T$ with eigenvalue $\alpha$.  If $\alpha$ is an eigenvalue of $T$,
then
\begin{equation}
	\{v \in V : T(v) = \alpha \, v\}
\end{equation}
is called the \emph{eigenspace}\index{eigenspace} associated to the
eigenvalue $\alpha$, and it is a vector subspace of $V$.

	Equivalently, $\alpha$ is not in the spectrum of $T$ exactly
when $T - \alpha I$ is an invertible linear operator on $V$, where $I$
denotes the identity transformation on $V$.  This uses the fact that a
linear operator on $V$ is invertible if and only if its kernel is
trivial.

	These notions also make sense for linear transformations on a
real vector space, but the next result is not true in general in that
case.  There are substitutes for this, a basic instance of which is
mentioned in Section \ref{Inner product spaces (and spectral theory)}.

\begintheorem
\label{the spectrum is nontrivial (over the complex numbers)}
If $V$ is a vector space of nonzero dimension and $T$ is a linear
operator on $V$, then the spectrum of $T$ is nonempty.
\end{theorem}

	To prove this, one can use the following characterization of
the spectrum of $T$: $\alpha$ lies in the spectrum of $T$ exactly when
the determinant of $T - \alpha I$ is $0$.  On the other hand, the
determinant of $T - \alpha I$ is a polynomial of $\alpha$, and the
degree of the polynomial is equal to the dimension of $V$.  The
``Fundamental Theorem of Algebra'' states that for every nonconstant
polynomial $P(z)$ on the complex numbers there is at least one $a \in
{\bf C}$ such that $P(a) = 0$.  (See \cite{Ru1} for a proof.)
Thus $\det(T - \alpha I)$ has at least one root, and the theorem
follows.

	Note that the number of elements of the spectrum of $T$ is at
most the dimension of $V$.  This can also be seen as a consequence of
the fact that the spectrum consists of the zeros of the polynomial
$\det (T - \alpha I)$, or one can observe that nonzero eigenvectors of
$T$ corresponding to distinct eigenvalues are automatically linearly
independent.

	Given a linear operator $T$ on a vector space $V$, define
the \emph{spectral radius}\index{spectral radius} $\Rad(T)$ by
\begin{equation}
\label{defn of Rad(T)}
	\Rad(T) = \max \{\alpha \in {\bf C}: \alpha 
				\hbox{ lies in the spectrum of } T\}.
\end{equation}
This is equivalent to
\begin{eqnarray}
\label{defn of Rad(T), 2}
\lefteqn{\quad \Rad(T) = }						\\
	&&  \inf\{r \ge 0 : T - \alpha I
			\hbox{ is invertible for all } 
			  \alpha \in {\bf C}
			\hbox{ such that } |\alpha| > r\}.	\nonumber
\end{eqnarray}

\section{Spectral radius and norms}
\label{Spectral radius and norms}

	Let $V$ be a vector space, and let $\|\cdot \|$ be a norm on
$V$.  If $T$ is a linear transformation on $V$, then we can define the
operator norm\index{operator norm} $\|T\|_{op}$ of $T$ using $\|\cdot
\|$ on $V$, as in Section \ref{Linear transformations and norms}.
More precisely, we are using $\| \cdot \|$ as the norm on $V$ in the
role of both $\|\cdot \|_1$ and $\|\cdot \|_2$ in the context of
Section \ref{Linear transformations and norms}, i.e., for both the
domain and range.

\beginlemma
\label{Rad(T) le ||T||_{op}}
$\Rad(T) \le \|T\|_{op}$.
\end{lemma}

	If we think of $\Rad(T)$ as being given by (\ref{defn of
Rad(T)}), then we can establish the lemma as follows.  Let $\alpha$
be an eigenvalue of $T$, and let $v \in V$ be a nonzero eigenvector 
corresponding to $\alpha$.  Thus $T(v) = \alpha \, v$, and from this
it is easy to see that $|\alpha| \le \|T\|_{op}$.  The lemma follows
easily from this.

	One can also approach the lemma through the following.

\beginlemma
\label{(I-A)^{-1} when ||A||_{op} < 1}
If $A : V \to V$ is a linear mapping such that $\|A\|_{op}
< 1$, then $I - A$ is invertible, and 
\begin{equation}
\label{||(I - A)^{-1}||_{op} le (1 - ||A||_{op})^{-1}}
	\|(I - A)^{-1}\|_{op} \le (1 - \|A\|_{op})^{-1}.
\end{equation}
\end{lemma}

	In the context of Lemma \ref{Rad(T) le ||T||_{op}}, one can
apply Lemma \ref{(I-A)^{-1} when ||A||_{op} < 1} with $A = \alpha^{-1}
T$ when $\alpha \in {\bf C}$, $|\alpha| > \|T\|_{op}$ (so that
$\|A\|_{op} < 1$), to get that $I - \alpha^{-1} T$ is invertible.
This is the same as saying that $\alpha I - T$ is invertible when
$|\alpha| > \|T\|_{op}$, so that $\Rad(T) \le \|T\|_{op}$ by
(\ref{defn of Rad(T), 2}).

	A famous technique for deriving Lemma \ref{(I-A)^{-1} when
||A||_{op} < 1} is the method of ``Neumann series''.\index{Neumann
series}  Specifically, one would like to obtain the inverse to $I - A$
by summing the series
\begin{equation}
	\sum_{n=0}^\infty A^n.
\end{equation}
This series converges ``absolutely'' because
\begin{equation}
	\|A^n\|_{op} \le \|A\|_{op}^n,
\end{equation}
as in (\ref{op norm UT le op norm U x op norm T}) in Section
\ref{Linear transformations and norms}, and because $\|A\|_{op} < 1$
by hypothesis.  By a standard computation the sum of the series is the
inverse of $I - A$.  The norm of the sum is less than or equal to the
sum of the norms, and the latter is bounded by a geometric series
whose sum is $(1-\|A\|_{op})^{-1}$.

	Alternatively, one can first assert that $I - A$ is invertible
because its kernel is trivial, which amounts to the same argument as
the initial one for Lemma \ref{Rad(T) le ||T||_{op}}.  The inequality
(\ref{||(I - A)^{-1}||_{op} le (1 - ||A||_{op})^{-1}}) is equivalent
to the statement that
\begin{equation}
	\|(I-A)^{-1}(v)\| \le (1-\|A\|_{op})^{-1} \, \|v\|
			\qquad\hbox{for all } v \in V.
\end{equation}
This is in turn equivalent to saying that
\begin{equation}
	\|w\| \le (1-\|A\|_{op})^{-1} \, \|(I-A)(w)\|
			\qquad\hbox{for all } w \in V.
\end{equation}
To get this, one can compute as follows:
\begin{equation}
	\|w\| \le \|(I-A)(w)\| + \|A(w)\| 
		\le \|(I-A)(w)\| + \|A\|_{op} \, \|w\|.
\end{equation}

\section{Spectral radius and norms, 2}
\label{Spectral radius and norms, 2}

	Let $V$ be a vector space with norm $\|\cdot \|$, and
let $\|\cdot \|_{op}$ denote the corresponding operator norm
for linear transformations on $V$ as in the previous section.

\beginlemma
\label{Rad(T) le ||T^n||_{op}^{1/n}}
If $T$ is a linear transformation on $V$ and $n$ is a positive
integer, then $\Rad(T) \le \|T^n\|_{op}^{1/n}$.
\end{lemma}

	Note that $\|T^n\|_{op}^{1/n}$ is always less than or equal to
$\|T\|_{op}$, because of (\ref{op norm UT le op norm U x op norm T})
in Section \ref{Linear transformations and norms}.

	Let us mention a couple of arguments for Lemma \ref{Rad(T) le
||T^n||_{op}^{1/n}}.  For the first, let $\alpha$ be any eigenvalue of
$T$.  Then $\alpha^n$ is an eigenvalue of $T^n$, with the same
eigenvector.  Hence $|\alpha|^n \le \Rad(T^n)$.  Applying Lemma
\ref{Rad(T) le ||T||_{op}} to $T^n$ instead of $T$ we obtain that
$|\alpha|^n \le \|T^n\|_{op}$, which is the same as $|\alpha| \le
\|T^n\|_{op}^{1/n}$.  Since this holds for any eigenvalue of $T$,
we may conclude that $\Rad(T) \le \|T^n\|_{op}^{1/n}$.

	Alternatively, one can observe that $\alpha^n I - T^n$ is
equal to the product of $\alpha I - T$ and another operator (just as
$x^n - y^n$ can be written as the product of $x-y$ and another
algebraic expression).  This implies that $\alpha I - T$ is invertible
whenever $\alpha^n I - T^n$ is, and one can employ the
characterization (\ref{defn of Rad(T), 2}) of the radius of
convergence.

\beginproposition
\label{spectral radius of T equal to limit of ||T^n||_{op}^{1/n}}
If $T$ is a linear transformation on the vector space $V$, then
the limit
\begin{equation}
\label{lim_{n to infty} ||T^n||_{op}^{1/n}}
	\lim_{n \to \infty} \|T^n\|_{op}^{1/n}
\end{equation}
exists and is equal to $\Rad(T)$.
\end{proposition}

	This may seem a bit odd at first, since the spectral radius
$\Rad(T)$ does not depend on the norm on $V$, while (\ref{lim_{n to
infty} ||T^n||_{op}^{1/n}}) does involve the norm, but it is not hard
to show that (\ref{lim_{n to infty} ||T^n||_{op}^{1/n}}) in fact also
does not depend on the norm.  Specifically, any two norms on a
finite-dimensional vector space are bounded by positive constant
multiples of each other, as in Section \ref{definitions, etc. (normed
vector spaces)}, and this leads to a similar relation for the operator
norms.  The effects of these constants washes out in (\ref{lim_{n to
infty} ||T^n||_{op}^{1/n}}), since $\lim_{n \to \infty} A^{1/n} = 1$
for any positive real number $A$ (as in \cite{Ru1}).

	One can use this type of remark together with Jordan canonical
forms to prove the proposition.  Let us now describe another argument.

	If $T$ is as in the proposition, then $(\alpha I - T)^{-1}$
exists for all $\alpha \in {\bf C}$ such that $|\alpha| > \Rad(T)$.
It will be convenient to reformulate this as saying that $(I - \beta
T)^{-1}$ exists for all $\beta \in {\bf C}$ such that $|\beta| <
\Rad(T)^{-1}$.  This includes $\beta = 0$, for which invertibility is
trivial.  If $\Rad(T) = 0$, then we take $\Rad(T)^{-1} = + \infty$,
and $(I - \beta T)^{-1}$ exists for all $\beta \in {\bf C}$.

	In fact, $(I - \beta T)^{-1}$ is a rational function on ${\bf
C}$ with poles lying outside the disk $\{\beta \in {\bf C} : |\beta| <
\Rad(T)^{-1}\}$.  More precisely, it is a rational function with
values in the vector space of linear mappings on $V$, but this is not
too serious, since the vector space of linear mappings on $V$ has
finite dimension (since $V$ is finite-dimensional).

	Recall that a complex-valued function of a complex variable
$z$ is said to be rational if it can be written as $P(z)/Q(z)$, where
$P(z)$ and $Q(z)$ are polynomials in $z$ (with complex coefficients)
and $Q(z)$ is not identically $0$.  The points at which $Q(z)$
vanishes are called the poles of the rational function, except to the
extent that they can be cancelled with roots of $z$ (in such a way
that the function can be rewritten as a quotient of polynomials in
which the denominator does not vanish at the point in question).

	To say that $(I - \beta T)^{-1}$ is a rational function of
$\beta$ means that it can be written as a finite sum of functions
which are each a complex-valued rational function of $\beta$ times a
fixed linear mapping on $V$.  This is true in this case because of
Cramer's rule, which allows $(I - \beta T)^{-1}$ to be expressed as
$\det(I - \beta T)^{-1}$ times a linear transformation on $V$ defined
in terms of minors of $I - \beta T$.  These minors are polynomials in
$\beta$, as is $\det(I - \beta T)$ (which is also not identically $0$
since it is equal to $1$ at $\beta = 0$).  The poles of this rational
function come from the zeros of $\det(I - \beta T)$, and they lie
outside of the disk $\{\beta \in {\bf C} : |\beta| < \Rad(T)^{-1}\}$
by assumption.

	Because the poles of this rational function lie outside of
this disk, this rational function of $\beta$ is equal to a convergent
power series in $\beta$ on this disk (centered at $0$).  This works
for any rational function, and more generally for any complex-analytic
function on a disk.

	The form of the power series for $(I - \beta T)^{-1}$ about
the origin is clear, and is given by $\sum_{n=0}^\infty \beta^n \,
T^n$.  We conclude that this series converges for all $\beta \in {\bf
C}$ such that $|\beta| < \Rad(T)^{-1}$.  This implies that the
coefficients of the series do not grow too fast in a certain sense.
This works in a general way, as discussed in \cite{Ru1}; for rational
functions, one has even more information about the coefficients, but
for the moment it is only the approximate size that is needed.  An
adequate bottom line is that for every real number $r > \Rad(T)$,
there is a positive constant $L$ (which may depend on $r$) such that
\begin{equation}
\label{||T^n||_{op} le L r^n}
	\|T^n\|_{op} \le L \, r^n
\end{equation}
for all $n \ge 0$.  

	We can rewrite this inequality as
\begin{equation}
	\|T^n\|_{op}^{1/n} \le L^{1/n} \, r.
\end{equation}
Using this and Lemma \ref{Rad(T) le ||T^n||_{op}^{1/n}}, it is not
hard to show that (\ref{lim_{n to infty} ||T^n||_{op}^{1/n}}) exists
and is equal to $\Rad(T)$.  Of course it is important too that $r$
above can be taken to be any real number larger than $\Rad(T)$.

	In the preceding discussion of upper bounds for
$\|T^n\|_{op}$, we are not using special features of the operator norm
$\|\cdot \|_{op}$ for linear transformations on $V$.  In effect, this
is accommodated by the constant $L$.  One can think in terms of
starting with computations for the complex-valued setting; one can
deal with different scalar components separately, and then combine
them to get results for linear transformations.  The final packaging
in terms of $\|\cdot \|_{op}$, instead of some other norm, only
affects the constant $L$.

\beginremark
{\rm If $a_n = \|T^n\|_{op}^{1/n}$, then
\begin{equation}
\label{a_{m+n} le a_m^{m/(m+n)} a_n^{n/(m+n)}}
	a_{m+n} \le a_m^{m/(m+n)} \, a_n^{n/(m+n)}
\end{equation}
for all positive integers $m$ and $n$.  This follows from the fact
that the operator norm of a product of operators is less than or equal 
to the product of the operator norms.  (Note that this applies equally
well to operators on real or complex normed vector spaces.)

	Now suppose that $\{a_n\}_{n=1}^\infty$ is simply any sequence
of nonnegative real numbers such that (\ref{a_{m+n} le a_m^{m/(m+n)}
a_n^{n/(m+n)}}) holds for all positive integers $m$ and $n$.  This is
a multiplicative convexity property; if $a_n > 0$ for all $n$, then
we can rewrite it as
\begin{equation}
	\log a_n \le \textstyle\frac{m}{m+n} \log a_m
			+ \textstyle\frac{n}{m+n} \log a_n,
\end{equation}
as in ordinary ``additive'' convexity.  (If $a_n = 0$ for some $n$,
then $a_j = 0$ for all $j \ge n$, and one might as well not worry
about these $a_j$'s too much anyway.)

	A well-known result states that $\lim_{n \to \infty} a_n$
automatically exists under this convexity condition.  To show this,
it is enough to verify that $\limsup_{n \to \infty} a_n \le
\liminf_{n \to \infty} a_n$.  Let us first observe that
\begin{equation}
	a_{jn} \le a_n
\end{equation}
for all positive integers $j$ and $n$.  This can be derived from
(\ref{a_{m+n} le a_m^{m/(m+n)} a_n^{n/(m+n)}}) using induction.
Next, if $1 < i < n$, then
\begin{equation}
	a_{jn + i} \le a_{jn}^{jn/(jn+i)} \, a_i^{i/(jn+i)}
		\le a_n^{jn/(jn+i)} \, a_i^{i/(jn+i)}.
\end{equation}
Once one has this, it is not difficult to finish the argument.

}
\end{remark}

\section{Inner product spaces}
\label{Inner product spaces (and spectral theory)}

	In this section, vector spaces are assumed to be
finite-dimensional, but both real and complex vector spaces are
allowed.

	Let $(V, \langle \cdot, \cdot \rangle)$ be an inner product
space, as in Section \ref{Inner product spaces}, real or complex.  If
$T : V \to V$ is a linear transformation, then there is a unique
linear transformation $T^* : V \to V$\index{$T^*$} such that
\begin{equation}
\label{characterization of T^*}
	\langle T(v), w \rangle = \langle v, T^*(w) \rangle
\end{equation}
for all $v, w \in V$.  This is called the \emph{adjoint}\index{adjoint
(of a linear transformation)} of $T$, and when $V$ is a real vector
space, $T^*$ is sometimes called the \emph{transpose}\index{transpose
(of a linear transformation)} of $T$.  Specifically, if we fix an
orthonormal basis for $V$, and express $T$ in terms of a matrix
relative to this orthonormal basis, then $T^*$ corresponds to the
transpose\index{transpose (of a matrix)} of the matrix for $T$ when
$V$ is a real vector space, and $T^*$ corresponds to the conjugate
transpose\index{conjugate transpose (of a matrix)} of the matrix for
$T$ in the complex case.  (Recall that if $\{t_{j,k}\}$ is a matrix,
then the matrix $\{s_{l, m}\}$ is the \emph{transpose} of
$\{t_{j,k}\}$ if $l$ and $m$ run through the same ranges of indices
for both $s_{l,m}$ and $t_{m,l}$, and if $s_{l,m} = t_{m,l}$ for all
$l$ and $m$.  The matrix $\{u_{l,m}\}$ is the \emph{conjugate
transpose} of $\{t_{j,k}\}$ if $l$ and $m$ run through the same ranges
of indices for both $u_{l,m}$ and $t_{m,l}$, and if $u_{l,m} =
\overline{t_{m,l}}$ for all $l$ and $m$, where $\overline{a}$ denotes
the complex conjugate of $a$.)  It is easy to see that $T^*$ is
uniquely determined by (\ref{characterization of T^*}), and in
particular different orthonormal bases for $V$ lead to the same linear
transformation $T^*$ when one computes $T^*$ in terms of matrices.

	Let $\|\cdot \|$ be the norm on $V$ associated to the inner
product $\langle \cdot, \cdot \rangle$, and let $\|\cdot \|_{op}$
denote the operator norm for linear transformations on $V$ defined
using $\|\cdot \|$.  It is not hard to check that
\begin{equation}
\label{formula for ||T||_{op} in terms of the inner product}
	\|T\|_{op} = \sup \{ |\langle T(v), w \rangle| : 
			    v, w \in V, \ \|v\| = \|w\| = 1\}
\end{equation}
for any linear transformation $T$ on $V$.  (This uses the
Cauchy-Schwarz inequality (\ref{Cauchy-Schwarz inequality in inner
product spaces}) to show that $\|T\|_{op}$ is greater than or equal to
the right side, and the choice of $w$ equal to a scalar multiple of
$T(v)$ (unless $T(v) = 0$) for the other inequality.)  One also has
that
\begin{equation}
\label{||T^*||_{op} = ||T||_{op}}
	\|T^*\|_{op} = \|T\|_{op},
\end{equation}
as can easily be derived from (\ref{formula for ||T||_{op} in terms of
the inner product}).

	Note that $I^* = I$, where $I$ denotes the identity
transformation.  If $S$ and $T$ are linear transformations on $V$ and
$a$, $b$ are scalars, then
\begin{equation}
	(a \, S + b \, T)^* = a \, S^* + b \, T^*
\end{equation}
when $V$ is a real vector space and
\begin{equation}
	(a \, S + b \, T)^* = 
			\overline{a} \, S^* + \overline{b} \, T^*
\end{equation}
when $V$ is a complex vector space.  Also,
\begin{equation}
	(T^*)^* = T
\end{equation}
and
\begin{equation}
\label{(S T)^* = T^* S^*}
	(S \, T)^* = T^* \, S^*
\end{equation}
(in both cases).  These properties can be verified in a
straightforward manner.

	From (\ref{(S T)^* = T^* S^*}) it follows that if $T$
is invertible, then $T^*$ is invertible as well, and
\begin{equation}
	(T^*)^{-1} = (T^{-1})^*.
\end{equation}
It follows that $T - \lambda \, I$ is invertible if and only if $T^* -
\lambda \, I$ is invertible in the real case, and if and only if $T^*
- \overline{\lambda} \, I$ is invertible in the complex case.  Thus
$\lambda$ is in the spectrum of $T$ if and only if
$\overline{\lambda}$ is in the spectrum of $T^*$.

	A linear transformation $T$ on $V$ is said to be
\emph{self-adjoint}\index{self-adjoint (linear transformation)} if
$T^* = T$.  When $V$ is a real vector space, one might call $T$
\emph{symmetric}\index{symmetric (linear transformation)} in this
event.  Sums of self-adjoint linear transformations are again
self-adjoint, as are products of self-adjoint linear transformations
by \emph{real} numbers.  The eigenvalues of a self-adjoint linear
transformation on a complex inner product space are always real
numbers.  This can be derived from the observations in the preceding
paragraph, and one can also compute as follows.  If $v$ is a nonzero
vector in $V$ which is an eigenvector of $T$ with eigenvalue $\lambda$,
and if $T$ is self-adjoint, then
\begin{eqnarray}
	\lambda \, \langle v, v \rangle & = & \langle \lambda \, v, v \rangle 
		= \langle T(v), v \rangle = \langle v, T^*(v) \rangle 	\\
	& = & \langle v, T(v) \rangle = \langle v, \lambda \, v \rangle
		= \overline{\lambda} \, \langle v, v \rangle,  	\nonumber
\end{eqnarray}
and hence $\overline{\lambda} = \lambda$.

	If $T$ is a self-adjoint linear transformation on $V$ and
$v_1$, $v_2$ are vectors in $V$ which are eigenvectors of $T$ with
eigenvalues $\lambda_1$, $\lambda_2$, respectively, and if $\lambda_1
\ne \lambda_2$, then $v_1$ and $v_2$ are orthogonal to each other.
This is because
\begin{eqnarray}
\label{computation about self-adjoint T and orthogonality of eigenvectors}
	\lambda_1 \, \langle v_1, v_2 \rangle & = &
	\langle \lambda_1 \, v_1, v_2 \rangle = \langle T(v_1), v_2 \rangle
								\\
  & = & \langle v_1, T(v_2) \rangle = \langle v_1, \lambda_2 \, v_2 \rangle 
	= \lambda_2 \, \langle v_1, v_2 \rangle,		\nonumber
\end{eqnarray}
so that $\langle v_1, v_2 \rangle$ is necessarily $0$.

	An important feature of self-adjoint linear transformations is
that they can be diagonalized in an orthonormal basis.  Let us recall
how this can be proved.  A first step is to show that there is at
least one nonzero eigenvector.  On a complex vector space, this is
true for any linear transformation, as in Theorem \ref{the spectrum is
nontrivial (over the complex numbers)} in Section \ref{The spectrum
and spectral radius}.  For a self-adjoint operator $T$ one can also
obtain this by maximizing or minimizing $\langle T(v), v \rangle$ over
the unit sphere ($\|v\| = 1$), and this works in both the real and the
complex cases.  More precisely, one can check that the eigenvectors of
$T$ with norm $1$ are exactly the critical points of the real-valued
function $\langle T(v), v \rangle$ on the unit sphere.
(Alternatively, one might prefer to phrase this in terms of Lagrange
multipliers.)  The existence of points in the unit sphere at which the
maximum or minimum of $\langle T(v), v \rangle$ is attained follows
from general results about continuity and compactness, and these
points are automatically critical points.

	Next, suppose that $T$ is a self-adjoint operator and that $v$
is an eigenvector of $T$ with eigenvalue $\lambda$.  The \emph{orthogonal
complement}\index{orthogonal complement} of $v$ in $V$ is the subspace
of vectors $w \in V$ such that $\langle v, w \rangle = 0$, as in
Section \ref{Inner product spaces, continued}.  The self-adjointness
of $T$ implies that $T$ maps the orthogonal complement of
$v$ to itself, since
\begin{equation}
	\langle v, T(w) \rangle = \langle T(v), w \rangle = 
			\langle \lambda \, v, w \rangle = 0
\end{equation}
when $\langle v, w \rangle = 0$.

	One can restrict $T$ to the orthogonal complement of $v$ and
repeat the process, to get a new eigenvector which is orthogonal to
the previous one, etc.  In this way one can show that there is an
orthonormal basis for the whole space which consists of eigenvectors
of $T$, as desired.

	Now let us assume for a moment that we are working with a
\emph{complex} inner product space.  If $T$ is a linear operator on
the vector space, then $T$ can be diagonalized in an orthogonal basis
if and only if $T$ is \emph{normal}\index{normal (linear
transformation)}, which means that $T$ and $T^*$ commute ($T \, T^* =
T^* \, T$).  Indeed, if $T$ can be so diagonalized, then $T^*$ is
diagonalized by the same basis, and $T$ and $T^*$ commute.
Conversely, any linear operator $T$ can be written as $A + i \, B$,
where $A$ and $B$ are self-adjoint, by taking $A = (T + T^*)/2$ and $B
= (T - T^*)/(2i)$.  Thus $A$ and $B$ can be diagonalized by
orthonormal bases, but in general the bases are not the same.
However, if $T$ is normal, then $A$ and $B$ commute, and one can check
that the eigenspaces of $A$ and $B$ are invariant under each other.
This permits one to choose an orthogonal basis in which both are
diagonalized, and hence in which $T$ is diagonalized.  (See also
Section \ref{Commuting families of operators} concerning this last
point.)

	As a consequence of the diagonalization, if $T$ is a
self-adjoint or normal linear transformation on a complex inner
product space, then the spectral radius of $T$ is equal to the norm of
$T$.  (See also Section \ref{The $C^*$-identity}.)  An analogous
statement holds for self-adjoint linear operators on real inner
product spaces (for which we have not officially defined the spectral
radius).

	On either a real or complex inner product space, suppose that
$T$ is a linear operator such that
\begin{equation}
	\langle T(v), T(w) \rangle = \langle v, w \rangle
\end{equation}
for all $v$, $w$ in the vector space.  This is equivalent to asking
that
\begin{equation}
	\|T(v)\| = \|v\|
\end{equation}
for all $v$, because of \emph{polarization},\index{polarization} as in
Remark \ref{polarization and the parallelogram law} in Section
\ref{Inner product spaces}.  This condition is also equivalent to
\begin{equation}
	T^* = T^{-1}.
\end{equation}
The linear operators which satisfy this property are called
\emph{orthogonal transformations}\index{orthogonal transformations} in
the case of real inner product spaces, and \emph{unitary
transformations}\index{unitary transformations} in the case of complex
inner product spaces.  A unitary transformation $T$ on a complex inner
product space is normal in the sense defined above, since $T$
automatically commutes with $T^{-1}$, and hence $T$ can be
diagonalized in an orthonormal basis.

\section{The $C^*$-identity}
\label{The $C^*$-identity}

	Let $(V, \langle \cdot, \cdot \rangle)$ be a finite-dimensional
inner product space, real or complex, and let $T$ be a linear operator
on $V$.  The \emph{$C^*$-identity}\index{C^*-identity@$C^*$-identity}
states that
\begin{equation}
\label{C^* identity, 1}
	\|T^* \, T\|_{op} = \|T\|_{op}^2,
\end{equation}
where $\|\cdot\|_{op}$ denotes the operator norm of a linear operator
on $V$, relative to the norm $\|\cdot\|$ on $V$ associated to the inner
product.  To see this, notice first that
\begin{equation}
	\|T^* \, T\|_{op} \le \|T^*\|_{op} \, \|T\|_{op} = \|T\|_{op}^2,
\end{equation}
using (\ref{||T^*||_{op} = ||T||_{op}}) in the second step.  On the
other hand,
\begin{equation}
	\|T^* \, T\|_{op} \ge \|T\|_{op}^2
\end{equation}
can be derived from (\ref{formula for ||T||_{op} in terms of the
inner product}) applied to $T^* \, T$.  Thus we get (\ref{C^*
identity, 1}).  Similarly, or by applying replacing $T$ with $T^*$ in
(\ref{C^* identity, 1}), one can obtain that
\begin{equation}
\label{C^* identity, 2}
	\|T \, T^*\|_{op} = \|T\|_{op}^2.
\end{equation}

	As a special case, if $T$ is self-adjoint, then
\begin{equation}
\label{||T^2||_{op} = ||T||_{op}^2}
	\|T^2\|_{op} = \|T\|_{op}^2.
\end{equation}
For any operator $T$ we have the inequality
\begin{equation}
	\|T^2\|_{op} \le \|T\|_{op}^2,
\end{equation}
but the reverse inequality is not true in general.  For instance,
$T^2$ could be the zero operator while $T$ is not.

	If $T$ is self-adjoint, then so is $T^j$ for any positive
integer $j$, since $(T^j)^* = (T^*)^j = T^j$.  Thus we can apply
(\ref{||T^2||_{op} = ||T||_{op}^2}) to $T^j$ to obtain
\begin{equation}
	\|T^{2j}\|_{op} = \|T^j\|_{op}^2.
\end{equation}
It follows that
\begin{equation}
	\|T^l\|_{op} = \|T\|_{op}^l
\end{equation}
when $l$ is a power of $2$, and from this one can verify equality
holds for all positive integers $l$.  (Exercise, using the fact that
the norm of a product is always less than or equal to the product of
the norms.)

	From here it follows in turn that if $V$ is a complex vector
space, then the spectral radius of $T$ is equal to $\|T\|_{op}$, by
Proposition \ref{spectral radius of T equal to limit of
||T^n||_{op}^{1/n}} in Section \ref{Spectral radius and norms, 2}.
This and the preceding identities are also simple consequences of a
diagonalization of $T$ in an orthonormal basis, as in the previous
section.

	Now suppose that $T$ is normal, i.e., $T^* \, T = T \, T^*$.
For any positive integer $l$ we have that
\begin{equation}
	\|T^l\|_{op}^2 = \|(T^l)^* \, T^l\|_{op}
\end{equation}
by the $C^*$-identity.  Because $T$ and $T^*$ commute,
\begin{equation}
	(T^l)^* \, T^l = (T^* \, T)^l,
\end{equation}
and hence
\begin{equation}
	\|(T^l)^* \, T^l\|_{op} = \|(T^* \, T)^l\|_{op} 
		= \|T^* \, T\|_{op}^l = \|T\|_{op}^{2l}.
\end{equation}
The second step uses the fact that $T^* \, T$ is self-adjoint, so that
the earlier formulae are applicable.  Therefore,
\begin{equation}
\label{||T^l||_{op} = ||T||_{op}^l}
	\|T^l\|_{op} = \|T\|_{op}^l
\end{equation}
for all positive integers $l$.

	As before, when $V$ is a complex vector space, this implies
that the spectral radius of $T$ is equal to the norm of $T$.  In this
case, this statement and (\ref{||T^l||_{op} = ||T||_{op}^l}) also
follow from a diagonalization of $T$ in an orthonormal basis.

\section{Projections}
\label{Projections}  
\index{projections (on a vector space)}

	In this section vector spaces are again permitted to be real
or complex.

\begindefinition
\label{def of projections and complementary subspaces}
Let $V$ be a vector space.  A linear operator $P : V \to V$ is said to
be a \emph{projection} if $P^2 = P$.  Two subspaces $V_1$ and $V_2$ of
$V$ are said to be \emph{complementary}\index{complementary subspaces}
or \emph{complements of each other} if $V_1 \cap V_2 = \{0\}$ and
$\span(V_1, V_2) = V$.
\end{definition}

	Observe that $I - P$ is automatically a projection when $P$
is, since $(I - P)^2 = I - P - P + P^2 = I - P$.  The kernel of $I -
P$ is equal to the image of $P$, and the image of $I - P$ is equal to
the kernel of $P$.

\beginlemma
\label{some properties of complementary subspaces and projections}
Two subspaces $V_1$, $V_2$ of a vector space $V$ are complementary if
and only if for each $v \in V$ there are unique vectors $v_1 \in V_1$
and $v_2 \in V_2$ such that $v = v_1 + v_2$.  In this event the
mappings $v \mapsto v_1$ and $v \mapsto v_2$ are linear, and they
define projections from $V$ onto $V_1$ along $V_2$ and from $V$ onto
$V_2$ along $V_1$, respectively.

	If $P : V \to V$ is a linear mapping which is a projection,
then the kernel and image of $P$ are complementary subspaces of $V$,
and $P$ is determined uniquely by its kernel and image.
\end{lemma}

	(Exercise.)

	Let us say that \emph{$P$ is the projection of $V$ onto $V_2$
along $V_1$} when $P$ is a projection on $V$ with kernel $V_1$ and
image $V_2$.  Of course $I - P$ is then the projection of $V$ onto
$V_1$ along $V_2$.

	If $P$ is a projection, then the nonzero vectors in the kernel
of $P$ are eigenvectors of $P$ with eigenvalue $0$, and the nonzero
vectors in the image of $P$ are eigenvectors of $P$ with eigenvalue
$1$.  The spectrum of $P$ consists exactly of $0$ and $1$, except in
the degenerate cases where $P = 0$ or $P = I$, for which the spectrum
consists of only $0$ or only $1$, respectively (unless $V$ has dimension
$0$).

	If $(V, \langle \cdot, \cdot \rangle)$ is an inner product
space, then for each subspace $W$ of $V$ there is a special projection
of $V$ onto $W$, namely, the \emph{orthogonal}
projection,\index{orthogonal projections} as in Section \ref{Inner
product spaces, continued}.  In the terminology above, this is the
projection of $V$ onto $W$ along the orthogonal complement $W^\perp$
of $W$.  One can rephrase this by saying that orthogonal projections
are exactly the same as projections in which the kernel and image of
the operator are orthogonal complements\index{orthogonal complements}
of each other.

\beginlemma
A projection P on an inner product space $(V, \langle \cdot, \cdot
\rangle)$ is an orthogonal projection if and only if it is
self-adjoint.
\end{lemma}

	Suppose first that $P$ is an orthogonal projection.  Let $u$
and $v$ be arbitrary elements of $V$.  Then $(I - P)(v)$ is orthogonal
to $P(u)$, i.e.,
\begin{equation}
	\langle P(u), (I - P)(v) \rangle = 0.
\end{equation}
This is the same as saying that
\begin{equation}
	\langle P(u), v \rangle = \langle P(u), P(v) \rangle.
\end{equation}
Similarly, $(I - P)(u)$ is orthogonal to $P(v)$, so that
\begin{equation}
	\langle P(u), P(v) \rangle = \langle u, P(v) \rangle.
\end{equation}
Hence $\langle P(u), v \rangle = \langle u, P(v) \rangle$, as
desired.

	Conversely, if $P$ is a self-adjoint projection on $V$, then
\begin{equation}
	\langle P(u), (I - P)(v) \rangle 
		= \langle u, P((I - P)(v)) \rangle = 0
\end{equation}
for all $u, v \in V$.  This shows that the image of $P$ is orthogonal
to the kernel of $P$, so that $P$ is an orthogonal projection.

\beginlemma
\label{a projection is orthogonal iff it has norm 1}
\index{projections with operator norm $1$}
A nonzero projection P on an inner product space $(V, \langle \cdot,
\cdot \rangle)$ is an orthogonal projection if and only if $\|P\|_{op}
= 1$.
\end{lemma}

	The ``only if'' part of this lemma was mentioned at the end of
Section \ref{Inner product spaces, continued}.  For the ``if'' part,
suppose that $P$ is a projection on $V$ with operator norm $1$.  The
operator norm of any nonzero projection is greater than or equal to
$1$, since there are nonzero vectors which are mapped to themselves
by the projection.  Thus the real content of the hypothesis that
$\|P\|_{op} = 1$ is that
\begin{equation}
\label{||P(v)|| le ||v|| for all v in V (assumption)}
	\|P(v)\| \le \|v\|  \qquad\hbox{for all } v \in V.
\end{equation}

	Let $u$, $w$ be elements of $V$ such that $u$ is in the kernel
of $P$ and $w$ is in the range of $P$.  We would like to show that $u$
and $w$ are orthogonal to each other.  Assume for the sake of a
contradiction that this is not the case, i.e.,
\begin{equation}
	\langle u, w \rangle \ne 0.
\end{equation}
We may as well assume that $\langle u, w \rangle$ is a positive real
number, since this can be arranged by multiplying $u$ or $w$ by a
scalar.

	Set $v = w - t \, u$, where $t$ is a positive real number,
to be specified in a moment.  Note that $P(v) = w$, independently
of $t$.  On the other hand,
\begin{eqnarray}
	\|v\|^2 = \langle v, v \rangle 
	    & = & \langle w, w \rangle - 2 \, t \, \langle u, w \rangle
				+ t^2 \, \langle u, u \rangle		\\
	    & = & \|w\|^2 - 2 \, t \, \langle u, w \rangle + t^2 \, \|u\|^2.
								\nonumber
\end{eqnarray}
Because $\langle u, w \rangle > 0$, it follows that
\begin{equation}
	\|v\|^2 < \|w\|^2 = \|P(v)\|^2
\end{equation}
when $t > 0$ is sufficiently small.  This contradicts (\ref{||P(v)||
le ||v|| for all v in V (assumption)}), and the lemma follows.

	Now suppose that $V$ is a vector space equipped with a norm
$\|\cdot\|$, which is not required to come from an inner product.  A
nonzero projection on $V$ has operator norm greater than or equal to
$1$, for the same reason as before, but in general it is not as easy
to obtain projections with norm $1$, or even with a good bound on the
norm, as in the case of inner product spaces.  If we take $V$ to be
${\bf R}^n$ or ${\bf C}^n$, with norm $\|\cdot\|_p$ as in Section
\ref{definitions, etc. (normed vector spaces)}, then the usual
``coordinate projections'' have norm $1$.  These are the projections
defined by replacing certain coordinates of a vector $v = (v_1, v_2,
\ldots, v_n)$ with $0$ while leaving the other coordinates unchanged.

	As another special case, let $(V, \|\cdot\|)$ be any normed
vector space, and let $z$ be any vector in $V$ with $\|z\| = 1$.  As
in Section \ref{Second duals}, there is a linear functional $\lambda$
on $V$ such that $\|\lambda\|^* = 1$ and $\lambda(z) = 1$.  Define
$P : V \to V$ by
\begin{equation}
	P(v) = \lambda(v) \, z.
\end{equation}
It is easy to verify that $P$ is a projection of $V$ onto
the $1$-dimensional subspace generated by $z$ with operator
norm $1$.

\section{Remarks about diagonalizable operators}
\label{Remarks about diagonalizable operators}

	Let $V$ be a real or complex vector space, and let $T$ be a
linear operator on $V$.  For each eigenvalue $\lambda$ of $T$, let
$E(T, \lambda)$ denote the corresponding eigenspace, so that
\begin{equation}
	E(T, \lambda) = \{v \in V : T(v) = \lambda \, v\}.
\end{equation}

	Suppose that $\lambda_1, \ldots, \lambda_l$ are distinct
eigenvalues of $T$.  If $v_1, \ldots, v_l$ are vectors in $V$ such
that $v_j \in E(T, \lambda_j)$ for each $j$ and
\begin{equation}
\label{sum_{j=1}^l v_j = 0}
	\sum_{j=1}^l v_j = 0,
\end{equation}
then $v_j = 0$ for each $j$.  This is not hard to prove using induction
on $l$, for instance.

	If $V$ is equipped with an inner product $\langle u, w
\rangle$ and $T$ is self-adjoint with respect to this inner product,
then eigenvectors of $T$ associated to distinct eigenvalues are
orthogonal.  This was mentioned in Section \ref{Inner product spaces
(and spectral theory)}, in the paragraph containing (\ref{computation
about self-adjoint T and orthogonality of eigenvectors}).  The linear
independence property described in the preceding paragraph provides a
version of this for linear operators in general, without inner
products or self-adjointness.

	Let us say that $T$ is diagonalizable if there is a basis of
$V$ consisting of eigenvectors of $T$.  Thus, if one represents $T$ by
a matrix using this basis, then the matrix would be a diagonal matrix
(with nonzero entries only along the diagonal).

	If $T$ is diagonalizable, then clearly $V$ is spanned by the
set of eigenvectors of $T$.  The converse is also true, but requires a
bit more care.  Specifically, one can get a basis of $V$ consisting of
eigenvectors of $T$ by choosing arbitrary bases for the eigenspaces of
$T$, and then combining them to get a basis for all of $V$.  That the
union of the bases of the eigenspaces is a linearly independent set of
vectors in $V$ can be derived from the observation above, in the
paragraph containing (\ref{sum_{j=1}^l v_j = 0}), and from the linear
independence of each of the bases of the eigenspaces separately.  The
fact that the union of the bases of the eigenspaces of $T$ spans $V$
exactly follows from the assumption that $V$ is spanned by the
eigenvectors of $T$.

\section{Commuting families of operators}
\label{Commuting families of operators}

	Let $V$ be a real or complex vector space, and let
$\mathcal{F}$ be a family of linear operators on $V$.  We assume that
every operator $T$ in $\mathcal{F}$ is diagonalizable, and that
$\mathcal{F}$ is a commuting family of operators,\index{commuting
families of operators} which means that
\begin{equation}
	S \circ T = T \circ S
\end{equation}
for all $S$, $T$ in $\mathcal{F}$.  Under these conditions, the
operators in $\mathcal{F}$ can be simultaneously diagonalized.  In
other words, there is a basis $v_1, \ldots, v_n$ for $V$ such that
each $v_i$ is an eigenvector for every operator in $\mathcal{F}$.

	To see this, let us begin with the following observation.
Suppose that $T$ is an operator in $\mathcal{F}$, and that $\lambda$
is an eigenvalue of $T$.  Let $E(T, \lambda)$ denote the corresponding
eigenspace of $T$, as in Section \ref{Remarks about diagonalizable
operators}.  If $S$ is any other operator in $\mathcal{F}$, then
\begin{equation}
\label{S(E(T, lambda)) subseteq E(T, lambda)}
	S(E(T, \lambda)) \subseteq E(T, \lambda).
\end{equation}
Indeed, if $v \in E(T, \lambda)$, then
\begin{equation}
	T(S(v)) = S(T(v)) = S(\lambda \, v) = \lambda \, S(v),
\end{equation}
and hence $S(v) \in E(T, \lambda)$.  Of course we used the information
that $S$ and $T$ commute in the first step in this computation.

	Now let $\lambda_1, \ldots, \lambda_k$ be a list of all of the
eigenvalues of $T$, without repetitions.  The corresponding
eigenspaces $E(T, \lambda_j)$ may have dimension larger than $1$,
because of multiplicities.  For each vector $v$ in $V$, there exist
unique vectors $v_j$ in $E(T, \lambda_j)$, $j = 1, \ldots, k$, such that
\begin{equation}
\label{v = sum_{j=1}^k v_j}
	v = \sum_{j=1}^k v_j.
\end{equation}
The existence of the $v_j$'s comes from the fact that $V$ is spanned
by the eigenvectors of $T$, since $T$ is diagonalizable, and the
uniqueness follows from the observation in the paragraph containing
(\ref{sum_{j=1}^l v_j = 0}).

	For each $j = 1, \ldots, k$, the correspondence $v \mapsto
v_j$ determines a well-defined mapping on $V$, which is in fact
linear, and a projection from $V$ onto $E(T, \lambda_j)$.  This uses
the existence and the uniqueness of the decomposition in (\ref{v =
sum_{j=1}^k v_j}).  Let us denote this projection by $P_j$.

	If $S$ is in $\mathcal{F}$, then
\begin{equation}
\label{S circ P_j = P_j circ S}
	S \circ P_j = P_j \circ S
\end{equation}
for each $j$.  To see this, we can compute as follows.  Given $v$ in
$V$, decompose $v$ into a sum of $v_j \in E(T, \lambda_j)$ as above.
Thus
\begin{equation}
	S(v) = \sum_{j=1}^k S(v_j),
\end{equation}
and each $S(v_j)$ lies in $E(T, \lambda_j)$, by (\ref{S(E(T, lambda))
subseteq E(T, lambda)}).  This implies (\ref{S circ P_j = P_j circ
S}), since the decompositions are unique.

	Let $\alpha$ be any eigenvalue of $S \in \mathcal{F}$, and
$E(S, \alpha)$ the corresponding eigenspace.  Notice that
\begin{equation}
\label{P_j(E(S, alpha)) subseteq E(S, alpha)}
	P_j(E(S, \alpha)) \subseteq E(S, \alpha).
\end{equation}
This is the analogue of (\ref{S(E(T, lambda)) subseteq E(T, lambda)}),
with $T$ replaced by $S$, $\lambda$ by $\alpha$, and $S$ by $P_j$.
The main point is that $S$ and $P_j$ commute, as in (\ref{S circ P_j =
P_j circ S}), which is all that we needed before.

	In other words, (\ref{P_j(E(S, alpha)) subseteq E(S, alpha)})
says that if $v$ is an eigenvector of $S$ with eigenvalue $\alpha$,
and if we decompose $v$ into eigenvectors of $T$, as in (\ref{v =
sum_{j=1}^k v_j}), then the components of this decomposition are
eigenvectors of $S$ with eigenvalue $\alpha$.

	We now leave the rest of the argument as an exercise.

	Let us consider a variant of this, where $V$ is an inner
product space, and $\mathcal{F}$ consists of self-adjoint operators
(which are then automatically diagonalizable).  In this case one can
get a simultaneous diagonalization by an orthogonal basis.  This can
be obtained using the orthogonality of the individual eigenspaces
involved.  There are also slightly different arguments that one can
consider, where one goes from (\ref{S(E(T, lambda)) subseteq E(T,
lambda)}) to saying that the restriction of an $S \in \mathcal{F}$ to
an eigenspace of $T$ is self-adjoint as an operator on that subspace.

	In the context of complex inner product spaces, one can look
at commuting families of operators $\mathcal{F}$ such that $T^*$ lies
in $\mathcal{F}$ whenever $T$ lies in $\mathcal{F}$, and again get a
simultaneous diagonalization by an orthogonal basis.  This situation
can be reduced to the previous one, by replacing $\mathcal{F}$ with
the collection of self-adjoint operators which arise as $T + T^*$ or
$i \, (T - T^*)$ for $T$ in $\mathcal{F}$.

\chapter{Linear operators between inner product spaces}

\section{Preliminary remarks}

	Suppose that $V_1$, $V_2$ are finite-dimensional vector
spaces, both real or both complex, and that $\langle \cdot, \cdot
\rangle_1$ and $\langle \cdot, \cdot \rangle_2$ are inner products
on $V_1$ and $V_2$, respectively.  If $T$ is a linear transformation
from $V_1$ to $V_2$, then there is a unique linear transformation
$T^* : V_2 \to V_1$ such that
\begin{equation}
	\langle T(v), w \rangle_2 = \langle v, T^*(w) \rangle_1
\end{equation}
for all $v \in V_1$ and $w \in V_2$, for essentially the same reasons
as before.  This is again called the \emph{adjoint}\index{adjoint (of
a linear transformation)} of $T$.

	Let us make the standing assumptions in this and the next two
sections that $(V_1, \langle \cdot, \cdot \rangle_1)$ and $(V_2,
\langle \cdot, \cdot \rangle_2)$ are as above, and that $\|\cdot\|_1$
and $\|\cdot\|_2$ denote the norms on $V_1$ and $V_2$ associated to
these inner products.  For $a, b = 1, 2$, if $S : V_a \to V_b$ is a
linear mapping, then $\|S\|_{op,ab}$ denotes the operator norm of $S$
defined using $\|\cdot\|_a$ on the domain and $\|\cdot\|_b$ on the
range.  This can be characterized in terms of inner products by
\begin{equation}
\label{||S||_{op,ab} in terms of inner products}
	\|S\|_{op,ab}  =  \sup \{ |\langle S(v), w \rangle_b| : 
			     v, w \in V, \|v\|_a = \|w\|_b = 1\},
\end{equation}
as in (\ref{formula for ||T||_{op} in terms of the inner product}) in
Section \ref{Inner product spaces (and spectral theory)}.

	This characterization of the operator norm implies that
\begin{equation}
\label{||T||_{op,12} = ||T^*||_{op,21}}
	\|T\|_{op,12} = \|T^*\|_{op,21}.
\end{equation}
Also, $(T^*)^* = T$ for all $T : V_1 \to V_2$, and sums, scalar
multiples, and compositions of adjoints behave in the same manner as
before, with the obvious allowances for the changes in the domains and
ranges.  The $C^*$-identities\index{C^*-identity@$C^*$-identity}
\begin{equation}
\label{||T^* T||_{op,11}  = ||T T^*||_{op,22} = ||T||_{op,12}^2}
	\|T^* \, T\|_{op,11} = \|T \, T^*\|_{op,22} = \|T\|_{op,12}^2
\end{equation}
can be verified as well.

\section{Schmidt decompositions}

	If $(V_1, \langle \cdot, \cdot \rangle_1)$ is different from
$(V_2, \langle \cdot, \cdot \rangle_2)$, then it does not make sense
to talk about an operator $T : V_1 \to V_2$ being self-adjoint, or
normal, or diagonalized in an orthonormal basis.  However, we do get
the linear mappings
\begin{equation}
	T^* \, T : V_1 \to V_1, \quad T \, T^* : V_2 \to V_2
\end{equation}
(with the ranges the same as the domains), and $T^* \, T$ is
self-adjoint as a linear mapping on $V_1$ equipped with $\langle
\cdot, \cdot \rangle_1$, while $T \, T^*$ is self-adjoint as a linear
mapping on $V_2$ equipped with $\langle \cdot, \cdot \rangle_2$.  We
can use this to derive the following \emph{Schmidt
decomposition}\index{Schmidt decompositions} for arbitrary linear
mappings from $V_1$ to $V_2$.

\beginproposition 
Let $T$ be a linear operator from $V_1$ to $V_2$.  There exists an
orthonormal basis $\{u_j\}_{j=1}^m$ of $V_1$ and an orthogonal set of
vectors $\{w_j\}_{j=1}^m$ in $V_2$ such that
\begin{equation}
\label{T(v) = sum_{j=1}^m langle v, u_j rangle_1 w_j}
	T(v) = \sum_{j=1}^m \langle v, u_j \rangle_1 \, w_j
\end{equation}
for all $v \in V_1$.
\end{proposition}

	Here the $w_j$'s can be equal to $0$, as when the kernel of
$T$ has positive dimension.  In particular, this necessarily happens
when the dimension of $V_2$ is less than the dimension of $V_1$.

	To prove the proposition, let $\{u_j\}_{j=1}^m$ be an
orthonormal basis of $V_1$ consisting of eigenvectors of $T^* \, T$.
Such a basis exists, since $T^* \, T : V_1 \to V_1$ is self-adjoint.
Set $w_j = T(u_j)$ for each $j$, so that (\ref{T(v) = sum_{j=1}^m
langle v, u_j rangle_1 w_j}) holds automatically.  It is not hard to
verify that the $w_j$'s are orthogonal to each other in this situation.

	Note that $T(u_j)$ is an eigenvector of $T \, T^*$ for each
$j$, with the same eigenvalue as $u_j$ has for $T^* \, T$.

\section{The Hilbert-Schmidt norm}
\index{Hilbert-Schmidt norm}

	Let $T$ be a linear mapping from $V_1$ to $V_2$.  Suppose that
$\{a_j\}_{j=1}^m$ and $\{b_k\}_{k=1}^n$ are orthonormal
bases for $V_1$ and $V_2$, respectively.  Thus
\begin{equation}
	v = \sum_{j=1}^m \langle v, a_j \rangle_1 \, a_j, \qquad
	  w = \sum_{k=1}^n \langle w, b_k \rangle_2 \, b_k
\end{equation}
for all $v \in V_1$ and $w \in V_2$, and we can express $T$ as
\begin{equation}
	T(v) = \sum_{j=1}^m \sum_{k=1}^n \langle v, a_j \rangle_1 \,
		\langle T(a_j), b_k \rangle_2 \, b_k.
\end{equation}

	Consider the sum
\begin{equation}
\label{sum_{j=1}^m sum_{k=1}^n |langle T(a_j), b_k rangle_2|^2}
	\sum_{j=1}^m \sum_{k=1}^n |\langle T(a_j), b_k \rangle_2|^2.
\end{equation}
By definition of the adjoint $T^*$, this is the same as
\begin{equation}
\label{sum_{j=1}^m sum_{k=1}^n |langle a_j, T^*(b_k) rangle_1|^2}
	\sum_{j=1}^m \sum_{k=1}^n |\langle a_j, T^*(b_k) \rangle_1|^2.
\end{equation}
Because $\{a_j\}_{j=1}^n$ and $\{b_k\}_{k=1}^n$ are orthonormal bases
for $V$, these sums are also equal to
\begin{equation}
\label{sum_{j=1}^m ||T(a_j)||_2^2}
	\sum_{j=1}^m \|T(a_j)\|_2^2
\end{equation}
and to
\begin{equation}
\label{sum_{k=1}^n ||T^*(b_k)||_1^2}
	\sum_{k=1}^n \|T^*(b_k)\|_1^2.
\end{equation}

	The \emph{Hilbert-Schmidt norm} of $T$, denoted $\|T\|_{HS}$,
is defined to be the square root of (\ref{sum_{j=1}^m sum_{k=1}^n
|langle T(a_j), b_k rangle_2|^2}).  This can also be described as the
$\|\cdot \|_2$ norm in Section \ref{definitions, etc. (normed vector
spaces)}, applied to the $n^2$ matrix entries $\langle T(a_j), b_k
\rangle_2$ of $T$ with respect to these choices of bases.  The
Hilbert-Schmidt norm does not depend on the particular choices of the
bases, however; namely, the equality with (\ref{sum_{j=1}^m
||T(a_j)||_2^2}) shows that the Hilbert-Schmidt norm does not depend on
the choice of $\{b_k\}_{k=1}^n$, and similarly the equality with
(\ref{sum_{k=1}^n ||T^*(b_k)||_1^2}) shows that the Hilbert-Schmidt norm
does not depend on the choice of $\{a_j\}_{j=1}^m$.

	The equality between (\ref{sum_{j=1}^m sum_{k=1}^n |langle
T(a_j), b_k rangle_2|^2}) and (\ref{sum_{j=1}^m sum_{k=1}^n |langle
a_j, T^*(b_k) rangle_1|^2}) implies that $\|T^*\|_{HS} = \|T\|_{HS}$.
More precisely, one can rewrite (\ref{sum_{j=1}^m sum_{k=1}^n
|langle a_j, T^*(b_k) rangle_1|^2}) as
\begin{equation}
	\sum_{j=1}^m \sum_{k=1}^n |\langle T^*(b_k), a_j \rangle_1|^2,
\end{equation}
and this is the same as (\ref{sum_{j=1}^m sum_{k=1}^n |langle T(a_j),
b_k rangle_2|^2}) with $T$ replaced by $T^*$, and with the roles of
$V_1$ and $V_2$ exchanged (and exchanging the roles of
$\{a_j\}_{j=1}^m$ and $\{b_k\}_{k=1}^n$ with them).

	Now suppose that we have two other inner product spaces $(V_0,
\langle \cdot, \cdot \rangle_0)$ and $(V_3, \langle \cdot, \cdot
\rangle_3)$, which are real if $V_1$, $V_2$ are real, and complex if
$V_1$, $V_2$ are complex.  If $A : V_0 \to V_1$ and $C : V_2 \to V_3$
are arbitrary linear mappings, then we have that
\begin{equation}
	\|C_1 \, T \|_{HS} \le \|C_1\|_{op,01} \, \|T\|_{HS}, \quad
		\|T \, C_2 \|_{HS} \le \|C_2\|_{op,23} \, \|T\|_{HS}
\end{equation}
(where we extend our earlier notation to mappings from $V_p$ to $V_q$
for $p, q = 0, 1, 2, 3$).  This follows from the formulas
(\ref{sum_{j=1}^m ||T(a_j)||_2^2}) and (\ref{sum_{k=1}^n
||T^*(b_k)||_1^2}) for the square of the Hilbert-Schmidt norm (allowing
also for linear mappings from $V_0$ to $V_2$ and from $V_1$ to $V_3$).

	For any linear mapping $T : V_1 \to V_2$ we have that
\begin{equation}
\label{||T||_{op,12} le ||T||_{HS} le min(m,n)^{1/2} ||T||_{op,12}}
	\|T\|_{op,12} \le \|T\|_{HS} \le \min(m,n)^{1/2} \, \|T\|_{op,12}.
\end{equation}
(Note that $m$ and $n$ are the dimensions of $V_1$ and $V_2$,
respectively.)  This follows easily from the equality of $\|T\|_{HS}$
with (\ref{sum_{j=1}^m ||T(a_j)||_2^2}) and (\ref{sum_{k=1}^n
||T^*(b_k)||_1^2}), and the fact that $\|T^*\|_{op,21} =
\|T\|_{op,12}$.  The first inequality in (\ref{||T||_{op,12} le
||T||_{HS} le min(m,n)^{1/2} ||T||_{op,12}}) becomes an equality when
$T$ has rank $1$, and equality holds in the second inequality for some
(nonzero) operators as well, e.g., when $V_1$ and $V_2$ are the same
and $T$ is a multiple of the identity.

\section{A numerical feature}

\beginlemma
Suppose that $V$ is a complex vector space with inner product $\langle
\cdot, \cdot \rangle$, and that $T$ is a linear operator on $V$.  Then
$T$ is self-adjoint if and only if the numbers
\begin{equation}
\label{langle T(v), v rangle}
	\langle T(v), v \rangle
\end{equation}
are real for all $v \in V$.
\end{lemma}

	If $T$ is self-adjoint, then we have that
\begin{equation}
	\overline{\langle T(v), v \rangle} = \langle v, T(v) \rangle
		= \langle v, T^*(v) \rangle = \langle T(v), v \rangle,
\end{equation}
by the basic properties of the inner product and the adjoint of a
linear transformation.  Thus the numbers (\ref{langle T(v), v rangle})
are real in this case.

	Conversely, suppose that the numbers (\ref{langle T(v), v
rangle}) are all real.  For any linear operator $T$ on $V$, we can
write
\begin{equation}
	T = A + i B
\end{equation}
where $A$ and $B$ are self-adjoint, by taking $A = (T + T^*)/2$ and
$B = (T - T^*)/(2i)$.  If the numbers (\ref{langle T(v), v rangle})
are all real, then the same is true of the numbers
\begin{equation}
	i \langle B(v), v \rangle,
\end{equation}
since $\langle A(v), v \rangle \in {\bf R}$ for all $v \in V$ by the
first part of the proof.  On the other hand, $\langle B(v), v \rangle
\in {\bf R}$ for all $v \in V$ as well, since $B$ is self-adjoint,
and hence
\begin{equation}
	\langle B(v), v \rangle = 0.
\end{equation}
This last implies that $B = 0$, using the diagonalizability of $B$.
Thus $T = A$, so that $T$ is self-adjoint.  This completes the proof
of the lemma.

\section{Numerical range}
\label{Numerical range}
\index{numerical range}

	Let $(V, \langle \cdot, \cdot \rangle)$ be a
(finite-dimensional) inner product space, and let $T$ be a linear
operator on $V$.  Consider the set
\begin{equation}
	W(T) = \{\langle T(v), v \rangle : v \in V, \ \|v\| = 1\}.
\end{equation}
This set is compact, for standard reasons of continuity and the
compactness of the unit sphere in $V$.

	Suppose first that $V$ is a real vector space, so that
$W(T)$ is a subset of the real line.  If we set
\begin{equation}
	A_1 = \frac{T + T^*}{2} \hbox{ and }
		A_2 = \frac{T - T^*}{2},
\end{equation}
then $T = A_1 + A_2$, $A_1^* = A_1$, and $A_2^* = - A_2$, i.e., $A_1$
is symmetric and $A_2$ is antisymmetric.  For $A_2$ we have that
\begin{equation}
	\langle A_2(v), v \rangle = - \langle v, A_2(v) \rangle
			= - \langle A_2(v), v \rangle,
\end{equation}
so that $\langle A_2(v), v \rangle = 0$ for all $v \in V$.  Thus
\begin{equation}
	W(T) = W(A_1).
\end{equation}
Because $A_1$ is symmetric, it can be diagonalized in an orthonormal
basis.  Using this one can check that
\begin{equation}
	W(A_1) \hbox{ is convex},
\end{equation}
and in fact that $W(A_1)$ is equal to the convex hull of the eigenvalues
of $A_1$.

	Alternatively, one might observe that $W(T)$ is connected,
because the image of a connected set under a continuous mapping is
connected.  The unit sphere in $V$ is connected as long as the dimension
of $V$ is at least $2$, and when the dimension of $V$ is $1$, $W(T)$
consists of a single point anyway.  It is also nice to observe that
\begin{equation}
	\langle T(v), v \rangle = \langle T(-v), -v \rangle,
\end{equation}
so that $v$ and $-v$ contribute the same value to the numerical range.
Of course connected subsets of ${\bf R}$ are automatically convex.

	Now let us turn to the case where $V$ is a complex vector
space, and $W(T)$ is a subset of ${\bf C}$.  It is no longer true that
the numerical range of $T$ is equal to the numerical range of the
self-adjoint part of $T$, which is to say that the anti-self-adjoint
part of $T$ can contribute to $W(T)$.

	If $T$ is normal, so that the self-adjoint and
anti-self-adjoint parts of $T$ commute, then $T$ can be diagonalized
in an orthonormal basis of $V$, as mentioned in Section \ref{Inner
product spaces (and spectral theory)}.  As in the real case, it is not
hard to use this to verify that $W(T)$ is equal to the convex hull of
the spectrum of $T$.

	In fact, $W(T)$ is convex for any linear operator $T$ on $V$.
See Chapter 17 of \cite{Halmos2}.

\chapter{Subspaces and quotient spaces}

	As before, it will be convenient to make the standing assumption
in this chapter that all vector spaces are finite-dimensional.

\section{Linear algebra}
\label{Linear algebra}

	Let $V$ be a vector space (real or complex), and let $W$ be a
subspace of $V$.  There is a standard notion of the \emph{quotient $V
/ W$ of $V$ by $W$},\index{quotient of a vector space} which is
obtained by identifying pairs of points in $V$ whose difference lies
in $W$.  These identifications are compatible with addition and scalar
multiplication of vectors, so that the quotient $V / W$ is a vector
space in a natural way.  

	One also gets the \emph{quotient mapping}\index{quotient mapping}
$q : V \to V / W$, which reflects the identifications.  This is a
linear mapping from $V$ \emph{onto} $V / W$ whose kernel is equal to $W$.

	Since $V / W$ is a vector space, it has a dual space $(V /
W)^*$, and the dual space can be identified with a subspace of the
dual $V^*$ of $V$ in a natural way.  Specifically, every linear
functional on $V / W$ determines a linear functional on $V$, by
composition with the quotient mapping $q$.  The linear functionals on
$V$ which arise in this manner are exactly the ones that vanish on $W$.

	In general, if $X$ is a subspace of $V$, one sometimes writes
$X^\perp$\index{$X^\perp$ (subspace of $V^*$ when $X \subseteq V$)}
for the subspace of $V^*$ consisting of linear functionals $\lambda$
on $V$ which are equal to $0$ on $X$.  Since $V^*$ is a vector space
in its own right, one can consider the subspace $(X^\perp)^\perp$ of
the second dual $V^{**}$.  The second dual $V^{**}$ is isomorphic to
$V$ in a natural way, as in Section \ref{Second duals}, and it is easy
to check that $(X^\perp)^\perp$ corresponds exactly to $X$ under this
isomorphism.  The previous statement about $(V / W)^*$ can be
rephrased as saying that $(V / W)^*$ is isomorphic in a natural way to
$W^\perp \subseteq V^*$.

	Similarly, if $X$ is any vector subspace of $V$, then every
linear functional on $V$ restricts to a linear functional on $X$.
Every linear functional on $X$ arises in this manner, and two linear
functionals $\lambda_1$, $\lambda_2$ on $V$ induce the same linear
functional on $X$ when $\lambda_1 - \lambda_2 \in X^\perp$.  In this
way, $X^*$ is isomorphic in a natural way to $V^* / X^\lambda$.

	This second scenario is essentially equivalent to the first
one, with $X \subseteq V$ corresponding to $W^\perp \subseteq V^*$,
and using the isomorphism between a vector space and its second dual
and between $X$ and $(X^\perp)^\perp$.

\beginremark 
\label{complementary subspaces and quotients}
{\rm A vector subspace $U$ of $V$ is
\emph{complementary}\index{complementary subspaces} to $W$ (Definition
\ref{def of projections and complementary subspaces} in Section
\ref{Projections}) if and only if the restriction of the quotient
mapping $q : V \to V / W$ to $U$ is one-to-one and maps $U$ onto $V /
W$.  }
\end{remark}

\section{Quotient spaces and norms}
\label{Quotient spaces and norms}

	Now suppose that $V$ is a vector space equipped with a norm
$\|\cdot\|$.  Let $W$ be a vector subspace of $V$, and consider the
quotient $V / W$ and the associated quotient mapping $q : V \to V /
W$.  The \emph{quotient norm}\index{quotient norm} $\|\cdot\|_{Q}$ on
$V / W$ can be defined as
\begin{equation}
\label{||z||_Q = inf {||v|| : v in V, q(v) = z}}
	\|z\|_Q = \inf \{\|v\| : v \in V, \, q(v) = z\}.
\end{equation}
This is the same as
\begin{equation}
	\|q(y)\|_Q = \inf \{\|y+w\| : w \in W\}.
\end{equation}
It is not hard to check that $\|\cdot\|_Q$ does indeed define a norm
on $V / W$.

	Using the norm $\|\cdot\|$ on $V$ and the norm $\|\cdot\|_Q$
on $V / W$ one can define the operator norm for any linear transformation
from $V$ to $V / W$.  The operator norm of the quotient mapping
$q : V \to V / W$ is equal to $1$ (except in the degenerate
case where $W = V$).

	Because of our standing assumption that $V$ be
finite-dimensional, a minimizing $v \in V$ for the infimum in
(\ref{||z||_Q = inf {||v|| : v in V, q(v) = z}}) exists.  It is not
unique, in general, but this is the case under assumptions of strict
convexity.  See Section \ref{Uniqueness of points of minimal
distance}.

	Let $U$ be a subspace of $V$ which is complementary to $W$,
and let $P$ be the projection of $V$ onto $U$ with kernel $W$, as in
Section \ref{Projections}.  Let us verify that
\begin{equation}
\label{||P||_{op}^{-1} ||u|| le ||q(u)||_Q le ||u|| for all u in U}
	\|P\|_{op}^{-1} \, \|u\| \le \|q(u)\|_Q \le \|u\|
		\qquad\hbox{for all } u \in U,
\end{equation}
where $\|P\|_{op}$ denotes the operator norm of $P$ as a linear
operator on $V$ (with respect to the norm $\|\cdot\|$ on $V$).  The
second inequality in (\ref{||P||_{op}^{-1} ||u|| le ||q(u)||_Q le
||u|| for all u in U}) holds automatically from the definition of the
quotient norm.  As for the first inequality, if $w$ is any element
of $W$, then $P(u+w) = u$, and hence
\begin{equation}
	\|u\| \le \|P\|_{op} \, \|u+w\|.
\end{equation}
One can take the infimum of the right side over $w \in W$ to get
$\|u\| \le \|P\|_{op} \, \|q(u)\|_Q$, as desired.

	In particular, if $\|P\|_{op}$ happens to be equal to $1$,
then the restriction of the quotient mapping $q$ to $U$ defines an
isometry from $U$ equipped with the norm $\|\cdot\|$ to $V / W$
equipped with the quotient norm $\|\cdot\|_Q$.  This includes the case
where the norm $\|\cdot\|$ on $V$ comes from an inner product and $U$
is the orthogonal complement\index{orthogonal complement} of $W$.

	Let us consider the quotient norm $\|\cdot\|_Q$ in terms of
duality.  As in Section \ref{Linear algebra}, the dual of $V / W$ as
a vector space can be identified with $W^\perp \subseteq V^*$.  Let
$\|\cdot\|^*$ be the dual norm of $\|\cdot\|$ on $V^*$, which can
be restricted to $W^\perp$.  Under the isomorphism between $(V / W)^*$
and $W^\perp$, the dual of the quotient norm $\|\cdot\|_Q$ corresponds
exactly to $\|\cdot\|^*$ on $W^\perp$.  Indeed, if $\lambda \in W^\perp$,
then $\lambda(w) = 0$ when $w \in W$, and
\begin{equation}
	|\lambda(v)| \le \|\lambda\|^* \, \|v\|
		\qquad\hbox{for all } v \in V,
\end{equation}
by the definition of the dual norm.  We can combine these two pieces of
information to obtain
\begin{equation}
	|\lambda(v)| \le \|\lambda\|^* \, \inf \{\|v + w\| : w \in W\}
		\qquad\hbox{for all } v \in V.
\end{equation}
This is the same as saying that the linear functional on $V / W$
induced by $\lambda$ has norm less than or equal to $\|\lambda\|^*$
with respect to the norm $\|\cdot\|_Q$ on $V / W$.  In the other
direction, if the linear functional on $V / W$ induced by $\lambda \in
W^\perp$ has norm less than or equal to $k$ with respect to
$\|\cdot\|_Q$, then
\begin{equation}
	|\lambda(v)| \le k \, \inf \{\|v+w\| : w \in W\}
		\qquad\hbox{for all } v \in V.
\end{equation}
This implies trivially that
\begin{equation}
	|\lambda(v)| \le k \, \|v\|
		\qquad\hbox{for all } v \in V,
\end{equation}
i.e., $\|\lambda\|^* \le k$.  Hence $\|\lambda\|^*$ is less than or
equal to the norm of the linear functional on $V / W$ induced by $\lambda$
with respect to the norm $\|\cdot\|_Q$ on $V / W$.  This shows that the
dual of the quotient norm $\|\cdot\|_Q$ corresponds exactly to the dual
norm $\|\cdot\|^*$ restricted to $W^\perp \subseteq V^*$.

	It is instructive to look directly at the case of subspaces as
well, even though this could be derived from the preceding discussion
using previous results about second duals.  Specifically, suppose that
$X$ is a vector subspace of $V$, to which we can restrict the norm
$\|\cdot\|$.  This leads to a norm on the dual space of $X$.  As a
vector space, the dual of $X$ is isomorphic in a natural way to $V^* /
X^\perp$, as discussed in Section \ref{Linear algebra}.  Under this
isomorphism, the dual norm of $X$ corresponds exactly to the quotient
norm on $V^* / X^\perp$ defined using the dual norm $\|\cdot\|^*$ on
$V^*$.  To see this, it is not hard to check that the dual norm of $X$
is less than or equal to the quotient norm on $V^* / X^\perp$ defined
using $\|\cdot\|^*$, just by the various definitions.  The other
inequality is less obvious, and it amounts to the statement that a
linear functional on $X$ can be extended to a linear functional on all
of $V$ with the same norm.  This holds by Theorem \ref{extension
theorem}.

\section{Mappings between vector spaces}

	Suppose that $V_1$ and $V_2$ are vector spaces, both real
or both complex, and that $\|\cdot\|_1$ and $\|\cdot\|_2$ are norms
on $V_1$ and $V_2$, respectively.  Let $Y$ be a vector subspace of
$V_1$, and assume that we have a linear mapping $T : Y \to V_2$.
Under what conditions can we find an extension of $T$ to a linear
mapping from $V_1$ to $V_2$ which has the same operator norm as
the original mapping, or perhaps an operator norm which is not too
much larger?  Here the operator norms use $\|\cdot\|_1$ on the
domain and $\|\cdot\|_2$ on the range.

	In general an extension with the same norm does not exist,
and it is not so easy to find an extension whose norm is not too
much larger than the original mapping.  Let us mention a few basic
points related to this, however.

	Suppose that $P$ is a projection from $V_1$ onto $Y$, and let
$k$ be the operator norm of $P$ (relative to $\|\cdot\|_1$ in both the
domain and range).  One way to get an extension of $T$ is to simply
take the composition $T \circ P$.  The operator norm of $T \circ P$
(using $\|\cdot\|_1$ on the domain and $\|\cdot\|_2$ on the range) is
bounded by $k$ times the operator norm of $T$, as in (\ref{op norm UT
le op norm U x op norm T}) in Section \ref{Linear transformations and
norms}.

	This works quite well when the norm $\|\cdot\|_1$ on $V_1$
comes from an inner product, so that there is a projection of norm $1$
onto any subspace of $V_1$, namely, the orthogonal projection.

	Note that one might take $V_2 = Y$, $\|\cdot\|_2 =
\|\cdot\|_1$ on $Y$, and $T$ to be the identity on $Y$.  Extensions of
$T$ to mappings from $V_1$ to $V_2 = Y$ are then exactly projections
of $V_1$ onto $Y$.  Thus the original question becomes one about
norms of projections onto $Y$.

	In another vein, if $V_2$ is $1$-dimensional, then the
question reduces to that of extending linear functionals, and this is
possible while keeping the same norm because of Theorem \ref{extension
theorem}.  

	This can be carried over to the case where $V_2$ is ${\bf
R}^n$ or ${\bf C}^n$ for some $n$ and one employs the norm
$\|\cdot\|_{\infty}$ from Section \ref{definitions, etc. (normed
vector spaces)}.  A linear mapping from a vector space $V_1$ into
${\bf R}^n$ or ${\bf C}^n$ (according to whether $V_1$ is real or
complex) can be described in terms of $n$ linear functionals on
$V_1$, corresponding to the $n$ coordinates for ${\bf R}^n$ or
${\bf C}^n$.  When we employ the norm $\|\cdot\|_{\infty}$ on
${\bf R}^n$ or ${\bf C}^n$, then the norm of a mapping from $V_1$
to ${\bf R}^n$ or ${\bf C}^n$ is the same as the maximum of the
norms of these $n$ linear functionals on $V_1$ (using the norm
that we have on $V_1$).  In this way the extension problem reduces
to its counterpart for linear functionals.

	Here is a ``dual'' version of the question at the beginning of
the section.  Suppose that $Z_1$ and $Z_2$ are vector spaces, both
real or both complex, and equipped with norms.  Assume also that we
are given a subspace $Y$ of $Z_2$, which we can use to define the
quotient $Z_2 / Y$.  Let $q : Z_2 \to Z_2 / Y$ be the corresponding
quotient mapping.  If $A$ is a linear mapping from $Z_1$ to $Z_2 / Y$,
under what conditions can we ``lift'' $A$ to a linear mapping
$\widehat{A} : Z_1 \to Z_2$ such that $A = q \circ \widehat{A}$ and
the norm of $\widehat{A}$ is the same as the norm of $A$, or not too
much larger?  For the norm of $A$, we use the quotient norm on $Z_2 /
Y$ coming from the given norm on $Z_2$.

	There are natural duals of the remarks above for the previous
question.  If $U$ is a subspace of $Z_2$ which is complementary to
$Y$, then the restriction of $q$ to $U$ is one-to-one and maps $U$
onto $Z_2 / Y$.  Thus we can obtain liftings $\widehat{A}$ of $A$
which take values in $A$.  There is then the question of the norm of
the lifting, which can be analyzed in terms of projections, through
the remarks in Section \ref{Quotient spaces and norms}.  If the norm
on $Z_2$ comes from an inner product, then one can take $U$ to be the
orthogonal complement of $Y$, and the lifting to $U$ has the same norm
as the original mapping $A$.

	If $Z_1$ has dimension $1$, then it is easy to do the lifting,
while keeping the norm fixed.  If instead $Z_1$ is ${\bf R}^n$ or
${\bf C}^n$, equipped with the norm $\|\cdot\|_1$ from Section
\ref{definitions, etc. (normed vector spaces)} (so that $\|\cdot\|_1$
is not now just a generic name for a norm, but rather a very specific
one), then the argument for $1$-dimensional domains can be carried
over in a simple manner.  More precisely, let $e_j$, $1 \le j \le n$,
be the standard basis vectors for ${\bf R}^n$ or ${\bf C}^n$
(according to whether one is working with real or complex vector
spaces).  The idea is to first choose $\widehat{A}(e_j)$ for
each $j$ so that $q(\widehat{A}(e_j)) = A(e_j)$ and the norm
of $\widehat{A}(e_j)$ is equal to the (quotient) norm of $A(e_j)$
for each $j$.  Once this is done, $\widehat{A}$ is determined on the
whole domain by linearity.  Because of the specific norm that we have
on the domain, the norm of $\widehat{A}$ as a linear mapping will be the
same as the norm of $A$.

\chapter{Variation seminorms}

\section{Basic definitions}
\label{Basic definitions}

	Let ${\bf Z}$\index{$Z$@${\bf Z}$} denote the set of integers, and
let ${\bf Z}^n$ denote the set of $n$-tuples of positive integers (for
a given positive integer $n$).

	If $x$, $y$ are elements of ${\bf Z}^n$, then we shall say
that $x$ and $y$ are \emph{neighbors}\index{neighbors} if there
is an integer $i$, $1 \le i \le n$, such that $y_i - x_i = \pm 1$
and $y_j = x_j$ for all $j \ne i$, $1 \le j \le n$, where $x_l$ and
$y_l$ denotes the components of $x$ and $y$.  This is clearly symmetric
in $x$ and $y$.  We shall denote by $N(x)$\index{$N(x)$} the set of
neighbors of $x$.  Note that $x$ is not considered to be a neighbor
of itself.

	If $U$ is a subset of ${\bf Z}^n$, then we define $\Int U,
\partial U \subseteq U$ by
\begin{equation}
\label{def of Int U}
\index{$Int U$@$\Int U$}
	\Int U = \{x \in U : N(x) \subseteq U\}
\end{equation}
and
\begin{equation}
\label{def of partial U}
\index{$bU$@$\partial U$}
	\partial U = U \backslash \partial U.
\end{equation}
This notation is somewhat at odds with that in topology, but these
definitions will play an analogous role.  

	In this chapter we shall make the standing assumption that
\begin{equation}
\label{standing assumption}
	U \subseteq {\bf Z}^n \hbox{ is finite and } \Int U \ne \emptyset.
\end{equation}

	Let $f$ be a function on $U$.  By ``function'' we mean one
that is real-valued or complex-valued, so that we have the real vector
space of real-valued functions on $U$, or the complex vector space of
complex-valued functions on $U$.  Sometimes it will be helpful to
specify one or the other, but often it will not matter.

	For each real number $p \ge 1$, define the
$p$-variation\index{p-variation@$p$-variation} $V_p(f)$\index{$V_p(f)$}
of $f$ (on $U$) by
\begin{equation}
\label{def of V_p(f)}
	V_p(f) = 
    \biggl(\sum_{x \in \Int U} \sum_{y \in N(x)} |f(x) - f(y)|^p \biggr)^{1/p}.
\end{equation}
This is a \emph{seminorm}\index{seminorm} on the vector space of functions
on $U$ (real or complex), which means that it satisfies the same properties
as a norm except that $V_p(f)$ might be $0$ even if $f$ is not the zero
function.  In particular, $V_p(f) = 0$ whenever $f$ is a constant function.
The fact that $V_p(\cdot)$ satisfies the triangle inequality
\begin{equation}
	V_p(f + g) \le V_p(f) + V_p(g)
\end{equation}
for all functions $f$ and $g$ on $U$ is a consequence of Minkowski's
inequality for sums (using the same choice of $p$).

	As in the case of norms, we can apply the triangle inequality
twice to obtain that
\begin{equation}
	V_p(f) \le V_p(g) + V_p(f - g) \quad\hbox{and}\quad
		V_p(g) \le V_p(f) + V_p(f-g),
\end{equation}
and hence that
\begin{equation}
	|V_p(f) - V_p(g)| \le V_p(f - g)
\end{equation}
for all functions $f$ and $g$ on $U$.  One can use this to show that
$V_p(\cdot)$ is continuous on the vector space of functions on $U$,
because $V_p(\cdot)$ is bounded by a constant multiple of a standard
norm on functions.

	Observe that elements of $U$ which are not neighbors of elements
of $\Int U$ are not used in (\ref{def of V_p(f)}), so that the values of
$f$ at such points does not affect $V_p(f)$.  One might as well take out
these points from $U$.

	It will sometimes be of interest to restrict our attention
to functions $f$ on $U$ which are equal to $0$ on $\partial U$.
In this case $V_p(f)$ can be rewritten as
\begin{equation}
\label{V_p(f) when f = 0 on partial U}
	V_p(f) =
    \biggl( \sum_{(x,y) \in D(\Int U)} |f(x) - f(y)|^p \biggr)^{1/p},
\end{equation}
where
\begin{equation}
	D(\Int U) = \{ (x,y) \in \Int U \times \Int U :
			x, y \hbox{ are neighbors}\}.
\end{equation}
On the vector space of functions on $U$ that are $0$ on $\partial U$,
$V_p(f)$ is a norm, i.e., $V_p(f) = 0$ implies that $f$ is the zero
function on $U$.  (Exercise.)

\section{The $p = 2$ and $n = 1$, $p = 1$ cases}
\label{The p = 2 and n = 1, p = 1 cases}

	When $p = 2$, $V_2(f)^2$ is a quadratic form (Hermitian
quadratic form in the complex case), and this leads to special
features.  Let $\langle f, g \rangle$ be the standard inner product
for functions on $U$, defined by
\begin{equation}
	\langle f, g \rangle = \sum_{x \in U} f(x) \, g(x)
\end{equation}
if we are working with real-valued functions, and by
\begin{equation}
	\langle f, g \rangle = \sum_{x \in U} f(x) \, \overline{g(x)}
\end{equation}
in the setting of complex-valued functions.  In either situation,
there is a unique linear operator $A$ on functions on $U$ which is
self-adjoint and satisfies
\begin{equation}
	V_2(f) = \langle A(f), f \rangle
\end{equation}
for all $f$.  This is not hard to see (and is an instance of a
standard way of writing quadratic forms in general).

	It is also natural to consider the restriction to the subspace
of functions which vanish on $\partial U$.  There is again a linear
operator $A_0$ mapping this subspace to itself which is self-adjoint
(with respect to the restriction of the inner product to the subspace)
and satisfies
\begin{equation}
	V_2(f) = \langle A_0(f), f \rangle
\end{equation}
for all $f$ in the subspace.  This operator $A_0$ is the
\emph{compression} of $A$ to the subspace, which means that $A_0$ is
obtained by restricting $A$ to the subspace, and then composing it
with the orthogonal projection from the whole space to the subspace.

	We can describe $A_0$ in terms of a matrix as follows.  For
each $x \in U$, let $e_x$ denote the function on $U$ which is $1$
at $x$ and $0$ at all other points.  The functions $e_x$, $x \in U$,
define an orthonormal basis for the space of functions on $U$, and
the functions $e_x$, $x \in \Int U$, define an orthonormal basis for
the subspace of functions which vanish on $\partial U$.

	For each $x, y \in \Int U$, $x \ne y$, we have that
\begin{eqnarray}
	\langle A_0(e_x), e_y \rangle & = & -1 
	   \quad\hbox{if $x$, $y$ are neighbors}		\\
				    & = & 0
	   \quad\enspace \hbox{ if $x$, $y$ are not neighbors},		
								\nonumber 
\end{eqnarray}
and $\langle A_0(e_x), e_x \rangle$ is equal to the number of
neighbors of $x \in \Int U$ which also lie in $\Int U$.  (Exercise.
Note that this is symmetric in $x$ and $y$, as it should be, for
self-adjointness.)  The matrix for $A$ can be determined analogously.

	Now suppose that $n = 1$ and $p = 1$.  Let $a$, $b$ be
integers with $a \le b - 2$, and let $U$ denote the set of integers in
the interval $[a,b]$.  In the notation of Section \ref{Basic
definitions}, $\Int U$ is the set of integers in $(a, b)$, and
$\partial U = \{a, b\}$.

	If $f$ is a function on $U$, then 
\begin{equation}
	V_1(f) \ge |f(b) - f(a)|.
\end{equation}
Let us restrict our attention to functions which are real-valued.  If
$f$ is monotone increasing on $U$, then
\begin{equation}
	V_1(f) = f(b) - f(a),
\end{equation}
and
\begin{equation}
	V_1(f) = f(a) - f(b)
\end{equation}
when $f$ is monotone decreasing on $U$.  Conversely,
\begin{equation}
	V_1(f) = |f(b) - f(a)|
\end{equation}
implies that $f$ is monotone (increasing or decreasing), as one can
check.

\section{Minimization}
\label{Minimization}

	Suppose that $U \subseteq {\bf Z}^n$ is as in Section
\ref{Basic definitions}, and fix a real number $p \ge 1$.  Let $b$ be
a function on $\partial U$.  We shall be interested in functions $f$
on $U$ such that $f = b$ on $\partial U$ and $V_p(f)$ is as small as
possible.

	For the usual reasons of continuity and compactness,
minimizing functions $f$ exist.  More precisely, although the set of
all functions $f$ on $U$ which agree with $b$ on $\partial U$ is not
bounded, one only needs to consider a bounded subset for the
minimization.  In other words, it is enough to look at such functions
$f$ for which $V_p(f)$ is bounded, since we are trying to minimize it.
Because we are fixing the boundary values of the functions, bounding
$V_p(f)$ leads to a bounded set of functions.  Otherwise, there is
the problem that nonzero functions $h$ can have $V_p(h) = 0$, such as
nonzero constant functions.

	When $p = 1$, minimizers with prescribed boundary values like
this do not have to be unique.  Examples of this can be seen from the
second part of Section \ref{The p = 2 and n = 1, p = 1 cases}.

	Let us look at the question of uniqueness in a somewhat
general way.  Assume that $f_1$ and $f_2$ are two minimizers of
$V_p(\cdot)$ on $U$ for the same set of boundary values $b$.  In
particular, $V_p(f_1) = V_p(f_2)$.  For each $t \in [0,1]$, $t \, f_1
+ (1-t) \, f_2$ is a function on $U$ with boundary values $b$, and
\begin{equation}
	V_p(t \, f_1 + (1-t) \, f_2) 
		\le t \, V_p(f_1) + (1-t) \, V_p(f_2).
\end{equation}
On the other hand, $V_p(t \, f_1 + (1-t) \, f_2) \ge V_p(f_1) = V_p(f_2)$,
by the minimality of $f_1$ and $f_2$, and hence
\begin{equation}
\label{V_p(t f_1 + (1-t) f_2) = V_p(f_1) = V_p(f_2)}
	V_p(t \, f_1 + (1-t) \, f_2) = V_p(f_1) = V_p(f_2),
\end{equation}
and $t \, f_1 + (1-t) \, f_2$.  (This would work as well for minimizing
any convex function.)

	We can rewrite (\ref{V_p(t f_1 + (1-t) f_2) = V_p(f_1) =
V_p(f_2)}) as
\begin{eqnarray}
\label{equality for sums with t f_1 + (1-t) f_2, x, y}
\lefteqn{\sum_{x \in \Int U} \sum_{y \in N(x)}
  \bigl|t \, f_1(x) + (1-t) \, f_2(x) - t \, f_1(y) - (1-t) \, f_2(y) \bigr|^p}
								\\
	&& = \sum_{x \in \Int U} \sum_{y \in N(x)} |f_1(x) - f_1(y)|^p
	   = \sum_{x \in \Int U} \sum_{y \in N(x)} |f_2(x) - f_2(y)|^p.
							\nonumber
\end{eqnarray}
One the other hand,
\begin{eqnarray}
\label{inequalities for |t f_1(x) + (1-t) f_2(x) - t f_1(y) - (1-t) f_2(y)|^p}
\lefteqn{\bigl|t \, f_1(x) + (1-t) \, f_2(x) 
		- t \, f_1(y) - (1-t) \, f_2(y) \bigr|^p}		\\
   && \le \bigl(t \, |f_1(x) - f_1(y)| + (1-t) \, |f_2(x) - f_2(y)| \bigr)^p
							\nonumber \\
	&& \le t \, |f_1(x) - f_1(y)|^p + (1-t) \, |f_2(x) - f_2(y)|^p,
							\nonumber
\end{eqnarray}
by the triangle inequality and the monotonicity and convexity of
the function $u^p$ on $[0,\infty)$.  Because of the equality
\begin{eqnarray}
\lefteqn{\sum_{x \in \Int U} \sum_{y \in N(x)}
 \bigl|t \, f_1(x) + (1-t) \, f_2(x) - t \, f_1(y) - (1-t) \, f_2(y) \bigr|^p} 
								\\
	&& = \sum_{x \in \Int U} \sum_{y \in N(x)} 
	   \Bigl(t \, |f_1(x) - f_1(y)|^p + (1-t) \, |f_2(x) - f_2(y)|^p \Bigr)
							\nonumber
\end{eqnarray}
(resulting from (\ref{equality for sums with t f_1 + (1-t) f_2, x,
y})), we obtain that
\begin{eqnarray}
\label{equalities for |t f_1(x) + (1-t) f_2(x) - tf_1(y) - (1-t) f_2(y)|^p}
\lefteqn{\bigl|t \, f_1(x) + (1-t) \, f_2(x) 
		- t \, f_1(y) - (1-t) \, f_2(y)\bigr|^p}		\\
   && = \bigl(t \, |f_1(x) - f_1(y)| + (1-t) \, |f_2(x) - f_2(y)| \bigr)^p
							\nonumber \\
	&& = t \, |f_1(x) - f_1(y)|^p + (1-t) \, |f_2(x) - f_2(y)|^p
							\nonumber
\end{eqnarray}
for all $x \in \Int U$ and $y \in N(x)$, i.e., we have equality in
both inequalities in (\ref{inequalities for |t f_1(x) + (1-t) f_2(x) - t
f_1(y) - (1-t) f_2(y)|^p}).

	The first equality in (\ref{equalities for |t f_1(x) + (1-t)
f_2(x) - tf_1(y) - (1-t) f_2(y)|^p}) is equivalent to
\begin{eqnarray}
\lefteqn{\bigl|t \, f_1(x) + (1-t) \, f_2(x) 
		- t \, f_1(y) - (1-t) \, f_2(y)\bigr|}		\\
   && = t \, |f_1(x) - f_1(y)| + (1-t) \, |f_2(x) - f_2(y)|, 
							\nonumber 
\end{eqnarray}
which is to say that we have equality for the triangle inequality.
As a consequence, for each $x \in \Int U$ and each $y \in N(x)$
there is a nonzero real or complex number $w$ (depending on whether
we are working with real or complex-valued functions) such that
$f_1(x) - f_1(y)$ and $f_2(x) - f_2(y)$ are both multiples of
$w$ by nonnegative real numbers.

	Now suppose that $p > 1$.  Now the second inequality in
(\ref{equalities for |t f_1(x) + (1-t) f_2(x) - tf_1(y) - (1-t) f_2(y)|^p})
and the strict convexity of the function $u^p$ on $[0,\infty)$ imply that
\begin{equation}
	|f_1(x) - f_1(y)| = |f_2(x) - f_2(y)|
\end{equation}
for all $x \in \Int U$ and all $y \in N(x)$.  Therefore,
\begin{equation}
	f_1(x) - f_1(y) = f_2(x) - f_2(y)
\end{equation}
for all $x \in \Int U$ and $y \in N(x)$, because of the information
derived in the previous paragraph.

	From this it follows that $f_1 = f_2$ on all of $U$, since
$f_1$ and $f_2$ agree on $\partial U$ by assumption.  In short,
we have uniqueness for minimizers when $p > 1$.

\section{Truncations}

	Given real numbers $c$, $d$, consider the functions
$\tau_1(s)$ and $\tau_2(s)$ on ${\bf R}$ given by 
\begin{equation}
	\tau_1(s) = \max(s,c) \quad\hbox{and}\quad \tau_2(s) = \min(s,d).
\end{equation}
It is easy to check that
\begin{equation}
\label{|tau_1(s) - tau_1(t)|, |tau_2(s) - tau_2(t)| le |s-t|}
	|\tau_1(s) - \tau_1(t)|, \ |\tau_2(s) - \tau_2(t)| \le |s-t|
		\quad\hbox{for all } s, t \in {\bf R}
\end{equation}
and
\begin{eqnarray}
\label{|tau_1(s) - tau_1(t)| < |s-t|, |tau_2(u) - tau_2(v)| < |u-v|}
	&& |\tau_1(s) - \tau_1(t)| < |s-t|, \quad 
		|\tau_2(u) - \tau_2(v)| < |u-v|			\\
	&& \hbox{when $s$ or } t \in (-\infty, c) 
		\hbox{ and $u$ or } v \in (d, \infty).	\nonumber
\end{eqnarray}
Of course
\begin{equation}
	\tau_1(s) = s \hbox{ and } \tau_2(u) = u \hbox{ when }
		s \ge c, \ u \le d.
\end{equation}

	Let $U \subseteq {\bf Z}^n$ be as before, and let $f$ be a
real-valued function on $U$.  For each $p \ge 1$, we have that
\begin{equation}
	V_p(\tau_1 \circ f), \ V_p(\tau_2 \circ f) \le V_p(f)
\end{equation}
because of (\ref{|tau_1(s) - tau_1(t)|, |tau_2(s) - tau_2(t)| le
|s-t|}).  

	Fix a $p \ge 1$, and suppose that $b$ is a real-valued
function on $\partial U$, and that $f$ is a real-valued function on
$U$ which agrees with $b$ on $\partial U$ and for which $V_p(f)$ is as
small as possible.  If $b(x) \ge c$ for all $x \in \partial U$, then
$f(x) \ge c$ for all $x \in U$.  Indeed, if $f(x)$ were strictly less
than $c$ for any $x \in \Int U$, then we would have
\begin{equation}
	V_p(\tau_1 \circ f) < V_p(f),
\end{equation}
by (\ref{|tau_1(s) - tau_1(t)| < |s-t|, |tau_2(u) - tau_2(v)| <
|u-v|}).  On the other hand, $\tau_1 \circ f = \tau_1 \circ b = b$ on
$\partial U$, since $b(x) \ge c$ for all $x \in \partial U$.  This
shows that $f$ would not minimize $V_p$ among real-valued functions on
$U$ that agree with $b$ on $\partial U$.  Thus $f(x) \ge c$ for all $x
\in U$.

	Similarly, if $b(x) \le d$ for all $x \in \partial U$, then
we may conclude that $f(x) \le d$ for all $x \in U$.

	In the complex case, one can consider analogous ``truncation''
mappings as follows.  Let $H$ be a closed half-plane in ${\bf C}$, and
let $L$ be the line which is the boundary of $H$.  Define $\theta : {\bf
C} \to {\bf C}$ by taking $\theta(z) = z$ when $z \in H$, and $\theta(z)$
to be the orthogonal projection of $z$ into $L$ when $z \in {\bf C}
\backslash H$.  As before one has
\begin{equation}
  |\theta(z) - \theta(w)| \le |z-w| \quad\hbox{for all } z, w \in {\bf C}.
\end{equation}
Hence
\begin{equation}
	V_p(\theta \circ f) \le V_p(f)
\end{equation}
for any $p \ge 1$ and any complex-valued function $f$ on $U$.

	The strict inequality (\ref{|tau_1(s) - tau_1(t)|, |tau_2(s) -
tau_2(t)| le |s-t|}) does not work now in the same way as before,
because $|\theta(z) - \theta(w)| = |z-w|$ when $z, w \in {\bf C}
\backslash H$ lie in a line parallel to $L$.  However,
\begin{equation}
	|\theta(z) - \theta(w)| < |z-w|
\end{equation}
does hold when one of $z$, $w$ lies in $H$ and the other lies in ${\bf
C} \backslash H$.  This is enough to show that if $b$ is a
complex-valued function on $\partial U$ with values in $H$, and if $f$
is a complex function on $U$ which is equal to $b$ on $\partial U$ and
for which $V_p(f)$ is as small as possible (for some fixed $p \ge 1$),
then $f$ also takes values in $H$.  That is, if $f$ did not take
values in $H$, then one could consider $\theta \circ f$, which would
still agree with $b$ on $\partial U$.  In this case one can again get
$V_p(\theta \circ f) < V_p(f)$, contradicting the minimality of
$V_p(f)$, but one has to be slightly more careful than before.  (The
main point is that if there is an $x \in U$ such that $f(x) \in {\bf
C} \backslash H$ in these circumstances, then there is in fact an
$x \in \Int U$ such that $f(x) \in {\bf C} \backslash H$ and $f(y) \in H$
for some $y \in N(x)$.  For this $x$ and $y$ one obtains 
\begin{equation}
	|\theta(f(x)) - \theta(f(y))| < |f(x) - f(y)|
\end{equation}
as above, and hence that $V_p(\theta \circ f) < V_p(f)$.)

	Since this works for arbitrary half-planes $H$ in ${\bf C}$,
one can use this to show that the values of a minimizing function $f$
lie in the convex hull of the boundary values.

	Here is a variant of this type of argument.  Let $r$ be a positive
real number, and define $\sigma : {\bf C} \to {\bf C}$ by
\begin{equation}
	\sigma(z) = z \hbox{ when } |z| \le r,
		\quad \sigma(z) = r \, z / |z| \hbox{ when } |z| > r.
\end{equation}
One can again show that
\begin{equation}
	|\sigma(z) - \sigma(w)| \le |z-w| 
		\quad\hbox{for all } z, w \in {\bf C},
\end{equation}
and that
\begin{equation}
	|\sigma(z) - \sigma(w)| < |z-w|
		\quad\hbox{when $|z|$ or } |w| > r.
\end{equation}
(One may prefer to look at these inequalities in terms of
differentials of $\sigma$, rather than computing directly.)  As
before, $V_p(\sigma \circ f) \le V_p(f)$ for any $p \ge 1$ and any
complex-valued function $f$ on $U$, and if $f$ minimizes $V_p$ among
functions with prescribed boundary values $b$ on $\partial U$, and if
$|b(x)| \le r$ for all $x \in \partial U$, then we obtain that $|f(x)|
\le r$ for all $x \in U$.

\chapter{Groups}
\index{groups}

\section{General notions}

	Let $G$ be a finite group.  Thus $G$ is a finite set with a
distinguished element $e$, called the \emph{identity element}, and a
binary operation $(g,h) \mapsto g \, h$ from $G \times G$ into $G$,
such that the following three conditions are satisfied.  First, the
operation is \emph{associative}, which means that
\begin{equation}
	(g \, h) \, k = g \, (h \, k)
\end{equation}
for all $g$, $h$, and $k$ in $G$.  Second,
\begin{equation}
	e \, g = g \, e = g
\end{equation}
for all $g$ in $G$, which makes precise the sense in which $e$ is
the identity element.  Third, for each $g$ in $G$ there is a
$g^{-1}$ in $G$ such that
\begin{equation}
	g \, g^{-1} = g^{-1} \, g = e.
\end{equation}
It is easy to see that $g^{-1}$ is uniquely determined by $g$, and it
is called the inverse of $g$.

	Let $\mathcal{F}(G)$ denote the vector space of functions on
$G$.  More precisely, one can consider the real vector space of
real-valued functions on $G$, or the complex vector space of
complex-valued functions on $G$.  We may write $\mathcal{F}^r(G)$
for the former and $\mathcal{F}^c(G)$ for the latter, to specify
one or the other.

	For each $p \in [1,\infty]$, we can define a norm
$\|\cdot\|_p$ on $\mathcal{F}(G)$, as in Section \ref{definitions,
etc. (normed vector spaces)}.  That is, we set
\begin{equation}
	\|f\|_p = \Bigr(\sum_{h \in G} |f(h)|^p \Bigr)^{1/p}
\end{equation}
when $1 \le p < \infty$, and
\begin{equation}
	\|f\|_\infty = \max \{ |f(h)| : h \in G\}.
\end{equation}

	If $g$ is an element of $G$, define the corresponding 
\emph{left translation operator}\index{translation operators} $T_g :
\mathcal{F}(G) \to \mathcal{F}(G)$ by
\begin{equation}
	L_g(f)(h) = f(g^{-1} h)
\end{equation}
for any function $f(h)$ in $\mathcal{F}(G)$.  Notice that
\begin{equation}
	\|L_g(f)\|_p = \|f\|_p
\end{equation}
for all $g$ in $G$, for all $f$ in $\mathcal{F}(G)$, and all $p$, $1
\le p \le \infty$.  This uses the fact that $h \mapsto g^{-1} \, h$
defines a permutation on $G$ for each $g$ in $G$.

	In general, a norm $\|\cdot\|$ on $\mathcal{F}(G)$ is
said to be \emph{invariant under left translations} if 
\begin{equation}
	\|L_g(f)\| = \|f\|
\end{equation}
for all $f$ in $\mathcal{F}(G)$.

	Similarly, for each $g$ in $G$, define the corresponding
\emph{right translation operator} $R_g : \mathcal{F}(G) \to
\mathcal{F}(G)$ by
\begin{equation}
	R_g(f)(h) = f(h \, g)
\end{equation}
for every function $f(h)$ in $\mathcal{F}(G)$.  Again
\begin{equation}
	\|R_g(f)\|_p = \|f\|_p
\end{equation}
for all $g$ in $G$, for all $f$ in $\mathcal{F}(G)$, and all $p$, $1
\le p \le \infty$.  A norm $\|\cdot\|$ on $\mathcal{F}(G)$ is said
to be \emph{invariant under right translations} if
\begin{equation}
	\|R_g(f)\| = \|f\|
\end{equation}
for all $f$ in $\mathcal{F}(G)$.  

	Observe that
\begin{equation}
	L_g \circ L_k = L_{g \, k} \quad\hbox{and}\quad
		R_g \circ R_k = R_{g \, k}
\end{equation}
for all $g$ and $k$ in $G$.  Also,
\begin{equation}
	L_g \circ R_k = R_k \circ L_g
\end{equation}
for all $g$, $k$.

\section{Some operators on $\mathcal{F}(G)$}
\label{Some operators on mathcal{F}(G)}

	Let $a$ be a function on $G$, and consider the operators
$S_a$, $T_a$ on $\mathcal{F}(G)$ defined by
\begin{equation}
	S_a = \sum_{h \in G} a(h) \, L_h, \qquad
		T_a = \sum_{h \in G} a(h) \, R_h.
\end{equation}
In other words,
\begin{equation}
	S_a(f)(g) = \sum_{h \in G} a(h) \, f(h^{-1} g), \quad
		T_a(f)(g) = \sum_{h \in G} a(h) \, f(gh).
\end{equation}
For $1 \le p \le \infty$, we have that
\begin{eqnarray}
	\|S_a(f)\|_p = \Bigl\| \sum_{h \in G} a(h) \, L_h(f) \Bigr\|_p
		& \le & \sum_{h \in G} |a(h)| \, \|L_h(f)\|_p		\\
		& =   & \|a\|_1 \, \|f\|_p,			\nonumber
\end{eqnarray}
and, similarly,
\begin{equation}
	\|T_a(f)\|_p \le \|a\|_1 \, \|f\|_p.
\end{equation}
Thus, as operators on $\mathcal{F}(G)$ equipped with the norm
$\|\cdot\|_p$, $S_a$ and $T_a$ have norm less than or equal to
$\|a\|_1$.

	Now suppose that $p = 1$.  Let us write $\delta_e$ for the
function on $G$ such that $\delta_e(e) = 1$ and $\delta_e(h) = 0$ when
$h \ne e$.  Then
\begin{equation}
	S_a(\delta_e)(k) = a(k)
\end{equation}
and
\begin{equation}
	T_a(\delta_e)(k) = a(k^{-1}),
\end{equation}
as one can verify.  It follows that the operator norms of $S_a$ and
$T_a$ on $\mathcal{F}(G)$ equipped with the norm $\|\cdot\|_1$ are
equal to $\|a\|_1$, since $\|\delta_e\|_1 = 1$.

\beginexercise
{\rm 
The Hilbert-Schmidt norms of $S_a$ and $T_a$ on
$\mathcal{F}(G)$ with respect to the norm $\|\cdot\|_2$
are equal to $\Order(G)^{1/2} \, \|a\|_2$, where $\Order(G)$
is the number of elements of $G$.
}
\end{exercise}

\section{Commutative groups}

	In this section we assume that $G$ is a finite commutative
group, so that
\begin{equation}
	g \, h = h \, g
\end{equation}
for all $g$, $h$ in $G$.  We shall also make the standing assumption
that $\mathcal{F}(G)$ is the complex vector space of complex-valued
functions on $G$.

	The norm $\|\cdot\|_2$ on $\mathcal{F}(G)$ is associated
to the inner product
\begin{equation}
\label{inner product on mathcal{F}(G)}
	\langle f_1, f_2 \rangle
		= \sum_{h \in G} f_1(h) \, \overline{f_2(h)}.
\end{equation}
This inner product is preserved by the translation operators, which is
to say that they define unitary operators on $\mathcal{F}(G)$ equipped
with $\langle f_1, f_2 \rangle$.  

	As in Section \ref{Inner product spaces (and spectral
theory)}, a unitary operator on a complex inner product space can be
diagonalized in an orthogonal basis.  The translation operators are
all unitary on $\mathcal{F}(G)$ (with this choice of inner product),
and so they can each be diagonalized in an orthonormal basis.  A key
point now is that the translation operators all commute with each
other, since the group $G$ is commutative.  This implies that
\begin{eqnarray}
  && \hbox{there is an orthogonal basis of $\mathcal{F}(G)$ in which all the}
									   \\
  && \hbox{translation operators are simultaneously diagonalized,}
								\nonumber
\end{eqnarray}
as in Section \ref{Commuting families of operators}.  The key point is
that an eigenspace of one of the operators is invariant under the
other operators (and one can make use of the information that the
operators are unitary here to adjust the argument somewhat).

	What do the elements of such an orthogonal basis look like?
In other words, if $f$ is an element of $\mathcal{F}(G)$ which is an
eigenvector of all of the translation operators, then what can we say
about $f$?  One can show that $f$ is an eigenvector of all of the
translation operators if and only if there is a complex number $c$ and
a nonzero complex-valued function $\chi$ on $G$ such that $f = c \,
\chi$ and
\begin{equation}
\label{chi(g h) = chi(g) chi(h)}
	\chi(g \, h) = \chi(g) \, \chi(h)
\end{equation}
for all $g$ and $h$ in $G$.  (Exercise.  Observe also that $|\chi(g)|
= 1$ for all $g$.)

	As part of the story of eigenvectors of the translation
operators, one gets that distinct functions $\chi$ on $G$ which
satisfy (\ref{chi(g h) = chi(g) chi(h)}) are orthogonal in
$\mathcal{F}(G)$, using the inner product above.  It is a good
exercise to extract a direct proof of this statement.

	In this situation, operators of the form $S_a$ and $T_a$ from
the preceding section are diagonalized by the same basis as the
translation operators, since $S_a$ and $T_a$ are, after all, linear
combinations of the translation operators.  Specifically, if $\chi$ is
a nonzero complex-valued function on $G$ which satisfies (\ref{chi(g
h) = chi(g) chi(h)}), then
\begin{equation}
	S_a(\chi) = \Bigl(\sum_{h \in G} a(h) \, \chi(h^{-1})\Bigr) \, \chi,
\end{equation}
and there is an analogous formula for $T_a(\chi)$.  As a consequence,
\begin{equation}
	\|S_a\|_{op, 22} = 
   \max_{\chi} \Bigl|\sum_{h \in G} a(h) \, \chi(h^{-1})\Bigr|,
\end{equation}
where $\|\cdot\|_{op, 22}$ denotes the operator norm of a linear
operator on $\mathcal{F}(G)$ with respect to the norm $\|\cdot\|_2$ on
$\mathcal{F}(G)$, and the maximum is taken over all nonzero
complex-valued functions $\chi$ on $G$ which satisfy (\ref{chi(g h) =
chi(g) chi(h)}).  This is consistent with the fact that
\begin{equation}
	\|S_a\|_{op, 22} \le \|a\|_1,
\end{equation}
mentioned in Section \ref{Some operators on mathcal{F}(G)}, because
$|\chi(h)| = 1$ for all $h$ in $G$.  Note that $\|S_a\|_{op, 22}$ can
be strictly less than $\|a\|_1$.

\section{Special cases}

	Let $n$ be a positive integer, and let $G$ be the group
consisting of $0, 1, \ldots, n-1$ with the group operation given
by ``modular addition''.  In other words, if $j$, $k$ are two
integers with $0 \le j, k \le n-1$, then the group operation
simply adds $j$ and $k$ when $j + k \le n-1$, and otherwise
it adds $j$ and $k$ and then subtracts $n$ to get back within
the range from $0$ to $n-1$.  It is not hard to check that this
is a commutative group.

	In this case one can check that the functions $\chi$ on $G$
as in the previous section are of the form
\begin{equation}
	\chi(j) = \alpha^j,
\end{equation}
where $\alpha$ is a complex number which is an $n$th root of unity,
i.e., $\alpha^n = 1$.  This includes $\alpha = 1$, and $\alpha$'s
which might be $m$th roots of unity for an $m$ which is smaller than
$n$ (and which divides $n$).  There are exactly $n$ such roots of unity,
which fits with the fact that $\mathcal{F}(G)$ has dimension $n$.
These are versions of versions of ``complex exponentials''.

	Now suppose that $G$ is the set of $l$-tuples $x = (x_1, x_2,
\ldots, x_l)$ where each $x_j$ is either $1$ or $-1$.  Given two
such $l$-tuples $x$, $y$, define $x \, y$ by
\begin{equation}
	x \, y = (x_1 \, y_1, \, x_2 \, y_2, \, \ldots, \, x_l \, y_l),
\end{equation}
using ordinary multiplication in each of the coordinates.  This also
defines a commutative group.

	For $i = 1, 2, \ldots, l$, define $\rho_i$ on $G$ by
$\rho_i(x) = x_i$.  Thus $\rho_i(x \, y) = \rho_i(x) \, \rho(y)$
for all $x, y \in G$.  These are some of the functions that we want,
but not all.

	If $I$ is a subset of the set $\{1, 2, \ldots, l\}$, define
$W_I$ on $G$ by 
\begin{equation}
	W_I = \prod_{i \in I} \rho_i.
\end{equation}
We interpret this as $W_I(x) = 1$ for all $x$ in $G$ when $I =
\emptyset$.  These functions $W_I$ are all distinct, and all satisfy
$W_I(x \, y) = W_I(x) \, W_I(y)$.  They account for all functions
of this type.

	The functions $W_I$ are versions of ``Walsh functions''.
The functions $\rho_i$ are versions of ``Rademacher functions''.

\section{Groups of matrices}

	Let $n$ be a positive integer.  The group $GL(n, {\bf R})$
consists of all invertible $n \times n$ matrices with real entries,
where the group operation is of course matrix multiplication.

	One can identify an $n \times n$ real matrix with a linear
mapping on ${\bf R}^n$.  The group $GL(n, {\bf R})$ acts on the vector
space of polynomials on ${\bf R}^n$, through
\begin{equation}
	P(x) \mapsto P(T^{-1}(x)).
\end{equation}
That is, if $P(x)$ is a polynomial on ${\bf R}^n$ and $T$ is an
invertible linear mapping on ${\bf R}^n$, then the action of $T$ on
the polynomial $P(x)$ is defined to be $P(T^{-1}(x))$.  If $S$ is
another invertible linear mapping on ${\bf R}^n$, then the action of
$S \, T$ on $P(x)$ is $P((S \, T)^{-1}(x)) = P(T^{-1}(S^{-1}(x)))$,
which is the same as first taking the action of $T$ on $P(x)$ and then
taking the action of $S$ on the polynomial resulting from the first
step.

	If $d$ is a nonnegative integer, then a polynomial $P(x)$ on
${\bf R}^n$ is said to be \emph{homogeneous of degree $d$} if
\begin{equation}
	P(t \, x) = t^d \, P(x)
\end{equation}
for all real numbers $t$ and all $x \in {\bf R}^n$.  This is the same
as saying that one can write $P(x)$ as a linear combination of
monomials $x_1^{j_1} x_2^{j_2} \cdots x_n^{j_n}$, where $\sum_{l=1}^n
j_l = d$ for each term.  Every polynomial on ${\bf R}^n$ can be
written as a sum of homogeneous polynomials, where the homogeneous
components are unique.  The subspaces of homogeneous polynomials are
mapped to themselves by the action of $GL(n, {\bf R})$.

	Now let us consider the group $O(n)$, which is the subgroup of
$GL(n, {\bf R})$ consisting of matrices whose inverses are equal to
their transposes.  These matrices correspond to \emph{orthogonal}
linear transformations on ${\bf R}^n$, i.e., linear transformations
which preserve the standard inner product on ${\bf R}^n$ (as in
Section \ref{Inner product spaces (and spectral theory)}).  For these
matrices, the action on polynomials can be refined.  Indeed, the
polynomials $|x|^{2k} = (\sum_{l=1}^n x_l^2)^k$ are invariant under
the action of $O(n)$.  This action also preserves the class of
\emph{harmonic} polynomials, which are the polynomials $P(x)$ such
that
\begin{equation}
	\sum_{l=1}^n \frac{\partial^2}{\partial x_l^2} \, P(x) = 0.
\end{equation}
These two properties of the action of $O(n)$ on polynomials on ${\bf
R}^n$ are actually ``dual'' to each other in a natural way, as in
Section IV.2 in \cite{SW-book} and Section III.3 of \cite{St2}.

	In a somewhat different direction, suppose that we consider
$GL(m, {\bf C})$, which is the group of $m \times m$ invertible
matrices with complex entries.  These matrices can be identified
with invertible linear transformations on ${\bf C}^m$, and one can
also view them as defining invertible real linear transformations on
${\bf R}^{2m}$.  In this way one can think of $GL(m, {\bf C})$ as
a subgroup of $GL(2m, {\bf R})$.

	For the action on polynomials, it will be helpful to consider
polynomials with complex coefficients.  A ``holomorphic'' polynomial
on ${\bf C}^m$ is one that can be written as a linear combination of
products of the complex coordinate functions $z_1, \ldots, z_m$, and
these are preserved by the action of $GL(m, {\bf C})$.  Arbitrary
polynomials, which are not necessarily holomorphic, can be expressed
as a linear combination of products of the complex coordinate
functions $z_1, \ldots, z_m$ and their complex conjugates
$\overline{z_1}, \ldots, \overline{z_m}$.  If we identity ${\bf C}^m$
with ${\bf R}^{2m}$, as a real vector space, then these are the usual
polynomials on ${\bf R}^{2m}$ with complex coefficients, written in a
slightly different way.  Normally one might write these polynomials in
terms of linear combinations of products of the real and imaginary
parts of the $z_j$'s, and this is equivalent to taking linear
combinations of products of the $z_j$'s and their complex conjugates.

	The group $GL(2m, {\bf R})$ acts on the full space of (not
necessarily holomorphic) polynomials on ${\bf C}^m \approx {\bf
R}^{2m}$, and the subspaces of homogeneous polynomials of various
degrees are preserved by this action, as before.  For $GL(m, {\bf
C})$, the subspace of holomorphic polynomials is also preserved, and,
more generally, the subspace spanned by monomials which are products
of $a$ $z_j$'s and $b$ $\overline{z_k}$'s is preserved for each $a$
and $b$.

	Inside $GL(m, {\bf C})$ there is the subgroup $U(m)$ of
matrices whose inverse is given by their adjoint (the complex
conjugate of the transpose).  These are the matrices corresponding
to linear transformations on ${\bf C}^m$ that are unitary, which is to
say that they preserve the standard Hermitian inner product on ${\bf
C}^n$, as in Section \ref{Inner product spaces (and spectral theory)}.
If we view $GL(m, {\bf C})$ as a subgroup of $GL(2m, {\bf R})$, as
above, then $U(m)$ is the same as the intersection of $GL(m, {\bf C})$
with $O(2m)$ (as one can verify).

	These groups and their actions on polynomials on ${\bf R}^n$,
${\bf C}^m$ (and some other vector spaces) are fundamental in
mathematics.  In particular, noncommutativity shows up.  This is
reflected in the limited number of $1$-dimensional spaces of
polynomials which are preserved by the various actions.

\chapter{Some special families of functions}

\section{Rademacher functions}
\label{Rademacher functions}\index{Rademacher functions}

	Given a positive integer $j$, the \emph{$j$th Rademacher
function $r_j(x)$} is the function on the unit interval $[0,1)$
defined by
\begin{equation}
	r_j(x) = (-1)^i \quad\hbox{when } x \in [i \, 2^{-j}, (i+1) \, 2^{-j}),
\end{equation}
where $i$ runs through all nonnegative integers such that $i+1 \le 2^{-j}$.

	Here are some basic properties of $r_j(x)$:
\begin{eqnarray}
\label{|r_j(x)| = 1}
    	&& \hbox{$|r_j(x)| = 1$ for all $x \in [0,1)$;}
									\\
\label{r_j(x) is constant on ....}
    	&& \hbox{$r_j(x)$ is constant on dyadic subintervals}
									\\
	&& \hbox{of $[0,1)$ of length $2^{-j}$;}		\nonumber\\
\label{int_L r_j = 0 when |L| > 2^{-j}}
	&& \int_L r_j(x) \, dx = 0 \hbox{ for all dyadic subintervals}   \\
	&& \hbox{$L$ of $[0,1)$ of length greater than $2^{-j}$.}	
								\nonumber
\end{eqnarray}
The latter comes from the way that $r_j(x)$ alternates between $1$ and $-1$
on dyadic intervals of length $2^{-j}$.

\beginlemma
\label{integrals of products of r_j's}
Suppose that $j_1, j_2, \ldots, j_n$ are distinct positive integers.  Then
\begin{equation}
\label{the integral of the product of r_j's}
	\int_{[0,1)} r_{j_1}(x) \, r_{j_2}(x) \cdots r_{j_n}(x) \, dx = 0.
\end{equation}
\end{lemma}

	To see this, we may as well assume that $j_n$ is the largest
of the $j_l$'s.  This implies that each of the other $r_{j_l}$'s is
constant on dyadic subintervals of $[0,1)$ of length $2^{-j_n + 1}$,
as in (\ref{r_j(x) is constant on ....}).  However, the integral of
$r_{j_n}$ over any dyadic subinterval of $[0,1)$ of length $2^{-j_n +
1}$ is $0$, by (\ref{int_L r_j = 0 when |L| > 2^{-j}}), and the same
must be true for the product of all of the $r_{j_l}$'s, since the others
are constant on these dyadic subintervals.  One can combine these integrals
for all dyadic subintervals of $[0,1)$ of length $2^{-j_n + 1}$ to get
(\ref{the integral of the product of r_j's}).

\beginproposition
\label{inequalities for p norms}
For each positive real number $p$, there is a constant $C(p) > 0$ so that
the following is true:
\begin{eqnarray}
\label{p-norm of sum comparable to 2-norm of coefficients}
	C(p)^{-1}\,  \Bigl(\sum_{j=1}^m |\alpha_j|^2\Bigr)^{1/2}
 & \le & \biggl(\int_{[0,1)} \, \Bigl|\sum_{j=1}^m \alpha_j \, r_j(x)\Bigr|^p 
					\, dx \biggr)^{1/p}		\\
	& \le & C(p) \, \Bigl(\sum_{j=1}^m |\alpha_j|^2\Bigr)^{1/2}
								\nonumber
\end{eqnarray}
for all positive integers $m$ and all choices of real numbers $\alpha_1,
\alpha_2, \ldots, \alpha_m$.
\end{proposition}

	When $p = 2$ one can take $C(p) = 1$, i.e.,
\begin{equation}
\label{2-norm of sum equal to 2-norm of coefficients}
	\int_{[0,1)} \, \Bigl|\sum_{j=1}^m \alpha_j \, r_j(x)\Bigr|^2 \, dx 
	=  \sum_{j=1}^m |\alpha_j|^2.
\end{equation}
This follows from the fact that the $r_j$'s are \emph{orthonormal
functions}, with respect to the usual integral inner product.  The
orthonormality of the $r_j$'s can be obtained from (\ref{|r_j(x)| =
1}) and Lemma \ref{integrals of products of r_j's}.

	Next, let us recall that if $p$ and $q$ are positive real
numbers such that $p \le q$, then
\begin{equation}
\label{p-norm le q-norm when p le q}
	\Bigl(\int_{[0,1)} |f(x)|^p \, dx \Bigr)^{1/p}
		\le \Bigl(\int_{[0,1)} |f(x)|^q \, dx \Bigr)^{1/q}
\end{equation}
for any function $f(x)$ on $[0,1)$.  This is a consequence of Jensen's
inequality (\ref{(|J|^{-1} int_J f(x) dx)^p le |J|^{-1} int_J f(x)^p dx}).
Thus we get that
\begin{equation}
  \biggl(\int_{[0,1)} \, \Bigl|\sum_{j=1}^m \alpha_j \, r_j(x)\Bigr|^p 
					\, dx \biggr)^{1/p}	
	 \le  \Bigl(\sum_{j=1}^m |\alpha_j|^2\Bigr)^{1/2}
\end{equation}
when $p \le 2$, and 
\begin{equation}
	\Bigl(\sum_{j=1}^m |\alpha_j|^2\Bigr)^{1/2}
  \le \biggl(\int_{[0,1)} \, \Bigl|\sum_{j=1}^m \alpha_j \, r_j(x)\Bigr|^p 
					\, dx \biggr)^{1/p}
\end{equation}
when $p \ge 2$, by (\ref{2-norm of sum equal to 2-norm of coefficients}).
		
	Now let us show that for each $p \ge 2$, there is a constant
$C(p) > 0$ so that
\begin{equation}
\label{p-norm of sum bounded by 2-norm of coefficients}
	\biggl(\int_{[0,1)} \, \Bigl|\sum_{j=1}^m \alpha_j \, r_j(x)\Bigr|^p 
					\, dx \biggr)^{1/p}	
	 \le  C(p) \, \Bigl(\sum_{j=1}^m |\alpha_j|^2\Bigr)^{1/2}
\end{equation}
for all $\alpha_1, \alpha_2, \ldots, \alpha_m \in {\bf R}$.  It is
enough to do this when $p$ is an even integer, because of (\ref{p-norm
le q-norm when p le q}) (which implies that (\ref{p-norm of sum
bounded by 2-norm of coefficients}) holds when $p$ is replaced by any
smaller number when it holds for $p$).

	When $p$ is an even integer, one can expand out
\begin{equation}
	\Bigl|\sum_{j=1}^m \alpha_j \, r_j(x)\Bigr|^p
\end{equation}
algebraically, into a $p$-fold sum
\begin{equation}
	\sum_{j_1 = 1}^m \sum_{j_2 = 1}^m \cdots \sum_{j_p = 1}^m
		\alpha_{j_1} \, \alpha_{j_2} \cdots \alpha_{j_p} \,
		   r_{j_1}(x) \, r_{j_2}(x) \cdots r_{j_p}(x).
\end{equation}
When this is integrated in $x$ over $[0,1)$, most of the terms drop
out, because of Lemma \ref{integrals of products of r_j's}.  In fact,
this will be true whenever there is an integer $j$ which occurs among
$j_1, j_2, \ldots, j_p$ an odd number of times.  Indeed, the product
\begin{equation}
\label{r_{j_1}(x) r_{j_2}(x) cdots r_{j_p}(x)}
	r_{j_1}(x) \, r_{j_2}(x) \cdots r_{j_p}(x)
\end{equation}
can be reduced to one in which no $r_l(x)$ occurs more than once,
since $r_l(x)^2 \equiv 1$ for all $l$.  If there is an integer $j$
which occurs an odd number of times among $j_1, j_2, \ldots, j_p$,
then the reduced product will be nontrivial (containing a factor of
$r_j(x)$), and the integral of it in $x$ over $[0,1)$ vanishes, by
Lemma \ref{integrals of products of r_j's}.

	Thus one is left with the terms in which every integer $j$
occurs among $j_1, j_2, \ldots, j_p$ an even number of times (zero
times for some $j$'s).  In this case the product (\ref{r_{j_1}(x)
r_{j_2}(x) cdots r_{j_p}(x)}) reduces to the function which is
identically $1$ on $[0,1)$, and whose integral over $[0,1)$ is simply
$1$.  The product of the coefficients
\begin{equation}
	\alpha_{j_1} \, \alpha_{j_2} \cdots \alpha_{j_p}
\end{equation}
becomes a product of squares of $\alpha_j$'s, for $p/2$ $j$'s (which
may include repetitions of some $j$'s).

	In the end, 
\begin{equation}
	\int_{[0,1)} \, \Bigl|\sum_{j=1}^m \alpha_j \, r_j(x)\Bigr|^p \, dx
\end{equation}
becomes a sum of products of $p/2$ $\alpha_j^2$'s.  It is not hard to
bound this by a constant times 
\begin{equation}
	\Bigl(\sum_{j=1}^m |\alpha_j|^2\Bigr)^{p/2},
\end{equation}
which is exactly what is needed for (\ref{p-norm of sum bounded by
2-norm of coefficients}).

	It remains to show that for each $p < 2$ there is a constant
$C(p) > 0$ so that
\begin{equation}
\label{2-norm of coefficients bound by p-norm of sum}
	C(p)^{-1}\,  \Bigl(\sum_{j=1}^m |\alpha_j|^2\Bigr)^{1/2}
 \le \biggl(\int_{[0,1)} \, \Bigl|\sum_{j=1}^m \alpha_j \, r_j(x)\Bigr|^p 
					\, dx \biggr)^{1/p}.
\end{equation}
We can derive this from the $p > 2$ case, as follows.

	Fix $p \in (0,2)$.  For any function $f(x)$ on $[0,1)$, we
have that
\begin{equation}
\label{2 norm bounded in terms of 4-norm and p-norm}
	\enspace \Bigl(\int_{[0,1)} |f(x)|^2 \, dx\Bigr)^{1/2}
	   \le \Bigl(\int_{[0,1)} |f(x)|^4 \, dx\Bigr)^{a/4}
		\, \Bigl(\int_{[0,1)} |f(x)|^p \, dx\Bigr)^{(1-a)/p},
\end{equation}
where $a \in (0,1)$ is chosen so that $1/2 = a/4 + (1-a)/p$.  This is
a standard consequence of H\"older's inequality, which one might
rewrite as
\begin{equation}
	\quad \int_{[0,1)} |f(x)|^2 \, dx
	    \le \Bigl(\int_{[0,1)} (|f(x)|^{2a})^{q} \, dx\Bigr)^{1/q}
 \, \Bigl(\int_{[0,1)} (|f(x)|^{2(1-a)})^{r} \, dx\Bigr)^{1/r},
\end{equation}
where $1/q = a/2$ and $1/r = 2 (1-a)/p$.

	Suppose now that $f(x) = \sum_{j=1}^m \alpha_j \, r_j(x)$, as
in (\ref{2-norm of coefficients bound by p-norm of sum}).  From
(\ref{p-norm of sum bounded by 2-norm of coefficients}) with $p$ replaced
by $4$ and (\ref{2-norm of sum equal to 2-norm of coefficients}) we have that
\begin{equation}
	\Bigl(\int_{[0,1)} \, |f(x)|^4 \, dx \Bigr)^{1/4}	
   \le  C(4) \, \Bigl(\int_{[0,1)} \, |f(x)|^2 \, dx \Bigr)^{1/2}.
\end{equation}
We can combine this with (\ref{2 norm bounded in terms of 4-norm and p-norm})
to obtain
\begin{equation}
	 \Bigl(\int_{[0,1)} |f(x)|^2 \, dx\Bigr)^{(1-a)/2}
	   \le C(4)^a \, \Bigl(\int_{[0,1)} |f(x)|^p \, dx\Bigr)^{(1-a)/p}.
\end{equation}
This and (\ref{2-norm of sum equal to 2-norm of coefficients}) imply
(\ref{2-norm of coefficients bound by p-norm of sum}), with $C(p) =
C(4)^{a/(1-a)}$.  This completes the proof of Proposition \ref{inequalities
for p norms}.

\section{Linear functions on spheres}
\label{linear functions on spheres}

	Let $n$ be a positive integer, and let ${\bf S}^{n-1}$ denote the
unit sphere in ${\bf R}^n$, i.e.,
\begin{equation}
	{\bf S}^{n-1} = \{x \in {\bf R}^n : |x| = 1\},
\end{equation}
where $|x|$ denotes the standard Euclidean norm of $x$, $|x| =
(\sum_{i=1}^n x_i^2)^{1/2}$.

	Let $\mathcal{L}({\bf S}^{n-1})$ denote the set of functions
on ${\bf S}^{n-1}$ which are the restrictions to ${\bf S}^{n-1}$ of
linear functions on ${\bf R}^n$.  Specifically, a linear function on
${\bf R}^n$ is a function $f(x)$ which can be given by $f(x) = \langle
x, v \rangle$, where $v \in {\bf R}^n$ and $\langle u, w \rangle$ is
the standard inner product on ${\bf R}^n$,
\begin{equation}
	\langle u, w \rangle = \sum_{i=1}^{n+1} u_i \, w_i,
				\quad u, w \in {\bf R}^n.
\end{equation}
Thus $\mathcal{L}({\bf S}^{n-1})$ is a real vector space of dimension $n$.

	If $p$ is a positive real number and $f(x) \in
\mathcal{L}({\bf S}^{n-1})$, consider the quantity
\begin{equation}
\label{(frac{1}{nu_{n-1}} int_{{bf S}^{n-1}} |f(x)|^p dx)^{1/p}}
   \Bigl(\frac{1}{\nu_{n-1}} \int_{{\bf S}^{n-1}} |f(x)|^p \, dx \Bigr)^{1/p},
\end{equation}
where $dx$ is the standard volume element of integration on ${\bf
S}^{n-1}$ and $\nu_{n-1}$ is the total volume of ${\bf S}^{n-1}$.
(When $n = 1$, ${\bf S}^0 = \{1, -1\}$, the integral above should be
replaced by the sum over these two points, and $\nu_0 = 2$.)

	We may suppose that $f(x)$ is given as $t \langle x, u
\rangle$, where $u \in {\bf R}^n$ satisfies $|u| = 1$ and $t \in {\bf
R}$.  Then (\ref{(frac{1}{nu_{n-1}} int_{{bf S}^{n-1}} |f(x)|^p
dx)^{1/p}}) reduces to
\begin{equation}
	|t| \,
  \Bigl(\frac{1}{\nu_{n-1}} \int_{{\bf S}^{n-1}} 
			|\langle x, u \rangle|^p \, dx \Bigr)^{1/p}.
\end{equation}
A key point now is that this integral does not depend on $u$.  In other
words, (\ref{(frac{1}{nu_{n-1}} int_{{bf S}^{n-1}} |f(x)|^p dx)^{1/p}})
is equal to
\begin{equation}
	|t| \,
  \Bigl(\frac{1}{\nu_{n-1}} \int_{{\bf S}^{n-1}} |x_1|^p \, dx \Bigr)^{1/p},
\end{equation}
where $x_1$ denotes the first coordinate of $x$.  This is because the
volume element of integration on ${\bf S}^{n-1}$ is invariant under
orthogonal transformations on ${\bf R}^n$, which can be used to change
from $u$ to any other unit vector.  Here we are using the first
standard basis vector (with first coordinate $1$ and the rest $0$).

	Thus we may conclude that if $p_1$, $p_2$ are positive real
numbers, then there is a positive constant $C(p_1, p_2, n)$ such that
\begin{eqnarray}
\label{comparison between integrals of powers of linear functions}
\lefteqn{\Bigl(\frac{1}{\nu_{n-1}} \int_{{\bf S}^{n-1}} |f(x)|^{p_1}
						\, dx \Bigr)^{1/p_1}}	\\
	& & = C(p_1, p_2, n) \,
   \Bigl(\frac{1}{\nu_{n-1}} \int_{{\bf S}^{n-1}} |f(x)|^{p_2} 
							 \, dx \Bigr)^{1/p_2}
								\nonumber
\end{eqnarray}
for all $f \in \mathcal{L}({\bf S}^{n-1})$.

	Concerning the $p=2$ case, let us note the orthogonality
conditions
\begin{equation}
	\int_{{\bf S}^{n-1}} x_i \, x_j \, dx = 0
		\quad\hbox{when } 1 \le i, j \le n, \ i \ne j.
\end{equation}

\section{Linear functions, continued}
\label{More on linear functions}

	Let $n$ be a positive integer, and consider expressions of the
form
\begin{equation}
\label{integral of |f(y)|^p against a Gaussian, etc.}
	\Bigl( \frac{1}{\mu_n} \int_{{\bf R}^n} 
				|f(y)|^p \, e^{-|y|^2} \, dy \Bigr)^{1/p},
\end{equation}
where $0 < p < \infty$, $f(y)$ is a function on ${\bf R}^n$, and
\begin{equation}
	\mu_n = \int_{{\bf R}^n} e^{-|y|^2} \, dy.
\end{equation}

	Suppose that $f(x)$ is a linear function on ${\bf R}^n$, $f(x) =
t \langle x, u \rangle$, where $t$ is a real number and $u$ is a unit vector
in ${\bf R}^n$.  Then (\ref{integral of |f(y)|^p against a Gaussian, etc.})
is equal to
\begin{equation}
	|t| \Bigl( \frac{1}{\mu_n} \int_{{\bf R}^n} 
		  |\langle y, u \rangle|^p \, e^{-|y|^2} \, dy \Bigr)^{1/p}.
\end{equation}
As in the previous section, one can apply rotation invariance to
obtain that this integral does not depend on the choice of unit vector
$u$.  Thus (\ref{integral of |f(y)|^p against a Gaussian, etc.}) is
the same as
\begin{equation}
	|t| \Bigl( \frac{1}{\mu_n} \int_{{\bf R}^n} 
		  	 |y_1|^p \, e^{-|y|^2} \, dy \Bigr)^{1/p}.
\end{equation}

	We also have that
\begin{equation}
	\frac{1}{\mu_n} \int_{{\bf R}^n} |y_1|^p \, e^{-|y|^2} \, dy
	= \frac{1}{\mu_1} \int_{\bf R} |z|^p \, e^{-|z|^2} \, dz.
\end{equation}
This implies that (\ref{integral of |f(y)|^p against a Gaussian,
etc.}) is equal to
\begin{equation}
|t| \Bigl( \frac{1}{\mu_1} \int_{\bf R} |z|^p \, e^{-|z|^2} \, dz \Bigr)^{1/p}.
\end{equation}
Hence we obtain that if $p_1$, $p_2$ are positive real numbers, then
there is a positive constant $C_1(p_1, p_2)$ such that
\begin{eqnarray}
\label{comparing integrals against a Gaussian}
\lefteqn{\Bigl( \frac{1}{\mu_n} \int_{{\bf R}^n} 
		|f(y)|^{p_1} \, e^{-|y|^2} \, dy \Bigr)^{1/p_1}}	\\
	& & = C_1(p_1, p_2) \,
	   \Bigl( \frac{1}{\mu_n} \int_{{\bf R}^n} 
		|f(y)|^{p_2} \, e^{-|y|^2} \, dy \Bigr)^{1/p_2}	
								\nonumber
\end{eqnarray}
for all positive integers $n$ and all linear functions $f(x)$ on ${\bf
R}^n$.

	Note that for each positive real number $p$ and each function
$h(y)$ on ${\bf R}^n$ which is homogeneous of degree $p$ --- so that
$h(ry) = r^p \, h(y)$ for all $y \in {\bf R}^n$ and $r > 0$ --- the
integral of $h$ on ${\bf S}^{n-1}$ is equal to a constant times the
integral of $h(y)$ times $e^{-|y|^2}$ on ${\bf R}^n$, where the
constant depends on $p$ and $n$ but not on $h(y)$.  (This works as
well for $p > -n$, and for other integrable radial functions besides
$e^{-|y|^2}$.)  In particular, this applies to $h(y) = |f(y)|^p$ when
$f(y)$ is linear in $y$.

	Let us also note the orthogonality conditions
\begin{equation}
	\int_{{\bf R}^n} y_i \, y_j \, e^{-|y|^2} \, dy = 0
			\quad\hbox{when } 1 \le i, j \le n, \ i \ne j.
\end{equation}

\section{Lacunary sums, $p = 4$}
\label{Lacunary sums, p = 4}

	Let $n$ be a positive integer and $a_{-n}, a_{-n+1}, \ldots,
a_{n-1}, a_n$ be complex numbers, and consider the function
\begin{equation}
	f(x) = \sum_{j=-n}^n a_j \, \exp(2 \, \pi \, i \, j \, x)
\end{equation}
on $[0,1]$ (or as a periodic function on ${\bf R}$, with period $1$).
Here $\exp z$ denotes the usual exponential function, also written $e^z$.
Recall that
\begin{equation}
\label{int_0^1 |f(x)|^2 dx = sum_{j=-n}^n |a_j|^2}
	\int_0^1 |f(x)|^2 \, dx = \sum_{j=-n}^n |a_j|^2.
\end{equation}
This follows from the orthonormality of the functions $\exp (2 \, \pi
\, i \, j \, x)$ with respect to the inner product
\begin{equation}
	\langle g, h \rangle = \int_0^1 g(x) \, \overline{h(x)} \, dx.
\end{equation}
This orthonormality itself reduces to the fact that
\begin{eqnarray}
	\int_0^1 \exp (2 \, \pi \, i \, k \, x) \, dx 
		& = & 1  \quad\hbox{when } k = 0 			\\
		& = & 0  \quad\hbox{when } k \ne 0.		\nonumber
\end{eqnarray}
(For $k \ne 0$, note that $\exp (2 \, \pi \, i \, k \, x)$ is the
derivative of $(2 \, \pi \, i \, k)^{-1} \, \exp (2 \, \pi \, i \, k
\, x)$.)

	Now let us look at a kind of \emph{lacunary}\index{lacunary
sums} sum, of the form
\begin{equation}
	\phi(x) = \sum_{j = 0}^m c_j \, \exp (2 \, \pi \, i \, 2^j \, x),
\end{equation}
where $m$ is a positive integer and $c_1, \ldots, c_m$ are complex numbers.
From (\ref{int_0^1 |f(x)|^2 dx = sum_{j=-n}^n |a_j|^2}) we have that
\begin{equation}
	\int_0^1 |\phi(x)|^2 \, dx = \sum_{j=0}^m |c_j|^2.
\end{equation}
Consider 
\begin{equation}
\label{int_0^1 |phi(x)|^4 dx}
	\int_0^1 |\phi(x)|^4 \, dx.
\end{equation}

	We can write
\begin{equation}
	|\phi(x)|^4 = \phi(x)^2 \, \overline{\phi(x)}^2
\end{equation}
and multiply out the sums to obtain
\begin{eqnarray}
\lefteqn{\quad |\phi(x)|^4 = }						\\
	& & \sum_{j_1 = 0}^m \sum_{j_2 = 0}^m \sum_{j_3 = 0}^m \sum_{j_4 = 0}^m
	    \, c_1 \, c_2 \, \overline{c_3} \, \overline{c_4} \,
    \exp (2 \, \pi \, i \, (2^{j_1} + 2^{j_2} - 2^{j_3} - 2^{j_4}) \, x).
								\nonumber
\end{eqnarray}
Hence 
\begin{eqnarray}
\lefteqn{\quad\enspace \int_0^1 |\phi(x)|^4 \, dx = }			\\
	& & \sum \{ c_1 \, c_2 \, \overline{c_3} \, \overline{c_4} :
		0 \le j_1, \, j_2, \, j_3, \, j_4 \le m, \ 
			2^{j_1} + 2^{j_2} - 2^{j_3} - 2^{j_4} = 0 \}.
								\nonumber
\end{eqnarray}

\beginlemma
Suppose that $j_1$, $j_2$, $j_3$, $j_4$ are nonnegative integers such that
$2^{j_1} + 2^{j_2} - 2^{j_3} - 2^{j_4} = 0$.  Then either $j_1 = j_3$ and
$j_2 = j_4$, or $j_1 = j_4$ and $j_2 = j_3$ (or both, so that $j_1 = j_2
= j_3 = j_4$).
\end{lemma}

	(Exercise.)

	Using the lemma we obtain that
\begin{equation}
	\int_0^1 |\phi(x)|^4 \, dx 
  \le 2 \, \Bigl(\sum_{j_1 = 0}^m |c_{j_1}|^2 \Bigr)
		\Bigl(\sum_{j_2 = 0}^m |c_{j_2}|^2 \Bigr)
	= 2 \, \Bigl(\sum_{j = 0}^m |c_j|^2 \Bigr)^2.
\end{equation}
This can be rewritten as
\begin{equation}
	\int_0^1 |\phi(x)|^4 \, dx
		\le 2 \, \Bigr(\int_0^1 |\phi(y)|^2 \, dy \Bigl)^2.
\end{equation}

	There are numerous extensions and variants of these results.

\chapter{Maximal functions}

\section{Definitions and basic properties}
\index{maximal functions}

	Let $f$ be a function on $[0,1)$.  Define the \emph{dyadic
maximal function}\index{dyadic maximal functions} $M(f)$ on $[0,1)$ by
\begin{equation}
\label{def of M(f)}\index{$M(f)$}
	M(f)(x) = \sup_{k \ge 0} |E_k(f)(x)|,
\end{equation}
where $E_k(f)$ is as in Section \ref{Functions on the unit interval}.
More precisely, the supremum in (\ref{def of M(f)}) is taken over
all nonnegative integers $k$.

	This definition of $M(f)$ is equivalent to setting
\begin{eqnarray}
\label{def of M(f), 2}
  M(f)(x) & = & \sup \biggl\{ \biggl| \frac{1}{|J|} \int_J f(y) \, dy \biggr|
		 : \hbox{$J$ is a dyadic subinterval}			\\
	   &&   \qquad\qquad\qquad\qquad\quad\enspace
		 \hbox{of $[0,1)$ and $x \in J$}\biggr\}.	\nonumber
\end{eqnarray}
The equivalence of (\ref{def of M(f)}) and (\ref{def of M(f), 2}) is
not hard to check, since $E_k(f)$ is defined in terms of averages over
dyadic intervals, as in (\ref{def of E_k(f)}), and all dyadic
intervals arise in this manner.  If $J$ is a dyadic subinterval of
$[0,1)$, then its length is of the form $2^{-k}$ for some nonnegative
integer $k$, and this gives the correspondence between the $J$'s in
(\ref{def of M(f), 2}) and the $k$'s in (\ref{def of M(f)}).

	Given a nonnegative integer $l$, define $M_l(f)$ on $[0,1)$ by
\begin{equation}
\label{def of M_l(f)}\index{$M_l(f)$}
	M_l(f)(x) = \sup_{0 \le k \le l} |E_k(f)(x)|.
\end{equation}
As above, we are implicitly restricting ourselves to $k$'s which are
integers here.  This is equivalent to taking
\begin{eqnarray}
\label{def of M_l(f), 2}
\quad
 M_l(f)(x) & = & \sup \biggl\{ \biggl| \frac{1}{|J|} \int_J f(y) \, dy \biggr|
		 : \hbox{$J$ is a dyadic subinterval of}	\\
	   &&   \qquad\qquad\qquad\qquad\quad\enspace
		    \hbox{$[0,1)$, $x \in J$, and $|J| \ge 2^{-l}$}\biggr\}.
								\nonumber
\end{eqnarray}
Observe that
\begin{equation}
\label{M_l(f) le M_p(f) when p ge l}
	M_l(f) \le M_p(f)  \qquad\hbox{when $p \ge l$}
\end{equation}
and
\begin{equation}
\label{M(f)(x) = sup_{l ge 0} M_l(f)(x) for all x in [0,1)}
	M(f)(x) = \sup_{l \ge 0} M_l(f)(x) \qquad\hbox{for all $x \in [0,1)$}.
\end{equation}

	For any pair of functions $f_1$, $f_2$ on $[0,1)$, 
\begin{equation}
\label{M(f_1 + f_2) le M(f_1) + M(f_2)}
	M(f_1 + f_2) \le M(f_1) + M(f_2)
\end{equation}
and
\begin{equation}
\label{M_l(f_1 + f_2) le M_l(f_1) + M_l(f_2)}
	M_l(f_1 + f_2) \le M_l(f_1) + M_l(f_2)
\end{equation}
for all $l \ge 0$.  Also,
\begin{equation}
\label{M(c f) = |c| M(f)}
	M(c \, f) = |c| \, M(f), \quad M_l(c \, f) = |c| \, M_l(f)
\end{equation}
when $c$ is a constant.  In other words, (\ref{M(f_1 + f_2) le M(f_1)
+ M(f_2)}), (\ref{M_l(f_1 + f_2) le M_l(f_1) + M_l(f_2)}), and
(\ref{M(c f) = |c| M(f)}) say that $M(f)$ and $M_l(f)$ are
\emph{sublinear}\index{sublinear (operator)} in $f$.

\beginlemma
\label{M(f), M_l(f), E_l(f), etc.}
If $f$ is constant on the dyadic subintervals of $[0,1)$ of length
$2^{-l}$, then $M(f)$ is also constant on the dyadic subintervals of
$[0,1)$ of length $2^{-l}$, and $M(f) = M_l(f)$.  For any function
$f$,
\begin{equation}
\label{M_l(f) = M(E_l(f))}
	M_l(f) = M(E_l(f)),
\end{equation}
and $M_l(f)$ is constant on dyadic intervals of length $2^{-l}$.
\end{lemma}

	(Exercise.  One can use part (c) of Lemma \ref{Some properties
of E_k(f)}.)

\begincorollary
\label{M(f) for dyadic step functions}
If $f$ is a dyadic step function, then so is $M(f)$, and $M(f) = M_j(f)$
for sufficiently large $j$.
\end{corollary}

\section{The size of the maximal function}
\label{The size of the maximal function}

	Let $f$ be a function on $[0,1)$, and let $M(f)$ be the
corresponding dyadic maximal function, as in the preceding section.

	What can we say about $M(f)$?  How does it compare to $f$,
in terms of overall size?  In other words, if $f$ is not too large,
can we say that $M(f)$ is not too large? 

\beginlemma [Supremum bound for M(f)]
\label{supremum bound for M(f)}
Suppose that $A$ is a nonnegative real number such that $|f(x)| \le A$
for all $x \in [0,1)$.  Then $M(f)(x) \le A$ for all $x \in [0,1)$.
\end{lemma}

	This is an easy consequence of the definitions.  If $|f(x)| \le A$
for all $x \in [0,1)$, then all averages of $f$ have absolute value less
than or equal to $A$.  This implies that $M(f) \le A$ at all points in 
$[0,1)$, since $M(f)$ is defined in terms of suprema of absolute values of
averages of $f$.

	This would also work if we assumed that $|f(x)| \le A$ holds
on $[0,1)$ except for a very small set, like a finite set, a countable
set, or, more generally, a set of measure $0$.  That is, the behavior
of $f$ on a small set like this would not affect any of the integrals
of $f$, and so one could just as well replace the values of $f$ with
$0$ on such a set.

\beginproposition [Weak-type estimate for M(f)]
\label{weak-type estimate for M(f)}\index{weak-type estimates}
For each positive real number $\lambda$,
\begin{equation}
\label{weak-type inequality for M(f)}
	|\{x \in [0,1) : M(f)(x) > \lambda\}| 
		\le \frac{1}{\lambda} \int_{[0,1)} |f(w)| \, dw.
\end{equation}
\end{proposition}

	The left-hand side of (\ref{weak-type inequality for M(f)})
refers to the \emph{measure} of the set in question, and the proof of the
lemma will give a very simple meaning for this.  

	Let $\lambda > 0$ be given.  Define $\mathcal{F}$ to be the
collection of dyadic intervals $L$ in $[0,1)$ such that
\begin{equation}
\label{|frac{1}{|L|} int_L f(y) dy| > lambda}
	\biggl| \frac{1}{|L|} \int_L f(y) \, dy \biggr| > \lambda,
\end{equation}
Let us check that 
\begin{equation}
\label{{x in [0,1) : M(f)(x) > lambda} = bigcup_{L in mathcal{F}} L}
	\{x \in [0,1) : M(f)(x) > \lambda\} = \bigcup_{L \in \mathcal{F}} L.
\end{equation}
If $x \in [0,1)$ and $M(f)(x) > \lambda$, then there is a dyadic
interval $L$ in $[0,1)$ such that $x \in L$ and $L$ satisfies
(\ref{|frac{1}{|L|} int_L f(y) dy| > lambda}), because of (\ref{def of
M(f), 2}).  This shows that the left side of (\ref{{x in [0,1) :
M(f)(x) > lambda} = bigcup_{L in mathcal{F}} L}) is contained in the
right side of (\ref{{x in [0,1) : M(f)(x) > lambda} = bigcup_{L in
mathcal{F}} L}).  Conversely, if $L \in \mathcal{F}$, then
\begin{equation}
	M(f)(x) \ge \biggl| \frac{1}{|L|} \int_L f(y) \, dy \biggr|
			> \lambda
\end{equation}
for all $x \in L$.  This shows that $L$ is contained in the left side
of (\ref{{x in [0,1) : M(f)(x) > lambda} = bigcup_{L in mathcal{F}}
L}).  Hence the right side of (\ref{{x in [0,1) : M(f)(x) > lambda} =
bigcup_{L in mathcal{F}} L}) is contained in the left side.  This proves
(\ref{{x in [0,1) : M(f)(x) > lambda} = bigcup_{L in mathcal{F}} L}).

	We may as well assume that the left side of (\ref{{x in [0,1)
: M(f)(x) > lambda} = bigcup_{L in mathcal{F}} L}) is not empty, since
otherwise (\ref{weak-type inequality for M(f)}) is automatic.  This
implies that $\mathcal{F}$ is nonempty as well.  Now we apply Lemma
\ref{Structure of unions of dyadic intervals} to obtain a subcollection
$\mathcal{F}_0$ of $\mathcal{F}$ such that
\begin{equation}
\label{bigcup_{L in mathcal{F}_0} L = bigcup_{L in mathcal{F}} L}
	\bigcup_{L \in \mathcal{F}_0} L = \bigcup_{L \in \mathcal{F}} L
\end{equation}
and the intervals in $\mathcal{F_0}$ are pairwise disjoint.  Combining
these properties with (\ref{{x in [0,1) : M(f)(x) > lambda} =
bigcup_{L in mathcal{F}} L}), we obtain that
\begin{equation}
\label{|{x in [0,1) : M(f)(x) > lambda}| = sum_{L in mathcal{F}_0} |L|}
	|\{x \in [0,1) : M(f)(x) > \lambda\}| = \sum_{L \in \mathcal{F}_0} |L|.
\end{equation}
For the proof of (\ref{weak-type inequality for M(f)}), we do not need
equality in (\ref{|{x in [0,1) : M(f)(x) > lambda}| = sum_{L in
mathcal{F}_0} |L|}), but only the inequality $\le$.  This inequality
does not require the disjointness of the intervals in $\mathcal{F}_0$,
but this will be needed in a moment.

\beginremark
{\rm If $f$ is a dyadic step function, then the set $\{x \in [0,1) :
M(f)(x) > \lambda\}$ can be given as a union of finitely many dyadic
intervals.  More precisely, one could use only intervals $L$ with size
$|L| \ge 2^{-j}$ for some $j$, because of Corollary \ref{M(f) for
dyadic step functions}.}
\end{remark}

	Now let us look at the right side of (\ref{|{x in [0,1) :
M(f)(x) > lambda}| = sum_{L in mathcal{F}_0} |L|}).  Each interval $L$
in the sum satisfies (\ref{|frac{1}{|L|} int_L f(y) dy| > lambda}),
which we can rewrite as 
\begin{equation}
\label{|L| < frac{1}{lambda} | int_L f(y) dy|}
	|L| < \frac{1}{\lambda} \Bigl| \int_L f(y) \, dy \Bigr|.
\end{equation}
Hence
\begin{eqnarray}
\label{sum_{L in mathcal{F}_0} |L| < ...}
	\sum_{L \in \mathcal{F}_0} |L| 
 < \sum_{L \in \mathcal{F}_0} \frac{1}{\lambda} \Bigl| \int_L f(y) \, dy \Bigr|
 & \le & \sum_{L \in \mathcal{F}_0} \frac{1}{\lambda}  \int_L |f(y)| \, dy
									\\
 & = & \frac{1}{\lambda} \int_{\bigcup_{L \in \mathcal{F}_0} L} |f(y)| \, dy.
								\nonumber
\end{eqnarray}
The last step uses the disjointness of the intervals $L$ in $\mathcal{F}_0$.

	The combination of (\ref{|{x in [0,1) : M(f)(x) > lambda}| =
sum_{L in mathcal{F}_0} |L|}) and (\ref{sum_{L in mathcal{F}_0} |L| <
...}) implies (\ref{weak-type inequality for M(f)}).  In fact, we get
the slightly more precise inequality that
\begin{equation}
	|\{x \in [0,1) : M(f)(x) > \lambda\}| 
\le \frac{1}{\lambda} \int_{\{x \in [0,1) : M(f)(x) > \lambda\}} |f(y)| \, dy,
\end{equation}
using (\ref{{x in [0,1) : M(f)(x) > lambda} = bigcup_{L in mathcal{F}}
L}) and (\ref{bigcup_{L in mathcal{F}_0} L = bigcup_{L in mathcal{F}} L}).

\section{Some variations}
\label{Some variations}

	Instead of dyadic intervals (and averages of functions over
them, etc.), one can also look at arbitrary intervals\index{nondyadic
intervals} in the real line.  For these, it is not as easy to have
disjointness, or reduce to that case.  However, there is a simple
substitute, indicated by the following lemmas.

\beginlemma [$3$ intervals to $2$]
\label{$3$ intervals to $2$}
Suppose that $I$, $J$, and $K$ are intervals in the real line such that
there is a point $x \in {\bf R}$ which lies in all three of them.  Then
one of the intervals is contained in the union of the other two.
\end{lemma}

	This is not hard to check.  One of the intervals will go as
far as possible to the left of $x$, another of the intervals (perhaps
the same one) will go as far as possible to the right of $x$, and the
union of these two intervals (which may be the same interval) will 
contain all of $I \cup J \cup K$.

\beginlemma [No point counted more than twice]
\label{No point counted more than twice}\index{nondyadic intervals}
Let $\mathcal{A}$ be an arbitrary finite collection of intervals in
${\bf R}$.  There is a subcollection $\mathcal{A}_1$ of $\mathcal{A}$
with the following properties:
\begin{eqnarray}
\label{bigcup_{J in mathcal{A}_1} J = bigcup_{J in mathcal{A}} J}
	&& \bigcup_{J \in \mathcal{A}_1} J = \bigcup_{J \in \mathcal{A}} J;
									\\
\label{no point in R belongs to more than two intervals in mathcal{A}_1}
	&& \hbox{no point in ${\bf R}$ belongs to more than two intervals in }
		\mathcal{A}_1.				
\end{eqnarray}
\end{lemma}

	This can be derived from Lemma \ref{$3$ intervals to $2$}.
Specifically, one starts with $\mathcal{A}$ itself, and one leaves it
alone if it already satisfies (\ref{no point in R belongs to more than
two intervals in mathcal{A}_1}).  If not, there is a point in ${\bf
R}$ which lies in $3$ intervals in $\mathcal{A}$, and one can throw
away one of the intervals without changing the total union, by Lemma
\ref{$3$ intervals to $2$}.  One repeats this process until there are
no points in ${\bf R}$ which lie in $3$ intervals remaining in the
collection (i.e., that have not been thrown away).  This has to happen
eventually, since we assumed that $\mathcal{A}$ is finite.  The resulting
collection can be used for $\mathcal{A}_1$.

	One can rephrase (\ref{no point in R belongs to more than two
intervals in mathcal{A}_1}) in terms of indicator functions, as
follows:
\begin{equation}
\label{sum_{J in mathcal{A}_1} {bf 1}_J(x) le 2 cdot ....}
	\sum_{J \in \mathcal{A}_1} {\bf 1}_J(x) 
	   \le 2 \cdot {\bf 1}_{\bigl(\bigcup_{J \in \mathcal{A}_1} J\bigr)}(x)
			\qquad\hbox{for all $x \in {\bf R}$}.
\end{equation}
If $f$ is a nonnegative function on ${\bf R}$, one can integrate
(\ref{sum_{J in mathcal{A}_1} {bf 1}_J(x) le 2 cdot ....}) with $f$ to
get
\begin{equation}
\label{sum_{J in mathcal{A}_1} int_J f(x) dx le 2 ....}
	\sum_{J \in \mathcal{A}_1} \int_J f(x) \, dx
	   \le 2 \int_{\bigcup_{J \in \mathcal{A}_1} J} f(x) \, dx.
\end{equation}
This could be used in a setting like that of the last step in (\ref{sum_{L
in mathcal{F}_0} |L| < ...}).  Of course
\begin{equation}
	\sum_{J \in \mathcal{A}_1} \int_J f(x) \, dx
	   \ge \int_{\bigcup_{J \in \mathcal{A}_1} J} f(x) \, dx
\end{equation}
holds automatically (so that the equality in the last step in (\ref{sum_{L
in mathcal{F}_0} |L| < ...}) corresponds to these two inequalities, going
in opposite directions).

	These observations are special to dimension $1$, although there
are well-known substitutes for other settings (like ${\bf R}^n$, $n > 1$).

	In ${\bf R}^n$, $n > 1$, one can also consider \emph{dyadic
cubes}\index{dyadic cubes}, which behave in much the same manner as
dyadic intervals in ${\bf R}$.  A \emph{dyadic cube} in ${\bf R}^n$ is
a Cartesian product of $n$ dyadic intervals in ${\bf R}$ of the same
size.  As above, one can focus on dyadic subcubes of the unit cube,
which is the Cartesian product of $n$ copies of the unit interval
$[0,1)$.  (This is not a serious restriction, however.)  With dyadic
cubes, one has much the same kind of disjointness and partitioning
properties as for dyadic intervals.

\section{More on the size of the maximal function}
\label{More on the size of the maximal function}

	Let $f$ be a function on $[0,1)$, as before.

\beginproposition
\label{p-integral inequalities for M(f)}
For each real number $p > 1$,
\begin{equation}
\label{the p-integral inequality}
	\int_{[0,1)} M(f)(x)^p \, dx 
		\le \frac{2^p \, p}{p-1} \int_{[0,1)} |f(y)|^p \, dy.
\end{equation}
\end{proposition}

	To prove this, we shall use Lemma \ref{supremum bound for
M(f)} and Proposition \ref{weak-type estimate for M(f)}.  In fact, we
shall only need the inequalities in these lemmas, together with the
fact that $M(f)$ is sublinear in $f$, as in (\ref{M(f_1 + f_2) le
M(f_1) + M(f_2)}) and (\ref{M(c f) = |c| M(f)}).  

	Let us first state and prove two lemmas that we shall use (and
which are of broader concern).

\beginlemma
\label{modified weak-type estimate for M(f)}
For each $\lambda > 0$,
\begin{equation}
\label{modified weak-type inequality}
	|\{x \in [0,1): M(f)(x) > 2 \, \lambda\}|
 \le \frac{1}{\lambda} \int_{\{u \in [0,1) : |f(u)| > \lambda\}}  |f(u)| \, du.
\end{equation}
\end{lemma}

	This is analogous to (\ref{weak-type inequality for M(f)}) in
Proposition \ref{weak-type estimate for M(f)}, except that we have
replaced $\lambda$ by $2 \, \lambda$ on the left side of the
inequality (making it smaller), and we have restricted the domain of
integration on the right side.  (The appearance of this factor of $2$
will come from simple arithmetic, rather than depending on the details
of the situation, e.g., as with other settings as in Section \ref{Some
variations}.)

	Let $\lambda > 0$ be given, and define the functions $f_1(x)$,
$f_2(x)$ on $[0,1)$ from the function $f(x)$ by setting
\begin{equation}
	f_1(x) = f(x) \enspace\hbox{when}\enspace |f(x)| \le \lambda,
	\quad f_1(x) = 0 \enspace\hbox{when}\enspace |f(x)| > \lambda,
\end{equation}
and
\begin{equation}
	f_2(x) = f(x) \enspace\hbox{when}\enspace |f(x)| > \lambda,
	\quad f_2(x) = 0 \enspace\hbox{when}\enspace |f(x)| \le \lambda.
\end{equation}
Thus $f = f_1 + f_2$.

	Because $|f_1(x)| \le \lambda$ for all $x$, we have that
$M(f_1)(x) \le \lambda$ for all $x$ as well, as in Lemma \ref{supremum
bound for M(f)}.  This implies that
\begin{equation}
	M(f)(x) \le \lambda + M(f_2)(x)
\end{equation}
for all $x \in [0,1)$, since $M(f) \le M(f_1) + M(f_2)$, as in
(\ref{M(f_1 + f_2) le M(f_1) + M(f_2)}).  Therefore
\begin{equation}
	\{x \in [0,1) : M(f)(x) > 2 \, \lambda\}
		\subseteq \{x \in [0,1) : M(f_2)(x) > \lambda\},
\end{equation}
and hence
\begin{equation}
\label{|{x : M(f)(x) > 2 lambda}| le |{x : M(f_2)(x) > lambda}|}
	|\{x \in [0,1) : M(f)(x) > 2 \, \lambda\}|
		\le |\{x \in [0,1) : M(f_2)(x) > \lambda\}|.
\end{equation}

	On the other hand, we can apply Proposition \ref{weak-type
estimate for M(f)} with $f$ replace by $f_2$ to get that
\begin{equation}
	|\{x \in [0,1) : M(f_2)(x) > \lambda\}| 
		\le \frac{1}{\lambda} \int_{[0,1)} |f_2(u)| \, du.
\end{equation}
The right side of this inequality is equal to the right side of
(\ref{modified weak-type inequality}), by the definition of $f_2$,
and so (\ref{modified weak-type inequality}) follows from this and
(\ref{|{x : M(f)(x) > 2 lambda}| le |{x : M(f_2)(x) > lambda}|}).
This proves Lemma \ref{modified weak-type estimate for M(f)}.
	
\beginlemma
\label{integrals from distribution functions}
Let $g(x)$ be a nonnegative real-valued function on $[0,1)$, and let
$p$ be a positive real number.  Then
\begin{equation}
\label{formula for integrals, distribution function}
	\int_{[0,1)} g(x)^p \, dx = 
\int_0^\infty p \, \lambda^{p-1} |\{x \in [0,1): g(x) > \lambda\}| \, d\lambda.
\end{equation}
\end{lemma}

	This is not at all special to the unit interval, but works in a
very general way.

	There is a nice (and well-known) geometric way to look at 
(\ref{formula for integrals, distribution function}).  Consider the
set 
\begin{equation}
	\{(x,\lambda) \in {\bf R}^2 : x \in [0,1), 0 < \lambda < g(x)\}.
\end{equation}
We want to integrate the function $p \, \lambda^{p-1}$ on this set.
(This function does not depend on $x$, but the region on which we are
integrating does.)  If we integrate first in $\lambda$, and afterwards
in $x$, we get the left side of (\ref{formula for integrals,
distribution function}).  If instead we integrate in $x$ first, and
afterwards in $\lambda$, then we obtain the right side of
(\ref{integrals from distribution functions}).  This proves Lemma
\ref{integrals from distribution functions}.

	Now let us apply this formula to the proof of Proposition
\ref{p-integral inequalities for M(f)}.  Thus we obtain that
\begin{equation}
	\int_{[0,1)} M(f)(x)^p \, dx 
	     =  \int_0^\infty p \, \lambda^{p-1} 
			|\{x \in [0,1): M(f)(x) > \lambda\}| \, d\lambda.
\end{equation}
If we apply (\ref{modified weak-type inequality}) to the right side
(with $\lambda$ replaced by $\lambda/2$), then we get
\begin{equation}
	\int_{[0,1)} M(f)(x)^p \, dx 
	     \le  \int_0^\infty p \, \lambda^{p-1} 
   \biggl(\frac{2}{\lambda} \int_{\{u \in [0,1) : |f(u)| > \lambda/2\}}  
			|f(u)| \, du \biggr)        
						\, d\lambda.
\end{equation}
Let us rewrite this as
\begin{equation}
	\int_{[0,1)} M(f)(x)^p \, dx 
	     \le  \int_0^\infty 2p \, \lambda^{p-2} 
    \int_{\{u \in [0,1) : |f(u)| > \lambda/2\}} |f(u)| \, du \, d\lambda.
\end{equation}
If we interchange the order of integration on the right side, we obtain
\begin{equation}
	\int_{[0,1)} M(f)(x)^p \, dx
  \le \int_{[0,1)} \int_0^{2 \, |f(u)|} 
	2p \, \lambda^{p-2} |f(u)| \, d\lambda \, du.
\end{equation}
Because $p > 1$, by assumption (in Proposition \ref{p-integral
inequalities for M(f)}), this reduces to
\begin{eqnarray}
	\int_{[0,1)} M(f)(x)^p \, dx					
  & \le & \int_{[0,1)} 2p \, (p-1)^{-1} (2 \, |f(u)|)^{p-1} |f(u)|  \, du
								         \\
  &  =  & \frac{2^p \, p}{p-1} \int_{[0,1)} |f(u)|^p \, du.	\nonumber
\end{eqnarray}
This gives (\ref{the p-integral inequality}), as desired.

	The argument used here is an instance of \emph{Marcinkeiwicz
interpolation of operators}.\index{Marcinkeiwicz interpolation of
operators}\index{interpolation of operators!Marcinkeiwicz} Variants of
this will be employed in Chapter \ref{Square functions}.

\chapter{Square functions}
\label{Square functions}

\section{$S$-functions}
\label{$S$-functions}

	Let $f$ be a function on $[0,1)$.  Define the associated square
function\index{square functions} $S(f)$ \index{$S(f)$} on $[0,1)$ by
\begin{equation}
\label{def of S(f)}
	S(f)(x) = 
	   \Bigl(|E_0(f)(x)|^2 + 
		\sum_{j=1}^\infty |E_j(f)(x) - E_{j-1}(f)(x)|^2 \Bigr)^{1/2},
\end{equation}
where $E_j(f)$ is as in Section \ref{Functions on the unit interval}.
For each nonnegative integer $l$, define $S_l(f)$ \index{$S_l(f)$} on
$[0,1)$ by
\begin{equation}
\label{def of S_l(f)}
	S_l(f)(x) = 
	   \Bigl(|E_0(f)(x)|^2 + 
		\sum_{j=1}^l |E_j(f)(x) - E_{j-1}(f)(x)|^2\Bigr)^{1/2}.
\end{equation}
If $l = 0$, then the sum in (\ref{def of S_l(f)}) is interpreted as
being $0$.  Clearly
\begin{equation}
	S_l(f)(x) \le S_p(f)(x) \qquad\hbox{when $l \le p$},
\end{equation}
and 
\begin{equation}
	S(f)(x) = \sup_{l \ge 0} S_l(f)(x),
\end{equation}
where the supremum is implicitly taken over nonnegative integers $l$. 

	These square functions are sublinear in $f$, i.e.,
\begin{eqnarray}
	S(f_1 + f_2)(x) & \le & S(f_1)(x) + S(f_2)(x), 		\\
	S_l(f_1 + f_2)(x) & \le & S_l(f_1)(x) + S_l(f_2)(x)
\end{eqnarray}
and
\begin{equation}
	S(c \, f)(x) = |c| \, S(f)(x), \quad
		S_l(c \, f) = |c| \, S_l(f)(x)
\end{equation}
for all functions $f$, $f_1$, $f_2$, all constants $c$, and all nonnegative
integers $l$.  This follows from the fact that $(\sum_j |a_j|^2)^{1/2}$
defines a norm on sequences $\{a_j\}$.

\beginlemma
\label{S(f), S_l(f), E_l(f), etc.}
If $f$ is constant on the dyadic subintervals of $[0,1)$ of length
$2^{-l}$, then $S_l(f)$ is also constant on the dyadic subintervals of
$[0,1)$ of length $2^{-l}$, and $S(f) = S_l(f)$.  For any function
$f$,
\begin{equation}
	S_l(f) = S(E_l(f)),
\end{equation}
and $S_l(f)$ is constant on dyadic intervals of length $2^{-l}$.
\end{lemma}

	(Exercise.  One can use Lemma \ref{Some properties of E_k(f)}.)

\begincorollary
If $f$ is a dyadic step function on $[0,1)$, then so is $S(f)$, and
$S(f) = S_j(f)$ for sufficiently large $j$.
\end{corollary}

\beginlemma
\label{2-norm of S(f)}
If $f$ is any function on $[0,1)$, then
\begin{equation}
\label{int S_l(f)^2 = int |E_l(f)|^2}
	\int_{[0,1)} S_l(f)(x)^2 \, dx = \int_{[0,1)} |E_l(f)(x)|^2 \, dx
\end{equation}
for all $l \ge 0$, and
\begin{equation}
\label{int S(f)^2 = int |f|^2}
	\int_{[0,1)} S(f)(x)^2 \, dx = \int_{[0,1)} |f(x)|^2 \, dx.
\end{equation}
\end{lemma}

	As usual, one should at least assume that $f$ is mildly well-behaved,
such as being integrable.

	The main point behind Lemma \ref{2-norm of S(f)} is simply
that the functions $E_0(f)$, $E_j(f) - E_{j-1}(f)$, $j \in {\bf Z}_+$,
are all orthogonal to each other (with respect to the usual integral
inner product).  This is not hard to check, using the fact that
$E_j(f) - E_{j-1}(f)$ has integral $0$ over every dyadic subinterval
of $[0,1)$ of length $2^{-j+1}$, and that $E_k(f)$ is constant on
dyadic intervals of length $2^{-k}$.  One also uses the formula
\begin{equation}
	E_l(f) = E_0(f) + \sum_{j=1}^l (E_j(f) - E_{j-1}(f)).
\end{equation}

\beginremark 
{\rm If $f(x)$ is a sum of constant multiples of Rademacher functions
(Section \ref{Rademacher functions}), then $S(f)(x)$ is
\emph{constant}, and $S(f)(x)^2$ is the sum of the squares of the
absolute values of the coefficients of the Rademacher functions.  Of
course, in general, $S(f)(x)$ is not constant.}
\end{remark}

\section{Estimates, 1}
\label{Estimates, 1}

\beginproposition
\label{p-norm of S(f) bounded by p-norm of M(f), p < 2}
Let $p$ be a positive real number, $p < 2$.  There is a constant
$C_1(p) > 0$ so that
\begin{equation}
\label{int S(f)^p le C_1(p) int M(f)^p}
	\int_{[0,1)} S(f)(x)^p \, dx
		\le C_1(p) \int_{[0,1)} M(f)(x)^p \, dx
\end{equation}
for any function $f$ on $[0,1)$.
\end{proposition}

	Let $p < 2$ and $f$ be given.  Also let $\lambda > 0$ be
given.  We shall try to get an estimate for
\begin{equation}
	|\{x \in [0,1) : S(f)(x) > \lambda\}|,
\end{equation}
and use that to analyze the integral of $S(f)(x)^p$.  

	Let $\mathcal{F}$ denote the set of dyadic subintervals $J$
of $[0,1)$ such that
\begin{equation}
	\frac{1}{|J|} \Bigl| \int_J f(y) \, dy \Bigr| > \lambda.
\end{equation}
If $\mathcal{F}$ is empty, then this will be fine for our eventual
conclusion, and so we shall assume that $\mathcal{F}$ is not empty.
Define $\mathcal{F}_0$ to be the set of maximal intervals in
$\mathcal{F}$.  As in Lemma \ref{Structure of unions of dyadic
intervals}, we have that
\begin{equation}
\label{union of J's in mathcal{F}_0 = union of J's in mathcal{F}}
	\bigcup_{J \in \mathcal{F}_0} J = \bigcup_{J \in \mathcal{F}} J
\end{equation}
and
\begin{equation}
	J_1 \cap J_2 = \emptyset 
		\quad\hbox{when } J_1, J_2 \in \mathcal{F}_0, J_1 \ne J_2.
\end{equation}

	Assume that $[0,1)$ is not an element of $\mathcal{F}_0$.  If
it is an element of $\mathcal{F}_0$, then this will also be fine for
our eventual conclusion.  Let $\mathcal{F}_1$ denote the set of dyadic
subintervals $L$ of $[0,1)$ such that there is a $J \in \mathcal{F}_0$
such that
\begin{equation}
	J \subseteq L \hbox{ and } |J| = |L|/2.
\end{equation}
Because $\mathcal{F}_0$ consists of maximal intervals in $\mathcal{F}$,
each $L$ in $\mathcal{F}_1$ does not lie in $\mathcal{F}$.  Therefore,
\begin{equation}
\label{frac{1}{|L|} |int_L f(y) dy| le lambda}
	\frac{1}{|L|} \Bigl| \int_L f(y) \, dy \Bigr| \le \lambda
\end{equation}
for all $L \in \mathcal{F}_1$.

	The elements of $\mathcal{F}_1$ need not be disjoint, and so
we let $\mathcal{F}_{10}$ be a subcollection of $\mathcal{F}_1$ as
in Lemma \ref{Structure of unions of dyadic
intervals}.  Thus
\begin{equation}
	\bigcup_{L \in \mathcal{F}_{10}} L = \bigcup_{L \in \mathcal{F}_1} L
\end{equation}
and
\begin{equation}
	L_1 \cap L_2 = \emptyset 
		\quad\hbox{when } L_1, L_2 \in \mathcal{F}_{10}, L_1 \ne L_2.
\end{equation}

	Define a function $f_\lambda(x)$ on $[0,1)$ by
\begin{eqnarray}
\label{def of f_lambda(x)}
	f_\lambda(x) & = & \frac{1}{|L|} \int_L f(y) \, dy
			    \qquad\hbox{when } x \in L, L \in \mathcal{F}_{10},
									\\
		     & = & f(x)  \qquad\qquad\qquad
 		\hbox{when } x \in [0,1) \backslash 
			\biggl(\bigcup_{L \in \mathcal{F}_{10}} L \biggr).
								\nonumber
\end{eqnarray}

\beginlemma 
\label{averages of f and averages of f_lambda}
Let $K$ be a dyadic subinterval of $[0,1)$ such that $K$ contains a
point $x$ in $[0,1) \backslash \biggl(\bigcup_{L \in \mathcal{F}_{10}} L
\biggr)$.  Then
\begin{equation}
	\frac{1}{|K|} \int_K f(u) \, du 
		= \frac{1}{|K|} \int_K f_\lambda(u) \, du.
\end{equation}
(The same is true if $K$ contains an element of $\mathcal{F}_{10}$
as a subinterval.)
\end{lemma}

	(Exercise.  Under these conditions, $K$ is equal to the union
of the intervals $L \in \mathcal{F}_0$ such that $L \subseteq K$ and
the set $K \, \backslash \bigcup_{L \in \mathcal{F}_0} L$, and these
sets are all disjoint.)

\begincorollary
\label{S(f)(x) = S(f_lambda)(x) on a certain set}
If $x \in [0,1) \backslash \biggl(\bigcup_{L \in \mathcal{F}_{10}} L
\biggr)$, then $S(f)(x) = S(f_\lambda)(x)$.
\end{corollary}

	Using this corollary we get that
\begin{eqnarray}
\label{{x : S(f)(x) > lambda} subseteq ....}
\lefteqn{\{x \in [0,1) : S(f)(x) > \lambda\}}	\\
	& & \qquad \subseteq \biggl(\bigcup_{L \in \mathcal{F}_{10}} L \biggr)
			\cup \{x \in [0,1) : S(f_\lambda)(x) > \lambda\}.
								\nonumber
\end{eqnarray}
Hence
\begin{eqnarray}
\label{|{x : S(f)(x) > lambda}| le ....}
\lefteqn{|\{x \in [0,1) : S(f)(x) > \lambda\}|}		\\
	& & \qquad \le \sum_{L \in \mathcal{F}_{10}} |L| 
			+ |\{x \in [0,1) : S(f_\lambda)(x) > \lambda\}|.
								\nonumber
\end{eqnarray}
For each $L \in \mathcal{F}_{10}$, there is an interval $J \in
\mathcal{F}_0$ such that $J \subseteq L$ and $|J| \ge |L|/2$.  The
intervals in $\mathcal{F}_{10}$ are disjoint, and this implies that each
$J \in \mathcal{F}_0$ is associated to at most one $L \in \mathcal{F}_{10}$.
This implies that
\begin{equation}
  \sum_{L \in \mathcal{F}_{10}} |L| \le 2 \sum_{J \in \mathcal{F}_0} |J|.
\end{equation}
The $J$'s are disjoint, so that
\begin{equation}
	\sum_{L \in \mathcal{F}_{10}} |L|
		 \le 2 \, \biggl| \bigcup_{J \in \mathcal{F}_0} J \biggr|
		  = 2 \, \biggl| \bigcup_{J \in \mathcal{F}} J \biggr|,
\end{equation}
where the last step uses (\ref{union of J's in mathcal{F}_0 = union of
J's in mathcal{F}}).

	We also have that
\begin{equation}
	\bigcup_{J \in \mathcal{F}} J = \{x \in [0,1) : M(f)(x) > \lambda\}.
\end{equation}
See (\ref{{x in [0,1) : M(f)(x) > lambda} = bigcup_{L in mathcal{F}} L});
at the moment, we really only need that the left side is contained in the
right side, which one can easily get from (\ref{def of M(f), 2}).
From this we obtain that
\begin{equation}
	\biggl|\bigcup_{J \in \mathcal{F}} J \biggr|
		 = |\{x \in [0,1) : M(f)(x) > \lambda\}|,
\end{equation}
and hence
\begin{equation}
	\sum_{L \in \mathcal{F}_{10}} |L| 
		\le 2 \, |\{x \in [0,1) : M(f)(x) > \lambda\}|.
\end{equation}
Plugging this into (\ref{|{x : S(f)(x) > lambda}| le ....}) leads to
\begin{eqnarray}
\label{|{x : S(f)(x) > lambda}| le ...,2}
\lefteqn{\quad |\{x \in [0,1) : S(f)(x) > \lambda\}|}		\\
	& & \le 2 \, |\{x \in [0,1) : M(f)(x) > \lambda\}|
			+ |\{x \in [0,1) : S(f_\lambda)(x) > \lambda\}|.
								\nonumber
\end{eqnarray}

	The first term on the right side of (\ref{|{x : S(f)(x) >
lambda}| le ...,2}) is fine for the goal of getting estimates in terms
of $M(f)$.  For the second we apply the Tchebychev inequality and
Lemma \ref{2-norm of S(f)}, to obtain
\begin{eqnarray}
	\lambda^2 \, |\{x \in [0,1) : S(f_\lambda)(x) > \lambda\}|
		& \le & \int_{[0,1)} S(f_\lambda)(x)^2 \, dx		\\
		& = &   \int_{[0,1)} |f_\lambda(x)|^2 \, dx.	\nonumber
\end{eqnarray}
Thus 
\begin{eqnarray}
\label{|{x : S(f)(x) > lambda}| le ...,3}
\lefteqn{|\{x \in [0,1) : S(f)(x) > \lambda\}|}		\\
	& & \le 2 \, |\{x \in [0,1) : M(f)(x) > \lambda\}|
			+ \lambda^{-2} \int_{[0,1)} |f_\lambda(x)|^2 \, dx.
								\nonumber
\end{eqnarray}

\beginlemma
\label{|f_lambda(x)| le min(lambda, M(f)(x))}
$|f_\lambda(x)| \le \min (\lambda, M(f)(x))$ (at least almost everywhere).
\end{lemma}

	(Exercise.  The main point is that
\begin{equation}
	\frac{1}{L} \Bigl|\int_L f(y) \, dy \Bigr|
\end{equation}
is always less than or equal to $M(f)(x)$ when $x \in L$, as in
(\ref{def of M(f), 2}), and that it is less than or equal to $\lambda$
if $L \in \mathcal{F}_1$ (and hence if $L \in \mathcal{F}_{10}$), or
if $x \in L$ and $x \in [0,1) \backslash \biggl(\bigcup_{L \in
\mathcal{F}_{10}} L \biggr)$.  Indeed, these conditions on $L$ imply
that $L \not\in \mathcal{F}$.)

	(Actually, one also has that $M(f_\lambda)(x) \le
\min(\lambda, M(f)(x))$ for all $x \in [0,1)$.  This is not hard to check.)

	Using the lemma, we may replace 
(\ref{|{x : S(f)(x) > lambda}| le ...,3}) with 
\begin{eqnarray}
\label{|{x : S(f)(x) > lambda}| le ...,4}
\lefteqn{\quad\enspace |\{x \in [0,1) : S(f)(x) > \lambda\}|}		\\
	& & \le 2 \, |\{x \in [0,1) : M(f)(x) > \lambda\}|
		+ \lambda^{-2} \int_{[0,1)} \min(\lambda, M(f)(x))^2 \, dx.
								\nonumber
\end{eqnarray}

	In deriving (\ref{|{x : S(f)(x) > lambda}| le ...,4}), we made
two technical assumptions.  The first was that $\mathcal{F} \ne
\emptyset$.  If $\mathcal{F} = \emptyset$, then the first term on the
right side of (\ref{|{x : S(f)(x) > lambda}| le ...,4}) is equal to
$0$.  In this case, we can take $f_\lambda = f$, and the argument
works in the same manner as before (and is simpler).  (Note that
$M(f)(x) \le \lambda$ for all $x \in [0,1)$ in this case.)  The second
assumption was that $[0,1)$ is not an element of $\mathcal{F}_0$.  If
$[0,1)$ is an element of $\mathcal{F}_0 \subseteq \mathcal{F}$, then
\begin{equation}
	\Bigl|\int_{[0,1)} f(y) \, dy \Bigr| > \lambda,
\end{equation}
and $\{x \in [0,1) : M(f)(x) > \lambda\} = [0,1)$.  In this situation,
(\ref{|{x : S(f)(x) > lambda}| le ...,4}) holds automatically, with
neither the second term on the right side nor the factor of two in the
first term on the right side being needed.  To summarize, (\ref{|{x :
S(f)(x) > lambda}| le ...,4}) in fact holds without these two
assumptions.

	By Lemma \ref{integrals from distribution functions}, we have that
\begin{equation}
\label{formula for int S(f)^p}
	\int_{[0,1)} S(f)(x)^p \, dx 
		= \int_0^\infty p \, \lambda^{p-1} 
			|\{x \in [0,1): S(f)(x) > \lambda\}| \, d\lambda
\end{equation}
and
\begin{equation}
\label{formula for int M(f)^p}
	\int_{[0,1)} M(f)(x)^p \, dx  
		= \int_0^\infty p \, \lambda^{p-1} 
			|\{x \in [0,1): M(f)(x) > \lambda\}| \, d\lambda.
\end{equation}
Integrating (\ref{|{x : S(f)(x) > lambda}| le ...,4}) in this way, we obtain
that
\begin{eqnarray}
\label{int_{[0,1)} S(f)(x)^p dx le ....}
\lefteqn{\quad \int_{[0,1)} S(f)(x)^p \, dx}	\\
	& & \le 2 \int_{[0,1)} M(f)(x)^p \, dx 
	      + \int_0^\infty p \, \lambda^{p-3} 
		     \int_{[0,1)} \min(\lambda, M(f)(x))^2 \, dx \, d\lambda.
								\nonumber
\end{eqnarray}

	Let us interchange the order of integrations in the second term
on the right side of (\ref{int_{[0,1)} S(f)(x)^p dx le ....}).  This leads
to
\begin{equation}
	\int_{[0,1)} \int_0^\infty p \, \lambda^{p-3}
		\min(\lambda, M(f)(x))^2 \, d\lambda \, dx.
\end{equation}
The integral in $\lambda$ here can be computed exactly.  Specifically,
we have that
\begin{equation}
	\int_{M(f)(x)}^\infty p \, \lambda^{p-3} \, M(f)(x)^2 \, d\lambda
		= \frac{p}{p-2} M(f)(x)^2.
\end{equation}
Note that it is important that $p < 2$ here, for the convergence of the
integral.  Also,
\begin{equation}
	\int_0^{M(f)(x)} p \, \lambda^{p-3} \, \lambda^2 \, d\lambda
		= M(f)(x)^p.
\end{equation}
Plugging these formulae into (\ref{int_{[0,1)} S(f)(x)^p dx le ....}),
we obtain (\ref{int S(f)^p le C_1(p) int M(f)^p}), as desired.

	The constant that we get here blows up as $p \to 2$, which is
silly, since we know that the $p=2$ case behaves well, as in Lemma
\ref{2-norm of S(f)}.  This is an artifact of the argument, and one
can get rid of it using interpolation.  For this it is convenient
to consider the analogue of (\ref{int S(f)^p le C_1(p) int M(f)^p})
with $M(f)$ replaced by $|f|$ (for $p > 1$, such as $p = 3/2$).

\section{Estimates, 2}
\label{Estimates, 2}

\beginproposition
\label{p-norm of M(f) bounded by p-norm of S(f), p < 2}
Let $p$ be a positive real number, $p < 2$.  There is a constant
$C_2(p) > 0$ so that
\begin{equation}
\label{int M(f)^p le C_2(p) int S(f)^p}
	\int_{[0,1)} M(f)(x)^p \, dx 
		\le C_2(p) \int_{[0,1)} S(f)(x)^p \, dx
\end{equation}
for any function $f$ on $[0,1)$.
\end{proposition}

	The argument for this is very much analogous to the one in
Section \ref{Estimates, 1}.

	Let $p < 2$ and $f$ be given, and let $\lambda > 0$ be given too.
As before, we would like to get a helpful upper bound for 
\begin{equation}
\label{|{x in [0,1) : M(f)(x) > lambda}|}
	|\{x \in [0,1) : M(f)(x) > \lambda\}|.
\end{equation}

\beginlemma
\label{{x in [0,1) : S(f)(x) > lambda} is a union of dyadic intervals}
The set	$\{x \in [0,1) : S(f)(x) > \lambda\}$ is a union of dyadic
subintervals of $[0,1)$.
\end{lemma}

	Indeed, suppose that $w \in [0,1)$ and $S(f)(w) > \lambda$.
This is the same as saying that
\begin{equation}
	\Bigl(|E_0(f)(w)|^2 + 
		\sum_{j=1}^\infty |E_j(f)(w) - E_{j-1}(f)(w)|^2 \Bigr)^{1/2}
			> \lambda,
\end{equation}
by (\ref{def of S(f)}).  From this we get automatically that
\begin{equation}
	\Bigl(|E_0(f)(w)|^2 + 
		\sum_{j=1}^l |E_j(f)(w) - E_{j-1}(f)(w)|^2 \Bigr)^{1/2}
			> \lambda
\end{equation}
for some $l$.  On the other hand, $E_j(f)$ is constant on dyadic
intervals of length $2^{-j}$, and hence $E_j(f)(w) = E_j(f)(y)$ for
all $j \le l$ and for all $y$ which lie in the dyadic subinterval of
$[0,1)$ of length $2^{-l}$ that contains $w$.  For these $y$'s, we
obtain that
\begin{equation}
	\Bigl(|E_0(f)(y)|^2 + 
		\sum_{j=1}^l |E_j(f)(y) - E_{j-1}(f)(y)|^2 \Bigr)^{1/2}
			> \lambda.
\end{equation}
Therefore $S(f)(y) > \lambda$ for all $y$ in this same dyadic interval.
Thus there is a dyadic interval which contains $w$ and which lies in our
set, and the lemma follows, since $w$ was an arbitrary element of the set.

	Let $\mathcal{G}_0$ denote the collection of maximal dyadic
subintervals of $[0,1)$ which are contained in $\{x \in [0,1) :
S(f)(x) > \lambda\}$.  We shall assume for the time being that $\{x
\in [0,1) : S(f)(x) > \lambda\} \ne \emptyset$, so that $\mathcal{G}_0
\ne \emptyset$ as well.  As usual, we have that
\begin{equation}
\label{bigcup_{J in mathcal{G}_0} J = {x in [0,1) : S(f)(x) > lambda}}
	\bigcup_{J \in \mathcal{G}_0} J = \{x \in [0,1) : S(f)(x) > \lambda\},
\end{equation}
and that the intervals in $\mathcal{G}_0$ are pairwise disjoint.  In
particular,
\begin{equation}
\label{sum_{J in mathcal{G}_0} |J| = |{x in [0,1) : S(f)(x) > lambda}|}
	\sum_{J \in \mathcal{G}_0} |J| = |\{x \in [0,1) : S(f)(x) > \lambda\}|.
\end{equation}

	Let us also assume that $\{x \in [0,1) : S(f)(x) > \lambda\}$
is not equal to the whole unit interval $[0,1)$.  Let $\mathcal{G}_1$
denote the set of dyadic subintervals $L$ of $[0,1)$ such that there
is a $J \in \mathcal{G}_0$ which satisfies
\begin{equation}
	J \subseteq L \hbox{ and } |J| = |L| / 2.
\end{equation}
Thus $\mathcal{G}_1 \ne \emptyset$ since $\mathcal{G}_0 \ne \emptyset$
and $\{x \in [0,1) : S(f)(x) > \lambda\} \ne [0,1)$.  Because the
elements of $\mathcal{G}_0$ are maximal dyadic intervals contained in
$\{x \in [0,1) : S(f)(x) > \lambda\}$, we obtain that each $L \in
\mathcal{G}_1$ is not a subset of $\{x \in [0,1) : S(f)(x) >
\lambda\}$, i.e., there is a point $u \in L$ such that
\begin{equation}
	S(f)(u) \le \lambda.
\end{equation}
The definition (\ref{def of S(f)}) of $S(f)$ then implies that
\begin{equation}
	\Bigl(|E_0(f)(u)|^2 + 
	    \sum_{j=1}^{\ell(L)} |E_j(f)(u) - E_{j-1}(f)(u)|^2 \Bigr)^{1/2}
			\le \lambda,
\end{equation}
where $\ell(L)$ is chosen so that $2^{-\ell(L)} = |L|$.  (If $\ell(L) = 0$,
then the sum on the left side of the inequality is interpreted as being $0$.)
We conclude that
\begin{eqnarray}
\label{S_{ell(L)}(f) le lambda on L, L in mathcal{G}_1}
   & &	\Bigl(|E_0(f)(w)|^2 + 
	    \sum_{j=1}^{\ell(L)} |E_j(f)(w) - E_{j-1}(f)(w)|^2 \Bigr)^{1/2}
			\le \lambda					\\
   & & \qquad\qquad\hbox{for all } w \in L,			\nonumber
\end{eqnarray}
since $E_j(f)(w)$ is constant on $L$ when $j \le \ell(L)$.

	The intervals in $\mathcal{G}_1$ need not be pairwise
disjoint, and, as usual, we can pass to a subset $\mathcal{G}_{10}$
(of maximal elements of $\mathcal{G}_1$, as in Lemma \ref{Structure of
unions of dyadic intervals}) so that
\begin{equation}
\label{bigcup_{L in mathcal{G}_{10}} L = bigcup_{L in mathcal{G}_1} L}
	\bigcup_{L \in \mathcal{G}_{10}} L
		= \bigcup_{L \in \mathcal{G}_1} L
\end{equation}
and 
\begin{equation}
	L_1 \cap L_2 = \emptyset 
		\quad\hbox{when } L_1, L_2 \in \mathcal{G}_1, L_1 \ne L_2.
\end{equation}

	The definition of $\mathcal{G}_1$ implies that
\begin{equation}
	\bigcup_{L \in \mathcal{G}_1} L 
		\supseteq \bigcup_{J \in \mathcal{G}_0} J.
\end{equation}
Hence 
\begin{equation}
\label{bigcup_{L in mathcal{G}_{10}} L supseteq {S(f) > lambda}}
	\bigcup_{L \in \mathcal{G}_{10}} L 
		\supseteq \{x \in [0,1) : S(f)(x) > \lambda\},
\end{equation}
by (\ref{bigcup_{L in mathcal{G}_{10}} L = bigcup_{L in mathcal{G}_1} L})
and (\ref{bigcup_{J in mathcal{G}_0} J = {x in [0,1) : S(f)(x) > lambda}}).
On the other hand,
\begin{equation}
\label{sum_{L in mathcal{G}_{10}} |L| le ....}
	\sum_{L \in \mathcal{G}_{10}} |L| 
		\le \sum_{J \in \mathcal{G}_0} 2 \, |J|
		= 2 \, |\{x \in [0,1) : S(f)(x) > \lambda\}|,
\end{equation}
because of the way that the intervals $L \in \mathcal{F}_1$ were
chosen and (\ref{sum_{J in mathcal{G}_0} |J| = |{x in [0,1) : S(f)(x)
> lambda}|}).

	Define a function $g_\lambda(x)$ on $[0,1)$ by
\begin{eqnarray}
\label{def of g_lambda}
	g_\lambda(x) & = & \frac{1}{|L|} \int_L f(y) \, dy
			\quad\hbox{when } x \in L, L \in \mathcal{G}_{10}  \\
		     & = & f(x)
			\qquad\qquad\quad\hbox{when } 
    x \in [0,1) \backslash \biggl(\bigcup_{L \in \mathcal{G}_{10}} L \biggr).
								\nonumber
\end{eqnarray}

\beginlemma
\label{frac{1}{|K|} int_K g_lambda = frac{1}{|K|} int_K g when ....}
If $K$ is a dyadic subinterval of $[0,1)$, and if either $K \supseteq
L$ for some $L \in \mathcal{G}_{10}$, or $x \in K$ for some $x \in
[0,1) \backslash \biggl(\bigcup_{L \in \mathcal{G}_{10}} L \biggr)$,
then
\begin{equation}
	\frac{1}{|K|} \int_K g_\lambda(u) \, du
		= \frac{1}{|K|} \int_K f(u) \, du.
\end{equation}
\end{lemma}

	This is not hard to check, by writing $K$ as the disjoint
union of the sets $K \backslash \biggl(\bigcup_{L \in
\mathcal{G}_{10}} L \biggr)$ and the intervals in $\mathcal{G}_{10}$
which are subsets of $K$.

\begincorollary
\label{M(f) = M(g_lambda) on the proper set}
If $x \in [0,1) \backslash \biggl(\bigcup_{L \in \mathcal{G}_{10}} L
\biggr)$, then $M(f)(x) = M(g_\lambda)(x)$.
\end{corollary}

	This follows immediately from Lemma \ref{frac{1}{|K|} int_K
g_lambda = frac{1}{|K|} int_K g when ....}.

\begincorollary
\label{S(f), S(g_lambda), etc.}
If $x \in [0,1) \backslash \biggl(\bigcup_{L \in \mathcal{G}_{10}} L
\biggr)$, then $S(f)(x) = S(g_\lambda)(x)$.  If $L \in
\mathcal{G}_{10}$, then
\begin{equation}
\label{S(g_lambda)(v) = ... on L, L in mathcal{G}_{10}}
	S(g_\lambda)(v) = \Bigl(|E_0(f)(v)|^2 + 
	    \sum_{j=1}^{\ell(L)} |E_j(f)(v) - E_{j-1}(f)(v)|^2 \Bigr)^{1/2}
\end{equation}
for all $v \in L$ (where $2^{-\ell(L)} = |L|$, as before).
\end{corollary}

	This can be derived from Lemma \ref{frac{1}{|K|} int_K
g_lambda = frac{1}{|K|} int_K g when ....}  and the definition
(\ref{def of g_lambda}) of $g_\lambda$.

\begincorollary
\label{S(g_lambda) le min(lambda, S(f))}
$S(g_\lambda)(x) \le \min(\lambda, S(f)(x))$ for all $x \in [0,1)$.
\end{corollary}

	That $S(g_\lambda)(x) \le S(f)(x)$ for all $x \in [0,1)$
follows easily from Corollary \ref{S(f), S(g_lambda), etc.}.  If $x
\in [0,1) \backslash \biggl(\bigcup_{L \in \mathcal{G}_{10}} L
\biggr)$, then $S(g_\lambda)(x) = S(f)(x) \le \lambda$ because of
(\ref{bigcup_{L in mathcal{G}_{10}} L supseteq {S(f) > lambda}}).  If
$x \in L$ for some $L \in \mathcal{G}_{10} \subseteq \mathcal{G}_1$,
then $S(g_\lambda)(x) \le \lambda$ by (\ref{S(g_lambda)(v) = ... on L,
L in mathcal{G}_{10}}) and (\ref{S_{ell(L)}(f) le lambda on L, L in
mathcal{G}_1}).

	Now let us apply $g_\lambda$ to the estimation of (\ref{|{x in
[0,1) : M(f)(x) > lambda}|}).  Because of Corollary \ref{M(f) =
M(g_lambda) on the proper set}, we have that
\begin{eqnarray}
\lefteqn{\{x \in [0,1) : M(f)(x) > \lambda\}}   \\
	 & & 	\subseteq \biggl(\bigcup_{L \in \mathcal{G}_{01}} L\biggr)
		  \cup \{x \in [0,1) : M(g_\lambda)(x) > \lambda\}.
							\nonumber
\end{eqnarray}
Hence
\begin{eqnarray}
\lefteqn{|\{x \in [0,1) : M(f)(x) > \lambda\}|}   \\
	 & & 	\le \sum_{L \in \mathcal{G}_{01}} |L|
		   + |\{x \in [0,1) : M(g_\lambda)(x) > \lambda\}|,
							\nonumber
\end{eqnarray}
and therefore
\begin{eqnarray}
\label{|{x in [0,1) : M(f)(x) > lambda}| le ....}
\lefteqn{\quad\enspace |\{x \in [0,1) : M(f)(x) > \lambda\}|}   \\
	 & & 	\le 2 \, |\{x \in [0,1) : S(f)(x) > \lambda\}|
		   + |\{x \in [0,1) : M(g_\lambda)(x) > \lambda\}|,
							\nonumber
\end{eqnarray}
by (\ref{sum_{L in mathcal{G}_{10}} |L| le ....}).

	Observe that
\begin{eqnarray}
	|\{x \in [0,1) : M(g_\lambda)(x) > \lambda\}|
		& \le & \lambda^{-2} \int_{[0,1)} M(g_\lambda)(u)^2 \, du  \\
		& \le & C \, \lambda^{-2} \int_{[0,1)} |g_\lambda(y)|^2 \, dy
								\nonumber
\end{eqnarray}
for some constant $C > 0$ (independent of $g$ and $\lambda$).  This uses
Proposition \ref{p-integral inequalities for M(f)} with $p = 2$.  On the
other hand,
\begin{eqnarray}
	\int_{[0,1)} |g_\lambda(y)|^2 \, dy 
		& = & \int_{[0,1)} S(g_\lambda)(w)^2 \, dw		\\
		& \le & \int_{[0,1)} \min(\lambda, S(f)(w))^2 \, dw,
								\nonumber
\end{eqnarray}
by Lemma \ref{2-norm of S(f)} and Corollary \ref{S(g_lambda) le
min(lambda, S(f))}.  Hence
\begin{equation}
	\quad |\{x \in [0,1) : M(g_\lambda)(x) > \lambda\}|
	    \le C \, \lambda^{-2} \int_{[0,1)} \min(\lambda, S(f)(w))^2 \, dw.
\end{equation}

	Putting this estimate into (\ref{|{x in [0,1) : M(f)(x) >
lambda}| le ....}), we obtain that
\begin{eqnarray}
\label{|{x in [0,1) : M(f)(x) > lambda}| le ...., 2}
\lefteqn{\qquad |\{x \in [0,1) : M(f)(x) > \lambda\}|}   \\
	 & & 	\le 2 \, |\{x \in [0,1) : S(f)(x) > \lambda\}|
	      + C \, \lambda^{-2} \int_{[0,1)} \min(\lambda, S(f)(w))^2 \, dw.
							\nonumber
\end{eqnarray}
Near the beginning of this argument, we made two assumptions, which
were that the set $\{x \in [0,1) : S(f)(x) > \lambda\}$ is neither
empty nor all of $[0,1)$.  If $\{x \in [0,1) : S(f)(x) > \lambda\}$ is
all of $[0,1)$, then (\ref{|{x in [0,1) : M(f)(x) > lambda}| le ....,
2}) holds automatically, without the second term on the right side, or
the factor of $2$ in the first term on the right side.  If $\{x \in
[0,1) : S(f)(x) > \lambda\} = \emptyset$, then the first term on the
right side of (\ref{|{x in [0,1) : M(f)(x) > lambda}| le ...., 2}) is
equal to $0$.  In this case we merely take $g_\lambda = f$ on all of
$[0,1)$, and the same argument works as before (with some
simplifications).  In particular, we have that $S(f) = S(g_\lambda)
\le \lambda$ everywhere on $[0,1)$.

	Thus we get (\ref{|{x in [0,1) : M(f)(x) > lambda}| le ....,
2}), and without these extra assumptions.  At this stage, the argument
can use the same kind of computations as in Section \ref{Estimates,
1}, starting from (\ref{|{x : S(f)(x) > lambda}| le ...,4}).  That is,
we apply Lemma \ref{integrals from distribution functions} and
integrate in $\lambda$.

	As before, we also get a constant which blows up as $p \to 2$,
and this is silly, because Lemma \ref{2-norm of S(f)} indicates that
$p = 2$ is fine.  This can be fixed.

\section{Duality}
\label{Duality}

	 If $f_1$, $f_2$ are functions on $[0,1)$, and $l$ is a nonnegative
integer, then
\begin{eqnarray}
	& & \quad \int_{[0,1)} E_l(f_1) \, E_l(f_2) \, dx =		\\
	& & \int_{[0,1)} \biggl(E_0(f_1) \, E_0(f_2) 
 + \sum_{j=1}^l (E_j(f_1) - E_{j-1}(f_1)) (E_j(f_2) - E_{j-1}(f_2))\biggr) 
				\, dx.				\nonumber
\end{eqnarray}
(As usual, the sum in $j$ on the right side is interpreted as being
$0$ when $l = 0$.)  This is a ``bilinear'' version of (\ref{int
S_l(f)^2 = int |E_l(f)|^2}) in Lemma \ref{2-norm of S(f)}.  The proof
is essentially the same, i.e., with respect to the usual integral
inner product, $E_0(f_1)$ is orthogonal to $E_k(f_2) - E_{k-1}(f_2)$
for $k \ge 1$, $E_0(f_2)$ is orthogonal to $E_j(f_1) - E_{j-1}(f_1)$
for $j \ge 1$, and $E_j(f_1) - E_{j-1}(f_1)$ is orthogonal to
$E_k(f_2) - E_{k-1}(f_2)$ when $j, k \ge 1$, $j \ne k$.  Also,
\begin{eqnarray}
\label{E_l(f_1) =  E_0(f_1) + sum_{j=1}^l (E_j(f_1) - E_{j-1}(f_1))}
	E_l(f_1) & = & E_0(f_1) + \sum_{j=1}^l (E_j(f_1) - E_{j-1}(f_1))   \\
\label{E_k(f_2) = E_0(f_2) + sum_{k=1}^l (E_k(f_2) - E_{k-1}(f_2))}
	E_k(f_2) & = & E_0(f_2) + \sum_{k=1}^l (E_k(f_2) - E_{k-1}(f_2)).
\end{eqnarray}
Thus one can multiply $E_l(f_1)$ and $E_l(f_2)$, expand out the
product using (\ref{E_l(f_1) = E_0(f_1) + sum_{j=1}^l (E_j(f_1) -
E_{j-1}(f_1))}) and (\ref{E_k(f_2) = E_0(f_2) + sum_{k=1}^l (E_k(f_2)
- E_{k-1}(f_2))}), integrate, and cancel out the cross terms using the
orthogonality properties.

	Under suitable conditions, one has that
\begin{eqnarray}
\label{int_{[0,1)} f_1 f_2 dx = ...}
	& & \quad \int_{[0,1)} f_1 \, f_2 \, dx =		\\
	& & \int_{[0,1)} \biggl(E_0(f_1) \, E_0(f_2) 
+ \sum_{j=1}^\infty (E_j(f_1) - E_{j-1}(f_1)) (E_j(f_2) - E_{j-1}(f_2))\biggr) 
				\, dx.				\nonumber
\end{eqnarray}
For instance, this is easy if at least one of $f_1$ and $f_2$ is a
dyadic step function.  (In this case
\begin{equation}
	\int_{[0,1)} f_1 \, f_2 \, dx = \int_{[0,1)} E_l(f_1) \, E_l(f_2) \, dx
\end{equation}
for sufficiently large $l$.)

	A consequence of (\ref{int_{[0,1)} f_1 f_2 dx = ...})  is that
\begin{equation}
\label{|int_{[0,1)} f_1(x) f_2(x) dx| le int_{[0,1)} S(f_1)(x) S(f_2)(x) dx}
	\Bigl|\int_{[0,1)} f_1(x) \, f_2(x) \, dx \Bigr|
		\le \int_{[0,1)} S(f_1)(x) \, S(f_2)(x) \, dx.
\end{equation}
This follows by applying the Cauchy-Schwarz inequality to the
integrand on the right side of (\ref{int_{[0,1)} f_1 f_2 dx = ...}).
That is,
\begin{eqnarray}
\lefteqn{\biggl|E_0(f_1)(x) \, E_0(f_2)(x) +}	\\
 	& &  \sum_{j=1}^\infty (E_j(f_1)(x) - E_{j-1}(f_1)(x)) 
			(E_j(f_2)(x) - E_{j-1}(f_2)(x))\biggr|		
								\nonumber
\end{eqnarray}
is less than or equal to the product of
\begin{equation}
\Bigl(|E_0(f_1)(x)|^2 
	+ \sum_{j=1}^\infty |E_j(f_1)(x) - E_{j-1}(f_1)(x)|^2 \Bigr)^{1/2}
\end{equation}
and
\begin{equation}
\Bigl(|E_0(f_2)(x)|^2 
	+ \sum_{j=1}^\infty |E_j(f_2)(x) - E_{j-1}(f_2)(x)|^2 \Bigr)^{1/2}
\end{equation}
for all $x \in [0,1)$, and this product is the same as
\begin{equation}
 	S(f_1)(x) \, S(f_2)(x).
\end{equation}

\beginproposition
\label{q-norm of f bounded by q-norm of S(f), q > 2}
For each real number $q > 2$, there is a constant $C_3(q) > 0$
such that
\begin{equation}
	\int_{[0,1)} |f(x)|^q \, dx \le C_3(q) \int_{[0,1)} S(f)(x)^q \, dx.
\end{equation}
\end{proposition}

	To prove this, we shall use a duality argument.  Let $q > 2$
be given, and let $p$, $1 < p < \infty$, denote the exponent dual to
$q$, so that $1/p + 1/q = 1$.  Thus $p <2$. 

	It suffices to show that there is a constant $C$ such that
\begin{equation}
\label{to show -- |int f f_2| le ...}
	 \Bigl|\int_{[0,1)} f(x) \, f_2(x) \, dx \Bigr|
		\le C \, \Bigl(\int_{[0,1)} S(f)^q \, dy \Bigr)^{1/q}
			\, \Bigl(\int_{[0,1)} |f_2|^p \, dw \Bigr)^{1/p}
\end{equation}
for all functions $f_2$ on $[0,1)$ which are dyadic step functions.
(The sufficiency of this comes from choosing $f_2$ so that $f \cdot
f_2$ is or approximates $|f|^q$, and using the fact that $p (q-1) = q$
and $1/p = 1 - 1/q$.  More precisely, given any positive integer $k$,
one can choose $f_2$ to be $\overline{E_k(f)} \cdot |E_k(f)|^{q-2}$,
where this is interpreted as being $0$ when $E_k(f)$ is $0$.  Then
$\int_{[0,1)} f \cdot f_2 \, dx = \int_{[0,1)} |E_k(f)|^q \, dx$ and
$|f_2| = |E_k(f)|^{q-1}$.)

	Recall that H\"older's inequality (with this choice of $q$,
$p$) says that
\begin{equation}
\label{Holder's inequality}
	\int_{[0,1)} g(x) \, h(x) \, dx
		\le \Bigl(\int_{[0,1)} g(y)^q \, dy \Bigr)^{1/q}
			\Bigl(\int_{[0,1)} h(w)^p \, dw \Bigl)^{1/p}
\end{equation}
for arbitrary nonnegative functions $g$, $h$ on $[0,1)$.  Applying
this to the right side of (\ref{|int_{[0,1)} f_1(x) f_2(x) dx| le
int_{[0,1)} S(f_1)(x) S(f_2)(x) dx}) (with $f_1 = f$), we obtain
\begin{equation}
	 \Bigl|\int_{[0,1)} f(x) \, f_2(x) \, dx \Bigr|
		\le \Bigl(\int_{[0,1)} S(f)^q \, dy \Bigr)^{1/q}
			\, \Bigl(\int_{[0,1)} S(f_2)^p \, dw \Bigr)^{1/p}.
\end{equation}
We can now get (\ref{to show -- |int f f_2| le ...}) from this and
Proposition \ref{p-norm of S(f) bounded by p-norm of M(f), p < 2},
since $p < 2$.  (We are also using Proposition \ref{p-integral
inequalities for M(f)} here.)

\beginremark
{\rm  This general method works without the restriction to $q > 2$,
i.e., inequalities like
\begin{equation}
	\Big(\int_{[0,1)} S(f_2)^p \, dw \Bigr)^{1/p}
		\le K \, \Bigl(\int_{[0,1)} |f_2|^p \, dw \Bigl)^{1/p}
\end{equation}
(for some constant $K$ and arbitrary functions $f$) lead to inequalities
of the form
\begin{equation}
	\Big(\int_{[0,1)} |f|^q \, du \Big)^{1/q}
		\le K \, \Bigl(\int_{[0,1)} S(f)^q \, du \Bigr)^{1/q}
\end{equation}
for the conjugate exponent $q$, $1/q + 1/p = 1$ (as long as $p, q \ge 1$).
}
\end{remark}

\section{Duality, continued}
\label{Duality, continued}

\beginproposition
\label{q-norm of S(f) bounded by q-norm of f, q > 2}
For each real number $q > 2$, there is a constant $C_4(q) > 0$ such
that 
\begin{equation}
\label{int_{[0,1)} S(f)(x)^q dx le C_4(q) int_{[0,1)} |f(x)|^q dx}
	\int_{[0,1)} S(f)(x)^q \, dx \le C_4(q) \int_{[0,1)} |f(x)|^q \, dx
\end{equation}
(for arbitrary $f$ on $[0,1)$).
\end{proposition}

	Let $q > 2$ be given.  It is enough to show that there is a
constant $C_4(q) > 0$ so that
\begin{equation}
\label{int_{[0,1)} S_l(f)(x)^q dx le C_4(q) int_{[0,1)} |f(x)|^q dx}
	\int_{[0,1)} S_l(f)(x)^q \, dx \le C_4(q) \int_{[0,1)} |f(x)|^q \, dx
\end{equation}
for all integers $l \ge 0$ and all $f$.

	Let $p$ be the conjugate exponent to $q$, so that $1/p + 1/q =
1$.  It suffices to show that there is a constant $C > 0$ so that
\begin{eqnarray}
\label{int (alpha_0 E_0(f) + sum_{j=1}^l alpha_j (E_j(f) - E_{j-1}(f))) le ...}
	& & \int_{[0,1)} \Bigl(\alpha_0(x) E_0(f)(x) 
	  + \sum_{j=1}^l \alpha_j(x) (E_j(f)(x) - E_{j-1}(f)(x))\Bigr) \, dx
									\\
	& & \qquad
	\le C \, \Bigl(\int_{[0,1)} |f(y)|^q \, dy\Bigr)^{1/q}  
    \, \Bigl(\int_{[0,1)} \Bigl(\sum_{j=0}^l |\alpha_j(w)|^2\Bigr)^{p/2} 
							\, dw \Bigl)^{1/p}
								\nonumber
\end{eqnarray}
for arbitrary functions $\alpha_0(x), \alpha_1(x), \ldots \alpha_l(x)$
on $[0,1)$.  (For the sufficiency of this one can choose the
$\alpha_i$'s to be given by
\begin{equation}
	\alpha_0  = \overline{E_0(f)} \cdot S_l(f)^{q-2}
\end{equation}
and
\begin{equation}
	\alpha_j = \overline{(E_j(f) - E_{j-1}(f))} \cdot S_l(f)^{q-2}
\end{equation}
when $j \ge 1$.  In this case
\begin{eqnarray}
\lefteqn{\alpha_0 \cdot E_0(f)
	 + \sum_{j=1}^l \alpha_j \cdot (E_j(f) - E_{j-1}(f))}		\\
	& & \qquad =  S_l(f)^2 \cdot S_l(f)^{q-2} = S_l(f)^q	\nonumber
\end{eqnarray}
and
\begin{equation}
	\Bigl(\sum_{j=0}^l |\alpha_j|^2\Bigr)^{1/2} 
		= S_l(f) \cdot S_l(f)^{q-2} = S_l(f)^{q-1}.
\end{equation}
Since $p$ is conjugate to $q$, $p \, (q-1) = q$ and $1/p = 1 - 1/q$.)

	It is easy to check that
\begin{equation}
\label{int_{[0,1)} g(x) E_j(h)(x) dx = int_{[0,1)} E_j(g)(x) h(x) dx}
	\int_{[0,1)} g(x) \, E_j(h)(x) \, dx 
		= \int_{[0,1)} E_j(g)(x) \, h(x) \, dx
\end{equation}
for any functions $g$, $h$ on $[0,1)$ and nonnegative integer $j$.
Thus (\ref{int (alpha_0 E_0(f) + sum_{j=1}^l alpha_j (E_j(f) -
E_{j-1}(f))) le ...}) is equivalent to
\begin{eqnarray}
	& & \int_{[0,1)} \Bigl(E_0(\alpha_0)(x) 
 + \sum_{j=1}^l (E_j(\alpha_j)(x) - E_{j-1}(\alpha_j)(x))\Bigr) \, f(x) \, dx
									\\
	& & \qquad
	\le C \, \Bigl(\int_{[0,1)} |f(y)|^q \, dy\Bigr)^{1/q}  
    \, \Bigl(\int_{[0,1)} \Bigl(\sum_{j=0}^l |\alpha_j(w)|^2\Bigr)^{p/2} 
							\, dw \Bigl)^{1/p}
								\nonumber
\end{eqnarray}
Because of H\"older's inequality, this inequality will follow if we
can prove that
\begin{eqnarray}
	& & \Bigl(\int_{[0,1)} \Bigl|E_0(\alpha_0)(x) +
		\sum_{j=1}^l (E_j(\alpha_j)(x) - E_{j-1}(\alpha_j)(x))\Bigr|^p 
						\, dx\Bigr)^{1/p}
									\\
	& & \qquad 
  \le C' \, \Bigl(\int_{[0,1)} \Bigl(\sum_{j=0}^l |\alpha_j(w)|^2\Bigr)^{p/2} 
							\, dw \Bigl)^{1/p}
								\nonumber
\end{eqnarray}
for some constant $C'$, all $l \ge 0$, and arbitrary functions
$\alpha_0(x), \alpha_1(x), \ldots \alpha_l(x)$.

	From Proposition \ref{p-norm of M(f) bounded by p-norm of
S(f), p < 2} we know that the left side of the above inequality is
bounded by a constant times
\begin{equation}
	\Bigl(\int_{[0,1)} S \Bigl(E_0(\alpha_0) +
		\sum_{j=1}^l (E_j(\alpha_j) - E_{j-1}(\alpha_j))\Bigr)(x)^p
						\, dx \Big)^{1/p}.
\end{equation}
Thus we would like to show that this is bounded by a constant times
\begin{equation}
\label{int_{[0,1)} (sum_{j=0}^l |alpha_j(w)|^2)^{p/2} dw)^{1/p}}
	\Bigl(\int_{[0,1)} \Bigl(\sum_{j=0}^l |\alpha_j(w)|^2\Bigr)^{p/2} 
							\, dw \Bigl)^{1/p}.
\end{equation}

	On the other hand, one can check that
\begin{eqnarray}
\lefteqn{S \Bigl(E_0(\alpha_0) +
		\sum_{j=1}^l (E_j(\alpha_j) - E_{j-1}(\alpha_j))\Bigr)(x)}
									\\
	 & & = \Bigl(|E_0(\alpha_0)(x)|^2 + 
   \sum_{j=1}^l |E_j(\alpha_j)(x) - E_{j-1}(\alpha_j)(x)|^2 \Bigr)^{1/2}.
								\nonumber
\end{eqnarray}
Hence we would like to show that
\begin{equation}
    \quad
	\Bigl(\int_{[0,1)} \Bigl(|E_0(\alpha_0)(x)|^2 + 
   \sum_{j=1}^l |E_j(\alpha_j)(x) - E_{j-1}(\alpha_j)(x)|^2 \Bigr)^{p/2}
							\, dx \Bigl)^{1/p}
\end{equation}
is bounded by a constant times (\ref{int_{[0,1)} (sum_{j=0}^l
|alpha_j(w)|^2)^{p/2} dw)^{1/p}}).  This will be handled in the next section.

\section{Some inequalities}
\label{Some inequalities}

	Let $l$ be a nonnegative integer, and let $\beta_0(x),
\beta_1(x), \ldots, \beta_l(x)$ be functions on $[0,1)$.  Given $p, r
\ge 1$, consider the problem of bounding
\begin{equation}
\label{(int_{[0,1)} (sum_{j=0}^l |E_j(beta_j)(x)|^r)^{p/r} dx)^{1/p}}
    \Bigl(\int_{[0,1)} \Bigl(\sum_{j=0}^l |E_j(\beta_j)(x)|^r \Bigr)^{p/r}
						\, dx \Bigr)^{1/p}
\end{equation}
by a constant times
\begin{equation}
\label{(int_{[0,1)} (sum_{j=0}^l |beta_j(x)|^r)^{p/r} dx)^{1/p}}
	\Bigl(\int_{[0,1)} \Bigl(\sum_{j=0}^l |\beta_j(x)|^r \Bigr)^{p/r}
						\, dx \Bigr)^{1/p}
\end{equation}
(where the constant does not depend on $l$ or $\beta_0(x), \beta_1(x),
\ldots, \beta_l(x)$).

	The situation at the end of the previous section corresponds
to $r = 2$, $1 < p < 2$.  At that stage, cancellation was not really
at issue, and so we pass to the variant here.  Furthermore, the
$\beta_j$'s might as well be nonnegative for these types of
inequalities.

	If $p$ happens to be equal to $r$, then the bound holds with
constant equal to $1$, and it reduces to the inequality
\begin{equation}
	\int_{[0,1)} |E_j(\beta)(x)|^p \, dx
		\le \int_{[0,1)} |\beta(x)|^p \, dx
\end{equation}
for any function $\beta$ on $[0,1)$.  This inequality follows from
that of Jensen.  More precisely, Jensen's inequality can be applied to
get 
\begin{equation}
\label{|frac{1}{|J|} int_J beta(y)|^p le frac{1}{|J|} int_J |beta(y)|^p}
	\biggl| \frac{1}{|J|} \int_J \beta(y) \, dy \biggr|^p
		\le \frac{1}{|J|} \int_J |\beta(y)|^p \, dy
\end{equation}
for any interval $J$.  The previous inequality uses this one for the
dyadic subintervals $J$ of $[0,1)$ of length $2^{-j}$, summing over
all these dyadic subintervals to get the integrals over $[0,1)$.

	If $r = \infty$, then we interpret
\begin{equation}
	\Bigl(\sum_{j=0}^l |\beta_j(x)|^r \Bigr)^{1/r},	\quad
		\Bigl(\sum_{j=0}^l |E_j(\beta_j)(x)|^r \Bigr)^{1/r}
\end{equation}
as being
\begin{equation}
	\sup_{0 \le j \le l} |\beta_j(x)|, \quad
		\sup_{0 \le j \le l} |E_j(\beta_j)(x)|
\end{equation}
(as usual).  In this case we may as well assume that all of the
$\beta_j$'s are the same, because we could replace them all with
$\sup_{0 \le j \le l} |\beta_j(x)|$ (which would not affect
(\ref{(int_{[0,1)} (sum_{j=0}^l |beta_j(x)|^r)^{p/r} dx)^{1/p}}), and
would increase (\ref{(int_{[0,1)} (sum_{j=0}^l
|E_j(beta_j)(x)|^r)^{p/r} dx)^{1/p}}) or keep it the same).  The
bounds in question would then be essentially ones about the maximal
function $M(f)$, and in fact we have such bounds for $1 < p \le
\infty$, as in Lemma \ref{supremum bound for M(f)} and Proposition
\ref{p-integral inequalities for M(f)}.

	Now suppose that $1 \le r < p < \infty$.  Let $s \in
(1,\infty)$ be conjugate to $p/r$, in the sense that
\begin{equation}
	1/s + r/p = 1.
\end{equation}
To show that (\ref{(int_{[0,1)} (sum_{j=0}^l |E_j(beta_j)(x)|^r)^{p/r}
dx)^{1/p}}) is bounded by a constant times (\ref{(int_{[0,1)}
(sum_{j=0}^l |beta_j(x)|^r)^{p/r} dx)^{1/p}}), it suffices to show
that there is a constant $C$ so that
\begin{eqnarray}
   & & \int_{[0,1)} \Bigl(\sum_{j=0}^l |E_j(\beta_j)(x)|^r \Bigr) \, h(x) \, dx
									\\
	& & \le C \,
\Bigl(\int_{[0,1)} \Bigl(\sum_{j=0}^l |\beta_j(y)|^r \Bigr)^{p/r} 
							\, dy \Bigr)^{r/p}
	\, \Bigl(\int_{[0,1)} h(w)^s \, dw \Bigr)^{1/s}
								\nonumber
\end{eqnarray} 
for all nonnegative functions $h$ on $[0,1)$, and all $l \ge 0$ and
$\beta_0, \beta_1, \ldots, \beta_l$.

	Since $r \ge 1$, we have that
\begin{equation}
	|E_j(\beta_j)(x)|^r \le E_j(|\beta_j|^r)(x)
\end{equation}
for all $x \in [0,1)$.  This follows from Jensen's inequality, as in
(\ref{|frac{1}{|J|} int_J beta(y)|^p le frac{1}{|J|} int_J
|beta(y)|^p}) (with the $p$ there replaced by $r$).  Thus
\begin{equation}
     \int_{[0,1)} \Bigl(\sum_{j=0}^l |E_j(\beta_j)(x)|^r \Bigr) \, h(x) \, dx
  \le \int_{[0,1)} \Bigl(\sum_{j=0}^l E_j(|\beta_j|^r)(x) \Bigr) \, h(x) \, dx.
\end{equation}
On the other hand, $\int_{[0,1)} E_j(g)(x) \, h(x) \, dx = 
\int_{[0,1)} g(x) \, E_j(h)(x) \, dx$ for arbitrary functions $g$ and $h$
on $[0,1)$, and so we get that
\begin{equation}
     \int_{[0,1)} \Bigl(\sum_{j=0}^l |E_j(\beta_j)(x)|^r \Bigr) \, h(x) \, dx
  \le \int_{[0,1)} \sum_{j=0}^l |\beta_j(x)|^r \, E_j(h)(x) \, dx.
\end{equation}
The maximal function $M(h) = \sup_j E_j(h)$ can be put into the right side,
to get
\begin{equation}
     \int_{[0,1)} \Bigl(\sum_{j=0}^l |E_j(\beta_j)(x)|^r \Bigr) \, h(x) \, dx
  \le \int_{[0,1)} \sum_{j=0}^l |\beta_j(x)|^r \, M(h)(x) \, dx.
\end{equation}
Since $M(h)$ does not depend on $j$, it can be separated from the sum,
and one can apply H\"older's inequality with exponents $p/r$, $s$ to
obtain
\begin{eqnarray}
   & & \int_{[0,1)} \Bigl(\sum_{j=0}^l |E_j(\beta_j)(x)|^r \Bigr) \, h(x) \, dx
									\\
	& & \le 
\Bigl(\int_{[0,1)} \Bigl(\sum_{j=0}^l |\beta_j(y)|^r \Bigr)^{p/r} 
							\, dy \Bigr)^{r/p}
	\, \Bigl(\int_{[0,1)} M(h)(w)^s \, dw \Bigr)^{1/s}
								\nonumber
\end{eqnarray}
This is exactly what we want, except that the last factor contains
$M(h)$ instead of $h$ by itself.  One can get rid of this, with an
extra constant factor on the right side, using Proposition
\ref{p-integral inequalities for M(f)}.  (Note that $s > 1$ when $p <
\infty$.)

	Thus (\ref{(int_{[0,1)} (sum_{j=0}^l |E_j(beta_j)(x)|^r)^{p/r}
dx)^{1/p}}) is bounded by a constant times (\ref{(int_{[0,1)}
(sum_{j=0}^l |beta_j(x)|^r)^{p/r} dx)^{1/p}}) when $1 \le r < p <
\infty$.  If $1 < p < r \le \infty$, then one can get a similar
inequality by duality. In other words, if $p'$, $r'$ are the exponents
conjugate to $p$, $r$, respectively, so that
\begin{equation}
	1/p + 1/p' = 1, \quad 1/r + 1/r' = 1,
\end{equation}
then one can get the same estimate for $(p,r)$ as for $(p',r')$.  

	To make this precise, let us note that (\ref{(int_{[0,1)}
(sum_{j=0}^l |E_j(beta_j)(x)|^r)^{p/r} dx)^{1/p}}) is bounded by a
constant $A$ times (\ref{(int_{[0,1)} (sum_{j=0}^l
|beta_j(x)|^r)^{p/r} dx)^{1/p}}) (for all $\beta_0(x), \beta_1(x),
\ldots, \beta_l(x)$) if and only if
\begin{equation}
\label{|int_{[0,1)} (sum_{j=0}^l E_j(beta_j)(x) gamma_j(x) dx|}
   \Bigl| \int_{[0,1)} 
	\Bigl(\sum_{j=0}^l E_j(\beta_j)(x) \, \gamma_j(x) \Bigr) \, dx \Bigr|
\end{equation}
is bounded by $A$ times
\begin{equation}
\label{products of integrals of sums with powers..}
	\Bigl(\int_{[0,1)} \Bigl(\sum_{j=0}^l |\beta_j(x)|^r \Bigr)^{p/r}
						\, dx \Bigr)^{1/p}
	\Bigl(\int_{[0,1)} \Bigl(\sum_{j=0}^l |\gamma_j(x)|^{r'} \Bigr)^{p'/r'}
						\, dx \Bigr)^{1/p'}
\end{equation}
for all $\beta_0(x), \beta_1(x), \ldots, \beta_l(x)$ and $\gamma_0(x),
\gamma_1(x), \ldots, \gamma_l(x)$.  The ``only if'' part of statement
can be verified using H\"older's inequality twice.  For the ``if''
part, one makes suitable choices of the $\gamma_j$'s in terms of the
$E_k(\beta_k)$'s to get back to the original integrals and sums.

	As in (\ref{int_{[0,1)} g(x) E_j(h)(x) dx = int_{[0,1)}
E_j(g)(x) h(x) dx}), (\ref{|int_{[0,1)} (sum_{j=0}^l E_j(beta_j)(x)
gamma_j(x) dx|}) is equal to
\begin{equation}
\label{|int_{[0,1)} (sum_{j=0}^l beta_j(x) E_j(gamma_j)(x) dx|}
   \Bigl| \int_{[0,1)} 
	\Bigl(\sum_{j=0}^l \beta_j(x) \, E_j(\gamma_j)(x) \Bigr) \, dx \Bigr|.
\end{equation}
This permits one to derive the statement mentioned above, about having
the same estimate for $(p,r)$ as for $(p',r')$.  The assumption that
$p < r$ exactly corresponds to $r' < p'$, which brings us back to the
previous case.

	In this kind of argument, one can allow for exponents which are
infinite, with suitable adjustments to the formulae.  In particular,
suprema would be used when necessary.

\section{Another inequality for $p=1$}
\label{Another inequality for p=1}

\beginproposition [Weak-type estimate for $S(f)$]
\label{weak-type estimate for S(f)}\index{weak-type estimates}
If $\lambda$ is a positive real number and $f(x)$ is an arbitrary
function on $[0,1)$, then
\begin{equation}
\label{|{ x in [0,1) : S(f) > lambda }| le ...}
	|\{ x \in [0,1) : S(f) > \lambda \}| 
		\le \frac{3}{\lambda} \int_{[0,1)} |f(x)| \, dx.
\end{equation}
\end{proposition}

	To prove this we can use much the same argument as in Section
\ref{Estimates, 1}.  More precisely, (\ref{|{x : S(f)(x) > lambda}| le ...,3})
and Lemma \ref{|f_lambda(x)| le min(lambda, M(f)(x))} in Section
\ref{Estimates, 1} imply that
\begin{eqnarray}
\label{|{x : S(f)(x) > lambda}| le ...,3, revised version}
\lefteqn{|\{x \in [0,1) : S(f)(x) > \lambda\}|}		\\
	& & \le 2 \, |\{x \in [0,1) : M(f)(x) > \lambda\}|
			+ \lambda^{-1} \int_{[0,1)} |f_\lambda(x)| \, dx,
								\nonumber
\end{eqnarray}
where $f_\lambda(x)$ is as in (\ref{def of f_lambda(x)}).  Proposition
\ref{weak-type estimate for M(f)} in Section \ref{The size of the
maximal function} can be applied to convert this to 
\begin{eqnarray}
\label{|{x : S(f)(x) > lambda}| le ...,3, revised version 2}
\lefteqn{|\{x \in [0,1) : S(f)(x) > \lambda\}|}		\\
	& & \le 2 \, \lambda^{-1} \int_{[0,1)} |f(x)| \, dx
			+ \lambda^{-1} \int_{[0,1)} |f_\lambda(x)| \, dx.
								\nonumber
\end{eqnarray}
On the other hand,
\begin{equation}
\label{int_{[0,1)} |f_lambda(x)| dx le int_{[0,1)} |f(x)| dx}
	\int_{[0,1)} |f_\lambda(x)| \, dx \le \int_{[0,1)} |f(x)| \, dx
\end{equation}
for all $\lambda > 0$.  This can be derived from the definition
(\ref{def of f_lambda(x)}) of $f_\lambda(x)$.  More precisely,
for any interval $L$ in $[0,1)$ we have that
\begin{equation}
	\Bigl| \frac{1}{|L|} \int_L f(x) \, dx \Bigr|
	   \le \frac{1}{|L|} \int_L |f(x)| \, dx,
\end{equation}
and this can be used in the decompositions involved in the definition
of $f_\lambda$.  Once one has (\ref{int_{[0,1)} |f_lambda(x)| dx le
int_{[0,1)} |f(x)| dx}), (\ref{|{ x in [0,1) : S(f) > lambda }| le
...})  follows easily from (\ref{|{x : S(f)(x) > lambda}| le ...,3,
revised version 2}).

\section{Variants}
\label{Variants}

	The techniques applied in this chapter have versions that are
used in numerous settings in analysis.  There are also some
alternatives, in this and other situations.

	A basic point is that one frequently considers classes of
linear operators (which are compatible with some underlying geometry
or structure).  One could employ the results of this chapter to the
study of linear operators, but one can also approach linear operators
directly, through similar arguments (i.e., with similar approximations
and comparisons).  For linear operators duality can be more
straightforward, especially with operators or classes of operators
which are self-dual.  Let us note that sublinear operators such as
those examined here can often be viewed in terms of linear operators
which take ordinary scalar-valued functions to functions with values
in a normed vector space.  This perspective is sometimes helpful, and
more generally one can look at operators which take vector-valued
functions to vector-valued functions.

	A very famous linear operator is the \emph{Hilbert transform}
or \emph{conjugation operator}, acting on functions on the unit
circle.  This operator can be described in terms of complex analysis
as the operator on boundary value functions which corresponds to
harmonic conjugation of harmonic functions on the unit disk in the
complex plane.  This operator also has a simple description in terms
of Fourier series, and from it one can obtain operators which truncate
a given Fourier series (so that bounds for the Hilbert transform lead
to bounds for the truncation operators).  There is a version of this
operator on the real line, with analogous features in terms of
harmonic conjugation and Fourier transforms (and explicit links to the
operator on the unit circle).

	Classically, connections with complex analysis were used as a
very strong tool in the study of certain special operators like the
Hilbert transform.  The results of this could then be applied to
Fourier series, for instance, where complex analysis did not seem to
be involved, and to other operators which could be treated in terms of
these building blocks.

	This is quite beautiful and remarkable, and has numerous
interesting aspects.  On the other hand, there was also the problem of
having more direct real-variable methods, which should both be
interesting in their own right and have further applications.

	The paper \cite{CZ} of Calder\'on and Zygmund was a
fundamental step in this development.  The arguments in this chapter
have followed the same rough plan (for which there are versions in a
number of settings).  In general, there are mild additional terms and
other adjustments involved, compared to the situation here.  Depending
on the circumstances, there can also be special structure which allows
simplification in different ways.  There are other circumstances in
which analogous questions come up and the basic methods are not
directly sufficient or applicable, and a variety of refinements and
additional tools have emerged.

\section{Some remarks concerning $p = 4$}
\label{Some remarks concerning $p = 4$}

	Suppose that $f(x)$ is a function on $[0,1)$, and consider
the integral
\begin{equation}
\label{int_{[0,1)} S(f)(x)^4 dx}
	\int_{[0,1)} S(f)(x)^4 \, dx.
\end{equation}
It turns out that this can be written as $A + 2 \, B$, where
\begin{equation}
\label{A = .....}
	A = \int_{[0,1)} \Bigl(|E_0(f)(x)|^4 + 
	      \sum_{j=1}^\infty |E_j(f)(x) - E_{j-1}(f)(x)|^4 \Bigr)
					\, dx
\end{equation}
and
\begin{eqnarray}
\label{B = .....}
\lefteqn{\qquad
	B = \int_{[0,1)} |E_0(f)(x)|^2 \, E_0(|f-E_0(f)|^2)(x) \, dx}  \\
	& & + \int_{[0,1)} 
		\sum_{j=1}^\infty |E_j(f)(x) - E_{j-1}(f)(x)|^2
				\, E_j(|f-E_j(f)|^2)(x)  \, dx.
							\nonumber
\end{eqnarray}
Let us first sketch the computations behind this, leaving the details
as an exercise, and then discuss some of its consequences.

	One can start with the definition (\ref{def of S(f)}) of
$S(f)(x)$, which expresses $S(f)(x)^2$ as a sum.  Thus $S(f)(x)^4$ can
be given as a square of a sum, which can be expanded into a double sum
of the products of the terms in the initial sum.  If $j$ and $k$ are
the indices in the double sum, then the double sum can be decomposed
into three parts, corresponding to $j = k$, $j < k$, and $j > k$.  The
integral of the $j = k$ part is given exactly by $A$ above, by
inspection.  The $j < k$ and $j > k$ parts are equal to each other, by
the symmetry of the sum.  It remains to show that the integral over
$[0,1)$ of the $j < k$ part is equal to $B$ above.  This equality does
not work directly at the level of the integrands, as for the $j = k$
part, but uses the integration.

	For each $j \ge 0$, let $B_j$ be the part of $B$ that
corresponds to $j$ in the obvious manner.  It is enough to show that
$B_j$ is equal to the integral over $[0,1)$ of the piece of the
original double sum that comes from restricting ourselves to this
particular $j$ and summing over $k > j$.

	Fix $j \ge 0$, and let $I$ be a dyadic subinterval of $[0,1)$
of length equal to $2^{-j}$.  A key point is that
\begin{eqnarray}
\label{2-norm formula for sum over k > j}
\lefteqn{\int_I \sum_{k=j+1}^\infty |E_k(f)(x) - E_{k-1}(f)(x)|^2 
		\, dx}						\\
	&& = \int_I |f(x) - E_j(f)(x)|^2 \, dx.		\nonumber
\end{eqnarray}
Note that $E_j(f)(x)$ is in fact constant on $I$, since $|I| =
2^{-j}$.  This equality comes from orthogonality of the $(E_k(f) -
E_{k-1}(f))$'s, $k > j$, on $I$, as in the situation of Lemma
\ref{2-norm of S(f)}.

	The formula for $B_j$ that we want says that one integral over
$[0,1)$ is the same as another integral over $[0,1)$.  To establish
this it suffices to show that the corresponding integrals over dyadic
intervals of length $2^{-j}$ are the same.  Fix such an interval $I$.
Inside both integrals one has the expression $|E_0(f)(x)|^2$ when $j =
0$ and $|E_j(f)(x) - E_{j-1}(f)(x)|^2$ when $j \ge 1$.  These
expressions are constant on $I$, since $|I| = 2^{-j}$.  Hence these
expressions can be pulled out of the integrals whose equality is to be
shown.  The integrals that are left reduce to the two sides of
(\ref{2-norm formula for sum over k > j}), and this yields the
desired equality.

	This completes the sketch of the proof that (\ref{int_{[0,1)}
S(f)(x)^4 dx}) is equal to $A + 2 \, B$, where $A$ and $B$ are as in
(\ref{A = .....}) and (\ref{B = .....}).  As a consequence of this,
it is not hard to see that there is a constant $C > 0$ so that
\begin{equation}
\label{int_{[0,1)} S(f)(x)^4 dx le ....}
	\int_{[0,1)} S(f)(x)^4 \, dx 
		\le C \, \int_{[0,1)} S(f)(x)^2 \, M(|f|^2)(x) \, dx
\end{equation}
for arbitrary functions $f(x)$ on $[0,1)$.  (Note that $M(f)(x)^2 \le
M(|f|^2)(x)$ for all $x$, by Jensen's inequality.)

	Using (\ref{int_{[0,1)} S(f)(x)^4 dx le ....}), one can obtain
that
\begin{equation}
	\int_{[0,1)} S(f)(x)^4 \, dx 
		\le C' \int_{[0,1)} |f(x)|^4 \,	dx
\end{equation}
for some constant $C' > 0$ and arbitrary $f$.  This also relies on the
Cauchy-Schwarz inequality, and on Proposition \ref{p-integral
inequalities for M(f)} in Section \ref{More on the size of the maximal
function} with $p = 2$ and $|f|^2$ instead of $f$ to estimate
$\int_{[0,1)} M(|f|^2)(x)^2 \, dx$ in terms of $\int_{[0,1)} |f(x)|^4
\, dx$.

	In other words, we get (\ref{int_{[0,1)} S(f)(x)^q dx le
C_4(q) int_{[0,1)} |f(x)|^q dx}) in Proposition \ref{q-norm of S(f)
bounded by q-norm of f, q > 2} in Section \ref{Duality, continued} for
$q = 4$.  Once one knows this, one can prove the analogue of
Proposition \ref{p-norm of S(f) bounded by p-norm of M(f), p < 2} in
Section \ref{Estimates, 1} for $p < 4$ instead of $p < 2$, through
essentially the same method as before.

	Similarly, the interpolation results in Chapter \ref{An
interpolation theorem} (especially Section \ref{A generalization
(interpolation)}) permit one to derive inequalities for $2 \le p \le
4$ from the cases of $p = 2, 4$ (and this works in a very general
manner).

	This type of approach works in a number of settings.  That is,
for special $p$'s there are special arguments or computations, and
then one can extend from there using other methods.  A famous instance
of this occurs in M.~Riesz's estimates \cite{MRiesz} for the Hilbert
transform (mentioned in the previous section).

\chapter{Interpolation of operators}
\label{An interpolation theorem}

	We follow here the method of \cite{MRiesz}.\index{Riesz
interpolation of operators}\index{interpolation of operators!Riesz}
(See \cite{BL, Peetre, Sadosky, SW-book, Z1} for further information.)

\section{The basic result}
\label{The basic result (interpolation)}

	Let $(a_{j, k})$ be an $n \times n$ matrix of complex numbers.
We associate to this a bilinear form $A(x,y)$ defined for $x, y \in {\bf C}^n$
by
\begin{equation}
\label{def of A(x,y)}
	A(x,y) = \sum_{j=1}^n \sum_{k=1}^n y_j \, a_{j, k} \, x_k.
\end{equation}

	For $1 \le p \le \infty$ define $M_p$ as follows.  If $1 < p < \infty$,
then we set
\begin{eqnarray}
\label{def of M_p, 1 < p < infty}
\lefteqn{\quad M_p =}		\\
	&&  \sup \biggl\{ |A(x,y)| : x, y \in {\bf C}^n, 
		\Bigl(\sum_{k=1}^n |x_k|^p \Bigr)^{1/p} \le 1,
		\Bigl(\sum_{j=1}^n |y_j|^{p'} \Bigr)^{1/p'} \le 1 \biggr\}.
								\nonumber
\end{eqnarray}
Throughout this chapter $p'$ denotes the conjugate exponent of $p$, so
that
\begin{equation}
\label{conjugate exponent condition}
	\frac{1}{p} + \frac{1}{p'} = 1.
\end{equation}
Thus $1 < p' < \infty$ when $1 < p < \infty$.  If $p = 1$, so that $p'
= \infty$, then we set
\begin{equation}
\label{def of M_p, p = 1}
	M_1 = \sup \biggl\{ |A(x,y)| : x, y \in {\bf C}^n, 
		\sum_{k=1}^n |x_k| \le 1,
		\max_{1 \le j \le n} |y_j| \le 1 \biggr\}.
\end{equation}
If $p = \infty$, so that $p' = 1$, then we put
\begin{equation}
\label{def of M_p, p = infty}
	M_\infty = \sup \biggl\{ |A(x,y)| : x, y \in {\bf C}^n, 
		\max_{1 \le k \le n} |x_k| \le 1,
		\sum_{j=1}^n |y_j| \le 1 \biggr\}.
\end{equation}
In each case, $M_p$ is a nonnegative real number.

\begintheorem
\label{M. Riesz' convexity theorem}
The function $\log M_p$ is convex as a function of $1/p \in [0,1]$.
That is, if $1 \le p < q \le \infty$ and $0 < t < 1$, and if $r$
satisfies
\begin{equation}
	\frac{1}{r} = \frac{t}{p} + \frac{1-t}{q},
\end{equation}
then
\begin{equation}
\label{M_r le M_p^t M_q^{1-t}}
	M_r \le M_p^t \, M_q^{1-t}.
\end{equation}
\end{theorem}

	The analogous result holds (and by essentially the same proof)
if the $a_{j, k}$'s are all real numbers, and we restrict our
attention to $x$, $y$ in ${\bf R}^n$ in the definition of $M_p$.

	Note that if $M_p = 0$ for some $p$, then $A \equiv 0$ and
$M_p = 0$ for all $p$.  One may as well assume that $A \not\equiv 0$
in the arguments that follow.

	We shall discuss the proof of this theorem in this and the next
sections.  We begin by observing that $M_p$ can also be given as
\begin{equation}
\label{second version of M_p, p < infty}
	M_p = \sup \biggl\{ 
 \Bigl(\sum_{j=1}^n \, \Bigl|\sum_{k=1}^n a_{j, k} \, x_k \Bigr|^p\Bigr)^{1/p}
	   : x \in {\bf C}^n, \Bigl(\sum_{k=1}^n |x_k|^p \Bigr)^{1/p} \le 1
								\biggr\}
\end{equation}
when $1 \le p < \infty$, and, for $p = \infty$, as
\begin{equation}
	M_\infty = \sup \biggl\{ 
 	\max_{1 \le j \le n} \, \Bigl|\sum_{k=1}^n a_{j, k} \, x_k \Bigr|
	   : x \in {\bf C}^n, \max_{1 \le k \le n} |x_k| \le 1 \biggr\}.
\end{equation}
Indeed, this definition of $M_p$ is greater than or equal to the
previous one because of H\"older's inequality.  Conversely, for each
$x \in {\bf C}^n$ one can choose $y \in {\bf C}^n$ so that the
inequality in H\"older's inequality becomes an equality, and so that
\begin{equation}
	\Bigl(\sum_{j=1} |y_j|^{p'} \Bigr)^{1/p'} = 1
\end{equation}
when $1 \le p' < \infty$, and 
\begin{equation}
	\max_{1 \le j \le n} |y_j| = 1
\end{equation}
when $p' = \infty$.  This is not hard to verify, and it implies that the
original definition of $M_p$ is greater than or equal to the version above.
Therefore the two are equal.

	Similarly,
\begin{equation}
\label{third version of M_p, p > 1}
	M_p = \sup \biggl\{ 
\Bigl(\sum_{k=1}^n \, 
	\Bigl|\sum_{j=1}^n y_j \, a_{j, k} \Bigr|^{p'}\Bigr)^{1/p'}
	   : y \in {\bf C}^n, \Bigl(\sum_{j=1}^n |y_j|^p \Bigr)^{1/p'} \le 1
								\biggr\}
\end{equation}
when $1 < p \le \infty$ (so that $p' < \infty$), and
\begin{equation}
	M_1 = \sup \biggl\{ 
 	\max_{1 \le k \le n} \, \Bigl|\sum_{j=1}^n y_j \, a_{j, k} \Bigr|
	   : y \in {\bf C}^n, \max_{1 \le j \le n} |y_j| \le 1 \biggr\}.
\end{equation}
The equivalence of these expressions for $M_p$ is a special case of
duality.

\beginremark
{\rm For $p = 1$ and $q = \infty$ there is an earlier and simpler
approach of Schur.  A key point is that
\begin{equation}
\label{formula for M_1}
	M_1 = \max_{1 \le k \le n} \sum_{j=1}^n |a_{j,k}|
\end{equation}
and
\begin{equation}
\label{formula for M_infty}
	M_\infty = \max_{1 \le j \le n} \sum_{k=1}^n |a_{j,k}|.
\end{equation}
(Exercise.)  Once one has these formulae, one can estimate $M_r$ for
$1 < r < \infty$ using Jensen's inequality or H\"older's inequality.
In general, $M_p$ does not admit such a nice expression, and Theorem
\ref{M. Riesz' convexity theorem} works around this problem.}
\end{remark}

\beginremark 
{\rm The quantity $M_p$ is the same as the operator norm of the linear
transformation associated to the matrix $a_{j,k}$ using the $p$-norm
$\|\cdot \|_p$ as in Section \ref{definitions, etc. (normed vector
spaces)} on both the domain and the range.  Instead of using the same
$p$-norm on the domain and the range, suppose that one considers the
operator norm defined using the $p_1$-norm on the domain and the
$p_2$-norm on the range.  In the context of Theorem \ref{M. Riesz'
convexity theorem}, one would also consider the operator norm using a
$q_1$-norm in the domain and a $q_2$-norm in the range, and one would
seek an interpolation inequality for the operator norm using an
$r_1$-norm on the domain and an $r_2$-norm in the range, where $r_1$
and $r_2$ are related to $p_1$, $q_1$ and $p_2$, $q_2$ in the same
manner as above (using the same $t$ for $r_1$, $p_1$, $q_1$ and for
$r_2$, $p_2$, $q_2$).

	Riesz's proof in \cite{MRiesz} allows for this more general
situation under the assumption that $p_1 \le p_2$ and $q_1 \le q_2$.
There are other methods which permit one to deal with arbitrary $p_1$,
$p_2$, $q_1$, $q_2$.  Note, however, that the condition $p_1 \le p_2$,
$q_1 \le q_2$ holds in many natural situations.}
\end{remark}

\section{A digression about convex functions}
\label{A digression about convex functions}

	Let $a$, $b$ be real numbers with $a < b$, and let $f(x)$
be a continuous real-valued function on $[a,b]$.  Suppose that
\begin{equation}
\label{f(frac{x+y}{2}) le frac{f(x) + f(y)}{2}}
	f\Bigl(\frac{x+y}{2}\Bigr) \le \frac{f(x) + f(y)}{2}
\end{equation}
for all $x, y \in [a,b]$.  Then $f$ satisfies the ``complete'' convexity
property 
\begin{equation}
\label{f(lambda x + (1-lambda) y) le lambda f(x) + (1-lambda) f(y)}
	f(\lambda \, x + (1-\lambda) \, y) 
		\le \lambda \, f(x) + (1-\lambda) \, f(y)
\end{equation}
for all $\lambda \in (0,1)$ and $x, y \in [a,b]$.  This is the same as
Exercise 24 on p101 of \cite{Ru1} (except for a minor change), and it
is not hard to prove.  For instance, one can repeat
(\ref{f(frac{x+y}{2}) le frac{f(x) + f(y)}{2}}) to obtain
(\ref{f(lambda x + (1-lambda) y) le lambda f(x) + (1-lambda) f(y)})
first with $\lambda$ equal to $1/4$ or $3/4$, in addition to $1/2$,
then with $\lambda$ of the form $j/8$, and so on.  Thus (\ref{f(lambda
x + (1-lambda) y) le lambda f(x) + (1-lambda) f(y)}) in fact holds when
$\lambda$ is of the form $k/2^l$, and this permits one to derive
(\ref{f(lambda x + (1-lambda) y) le lambda f(x) + (1-lambda) f(y)})
for all $\lambda \in (0,1)$, by continuity.

	Here is a more general version of this observation.

\beginlemma
\label{lemma about convexity conditions}
Suppose that $a$, $b$ are real numbers with $a < b$, and that $f(x)$
is a real-valued function on $[a,b]$ which is continuous.  Assume that
for each $x, y \in [a,b]$ there is a $\lambda_{x,y} \in (0,1)$ such
that
\begin{equation}
	f(\lambda_{x,y} \, x + (1-\lambda_{x,y}) \, y)
		\le \lambda_{x,y} \, f(x) + (1-\lambda_{x,y}) \, f(y).
\end{equation}
Then (\ref{f(lambda x + (1-lambda) y) le lambda f(x) + (1-lambda) f(y)})
holds for all $\lambda \in (0,1)$ and $x, y \in [a,b]$.
\end{lemma}

	To prove this, define $L_{x,y}$ for all $x, y \in [a,b]$ by
\begin{equation}
	\quad L_{x,y} = \{\lambda \in [0,1] : 
			f(\lambda \, x + (1-\lambda) \, y)
				\le \lambda \, f(x) + (1-\lambda) \, f(y)\}.
\end{equation}
Thus $0, 1 \in L_{x,y}$ automatically, and $\lambda_{x,y} \in L_{x,y}
\cap (0,1)$ by hypothesis.  The continuity of $f$ implies that
$L_{x,y}$ is always a closed subset of $[0,1]$.  We want to show that
$L_{x,y} = [0,1]$ for all $x, y \in [a,b]$.

	Fix $x, y \in [a,b]$, and suppose to the contrary that
$L_{x,y} \ne [0,1]$.  Because $(0,1) \backslash L_{x,y}$ is an open
set, it follows that there are $\lambda_1, \lambda_2 \in L_{x,y}$ such
that $\lambda_1 < \lambda_2$ and
\begin{equation}
\label{(lambda_1, lambda_2) subseteq (0,1) backslash L_{x,y}}
	(\lambda_1, \lambda_2) \subseteq (0,1) \backslash L_{x,y}.
\end{equation}

	Define $u, v \in [a,b]$ by
\begin{equation}
	u = \lambda_1 \, x + (1-\lambda_1) \, y, \quad
		v = \lambda_2 \, x + (1-\lambda_2) \, y.
\end{equation}
By hypothesis, there is a $\lambda_{u,v} \in (0,1)$ which lies
in $L_{u,v}$, so that
\begin{equation}
	f(\lambda_{u,v} \, u + (1-\lambda_{u,v}) \, v)
		\le \lambda_{u,v} \, f(u) + (1-\lambda_{u,v}) \, f(v).
\end{equation}
We also know that $\lambda_1$, $\lambda_2$ lie in $L_{x,y}$, and hence
\begin{equation}
	\qquad f(u) \le \lambda_1 \, f(x) + (1-\lambda_1) \, f(y), \quad
		f(v) \le \lambda_2 \, f(x) + (1-\lambda_2) \, f(y).
\end{equation}
Therefore,
\begin{eqnarray}
\lefteqn{\qquad\qquad f(\lambda_{u,v} \, u + (1-\lambda_{u,v}) \, v) \le} 
									\\
 & & 
\Bigl(\lambda_{u,v} \, \lambda_1 + (1-\lambda_{u,v}) \, \lambda_2\Bigr) 
	\, f(x) + 
\Bigl(\lambda_{u,v} \, (1-\lambda_1) + (1-\lambda_{u,v}) \, (1-\lambda_2)\Bigr)
	\, f(y).
								\nonumber
\end{eqnarray}
Set $\alpha = \lambda_{u,v} \, \lambda_1 + (1-\lambda_{u,v}) \lambda_2$.
Then $\alpha \in (\lambda_1, \lambda_2)$, and the preceding inequality
can be rewritten as
\begin{equation}
	f(\lambda_{u,v} \, u + (1-\lambda_{u,v}) \, v)
		\le \alpha \, f(x) + (1-\alpha) \, f(y).
\end{equation}
Similarly,
\begin{equation}
	\lambda_{u,v} \, u + (1-\lambda_{u,v}) \, v
		= \alpha \, x + (1-\alpha) \, y,
\end{equation}
and hence
\begin{equation}
	f(\alpha \, x + (1-\alpha) \, y) 
		\le \alpha \, f(x) + (1-\alpha) \, f(y).
\end{equation}
In other words, $\alpha \in L_{x,y}$.  This contradicts
(\ref{(lambda_1, lambda_2) subseteq (0,1) backslash L_{x,y}}),
and the lemma follows.

	As a consequence of the lemma, in order to establish Theorem
\ref{M. Riesz' convexity theorem}, it suffices to show that for each
$p, q \in [1,\infty]$ there is a $t \in (0,1)$ so that (\ref{M_r le
M_p^t M_q^{1-t}}) holds.  This uses the observation that $M_s$ is a
continuous function of $1/s$, $1/s \in [0,1]$, which is not difficult
to show.

\section{A place where the maximum is attained}
\label{A place where the maximum is attained}

	We continue with the same notation and assumptions from 
Section \ref{The basic result (interpolation)}.

	Fix a real number $r$, $1 < r < \infty$, and let $r'$ be its
conjugate exponent.  There exist $x^0, y^0 \in {\bf C}^n$ such that
\begin{equation}
\label{maximizing A(x,y)}
	|A(x^0, y^0)| = M_r
\end{equation}
and
\begin{eqnarray}
\label{x^0 has norm 1}
	\Bigl(\sum_{k=1}^n |x^0_k|^r \Bigr)^{1/r} = 1, 			\\
\label{y^0 has norm 1}
	\Bigl(\sum_{j=1}^n |y^0_j|^{r'} \Bigr)^{1/r'} = 1.
\end{eqnarray}
In other words, the supremum in the definition (\ref{def of M_p, 1 < p
< infty}) of $M_r$ is attained at $x^0$, $y^0$.  We are lead to the
conditions (\ref{x^0 has norm 1}), (\ref{y^0 has norm 1}) because
otherwise we could multiply $x^0$ or $y^0$ by a real number greater
than $1$, and $A(x^0, y^0)$ would be similarly multiplied.

	For this choice of $x^0$ and $y^0$ we have that
\begin{eqnarray}
\label{|A(x^0, y^0)| = ..., 1}
	|A(x^0, y^0)| 
  & = & \Bigl| \sum_{j=1}^n \sum_{k=1}^n y^0_j \, a_{j, k} \, x^0_k \Bigr|
									\\
  & = & \Bigl(\sum_{i=1}^n |y^0_i|^{r'} \Bigr)^{1/r'}
\Bigl(\sum_{j=1}^n \Bigl| \sum_{k=1}^n a_{j, k} \, x^0_k \Bigr|^r \Bigr)^{1/r}.
								\nonumber
\end{eqnarray}
More precisely, if the second equality is replaced by $\le$, then it
always holds, by H\"older's inequality.  We have equality in this case
because otherwise we could replace $y^0$ with an element of ${\bf
C}^n$ for which equality does occur (and with (\ref{y^0 has norm 1})),
and the value of $|A(x^0, y^0)|$ would increase, i.e., it would not be
the maximum.

	For the same reason, 
\begin{eqnarray}
\label{|A(x^0, y^0)| = ..., 2}
	|A(x^0, y^0)| 
  & = & \Bigl| \sum_{j=1}^n \sum_{k=1}^n y^0_j \, a_{j, k} \, x^0_k \Bigr|
									\\
  & = &   \Bigl(\sum_{k=1}^n \Bigl| \sum_{j=1}^n y^0_j \, a_{j, k} \Bigr|^{r'}
								  \Bigr)^{1/r'}
	\Bigl(\sum_{l=1}^n |x^0_l|^r \Bigr)^{1/r}.
								\nonumber
\end{eqnarray}

\section{The rest of the argument}
\label{The rest of the argument}

	Because of the second equality in (\ref{|A(x^0, y^0)| = ..., 1}),
we have that
\begin{equation}
\label{|sum_{k=1}^n a_{j, k} x^0_k| = mu |y^0_j|^{r' - 1}}
	\Bigl|\sum_{k=1}^n a_{j, k} \, x^0_k \Bigr| = \mu \, |y^0_j|^{r' - 1}
\end{equation}
for $j = 1, 2, \ldots, n$, where $\mu$ is a constant.  This can be
obtained from the proof of H\"older's inequality (and the conditions
for equality in the arithmetic-geometric mean inequalities).  Similarly,
\begin{equation}
\label{|sum_{j=1}^n y^0_j a_{j, k}| = nu |x^0_k|^{r-1}}
	\Bigl|\sum_{j=1}^n y^0_j \, a_{j, k} \Bigr| = \nu \, |x^0_k|^{r-1}
\end{equation}
for $k = 1, 2, \ldots, n$ and some constant $\nu$.

	From (\ref{|A(x^0, y^0)| = ..., 1}) and (\ref{|A(x^0, y^0)| = ..., 2})
we have that
\begin{equation}
	M_r = 
\Bigl(\sum_{j=1}^n \Bigl| \sum_{k=1}^n a_{j, k} \, x^0_k \Bigr|^r \Bigr)^{1/r}
	    =
\Bigl(\sum_{k=1}^n \Bigl| \sum_{j=1}^n y^0_j \, a_{j, k} \Bigr|^{r'}
								 \Bigr)^{1/r'}.
\end{equation}
This uses the fact that $M_r = |A(x^0, y^0)|$, by definitions, and
(\ref{x^0 has norm 1}), (\ref{y^0 has norm 1}).  Combining the first
equality with (\ref{|sum_{k=1}^n a_{j, k} x^0_k| = mu |y^0_j|^{r' - 1}})
we obtain
\begin{equation}
	M_r = \mu \, \Bigl(\sum_{j=1}^n |y^0_j|^{r (r'-1)}\Bigr)^{1/r}.
\end{equation}
On the other hand, $r (r' - 1) = r'$, since $1/r + 1/r' = 1$.  Thus we
may apply (\ref{y^0 has norm 1}) to get 
\begin{equation}
	M_r = \mu.
\end{equation}
For the same reasons,
\begin{equation}
	M_r = \nu.
\end{equation}

	Now suppose that $p$, $q$, and $t$ are real numbers such that
$1 \le p < q \le \infty$, $0 < t < 1$, and $1/r = t/p + (1-t)/q$.
Then
\begin{equation}
\label{(sum_{j=1}^n |sum_{k=1}^n a_{j, k} x^0_k|^p)^{1/p} le ..}
\Bigl(\sum_{j=1}^n \Bigl| \sum_{k=1}^n a_{j, k} \, x^0_k \Bigr|^p \Bigr)^{1/p}
	\le M_p \, \Bigl(\sum_{k=1}^n |x^0_k|^p \Bigr)^{1/p}
\end{equation}
and
\begin{equation}
\label{(sum_{k=1}^n |sum_{j=1}^n y^0_j a_{j, k}|^{q'})^{1/q'} le ..}
\Bigl(\sum_{k=1}^n \Bigl| \sum_{j=1}^n y^0_j \, a_{j, k} \Bigr|^{q'}
								 \Bigr)^{1/q'}
	\le M_q \, \Bigl(\sum_{j=1}^n |y^0_j|^{q'} \Bigr)^{1/q'},
\end{equation}
where $q'$ is the conjugate exponent of $q$.  (Note that $p, q' <
\infty$ under the present conditions.)  These inequalities follow from
(\ref{second version of M_p, p < infty}) and (\ref{third version of
M_p, p > 1}).

	Applying our earlier computations we get that
\begin{equation}
\label{M_r -- M_p inequality}
	M_r \, \Bigl(\sum_{j=1}^n |y^0_j|^{p (r'-1)}\Bigr)^{1/p}
	   \le M_p \, \Bigl(\sum_{k=1}^n |x^0_k|^p \Bigr)^{1/p}
\end{equation}
and
\begin{equation}
\label{M_r -- M_q inequality}
	M_r \, \Bigl(\sum_{k=1}^n |x^0_k|^{q' (r-1)} \Bigr)^{1/q'}
	   \le M_q \, \Bigl(\sum_{j=1}^n |y^0_j|^{q'} \Bigr)^{1/q'}.
\end{equation}
We can take the $t$th and $(1-t)$th powers of these inequalities and
then multiply to get
\begin{eqnarray}
\label{M_r -- M_p, M_q inequality}
\lefteqn{M_r \, \Bigl(\sum_{j=1}^n |y^0_j|^{p (r'-1)}\Bigr)^{t/p}
	\Bigl(\sum_{k=1}^n |x^0_k|^{q' (r-1)} \Bigr)^{(1-t)/q'}} 	\\
	& & \le M_p^t \, M_q^{1-t} \, \Bigl(\sum_{k=1}^n |x^0_k|^p \Bigr)^{t/p}
		\Bigl(\sum_{j=1}^n |y^0_j|^{q'} \Bigr)^{(1-t)/q'}.
								\nonumber
\end{eqnarray}

	We are going to need some identities with indices.  Let us
first check that
\begin{equation}
\label{t (frac{1}{r} - frac{1}{p}) = (1-t) (frac{1}{r'} - frac{1}{q'})}
	t \Bigl(\frac{1}{r} - \frac{1}{p}\Bigr)
		= (1-t) \Bigl(\frac{1}{r'} - \frac{1}{q'}\Bigr).
\end{equation}
Because $1/r = t/p + (1-t)/q$, we have that
\begin{equation}
	t \Bigl(\frac{1}{r} - \frac{1}{p}\Bigr)
		= t \Bigl(\frac{t-1}{p} + \frac{1-t}{q}\Bigr).
\end{equation}
Similarly, $1/r' = t/p' + (1-t)/q'$ (simply by using the formula for
$1/r$ in $1/r' = 1 - 1/r$), and this leads to
\begin{equation}
	(1-t) \Bigl(\frac{1}{r'} - \frac{1}{q'}\Bigr)
		= (1-t) \Bigl(\frac{t}{p'} - \frac{t}{q'}\Bigr)
		= (1-t) \Bigl(\frac{-t}{p} + \frac{t}{q}\Bigr).
\end{equation}
The last step uses $1/p' = 1 - 1/p$, $1/q' = 1 - 1/q$.  The identity
(\ref{t (frac{1}{r} - frac{1}{p}) = (1-t) (frac{1}{r'} -
frac{1}{q'})}) follows directly from these computations.

	Suppose that
\begin{equation}
\label{frac{t}{r} = frac{1-t}{r'}}
	\frac{t}{r} = \frac{1-t}{r'}.
\end{equation}
Then (\ref{t (frac{1}{r} - frac{1}{p}) = (1-t) (frac{1}{r'} -
frac{1}{q'})}) implies that $t/p = (1-t)/q'$.  This implies in turn that
\begin{equation}
\label{p (r'-1) = q', q' (r-1) = p}
	p (r'-1) = q', \qquad q' (r-1) = p.
\end{equation}
Indeed, $r' - 1 = r'/r$ and $r - 1 = r/r'$, since $1/r + 1/r' = 1$.
Thus $r' - 1 = (1-t)/t$ and $r - 1 = t/(1-t)$, by (\ref{frac{t}{r} =
frac{1-t}{r'}}).  This leads to (\ref{p (r'-1) = q', q' (r-1) = p}).

	Thus, under the assumption (\ref{frac{t}{r} = frac{1-t}{r'}}),
we obtain
\begin{equation}
\label{M_r le M_p^t M_q^{1-t}, assuming frac{t}{r} = frac{1-t}{r'}}
	M_r \le M_p^t \, M_q^{1-t}
\end{equation}
from (\ref{M_r -- M_p, M_q inequality}).  That is, the sums involving
$x^0$ and $y^0$ on the left and right sides of (\ref{M_r -- M_p, M_q
inequality}) exactly match up under these conditions, by the computations
in the previous paragraph.

	Note that for each $p$, $q$ such that $1 \le p < q \le
\infty$, there always is a $t \in (0,1)$ so that if $r$ is given by
$1/r = t/p + (1-t)/q$, then (\ref{frac{t}{r} = frac{1-t}{r'}}) is
satisfied.  Theorem \ref{M. Riesz' convexity theorem} follows, as
indicated at the end of Section \ref{A digression about convex
functions}.

\section{A reformulation}
\label{A reformulation (of Riesz interpolation)}

	Let $T$ be a linear mapping from the vector space of dyadic step
functions on $[0,1)$ into itself.  Suppose that $p$ and $q$ are real numbers
with $1 \le p < q \le \infty$, and that there are nonnegative real numbers
$N_p$, $N_q$ such that
\begin{equation}
\label{(int_{[0,1)} |T(f)(x)|^p dx)^{1/p} le ..}
	\Bigl(\int_{[0,1)} |T(f)(x)|^p \, dx \Bigr)^{1/p}
		\le N_p \, \Big(\int_{[0,1)} |f(x)|^p \, dx \Bigr)^{1/p}
\end{equation}
and
\begin{equation}
\label{(int_{[0,1)} |T(f)(x)|^q dx)^{1/q} le ..}
	\Bigl(\int_{[0,1)} |T(f)(x)|^q \, dx \Bigr)^{1/q}
		\le N_q \, \Big(\int_{[0,1)} |f(x)|^q \, dx \Bigr)^{1/q}
\end{equation}
for all $f$.  If $q = \infty$, then the latter should be replaced with
\begin{equation}
\label{sup_{x in [0,1)} |T(f)(x)| le N_infty sup_{[0,1)} |f(x)|}
	\sup_{x \in [0,1)} |T(f)(x)| \le N_\infty \, \sup_{[0,1)} |f(x)|,
\end{equation}
as usual.

	Let $t \in (0,1)$ be given, and define $r$ by $1/r = t/p + (1-t)/q$.
Then
\begin{equation}
\label{estimate for int |T(f)(x)|^r}
	\Bigl(\int_{[0,1)} |T(f)(x)|^r \, dx \Bigr)^{1/r}
   \le N_p^t \, N_q^{1-t} \, \Big(\int_{[0,1)} |f(x)|^r \, dx \Bigr)^{1/r}.
\end{equation}
for all dyadic step functions $f$.

	This assertion can be derived from Theorem \ref{M. Riesz'
convexity theorem}, as follows.  Let $m$ be a positive integer.  We
have that
\begin{equation}
	\Bigl(\int_{[0,1)} |E_m(T(f))(x)|^p \, dx \Bigr)^{1/p}
		\le N_p \, \Big(\int_{[0,1)} |f(x)|^p \, dx \Bigr)^{1/p}
\end{equation}
for all $f$, and analogously for $q$ instead of $p$.  Indeed,
\begin{equation}
	\int_{[0,1)} |E_m(g)(x)|^p \, dx \le \int_{[0,1)} |g(x)|^p \, dx
\end{equation}
for all $g$, as can be derived from Jensen's inequality.  We also have
that $\sup_{[0,1)} |E_m(g)| \le \sup_{[0,1)} |g|$.  This implies that
our hypotheses for $T$ hold as well for $E_m \circ T$.

	Theorem \ref{M. Riesz' convexity theorem} can be applied to obtain
that
\begin{equation}
\label{estimate for int |E_m(T(f))(x)|^r}
	\Bigl(\int_{[0,1)} |E_m(T(f))(x)|^r \, dx \Bigr)^{1/r}
   \le N_p^t \, N_q^{1-t} \, \Big(\int_{[0,1)} |f(x)|^r \, dx \Bigr)^{1/r}
\end{equation}
for all step functions $f$ on $[0,1)$ which are constant on dyadic
intervals of length $2^{-m}$.  In other words, we think of $E_m \circ
T$ as a linear transformation from the space of step functions on
$[0,1)$ which are constant on dyadic intervals of length $2^{-m}$ to
itself.  This space can be identified with ${\bf C}^n$, $n = 2^m$ in a
simple way, with the components of vectors corresponding to the values
of functions on the dyadic subintervals of $[0,1)$ of length $2^{-m}$.
(One can work with real-valued functions instead of complex-valued
functions, as long as $T$ maps real functions to real functions.)  In
this identification, the integrals of powers of absolute values of
functions that we have here correspond to the sums that we considered
before, except for constant factors, which can easily be handled.
This permits us to derive (\ref{estimate for int |E_m(T(f))(x)|^r})
from the earlier result for linear mappings on ${\bf C}^n$.  (One might
note that this process is reversible, i.e., one could go from vectors
in ${\bf C}^n$ to functions.)

	Once one has (\ref{estimate for int |E_m(T(f))(x)|^r}) for
step functions which are constant on dyadic intervals of length
$2^{-m}$, where $m$ is arbitrary, it is easy to derive (\ref{estimate
for int |T(f)(x)|^r}) for dyadic step functions in general.  (Of course
one could extend this to other kinds of functions too.)

\section{A generalization}
\label{A generalization (interpolation)}

	The maximal and square function operators $M$, $M_l$, $S$, $S_l$
are not linear, but the previous results also work for them.  The reason
for this is that although these operators are not linear, they are
given by suprema of absolute values of other linear mappings.

	More precisely, suppose that we have a linear operator $T$
which is not linear, but which can be described in terms of a family
$\{T_\alpha\}_\alpha$ of linear operators in the following manner:
$|T_\alpha(f)(x)| \le T(f)(x)$ for all $\alpha$, $f$, and $x$; and for
any given function $f$ there is an $\alpha$ such that $T(f) =
|T_\alpha(f)|$ (i.e., $T(f)(x) = |T_\alpha(f)(x)|$ for all $x$).  (One
could adjust this to allow for suitable approximations.)  Then there
are interpolation results for $T$ which are just like the ones for
linear operators.  That is, one would assume conditions like
(\ref{(int_{[0,1)} |T(f)(x)|^p dx)^{1/p} le ..}), (\ref{(int_{[0,1)}
|T(f)(x)|^q dx)^{1/q} le ..}), (\ref{sup_{x in [0,1)} |T(f)(x)| le
N_infty sup_{[0,1)} |f(x)|}) for $T$, and these imply analogous
conditions for all of the linear operators $T_\alpha$, uniformly in
$\alpha$.  One could then apply the interpolation result to get
(\ref{estimate for int |T(f)(x)|^r}) for $T_\alpha$, uniformly in
$\alpha$, and the analogous assertion for $T$ would be a consequence
because of the way that $T(f)$ can be given in terms of $T_\alpha$ for
some $\alpha$, depending on $f$.

	This is the situation for the maximal and square function
operators mentioned above.  For maximal functions, one lets $\alpha$
run through non-negative integer-valued functions on $[0,1)$, and for
$T_\alpha$ one puts
\begin{equation}
	T_\alpha(f)(x) = E_{\alpha(x)}(f)(x).
\end{equation}
Thus $T_\alpha(f)$ is linear in $f$, $T_\alpha$ can be bounded in
terms of maximal functions, and one can have inequalities in the other
direction by choosing $\alpha$ properly for a given $f$.  In the
context of square functions, one lets $\alpha(x)$ be a sequence-valued
function such that
\begin{equation}
	\Bigl(\sum_{i=0}^\infty |\alpha_i(x)|^2 \Bigr)^{1/2} = 1
\end{equation}
for all $x$.  For the counterpart of $T_\alpha$, one considers linear
operators of the form
\begin{equation}
	\alpha_0(x) \, E_0(f)(x) 
		+ \sum_{i=1}^\infty \alpha_i(x) (E_i(f)(x) - E_{i-1}(f)(x)).
\end{equation}
These expressions are linear in $f$ and bounded by square functions,
because of the Cauchy-Schwarz inequality.  One can go the other way by
choosing $\alpha$ properly for a given $f$.

	In this manner, one can get the same kind of interpolation
estimates for $M$, $M_l$, $S$, $S_l$ as for linear operators.

\chapter{Quasisymmetric mappings}
\index{quasisymmetric mappings}

\section{Basic notions}

	Let $(M, d(x,y))$ and $(N, \rho(u,v))$ be metric spaces.  Thus
$M$ is a nonempty set and $d(x,y)$ is a nonnegative real-valued
function on $M \times M$ such that $d(x,y) = 0$ exactly when $x = y$,
$d(x,y) = d(y,x)$ for all $x, y \in M$, and 
\begin{equation}
	d(x,z) \le d(x,y) + d(y,z)  \qquad\hbox{for all } x, y, z \in M
\end{equation}
(the triangle inequality),\index{triangle inequality} and similarly
for $(N, \rho(u,v))$.

	Suppose that $\eta(t)$ is a nonnegative real-valued function
on $[0,\infty)$.  Let us say that $\eta$ is
\emph{admissible}\index{admissible functions on $[0,\infty)$} if
$\eta(0) = 0$ and
\begin{equation}
	\lim_{t \to 0} \eta(0) = 0.
\end{equation}

	A mapping $f : M \to N$ is said to be \emph{quasisymmetric} if
$f$ is not constant (unless $M$ contains only one element) and if
there is an admissible function $\eta : [0,\infty) \to [0,\infty)$
such that
\begin{equation}
\label{rho(f(y), f(x)) le eta(t) rho(f(z), f(x))}
	\rho(f(y), f(x)) \le \eta(t) \, \rho(f(z), f(x))
\end{equation}
whenever $x$, $y$, and $z$ are elements of $M$ and $t$ is a positive
real number such that
\begin{equation}
\label{d(y,x) le t d(z,x)}
	d(y,x) \le t \, d(z,x).
\end{equation}

	Roughly speaking, this condition asks that \emph{relative}
distances be approximately preserved by $f$.  For instance, if $y$ is
much closer to $x$ than $z$ is to $x$, then (\ref{d(y,x) le t d(z,x)})
holds with $t$ small, and (\ref{rho(f(y), f(x)) le eta(t) rho(f(z),
f(x))}) implies that $f(y)$ is relatively close to $f(x)$, compared to
the distance between $f(z)$ and $f(x)$.  Similarly, if $y$ is not too
close to $x$ compared to the distance from $z$ to $x$, but still $y$
is not too far from $x$ compared to the distance from $z$ to $x$, then
one can have (\ref{d(y,x) le t d(z,x)}) with a $t$ which is not too
small and not too large, and (\ref{rho(f(y), f(x)) le eta(t) rho(f(z),
f(x))}) implies that $f(y)$ is not too far from $f(x)$ compared to the
distance from $f(z)$ to $f(x)$.

	We shall describe some examples in the next two sections, but
first let us mention a few simple facts pertaining to the definition.

	A quasisymmetric mapping $f : M \to N$ is one-to-one.  Indeed,
if $x$ and $z$ are two distinct points in $M$ such that $f(x) = f(z)$,
then one can use the condition above to obtain that $f(x) = f(y)$ for
all $y \in M$.  Constant mappings are excluded from the definition of
quasisymmetric mappings, except in the case where $M$ has only one
element, and we conclude that $f$ is one-to-one.

	If $\eta : [0,\infty) \to [0,\infty)$ is admissible, then the
function $\widehat{\eta} : [0,\infty) \to [0,\infty)$ given by
\begin{equation}
	\widehat{\eta}(t) = \inf \{\eta(s) : s \ge t \}
\end{equation}
satisfies
\begin{equation}
	\widehat{\eta}(t) \le \eta(t) \qquad\hbox{for all } t \in [0,\infty)
\end{equation}
and is admissible in its own right.  It is easy to see that
$\widehat{\eta}$ is monotone increasing, and is hence a kind of
regularization of $\eta$.  On the other hand, if $f$ is any
quasisymmetric mapping relative to the function $\eta$, then $f$ is
also quasisymmetric relative to $\widehat{\eta}$.  This follows easily
from the definition.

	On the other hand, if $f$ is a quasisymmetric mapping relative
to a function $\eta$, then $f$ is also quasisymmetric relative to any
function $\theta : [0,\infty) \to [0,\infty)$ such that
\begin{equation}
	\eta(t) \le \theta(t) \qquad\hbox{for all } t \in [0,\infty).
\end{equation}
In particular, one can replace $\eta(t)$ with $\theta(t) = \eta(t) +
\epsilon \, t$ for any fixed $\epsilon > 0$.  If $\eta(t)$ is already
monotone increasing on $[0,\infty)$, which can always be arranged as
in the previous paragraph, then this choice of $\theta(t)$ will be
strictly increasing.  It will also tend to $\infty$ as $t \to \infty$,
so that $\theta$ maps $[0,\infty)$ onto $[0,\infty)$.  

	Now assume that $f : M \to N$ is quasisymmetric relative to a
function $\eta$ and that $f(M) = N$.  In this case the inverse mapping
$f^{-1} : N \to M$ is defined, since $f$ is one-to-one, and let us
show that it is quasisymmetric.  By the remarks in the previous
paragraph, we may assume that $\eta : [0,\infty) \to [0,\infty)$ is
invertible.  We would like to show that $f^{-1}$ is quasisymmetric
relative to $\alpha : [0,\infty) \to [0,\infty)$ defined by
\begin{equation}
	\alpha(t) = \frac{1}{\eta^{-1}(t^{-1})} 
			\quad\hbox{when } t > 0, \quad \alpha(0) = 0.
\end{equation}
This is equivalent to the assertion that 
\begin{equation}
	d(y,x) \le \alpha(t) \, d(z,x)
\end{equation}
whenever $x$, $y$, and $z$ are elements of $M$ and $t$ is a
positive real number such that
\begin{equation}
	\rho(f(y), f(x)) \le t \, \rho(f(z), f(x)).
\end{equation}
It is convenient to rephrase this as the statement that
\begin{equation}
\label{d(y,x) le t d(z,x), 2}
	d(y,x) \le t \, d(z,x)
\end{equation}
whenever $x$, $y$, and $z$ are elements of $M$ and $t$ is a
positive real number such that
\begin{equation}
\label{rho(f(y), f(x)) le alpha^{-1}(t) rho(f(z), f(x))}
	\rho(f(y), f(x)) \le \alpha^{-1}(t) \, \rho(f(z), f(x)).
\end{equation}

	We may as well assume that $y \ne x$, since otherwise
(\ref{d(y,x) le t d(z,x), 2}) is trivial.  If (\ref{d(y,x) le
t d(z,x), 2}) does not hold, then
\begin{equation}
\label{d(y,x) > t d(z,x)}
	d(y,x) > t \, d(z,x),
\end{equation}
or
\begin{equation}
	d(z,x) < t^{-1} \, d(y,x),
\end{equation}
and hence
\begin{equation}
	d(z,x) \le s \, d(y,x)
\end{equation}
for some $s < t^{-1}$.  It follows that
\begin{equation}
	\rho(f(z), f(x)) \le \eta(s) \, \rho(f(y), f(x)),
\end{equation}
using the original quasisymmetry condition for $f$ with the
roles of $y$ and $z$ reversed, and with $t$ replaced with $s$.
This implies in turn that
\begin{equation}
	\rho(f(z), f(x)) < \eta(t^{-1}) \, \rho(f(y), f(x)),
\end{equation}
or
\begin{equation}
	\frac{1}{\eta(t^{-1})} \, \rho(f(z), f(x)) < \rho(f(y), f(x)).
\end{equation}
This is the same as saying that (\ref{rho(f(y), f(x)) le alpha^{-1}(t)
rho(f(z), f(x))}) does not hold, which is what we want.

	The composition of two quasisymmetric mappings is
quasisymmetric.  For the associated function $\eta$ for the
composition, one can use the composition of the functions $\eta$
associated to the two mappings being composed.

\section{Examples}

	For any metric space $(M, d(x,y))$, the identity mapping on $M$
is quasisymmetric, with respect to the function $\eta(t) \equiv t$.
More generally, let $(M, d(x,y))$ and $(N, \rho(u,v))$ be metric spaces,
and suppose that $f : M \to N$ is an \emph{isometry},\index{isometry}
so that
\begin{equation}
	\rho(f(x), f(y)) = d(x,y) \qquad\hbox{for all } x, y \in M.
\end{equation}
Then $f$ is quasisymmetric with $\eta(t) \equiv t$.

	A mapping $f : M \to N$ is said to be
\emph{bilipschitz\index{bilipschitz mappings} with constant $C > 0$}
if
\begin{equation}
	C^{-1} \, d(x,y) \le \rho(f(x), f(y)) \le C \, d(x,y)
\end{equation}
for all $x, y \in M$.  If $f$ is bilipschitz with constant $C$, then
$f$ is quasisymmetric with $\eta(t) = C^2 \, t$.  Note that $f$ is
an isometry exactly when it is bilipschitz with constant $1$.

	Let us focus for the moment on the case where $M = N = {\bf
R}^n$ for some positive integer $n$, with the usual metric $|x-y|$,
where $|x|$ denotes the Euclidean norm of $x$.  Suppose that $f : {\bf
R}^n \to {\bf R}^n$ is an isometry.  We can write $f$ as $T(x) + x_0$,
where $x_0$ is a fixed element of ${\bf R}^n$ and $T : {\bf R}^n \to
{\bf R}^n$ is an isometry which fixes the origin, by taking $x_0 =
f(0)$.  A well-known result (which is not too hard to prove) states
that $T$ is then linear.  Thus $T$ is an orthogonal transformation, as
discussed near the end of Section \ref{Inner product spaces (and
spectral theory)}.

	A mapping $f : {\bf R}^n \to {\bf R}^n$ is called a
\emph{similarity}\index{similarity} if it can be written as $f(x) = a
\, T(x) + x_0$, where $a$ is a positive real number, $T$ is an
orthogonal transformation on ${\bf R}^n$, and $x_0$ is a fixed element
of ${\bf R}^n$.  It is easy to see that a similarity is quasisymmetric
on ${\bf R}^n$, with $\eta(t) \equiv t$.  A similarity $f(x) = a \,
T(x) + x_0$ is also bilipschitz with constant $C = \max(a, a^{-1})$.
In particular, the bilipschitz constant depends on the scale factor
$a$, while the quasisymmetry function $\eta$ does not.

	Conversely, if $f : {\bf R}^n \to {\bf R}^n$ is quasisymmetric
with $\eta(t) \equiv t$, then $f$ is a similarity.  Here is an outline of
the argument.  First, if $f$ has this property, then $f^{-1}$ does too,
as in the previous section.  Hence
\begin{equation}
	|f(y) - f(x)| = t \, |f(z) - f(x)|
\end{equation}
when $x, y, z \in {\bf R}^n$ and $t > 0$ satisfy 
\begin{equation}
	|y-x| = t \, |z-x|.
\end{equation}
Using this, one can verify that there is an $a > 0$ such that
\begin{equation}
	|f(u) - f(v)| = a \, |u-v|
\end{equation}
for all $u, v \in {\bf R}^n$.  In other words, $a^{-1} f(x)$ is an
isometry on ${\bf R}^n$, and hence $f$ is a similarity.

	Now suppose that $b$ is a positive real number, and consider
the mapping $f$ on ${\bf R}^n$ defined by
\begin{equation}
\label{f(x) = |x|^{b-1} x}
	f(x) = |x|^{b-1} \, x.
\end{equation}
At $x = 0$ we interpret this as meaning that $f(0) = 0$ (which would
really only be in question for $b < 1$, and even then this is the
natural choice in terms of continuity).  We leave it as an exercise to
show that these mappings are quasisymmetric.

	As an extension of this example, suppose that for each real
number $s$ one chooses an orthogonal transformation $R_s$ on ${\bf R}^n$,
in such a way that
\begin{equation}
\label{||R_s - R_t|| le C |s-t|}
	\|R_s - R_t\| \le C \, |s-t|
\end{equation}
for some constant $C$ and all $s, t \in {\bf R}$.  Here $\|\cdot\|$ denotes
the operator norm for linear transformations on ${\bf R}^n$ (and any
other norm would work practically as well, affecting only the constant
$C$ and not this Lipschitz condition otherwise).  Then the mapping
on ${\bf R}^n$ given by
\begin{equation}
\label{x mapsto |x|^{b-1} R_{log |x|}(x) when x ne 0, 0 mapsto 0}
	x \mapsto |x|^{b-1} \, R_{\log |x|}(x) \quad\hbox{when } x \ne 0,
		\quad 0 \mapsto 0
\end{equation}
is also quasisymmetric, for all $b > 0$.

	This mapping is actually the composition of (\ref{f(x) =
|x|^{b-1} x}) and the mapping defined by
\begin{equation} 
\label{x mapsto R_{log |x|}(x) when x ne 0, 0 mapsto 0}
	x \mapsto R_{\log |x|}(x) \quad\hbox{when } x \ne 0,
		\quad 0 \mapsto 0.
\end{equation}
One can show that (\ref{x mapsto R_{log |x|}(x) when x ne 0, 0 mapsto
0}) is bilipschitz on ${\bf R}^n$, under the assumption (\ref{||R_s -
R_t|| le C |s-t|}).  As a result, the quasisymmetry of (\ref{x mapsto
|x|^{b-1} R_{log |x|}(x) when x ne 0, 0 mapsto 0}) follows from that
of (\ref{f(x) = |x|^{b-1} x}).

\section{Cantor sets}
\label{Cantor sets}

	Let $K$ denote the usual Cantor ``middle-thirds'' set
contained in $[0,1]$, as in \cite{Ru1}.  Thus $K$ is a closed set
obtained by removing the (open) middle-third $(1/3, 2/3)$ from
$[0,1]$, then removing the middle-thirds of the two intervals that
result from that, then removing the middle-thirds of the four
intervals coming from the previous step, and so on.  One can write $K$
as $\bigcap_{j=0}^\infty E_j$, where $E_j$ is the union of $2^j$
disjoint closed intervals of length $3^{-j}$ produced after $j$ steps
of the construction, with $E_0 = [0,1]$ and $E_{j+1} \subseteq E_j$
for all $j$.

	If $r$ is a real number such that $0 < r < 1$, we can define
an analogous set $K(r) \subseteq [0,1]$ using $r$ instead of $1/3$.
In the first step one removes $(1/2 - r/2, 1/2 + r/2)$ from $[0,1]$ to
get two closed intervals of length $(1-r)/2$, at the ends of $[0,1]$.
After $j$ steps, one gets a set $E_j(r)$ which is the union of $2^j$
disjoint closed intervals of length $((1-r)/2)^j$, where $E_0(r) =
[0,1]$ and $E_{j+1}(r) \subseteq E_j(r)$ for all $j$.  As before, we
can define $K(r)$ to be $\bigcap_{j=0}^\infty E_j(r)$, and $K(1/3)$ is
equal to the set $K$ from the previous paragraph.

	Fix $r \in (0,1)$.  There is a natural one-to-one mapping
$h_r$ from $K$ onto $K_r$ that one can define.  The basic idea is
that for each nonnegative integer $j$ there are $2^j$ intervals
in $E_j$ and in $E_j(r)$, and for each positive integer $m$ such
that $m \le 2^j$ we want $h_r$ to map the intersection of $K$ with
the $m$th interval of length $3^{-j}$ in $E_j$ to the intersection of
$K(r)$ with the $m$th interval of length $((1-r)/2)^j$ in $E_j(r)$.
When we refer to the $m$th interval in $E_j$ or $E_j(r)$, we go
from left to right.  It is easy to see that these requirements for
different $j$'s are compatible with each other.

	Note that $0$ and $1$ lie in both $K$ and $K(r)$, and $h_r$ is
defined so that $h_r(0) = 0$ and $h_r(1) = 1$.  Similarly, for $j \ge
1$ the endpoints of the intervals in $E_j$ all lie in $K$ (because
they remain in $E_i$ for all $i > j$), and the endpoints of the
intervals in $E_j(r)$ lie in $K(r)$ for the same reason.  The
endpoints of the intervals in $E_j$ are sent to the corresponding
endpoints of the corresponding intervals in $E_j(r)$ by $h_r$.  These
assignments for different $j$'s are compatible with each other because
of the way that the various intervals fit together.

	Notice also that $h_r$ is not a linear function (of the form
$c \, x + d$) on the intersection of $K$ with any open interval that
contains an element of $K$.

	One can view $K$ and $K_r$ as metric spaces, using the
ordinary Euclidean metric $d(x,y)$ on each.  It is not too difficult
to show that $h_r : K \to K_r$ is a quasisymmetric mapping.

	Here is a variant of this which is somewhat more complicated.
Define a set $J \subseteq [0,1]$ by removing two-fifths at each step,
in separate pieces, rather than one-third in one piece.  Specifically,
we start with $[0,1]$, as before, and in the first step we remove the
two open intervals $(1/5, 2/5)$ and $(3/5, 4/5)$.  This leaves the
three closed intervals $[0, 1/5]$, $[2/5, 3/5]$, and $[4/5, 1]$.  In
each of these one removes two ``fifths'' again, to get a total of $9$
intervals of length $5^{-2}$.  In general, after $l$ steps, there are
$3^l$ closed intervals of length $5^{-l}$.  If $F_l$ denotes
the union of these $3^l$ intervals, then $F_{l+1} \subseteq F_l$
for all $l$ and $J$ is defined to be $\bigcap_{l=0}^\infty F_l$.

	To define a mapping from $K$ to $J$, it is convenient to
give a different description of $K$.  Let us combine the first step
in the construction of $K$ with half of the second step, as follows.
In the new first step, we start with the interval $[0,1]$, and
remove two subintervals $(1/3, 2/3)$ and $(7/9, 8/9)$ to get the
three closed intervals $[0, 1/3]$, $[2/3, 7/9]$, and $[8/9, 1]$.
Now the intervals no longer have the same length, but we shall not
be too concerned about that.  

	One can repeat the process to each of these three intervals,
obtaining a total of $9$ closed subintervals of various lengths.
In general, after $l$ steps, one gets $3^l$ closed intervals of various
lengths.  Let us write $\widehat{E}_l$ for the set which is the union
of these $3^l$ intervals.

	Each of these $l$ intervals in $\widehat{E}_l$ also occurs
in the $E_j$'s.  More precisely, each of these subintervals occurs
in an $E_j$ with $l \le j \le 2 \, l$.  It is not hard to check that
\begin{equation}
	E_{2l} \subseteq \widehat{E}_l \subseteq E_l
\end{equation}
for all $l \ge 0$.  As a result, 
\begin{equation}
	K = \bigcap_{l=0}^\infty \widehat{E}_l.
\end{equation}

	This description of $K$ can be used to define a one-to-one
mapping from $K$ onto $J$.  As before, this mapping takes the
intersection of $K$ with the $n$th interval in $\widehat{E}_l$ to the
intersection of $J$ with the $n$th interval in $F_l$, $1 \le n \le
3^l$, $l \ge 0$.  It also sends $0$ to $0$, $1$ to $1$, and, in general,
the endpoints of the intervals in $\widehat{E}_l$ to the corresponding
endpoints of the corresponding intervals in $F_l$.

	It is not too difficult to show that this mapping from $K$
onto $J$ is also quasisymmetric, using the Euclidean metric on both
$K$ and $J$.

\section{Bounds in terms of $C \, t^a$}

	Suppose that $f : {\bf R}^n \to {\bf R}^n$ is a quasisymmetric
mapping (with the standard Euclidean metric is used on ${\bf R}^n$).
Then there are positive real numbers $C$, $a_1 \le 1$, and $a_2 \ge 1$
such that $f$ is quasisymmetric with respect to
\begin{equation}
	\eta(t) = C \, t^{a_1} \quad\hbox{for } 0 \le t \le 1,
		\quad \eta(t) = C \, t^{a_2} \quad\hbox{for } t \ge 1.
\end{equation}

	Here is the basic idea.  Because $f$ is quasisymmetric, there
are positive real numbers $t_1 \le 1$ and $L$ such that
\begin{eqnarray}
\label{|f(y) - f(x)| le frac{1}{2} |f(z) - f(x)| when |y-x| le t_1 |z-x|}
	&& |f(y) - f(x)| \le \textstyle\frac{1}{2} \, |f(z) - f(x)|	
		\hbox{ for all } 				\\
	&& x, y, z \in {\bf R}^n \hbox{ such that }
		|y-x| \le t_1 \, |z-x|				\nonumber
\end{eqnarray}
and
\begin{eqnarray}
\label{|f(y) - f(x)| le L |f(z) - f(x)| when |y-x| le 2 |z-x|}
	&& |f(y) - f(x)| \le L \, |f(z) - f(x)|	
		\hbox{ for all } 				\\
	&& x, y, z \in {\bf R}^n \hbox{ such that }
		|y-x| \le 2 \, |z-x|.				\nonumber
\end{eqnarray}
One can ``iterate'' these statements to obtain that
\begin{eqnarray}
	&& |f(y) - f(x)| \le \textstyle\frac{1}{2^k} \, |f(z) - f(x)|	
		\hbox{ for all } 				\\
	&& x, y, z \in {\bf R}^n \hbox{ such that }
		|y-x| \le t_1^k \, |z-x|			\nonumber
\end{eqnarray}
and
\begin{eqnarray}
	&& |f(y) - f(x)| \le L^k \, |f(z) - f(x)|	
		\hbox{ for all } 				\\
	&& x, y, z \in {\bf R}^n \hbox{ such that }
		|y-x| \le 2^k \, |z-x|				\nonumber
\end{eqnarray}
for all positive integers $k$.

	To be explicit, let us consider the case where $k = 2$.
Suppose that $x, y, z \in {\bf R}^n$ satisfy $|y-x| \le t_1^2 \, |z-x|$.
Let $\widetilde{z}$ be an element of ${\bf R}^n$ such that
\begin{equation}
	|\widetilde{z} - x| = t_1 \, |z-x|.
\end{equation}
Then
\begin{equation}
	|y-x| \le t_1 \, |\widetilde{z} - x|,
\end{equation}
and hence
\begin{equation}
\label{|f(y) - f(x)| le frac{1}{2} |f(widetilde{z}) - f(x)|}
	|f(y) - f(x)| \le \textstyle\frac{1}{2} \, |f(\widetilde{z}) - f(x)|,
\end{equation}
by (\ref{|f(y) - f(x)| le frac{1}{2} |f(z) - f(x)| when |y-x| le t_1
|z-x|}) with $z$ replaced by $\widetilde{z}$.  Similarly,
\begin{equation}
\label{|f(widetilde{z}) - f(x)| le frac{1}{2} |f(z) - f(x)|}
	|f(\widetilde{z}) - f(x)| \le \textstyle\frac{1}{2} \, |f(z) - f(x)|,
\end{equation}
by (\ref{|f(y) - f(x)| le frac{1}{2} |f(z) - f(x)| when |y-x| le t_1
|z-x|}) with $y$ replaced by $\widetilde{z}$ (and $z$ kept as it is
this time).  Combining (\ref{|f(y) - f(x)| le frac{1}{2}
|f(widetilde{z}) - f(x)|}) and (\ref{|f(widetilde{z}) - f(x)| le
frac{1}{2} |f(z) - f(x)|}) we obtain
\begin{equation}
	|f(y) - f(x)| \le \textstyle\frac{1}{2^2} \, |f(z) - f(x)|,
\end{equation}
as desired.  One can treat the case where $|y-x| \le 2^2 \, |z-x|$
similarly, using an element $\widehat{z}$ of ${\bf R}^n$ such that
\begin{equation}
	|\widehat{z} - x| = 2 \, |z-x|.
\end{equation}

	This argument works for quasisymmetric mappings between metric
spaces in general under modest assumptions on the domain.  To be more
precise, suppose that we have a quasisymmetric mapping whose domain is
a metric space $(M, d(x,y))$.  If $(M, d(x,y))$ has the property that
for each $x \in M$ and each positive real number $r$ there is a point
$w \in M$ such that $d(x,w) = r$, then exactly the same argument as
above goes through.  One can also work with bounded metric spaces, by
restricting the $r$'s to suitable ranges.  In particular, these
conditions hold when $M$ is connected.  

	In the case of the kind of Cantor sets considered in the
previous section, some adjustments are needed to take care of the
``gaps'', and this is not too hard to do.  Basically, one should be a
bit more careful with the points $\widetilde{z}$, $\widehat{z}$ as
above.  In this regard, it be preferable to use another number instead
of $2$ in (\ref{|f(y) - f(x)| le L |f(z) - f(x)| when |y-x| le 2
|z-x|}), and to modify the choice of $t_1$.

	Although it is convenient to think in terms of the domain
here, there is a natural symmetry between the domain and the range,
reflected in the fact that a surjective mapping is quasisymmetric
if and only if its inverse is quasisymmetric.

\backmatter

\newpage

\addcontentsline{toc}{chapter}{Index}
\printindex


\addcontentsline{toc}{chapter}{Bibliography}

\begin{thebibliography}{FraJW}


\bibitem [Ash] {Ash} J.~Ash, editor, {\it Studies in Harmonic
Analysis}, Proceedings of the expository conference ``A Survey of
Harmonic Analysis'', Mathematical Association of America, 1976.

\bibitem [BerL] {BL} J.~Bergh and J.~L\"oftstr\"om, {\it Interpolation
Spaces: An Introduction}, Springer-Verlag, 1976.

\bibitem [Bur1] {Bur1} D.~Burkholder, {\it Boundary value problems
and sharp inequalities for martingale transforms}, Annals of Probability
{\bf 12} (1984), 647--702.

\bibitem [Bur2] {Bur2} D.~Burkholder, {\it Some extremal problems in
martingale theory and harmonic analysis}, in {\it Harmonic Analysis
and Partial Differential Equations: Essays in Honor of Alberto
P.~Calder\'on}, 99--115, University of Chicago Press, 1999.

\bibitem [CalZ] {CZ} A.~Calder\'on and A.~Zygmund, {\it On the
existence of certain singular integrals}, Acta Mathematica {\bf 88}
(1952), 85--139.

\bibitem [CoiM] {CM} R.~Coifman and Y.~Meyer, {\it Au-Del\`a des
Op\'erateurs Pseudo-Diff\'erentiels}, Ast\'erisque {\bf 57},
Soci\'et\'e Math\'ematique de France, 1978.

\bibitem [CoiW1] {CW1} R.~Coifman and G.~Weiss, {\it Analyse Harmonique
Non-Commutative sur Certains Espaces Homog\`enes}, Lecture Notes in
Mathematics {\bf 242}, Springer-Verlag, 1971.

\bibitem [CoiW2] {CW2} R.~Coifman and G.~Weiss, {\it Extensions of Hardy
spaces and their use in analysis}, Bulletin of the American
Mathematical Society {\bf 83} (1977), 569--645.

\bibitem [CoiW3] {CW3} R.~Coifman and G.~Weiss, {\it Transference
Methods in Analysis}, Conference Board of the Mathematical Sciences
Regional Conference Series in Mathematics {\bf 31}, American
Mathematical Society, 1977.

\bibitem [Die] {Dieudonne} J.~Dieudonn\'e, {\it Special Functions and
Linear Representations of Lie Groups}, Conference Board of the
Mathematical Sciences Regional Conference Series in Mathematics {\bf
42}, American Mathematical Society, 1980.

\bibitem [Duo] {Duo} J.~Duoandikoetxea, {\it Fourier Analysis},
translated from the Spanish edition and revised by D.~Cruz-Uribe, SFO,
American Mathematical Society, 2001.

\bibitem [Dur] {Duren} P.~Duren, {\it Theory of $H^p$ Spaces}, Academic
Press, 1970.

\bibitem [FraJW] {FJW} M.~Frazier, B.~Jawerth, and G.~Weiss, {\it
Littlewood--Paley Theory and the Study of Function Spaces}, Conference
Board of the Mathematical Sciences Regional Conference Series in
Mathematics {\bf 79}, American Mathematical Society, 1991.

\bibitem [GarcR] {GCRF} J.~Garc\'{\i}a-Cuerva and J.~Rubio de Francia,
{\it Weighted Norm Inequalities and Related Topics}, North-Holland,
1985.

\bibitem [Garn] {Ga} J.~Garnett, {\it Bounded Analytic Functions},
Academic Press, 1981.

\bibitem [Hal1] {Halmos1} P.~Halmos, {\it Introduction to Hilbert Space},
Chelsea Publishing Company, 1957.

\bibitem [Hal2] {Halmos2} P.~Halmos, {\it A Hilbert Space Problem Book},
Van Nostrand, 1967.

\bibitem [Hei] {Juha} J.~Heinonen, {\it Lectures on Analysis on Metric
Spaces}, Springer-Verlag, 2001.

\bibitem [HerW] {HerWei} E.~Hern\'andez and G.~Weiss, {\it A First Course
on Wavelets}, CRC Press, 1996.

\bibitem [HewS] {HS} E.~Hewitt and K.~Stromberg, {\it Real and
Abstract Analysis}, Springer-Verlag, 1965.

\bibitem [Kat] {Kat} Y.~Katznelson, {\it An Introduction to Harmonic
Analysis}, Dover Publications, 1976.

\bibitem [Koo] {Koosis} P.~Koosis, {\it Introduction to $H_p$ Spaces},
Cambridge University Press, 1980.

\bibitem [Jou] {Journe} J.-L.~Journ\'e, {\it Calder\'on--Zygmund
Operators, Pseudodifferential Operators, and the Cauchy Integral of
Calder\'on}, Lecture Notes in Mathematics {\bf 994}, Springer-Verlag,
1983.

\bibitem [Mey1] {Meyer1} Y.~Meyer (Volume III with R.~Coifman), {\it
Ondelettes et Op\'erateurs I, II, III: Ondelettes, Op\'erateurs de
Calder\'on--Zygmund, Op\'erateurs Multilin\'eaires}, Hermann, 1990,
1991; English translations, {\it Wavelets and Operators}, {\it
Wavelets: Calder\'on--Zygmund and Multilinear Operators}, translated
by D.~Salinger, Cambridge University Press, 1992, 1997.

\bibitem [Mey2] {Meyer2} Y.~Meyer, {\it Ondelettes et Algortihmes
Concurrents}, Hermann, 1992.

\bibitem [Ner] {Neri} U.~Neri, {\it Singular Integrals}, Lecture Notes
in Mathematics {\bf 200}, Springer-Verlag, 1971.

\bibitem [Pe] {Peetre} J.~Peetre, {\it Marcel Riesz in Lund}, in {\it
Function Spaces and Applications}, Lecture Notes in Mathematics {\bf
1302}, 1--10, Springer-Verlag, 1988.

\bibitem [Pis] {Pisier} G.~Pisier, {\it Factorization of Linear
Operators and Geometry of Banach Spaces}, Conference Board of the
Mathematical Sciences Regional Conference Series in Mathematics {\bf
60}, American Mathematical Society, 1985.

\bibitem [Rie] {MRiesz} M.~Riesz, {\it Sur les maxima des formes
bilin\'eaires et sur les fonctionnelles lin\'eaires}, Acta Mathematica
{\bf 49} (1926), 465-497.

\bibitem [RieN] {RieN} F.~Riesz and B.~Sz.-Nagy, {\it Functional Analysis},
Dover Publications, 1990.

\bibitem [Rud1] {Ru1} W.~Rudin, {\it Principles of Mathematical Analysis},
McGraw-Hill, 1976.

\bibitem [Rud2]{Ru2} W.~Rudin, {\it Real and Complex Analysis},
McGraw-Hill, 1987.

\bibitem [Rud3] {Ru3} W.~Rudin, {\it Functional Analysis}, McGraw-Hill, 1973.

\bibitem [Sad] {Sadosky} C.~Sadosky, {\it Interpolation of Operators
and Singular Integrals: An Introduction to Harmonic Analysis}, Marcel
Dekker, 1979.

\bibitem [Sar] {Sar} D.~Sarason, {\it Function Theory on the Unit
Circle}, Virginia Polytechnic Institute and State University, 1978.

\bibitem [Sem1] {Se1} S.~Semmes, {\it A primer on Hardy spaces, and
some remarks on a theorem of Evans and M\"uller}, Communications in
Partial Differential Equations {\bf 19} (1994), 277--319.

\bibitem [Sem2] {Se2} S.~Semmes, {\it Metric spaces and mappings
seen at many scales}, appendix in {\it Metric Structures for
Riemannian and Non-Riemannian Spaces}, M.~Gromov et al., Birkh\"auser,
1999.

\bibitem [Ste1] {St1} E.~Stein, {\it Topics in Harmonic Analysis
Related to the Littlewood--Paley Theory}, Annals of Mathematics
Studies {\bf 63}, Princeton University Press, 1970.

\bibitem [Ste2] {St2} E.~Stein, {\it Singular Integrals and
Differentiability Properties of Functions}, Princeton University
Press, 1970.

\bibitem [Ste3] {St3} E.~Stein, {\it Harmonic Analysis: Real-Variable
Methods, Orthogonality, and Oscillatory Integrals}, Princeton
University Press, 1993.

\bibitem [SteW1] {SW1} E.~Stein and G.~Weiss, {\it Interpolation of
operators with change of measures}, Transactions of the American
Mathematical Society {\bf 87} (1958), 159--172.

\bibitem [SteW2] {SW-book} E.~Stein and G.~Weiss, {\it Introduction to
Fourier Analysis on Euclidean Spaces}, Princeton University Press,
1971.

\bibitem [StrT] {ST} J.-O.~Str\"omberg and A.~Torchinsky, {\it
Weighted Hardy Spaces}, Lecture Notes in Mathematics {\bf 1381},
Springer-Verlag, 1989.

\bibitem [Tor] {To} A.~Torchinsky, {\it Real-Variable Methods in
Harmonic Analysis}, Academic Press, 1986.

\bibitem [TroV] {TroV} D.~Trotsenko and J.~V\"ais\"al\"a, {\it
Upper sets and quasisymmetric maps}, Annales Academi{\ae}
Scientiarum Fennic{\ae} Mathematica {\bf 24} (1999), 465--488.

\bibitem [TukV] {TV} P.~Tukia and J.~V\"ais\"al\"a, {\it
Quasisymmetric embeddings of metric spaces}, Annales Academi{\ae}
Scientiarum Fennic{\ae} Mathematica {\bf 5} (1980), 97--114.

\bibitem [V\"ai] {Vai} J.~V\"ais\"al\"a, {\it The free quasiworld:
Freely quasiconformal and related mappings in Banach spaces}, in
{\it Quasiconformal Geometry and Dynamics}, 55--118, Banach Center
Publications {\bf 48}, Polish Academy of Sciences, 1999.

\bibitem [Wei1] {Weiss1} G.~Weiss, {\it Harmonic analysis}, in {\it
Studies in Real and Complex Analysis}, 124--178, Mathematical
Association of America, 1965.

\bibitem [Wei2] {Weiss2} G.~Weiss, {\it Harmonic analysis on compact
groups}, in \cite{Ash}, 198--223.

\bibitem [Woj] {Woj} P.~Wojtaszczyk, {\it Banach Spaces for Analysts},
Cambridge University Press, 1991.

\bibitem [Zyg1] {Z1} A.~Zygmund, {\it Trigonometric Series}, Volumes I and
II, Cambridge University Press, 1979.

\bibitem [Zyg2] {Z2} A.~Zygmund, {\it Int\'egrales Singuli\`eres},
Lecture Notes in Mathematics {\bf 204}, Springer-Verlag, 1971.



\end{thebibliography}
\end{document}